\documentclass[a4paper,12pt,leqno]{article}

\usepackage{amsmath}
\usepackage{amscd}
\usepackage[psamsfonts]{amssymb}
\usepackage{graphicx}
\usepackage{lscape}
\usepackage[french]{babel}
\usepackage[utf8]{inputenc}
\usepackage[all]{xy}
\usepackage{color}

\newcounter{sect}
\newcounter{subsect}[sect]
\newcommand{\sect}[1]
{\refstepcounter{subsect}\refstepcounter{sect}
\renewcommand{\thesubsect}{\thesect}
{\large{\bf\rule{0ex}{10ex}\thesect}
~~#1\rule[-2ex]{0ex}{2ex}\hfill}}

\newcommand{\subsect}[1]{\refstepcounter{subsect}
\renewcommand{\thesubsect}{\thesect.\arabic{subsect}} 
{{\bf\rule{0ex}{4ex}\thesubsect}
~~#1\rule[-2ex]{0ex}{2ex}\hfill}}

\newcounter{th}[subsect]
\renewcommand{\theth}{\thesubsect.\arabic{th}}
\newenvironment{theo}[1][]{\refstepcounter{th}\rule{0ex}{4ex}%
\textsc{Th\'{e}or\`{e}me}~\theth#1.~~\em}{\em%
\rule[-1.5ex]{0ex}{1.5ex}}

\newenvironment{pro}[1][]{\refstepcounter{th}\rule{0ex}{4ex}%
\textsc{Proposition}~\theth#1.~~\em}{\em\rule[-1.5ex]{0ex}{1.5ex}}

\newenvironment{lem}[1][]{\refstepcounter{th}\rule{0ex}{4ex}%
\textsc{Lemme}~\theth#1.~~\em}{\em\rule[-1.5ex]{0ex}{1.5ex}}

\newenvironment{cor}[1][]{\refstepcounter{th}\rule{0ex}{4ex}%
\textsc{Corollaire}~\theth#1.~~\em}{\em\rule[-1.5ex]{0ex}{1.5ex}}

\newenvironment{pro-def}[1][]{\refstepcounter{th}\rule{0ex}{4ex}%
\textsc{Proposition-D\'{e}finition}~\theth#1.~~\em}{\em\rule[-1.5ex]{0ex}{1.5ex}}

\newenvironment{cor-def}[1][]{\refstepcounter{th}\rule{0ex}{4ex}%
\textsc{Corollaire-D\'{e}finition}~\theth#1.~~\em}{\em\rule[-1.5ex]{0ex}{1.5ex}}

\newenvironment{lem-def}[1][]{\refstepcounter{th}\rule{0ex}{4ex}%
\textsc{Lemme-D\'{e}finition}~\theth#1.~~\em}{\em\rule[-1.5ex]{0ex}{1.5ex}}

\newenvironment{scho}[1][]{\refstepcounter{th}\rule{0ex}{4ex}%
\textsc{Scholie}~\theth#1.~~\em}{\em\rule[-1.5ex]{0ex}{1.5ex}}

\newenvironment{scho-def}[1][]{\refstepcounter{th}\rule{0ex}{4ex}%
\textsc{Scholie-D\'{e}finition}~\theth#1.~~\em}{\em\rule[-1.5ex]{0ex}{1.5ex}}

\def\numero{$\mathrm{n}^{\text{o}}$}

\begin{document}

\begin{center}

\large
\textbf{Complexes de modules équivariants \\ sur l'algèbre de Steenrod \\ associés à un $(\mathbb{Z}/2)^{n}$-CW-complexe fini}

\normalsize
\vspace{1cm}
\textsc{Dorra Bourguiba, Jean Lannes,\\ Lionel Schwartz} \textsl{et} \textsc{Saïd Zarati}

\end{center}

\sect{Introduction}

Soit $V$ un $2$-groupe abélien élémentaire~; $V$ est donc isomorphe à $(\mathbb{Z}/2)^{n}$ pour un certain entier $n$. Nous verrons le plus souvent $V$  comme un $\mathbb{Z}/2$-espace vectoriel dont la dimension $n$ sera notée $\dim V$.

\medskip
 Soit $X$ un $V$-CW-complexe (que nous supposerons le plus souvent fini, pour une référence sur la notion de CW-complexe équivariant, voir par exemple \cite[Chap. II, Sect. 1]{tD}). Nous étudions dans ce mémoire deux complexes de ``modules équivariants sur l'algèbre de Steenrod modulo $2$'' (cette notion sera bien sûr précisée dans la suite de cette introduction) associés à~$X$. Le premier est défini à l'aide de la filtration par les orbites de $X$ (voir ci-après), nous l'appellerons le ``complexe topologique''. Le second est défini de façon plus algébrique~; nous l'appellerons pour cela le ``complexe algébrique''.
 
\vspace{0.625cm}
\textsc{Le complexe topologique}
 
\medskip
Soit $W\subset V$ un sous-groupe~; la codimension de $W$, vu comme un sous-espace vectoriel de $V$, est notée $\mathop{\mathrm{codim}}W$. 
 
\medskip
Soit $p$ un entier avec $-1\leq p\leq n:=\dim V$~; on pose
$$
\hspace{24pt}
\mathrm{F}_{p}\hspace{1pt}X
\hspace{4pt}:=\hspace{4pt}
\bigcup_{\mathop{\mathrm{codim}}W\hspace{1pt}\leq\hspace{1pt}p} X^{W}
\hspace{24pt},
$$
$W$ décrivant l'ensemble des sous-groupes de $V$ de codimension inférieure ou égale à $p$ et $X^{W}$ désignant le sous-$V$-CW-complexe de $X$ constitué des points fixes par $W$. On a donc une filtration croissante de $X$ par des sous-$V$-CW-complexes~:
$$
\hspace{24pt}
\emptyset=\mathrm{F}_{-1}\hspace{1pt}X\subset\mathrm{F}_{0}\hspace{1pt}X\subset\mathrm{F}_{1}\hspace{1pt}X\subset\ldots\subset\mathrm{F}_{n-1}\hspace{1pt}X\subset\mathrm{F}_{n}\hspace{1pt}X=X
\hspace{24pt}.
$$
On observera que l'on a $\mathrm{F}_{0}\hspace{1pt}X=X^{V}$ et que $\mathrm{F}_{n-1}\hspace{1pt}X$ est le sous-$V$-CW-complexe de $X$ constitué des points dont le groupe d'isotropie est non trivial. Le sous-$V$-CW-complexe $\mathrm{F}_{n-1}\hspace{1pt}X$ sera aussi noté $\mathrm{Sing}_{V}X$ (``partie singulière'' de l'action de $V$ sur $X$)~; l'ouvert complémentaire $X-\mathrm{Sing}_{V}X$ est la réunion des orbites libres (``partie régulière'' de l'action de $V$ sur $X$). On définit un complexe de cochaînes
$$
\mathrm{C}_{\mathrm{top}}^{0}\hspace{1pt}X\to\mathrm{C}_{\mathrm{top}}^{1}\hspace{1pt}X\to\mathrm{C}_{\mathrm{top}}^{2}\hspace{1pt}X\to\ldots\to\mathrm{C}_{\mathrm{top}}^{p}\hspace{1pt}X\to\ldots\to\mathrm{C}_{\mathrm{top}}^{n}\hspace{1pt}X\to 0\to 0\to\ldots
$$
en posant
\begin{multline*}
\mathrm{C}_{\mathrm{top}}^{p}\hspace{1pt}X
\hspace{4pt}:=\hspace{4pt}
\Sigma^{-p}\hspace{2pt}\mathrm{H}_{V}^{*}(\mathrm{F}_{p}\hspace{1pt}X,\mathrm{F}_{p-1}\hspace{1pt}X;\mathbb{F}_{2})
\hspace{4pt}=\hspace{4pt}
\\
\bigoplus_{\mathop{\mathrm{codim}}W=p}
\Sigma^{-p}\hspace{2pt}
\mathrm{H}^{*}(V;\mathbb{F}_{2})
\otimes_{\mathrm{H}^{*}(V/W;\mathbb{F}_{2})}
\mathrm{H}_{V/W}^{*}(X^{W},\mathrm{Sing}_{V/W}X^{W};\mathbb{F}_{2})
\end{multline*}
et en prenant pour cobord le connectant de la triade $(\mathrm{F}_{p+1}\hspace{1pt}X,\mathrm{F}_{p}\hspace{1pt}X,\mathrm{F}_{p-1}\hspace{1pt}X)$. Décryptons un peu la notation. Les notations $\mathrm{H}_{V}^{*}( -;\mathbb{F}_{2})$ et $\mathrm{H}^{*}(V;\mathbb{F}_{2})$ dési\-gnent respectivement la cohomologie $V$-équivariante modulo $2$ et la cohomologie modulo $2$ du groupe $V$. Si $E=(E^{m})_{m\in\mathbb{Z}}$ est un ``objet'' gradué et $s$ un entier relatif alors $\Sigma^{s}E$ désigne l'objet gradué $(E^{m-s})_{m\in\mathbb{Z}}$. Ci-dessus $X^{W}$ est considéré comme un $V/W$-espace et $\mathrm{H}^{*}(V/W;\mathbb{F}_{2})$ est identifiée à une sous-algèbre de $\mathrm{H}^{*}(V;\mathbb{F}_{2})$.

\medskip
Le complexe $\mathrm{C}_{\mathrm{top}}^{\bullet}\hspace{1pt}X$ est muni d'une coaugmentation naturelle $\mathrm{H}_{V}^{*}(X;\mathbb{F}_{2})=:\mathrm{C}_{\mathrm{top}}^{-1}\hspace{1pt}X\to\mathrm{C}_{\mathrm{top}}^{0}\hspace{1pt}X$ ; le complexe coaugmenté associé est noté $\widetilde{\mathrm{C}}_{\mathrm{top}}^{\bullet}\hspace{1pt}X$.

\medskip
Notre premier résultat est le suivant~:

\begin{theo}\label{introCtop}
Soient $V$ un $2$-groupe abélien élémentaire et $X$ un $V$-CW-complexe fini. Les deux propriétés suivantes sont équivalentes~:
\begin{itemize}
\item[(i)] Le $\mathrm{H}^{*}(V;\mathbb{F}_{2})$-module $\mathrm{H}_{V}^{*}(X;\mathbb{F}_{2})$ est libre.
\item[(ii)] Le complexe coaugmenté $\widetilde{\mathrm{C}}_{\mathrm{top}}^{\bullet}\hspace{1pt}X$ est acyclique.
\end{itemize}
\end{theo}

Ici deux commentaires s'imposent~: 

\medskip
1) Le problème analogue, le $2$-groupe abélien élémentaire $V$ étant remplacé par un tore $T$ et $\mathrm{H}_{V}^{*}(-;\mathbb{F}_{2})$ par $\mathrm{H}_{T}^{*}(-;\mathbb{Q})$, est étudié en détail dans \cite{AFP1}.

\medskip
2) Peu de temps avant l'achèvement de la rédaction de ce mémoire est apparu le preprint de \cite{AFP2} qui traite notamment de $V$-espaces. Le théorème ci-dessus est un cas particulier du théorème 10.2 de cette référence (dans lequel $\widetilde{\mathrm{C}}_{\mathrm{top}}^{\bullet}\hspace{1pt}X$ est appelé {\em the augmented Atiyah-Bredon sequence for $X$}). Nous décrivons ci-après les techniques que nous employons pour démontrer le théorème \ref{introCtop}. Nous reviendrons sur \cite[Theorem 10.2]{AFP2} dans la dernière partie de cette introduction.

\pagebreak

\bigskip
Nous démontrons l'énoncé \ref{introCtop} en exploitant les structures ``linéaires'' de $\mathrm{H}_{V}^{*}(-;\mathbb{F}_{2})$ que nous passons en revue ci-après. La cohomologie que nous considérons dans ce mémoire (en particulier dans la suite de cette introduction) est  à coefficients dans $\mathbb{F}_{2}$, à l'exception de la section 1 dans laquelle elle est à coefficients dans $\mathbb{F}_{\ell}$, $\ell$ premier~; dans tous les cas  nous allègerons $\mathrm{H}^{*}(-;\mathbb{F}_{\ell})$ en $\mathrm{H}^{*}(-)$ et $\mathrm{H}_{V}^{*}(-;\mathbb{F}_{\ell})$ en $\mathrm{H}_{V}^{*}(-)$~; $\mathrm{H}^{*}(V;\mathbb{F}_{\ell})=\mathrm{H}_{V}^{*}(\mathrm{point};\mathbb{F}_{\ell})$ sera simplement notée $\mathrm{H}^{*}V$.

\medskip
Soit $(X,Y)$ une paire de $V$-espaces. La cohomologie équivariante  $\mathrm{H}_{V}^{*}(X,Y)$ est un $\mathbb{F}_{2}$-espace gradué. C'est un $\mathrm{H}^{*}V$-module (au sens gradué)~; on rappelle que $\mathrm{H}^{*}V$ est isomorphe à une algèbre de polynômes, à coefficients dans $\mathbb{F}_{2}$, en $n$ indéterminées, chacune de degré $1$. Soit $\mathrm{A}$ l'algèbre de Steenrod modulo $2$. La cohomologie équivariante  $\mathrm{H}_{V}^{*}(X,Y)$ est un $\mathrm{A}$-module et cet $\mathrm{A}$-module est instable ($\mathrm{Sq}^{i}:\mathrm{H}_{V}^{m}(X,Y)\to\mathrm{H}_{V}^{m+i}(X,Y)$ est trivial pour $i>m$)~; puisque l'on a $\mathrm{H}^{*}V=\mathrm{H}_{V}^{*}(\mathrm{point})$ il en est de même pour $\mathrm{H}^{*}V$. L'application de structure
$$
\mathrm{H}^{*}V\otimes\mathrm{H}_{V}^{*}(X,Y)\to\mathrm{H}_{V}^{*}(X,Y)
$$
est $\mathrm{A}$-linéaire ; on rappelle que la structure d'algèbre de Hopf de $\mathrm{A}$ fait du produit tensoriel de deux $\mathrm{A}$-modules (resp. $\mathrm{A}$-modules instables) un $\mathrm{A}$-module (resp. $\mathrm{A}$-module instable). On dit que $\mathrm{H}_{V}^{*}(X,Y)$ est un $\mathrm{H}^{*}V$-$\mathrm{A}$-module instable~; la catégorie de ces objets est notée $V$-$\mathcal{U}$ (la catégorie des $\mathrm{A}$-modules instables est elle notée $\mathcal{U}$). Si $(X,Y)$ est une paire de $V$-CW-complexes finis alors $\mathrm{H}_{V}^{*}(X,Y)$ est de type fini comme $\mathrm{H}^{*}V$-module~; la sous-catégorie pleine de $V$-$\mathcal{U}$ dont les objets sont les  $\mathrm{H}^{*}V$-$\mathrm{A}$-modules instables qui sont de type fini comme  $\mathrm{H}^{*}V$-modules est notée $V_{\mathrm{tf}}$-$\mathcal{U}$.

\vspace{0.625cm}
\textsc{Le complexe algébrique}

\medskip
Nous décrivons  maintenant un complexe de cochaînes uniquement défini en fonction du $\mathrm{H}^{*}V$-$\mathrm{A}$-module instable $\mathrm{H}_{V}^{*}X$.

\medskip
Le complexe $\mathrm{C}_{\mathrm{top}}^{\bullet}\hspace{1pt}X$ n'est rien d'autre que le terme $\mathrm{E}_{1}$ de la suite spectrale  convergeant vers $\mathrm{H}_{V}^{*}X$ définie par la filtration $V$-équivariante de $X$ introduite plus haut. La filtration décroissante de $\mathrm{H}_{V}^{*}X$ associée, disons $\mathrm{H}_{V}^{*}X=\mathrm{F}^{0}_{(1)}\supset\mathrm{F}^{1}_{(1)}\supset\ldots\supset\mathrm{F}^{n}_{(1)}\supset\mathrm{F}^{n+1}_{(1)}=0$, est définie par
$$
\mathrm{F}^{p}_{(1)}\mathrm{H}_{V}^{*}X
\hspace{4pt}:=\hspace{4pt}
\ker\hspace{1pt}(\mathrm{H}_{V}^{*}X\to\mathrm{H}_{V}^{*}\mathrm{F}_{p-1}\hspace{1pt}X)
$$
pour $0\leq p\leq n+1$. On peut définir une seconde filtration de $\mathrm{H}_{V}^{*}X$, variante de la précédente,  en posant
$$
\hspace{24pt}
\mathrm{F}^{p}_{(2)}\mathrm{H}_{V}^{*}X
\hspace{4pt}:=\hspace{4pt}
\bigcap_{\mathop{\mathrm{codim}}W<p}
\ker\hspace{1pt}(\mathrm{H}_{V}^{*}X\to\mathrm{H}_{V}^{*}X^{W})
\hspace{24pt};
$$
on observera que l'on a $\mathrm{F}^{p}_{(1)}\mathrm{H}_{V}^{*}X\subset\mathrm{F}^{p}_{(2)}\mathrm{H}_{V}^{*}X$ pour tout $p$.

\pagebreak

\medskip
Supposons à présent que $X$ est un $V$-CW-complexe fini. Dans ce cas cette seconde filtration peut être définie ``algébriquement'' en termes de la structure de $\mathrm{H}^{*}V$-$\mathrm{A}$-module instable de $\mathrm{H}_{V}^{*}X$. Précisons un peu. Soient $W\subset V$ un sous-groupe et $\mathrm{Fix}_{(V,W)}:V\text{-}\mathcal{U}\to V/W\text{-}\mathcal{U}$ ``le'' fonteur adjoint à gauche du foncteur ``extension des scalaires'' $V/W\text{-}\mathcal{U}\to V\text{-}\mathcal{U},N\mapsto\mathrm{H}^{*}V\otimes_{\mathrm{H}^{*}V/W}N$ (pour un ``primer'' sur les foncteurs $\mathrm{Fix}$ voir la section 1)~; l'homomorphisme $\mathrm{Fix}_{(V,W)}\mathrm{H}_{V}^{*}X\to\mathrm{H}_{V/W}^{*}X^{\mathrm{W}}$ adjoint de l'homomorphisme de $\mathrm{H}^{*}V$-$\mathrm{A}$-modules instables $\mathrm{H}_{V}^{*}X\to\mathrm{H}_{V}^{*}X^{\mathrm{W}}=\mathrm{H}^{*}V\otimes_{\mathrm{H}^{*}V/W}\mathrm{H}_{V/W}^{*}X^{\mathrm{W}}$ induit par l'inclusion $i:X^{W}\to X$ est un isomorphisme et l'unité d'adjonction 
$$
\eta_{(V,W)}:\mathrm{H}_{V}^{*}X\to\mathrm{H}^{*}V\otimes_{\mathrm{H}^{*}V/W}\mathrm{Fix}_{(V,W)}\mathrm{H}_{V}^{*}X
$$
s'identifie à $\mathrm{H}_{V}^{*}i$ (voir \ref{cohFix}). On est donc amené à introduire, pour tout $\mathrm{H}^{*}V$-$\mathrm{A}$-module instable $M$, une filtration décroissante (par des sous-$\mathrm{H}^{*}V$-$\mathrm{A}$-modules instables)
$$
M=\mathrm{F}^{0}M\supset\mathrm{F}^{1}M\supset\ldots\supset\mathrm{F}^{n}M\supset\mathrm{F}^{n+1}M=0
$$
en posant
$$
\mathrm{F}^{p}M
\hspace{4pt}:=\hspace{4pt}
\bigcap_{\mathop{\mathrm{codim}}W<p}
\ker\hspace{1pt}(\eta_{(V,W)}:M\to\mathrm{H}^{*}V\otimes_{\mathrm{H}^{*}V/W}\mathrm{Fix}_{(V,W)}M)
$$
telle que l'on a $\mathrm{F}^{p}_{(2)}\mathrm{H}_{V}^{*}X=\mathrm{F}^{p}\mathrm{H}_{V}^{*}X$. On constate que $\mathrm{F}^{1}M$ et $\mathrm{F}^{n}M$ peuvent être définis en termes de la structure de $\mathrm{H}^{*}V$-module de $M$~:  $\mathrm{F}^{1}M$ est la $\mathrm{H}^{*}V$-torsion et $\mathrm{F}^{n}M$ est constitué des éléments annulés par une puissance de l'idéal d'augmentation $\widetilde{\mathrm{H}}^{*}V$. Nous montrerons que toute la filtration  $(\mathrm{F}^{p}M)_{0\leq p\leq n+1}$ peut être en fait définie en termes de la structure de $\mathrm{H}^{*}V$-module de $M$ : c'est la ``filtration par la codimension du support''. On observera que si $M$ est de type fini comme $\mathrm{H}^{*}V$-module alors $\mathrm{F}^{n}M$ est le plus grand sous-$\mathrm{H}^{*}V$-module fini de $M$ ; nous le noterons $\mathrm{Pf}_{V}M$ ($\mathrm{Pf}$ pour ``partie finie'').

\medskip
Disposant de la filtration décrite ci-dessus, on définit un complexe de cochaines $\mathrm{C}^{\bullet}M$ de la façon suivante. Soit $M\to I^{\bullet}$ une résolution de $M$ dans la catégorie $V\text{-}\mathcal{U}$~; on pose $\mathrm{C}^{p}M=\mathrm{H}^{p}(\mathrm{F}^{p}I^{\bullet}/\mathrm{F}^{p+1}I^{\bullet})$ et on prend pour cobord le connectant évident. \`{A} nouveau $\mathrm{C}^{\bullet}M$ est muni d'une coaugmentation naturelle $M:=\mathrm{C}^{-1}M\to\mathrm{C}^{0}M$ ; le complexe coaugmenté associé est noté $\widetilde{\mathrm{C}}^{\bullet}M$. Supposons $M$ de type fini comme $\mathrm{H}^{*}V$-module~; on constate que $\mathrm{C}^{p}M$ vérifie alors une formule analogue à celle que nous avons donnée plus haut pour $\mathrm{C}_{\mathrm{top}}^{p}\hspace{1pt}X$~:
$$
\hspace{24pt}
\mathrm{C}^{p}M
\hspace{4pt}=\hspace{4pt}
\bigoplus_{\mathop{\mathrm{codim}}W=p}\mathrm{H}^{*}V
\otimes_{\mathrm{H}^{*}V/W}
\mathrm{R}^{p}\mathrm{Pf}_{V/W}\hspace{1pt}(\mathrm{Fix}_{(V,W)}M)
\hspace{24pt}.
$$
Expliquons la notation. Si $M$ est de type fini comme $\mathrm{H}^{*}V$-module alors il en est de même pour $\mathrm{Fix}_{(V,W)}M$ comme $\mathrm{H}^{*}V/W$-module ; $\mathrm{R}^{p}\mathrm{Pf}_{V/W}$ désigne le $p$-ième foncteur dérivé à droite de l'endofoncteur $\mathrm{Pf}_{V/W}$ de $V/W_{\mathrm{tf}}\text{-}\mathcal{U}$.

\medskip
Notre deuxième résultat est le suivant~:

\begin{theo}\label{introCalg}
Soient $V$ un $2$-groupe abélien élémentaire et $M$ un $\mathrm{H}^{*}V$-$\mathrm{A}$-module instable qui est de type fini comme $\mathrm{H}^{*}V$-module. Les deux propriétés suivantes sont équivalentes~:
\begin{itemize}
\item[(i)] Le $\mathrm{H}^{*}V$-module $M$ est libre.
\item[(ii)] Le complexe coaugmenté $\widetilde{\mathrm{C}}^{\bullet}M$ est acyclique.
\end{itemize}
\end{theo}

Posons $\widetilde{\mathrm{C}}_{\mathrm{alg}}^{\bullet}\hspace{1pt}X:=\widetilde{\mathrm{C}}^{\bullet}\mathrm{H}^{*}_{V}X$~; les deux complexes $\widetilde{\mathrm{C}}_{\mathrm{top}}^{\bullet}\hspace{1pt}X$ et $\widetilde{\mathrm{C}}_{\mathrm{alg}}^{\bullet}\hspace{1pt}X$ sont reliés~:

\begin{pro}\label{introCalgCtop}
Soient $V$ un $2$-groupe abélien élémentaire et $X$ un $V$-CW-complexe fini. 

\smallskip
{\em (a)} Il existe un homomorphisme de complexes de cochaînes coaugmentés
$$
\varkappa:
\hspace{4pt}\widetilde{\mathrm{C}}_{\mathrm{alg}}^{\bullet}\hspace{1pt}X
\longrightarrow
\widetilde{\mathrm{C}}_{\mathrm{top}}^{\bullet}\hspace{1pt}X
$$
tel que $\varkappa^{\hspace{1pt}p}$ est un homomorphisme de $\mathrm{H}^{*}V$-$\mathrm{A}$-modules pour $-1\leq p\leq n$ (et~l'identité pour $p=-1$).

\smallskip
{\em (b)} Si le $\mathrm{H}^{*}V$-module $\mathrm{H}^{*}_{V}X$ est libre alors $\varkappa$ est un isomorphisme.
\end{pro}

\medskip
Le lecteur notera que l'adjectif ``instable'' n'apparaît pas à la fin de l'énoncé du point (a) ci-dessus~; en effet le $\mathrm{A}$-module $\mathrm{C}_{\mathrm{top}}^{p}\hspace{1pt}X$ n'est pas \textit{a priori} instable pour $p>0$. Cependant le point (b) implique~:

\begin{scho}\label{introInstable}
Soient $V$ un $2$-groupe abélien élémentaire et $X$ un $V$-CW-complexe fini.  Si $\mathrm{H}^{*}_{V}X$ est libre comme $\mathrm{H}^{*}V$-module alors le $\mathrm{A}$-module $\mathrm{C}_{\mathrm{top}}^{p}\hspace{1pt}X$ est instable pour tout $p$.
\end{scho}

\vspace{0.625cm}
\textsc{Illustrations}

\medskip
Nous illustrons les énoncés \ref{introCtop}, \ref{introCalg} et \ref{introCalgCtop} en considérant des représentations linéaires réelles d'un $2$-groupe abélien élémentaire $V\simeq(\mathbb{Z}/2)^{n}$.

\medskip
Soit $E$ un $\mathbb{R}$-espace vectoriel de dimension finie $m$ muni d'une action linéaire de $V$. Soient $f$ la dimension du sous-espace invariant $E^{V}$ et $\mathrm{w}_{m-f}(E)$ la $(m-f)$-ième classe de Stiefel-Whitney de $E$ ($\mathrm{w}_{m-f}(E)$ appartient à $\mathrm{H}^{m-f}\hspace{1pt}V$ et est un produit d'éléments de  $\mathrm{H}^{1}V-\{0\}$).

\smallskip
L'isomorphisme de Thom fournit un isomorphisme de $\mathrm{H}^{*}V$-$\mathrm{A}$-modules instables $\mathrm{H}^{*}_{V}(E,E-\{0\})
\cong\Sigma^{f}\hspace{1pt}\mathrm{w}_{m-f}(E)\hspace{1pt}\mathrm{H}^{*}V$~; $\mathrm{H}^{*}_{V}(E,E-\{0\})$ est donc un exemple de $\mathrm{H}^{*}V$-$\mathrm{A}$-module instable qui est libre de dimension de dimension~$1$ comme $\mathrm{H}^{*}V$-module. En fait, un résultat de J.-P. Serre \cite[\S 2, Corollaire]{Ser} (voir \ref{Serre4}) implique que la correspondance $E\mapsto\mathrm{H}^{*}_{V}(E,E-\{0\})$
induit une bijection entre classes d'isomorphisme de représentations linéaires réelles de~$V$ et classes d'isomorphisme de $\mathrm{H}^{*}V$-$\mathrm{A}$-modules instables qui sont libres de dimension $1$ comme $\mathrm{H}^{*}V$-modules (voir \ref{libredim1bis}).

\medskip
Supposons que $E$ est un espace euclidien et que $V$ agit par isométries~; on~a d'après ce qui précède $\mathrm{H}^{*}_{V}(\mathrm{D}(E),\mathrm{S}(E))\cong\Sigma^{f}\hspace{1pt}\mathrm{w}_{m-f}(E)\hspace{1pt}\mathrm{H}^{*}V$ (les notations $\mathrm{D}(-)$ et $\mathrm{S}(-)$ désignent respectivement la boule et la sphère unité d'un espace euclidien), soit encore $\mathrm{H}^{*}_{V}\mathrm{S}(E\oplus\mathbb{R}_{0})\cong\Sigma^{f}\hspace{1pt}\mathrm{w}_{m-f}(E)\hspace{1pt}\mathrm{H}^{*}V\oplus\mathrm{H}^{*}V$ ($\mathbb{R}_{0}$ désigne l'espace euclidien $\mathbb{R}$ muni de l'action triviale de $V$).

\medskip
Le théorème \ref{introCtop} (en prenant $X=\mathrm{S}(E\oplus\mathbb{R}_{0}))$, le théorème \ref{introCalg} et la proposition~\ref{introCalgCtop} conduisent par exemple à l'énoncé suivant~:

\medskip
\begin{pro}\label{representation} Soit  $E_{\mathrm{r\acute{e}g}}$ le plus grand ouvert de $E$ sur lequel l'action de~$V$ est libre et $V\backslash E_{\mathrm{r\acute{e}g}}$ le quotient de cette action ($V\backslash E_{\mathrm{r\acute{e}g}}$ est une vari\'et\'e de classe $\mathrm{C}^{\infty}$).

\medskip
On a un isomorphisme canonique de $\mathrm{H}^{*}V$-$\mathrm{A}$-modules instables
$$
\hspace{24pt}
\mathrm{H}^{*}_{\mathrm{c}}(V\backslash E_{\mathrm{r\acute{e}g}})
\hspace{4pt}\cong\hspace{4pt}
\Sigma^{f+n}\hspace{2pt}\mathrm{R}^{n}\mathrm{Pf}_{V}\hspace{1pt}(\hspace{1pt}\mathrm{w}_{m-f}(E)\hspace{1pt}\mathrm{H}^{*}V\hspace{1pt})
\hspace{24pt}.
$$
	
\smallskip
{\em (La notation $\mathrm{H}^{*}_{\mathrm{c}}$ désigne ici  la cohomologie modulo $2$ à support compact.)}
\end{pro}

L'ouvert $E_{\mathrm{r\acute{e}g}}$ est le complémentaire d'un ensemble d'hyperplans et l'étude de la cohomologie des espaces du type $V\backslash E_{\mathrm{r\acute{e}g}}$ a fait l'objet de nombreuses recherches~;  Nguyen Dang Ho Hai nous a signalé à ce sujet les références \cite{DJ} et \cite{BP}.

\medskip
Nous explorons plus en détails le cas $E=\widetilde{\mathrm{R}}[V]^{\oplus\hspace{1pt}h}$, en clair le cas où $E$ est  la somme directe de $h$ copies de la représentation régulière réelle réduite de~$V$~; on a alors $E^{V}=0$, $\dim E=h\hspace{1pt}(2^{n}-1)$ et $\mathrm{w}_{h\hspace{1pt}(2^{n}-1)}=\mathrm{c}_{V}^{h}$, ($\mathrm{c}_{V}$~désignant le produit des éléments de $\mathrm{H}^{1}V-\{0\}$). On a~:
$$
\hspace{24pt}
\widetilde{\mathrm{C}}^{\bullet}_{\mathrm{top}}\hspace{1pt}\mathrm{S}(E\oplus\mathbb{R}_{0})\hspace{2pt}\cong\hspace{2pt}\hspace{2pt}\widetilde{\mathrm{C}}^{\bullet}_{\mathrm{alg}}\hspace{1pt}\mathrm{S}(E\oplus\mathbb{R}_{0})\hspace{2pt}\cong\hspace{2pt}\widetilde{\mathrm{C}}^{\bullet}(\mathrm{c}_{V}^{h}\mathrm{H}^{*}V)\oplus\widetilde{\mathrm{C}}^{\bullet}(\mathrm{H}^{*}V)
\hspace{24pt};
$$
les complexes ci-dessus (qui sont des complexes de cochaînes dans la catégorie $V_{\mathrm{tf}}$-$\mathcal{U}$) sont acycliques d'après \ref{introCtop} ou \ref{introCalg}, le premier isomorphisme est donné par \ref{introCalgCtop}. La spécificité de ces exemples tient à ce que les deux complexes de gauche  sont naturellement munis d'une action (à droite) du groupe $\mathrm{GL}(V)$, que l'isomorphisme $\varkappa$ de \ref{introCalgCtop} est équivariant et  que l'action respecte la décom\-position en somme directe de droite. Posons
$$
\hspace{12pt}
\mathrm{M}(V,h):=\mathrm{R}^{n}\mathrm{Pf}_{V}\hspace{1pt}(\hspace{1pt}\mathrm{c}_{V}^{h}\mathrm{H}^{*}V\hspace{1pt})\hspace{8pt}(\cong\mathrm{C}^{n}(\mathrm{c}_{V}^{h}\mathrm{H}^{*}V)\cong\Sigma^{-n}\mathrm{H}^{*}_{\mathrm{c}}(V\backslash \widetilde{\mathbb{R}}[V]^{\oplus\hspace{1pt}h}_{\mathrm{r\acute{e}g}}))
\hspace{12pt};
$$
nous montrons que $(\mathrm{M}(V,h))^{0}$ (le sous-espace des éléments de degré $0$ de $\mathrm{M}(V,h)$) est isomorphe comme $\mathbb{F}_{2}[\mathrm{GL}(V)]$-module à droite à la duale $\mathrm{St}_{V}^{*}$ de la représentation de Steinberg modulo $2$ de $\mathrm{GL}(V)$ (voir \ref{rapSt}) et que $\mathrm{M}(V,h)$ est engendré comme $\mathrm{H}^{*}V$-module par $(\mathrm{M}(V,h))^{0}$.

\medskip
Comme $\mathrm{St}_{V}^{*}$ est un objet projectif de la catégorie des $\mathbb{F}_{2}[\mathrm{GL}(V)]$-modules à droite,  le foncteur $\mathrm{Hom}_{(\mathbb{F}_{2}[\mathrm{GL}(V)])^{\mathrm{op}}}(\mathrm{St}_{V}^{*},-)=:\mathrm{e}_{V}^{\mathrm{St}}(-)$ est exact, si bien que l'on peut obtenir de  nouveaux complexes acycliques (dans la catégorie $\mathcal{U}$), en appliquant ce foncteur. C'est ce que vient de faire Nguyen Dang Ho Hai (NDHH) dans le cas $h=2$~; nous donnons ci-après (en petits caractères) quelques informations sur ses travaux.

\medskip
\footnotesize
Posons pour alléger $\mathrm{e}_{n}^{\mathrm{St}}:=\mathrm{e}_{(\mathbb{Z}/2)^{n}}^{\mathrm{St}}$ et $\omega_{n}:=\mathrm{c}_{(\mathbb{Z}/2)^{n}}$~; on a $\mathrm{e}_{n}^{\mathrm{St}}(\omega_{n}^{2}\hspace{1pt}\mathrm{H}^{*}(\mathbb{Z}/2)^{n})=\omega_{n}\mathrm{L}_{n}$, $\mathrm{L}_{n}$~désignant le sous-$\mathrm{A}$-module de $\mathrm{H}^{*}(\mathbb{Z}/2)^{n}$ défini par $\mathrm{L}_{n}:=\mathrm{e}_{n}^{\mathrm{St}}(\omega_{n}\hspace{1pt}\mathrm{H}^{*}(\mathbb{Z}/2)^{n}$) \cite{MP} (on~rappelle que l'on a $\mathrm{e}_{n}^{\mathrm{St}}(\mathrm{H}^{*}(\mathbb{Z}/2)^{n})=\mathrm{L}_{n}\oplus\mathrm{L}_{n-1}$). Le complexe $\mathrm{e}_{n}^{\mathrm{St}}(\widetilde{\mathrm{C}}^{\bullet}(\omega_{n}^{2}\hspace{1pt}\mathrm{H}^{*}(\mathbb{Z}/2)^{n})$ peut donc être vu comme une résolution, dans la catégorie $\mathcal{U}$, de $\omega_{n}\mathrm{L}_{n}$. NDHH montre que ce complexe est une résolution injective, somme directe de la résolution injective minimale
$$
0\to\omega_{n}\mathrm{L}_{n}\to\mathrm{L}_{n}\to\ldots\to\mathrm{L}_{n-k}\otimes\mathrm{J}(2^{p}-1)\to\ldots\to\ldots\to\mathrm{J}(2^{n}-1)
$$
de \cite{NST} et de ``décalés''  de complexes  de la forme $I\overset{\mathrm{id}}{\to}I$, avec $I$ injectif~; on rappelle que la notation $\mathrm{J}(k)$ désigne le $k$-ième $\mathcal{U}$-injectif de Brown-Gitler (voir le début de la section~2, prendre $V=0$).

\smallskip
Par ailleurs NDHH montre à l'aide d'un résultat d'Inoue \cite{In}, mais sans utiliser les résultats de ce mémoire, que l'on a un isomorphisme de $\mathrm{A}$-modules instables
$$
\hspace{24pt}
\mathrm{e}_{V}^{\mathrm{St}}\hspace{2pt}
(\Sigma^{-n}\hspace{1pt}\mathrm{H}^{*}_{\mathrm{c}}
(V\backslash\widetilde{\mathbb{R}}[V]^{\oplus 2}_{\mathrm{r\acute{e}g}}))
\hspace{4pt}\cong\hspace{4pt}\bigoplus_{0\leq i\leq n} \mathrm{J}(2^{i}-1)
\hspace{24pt}.
$$
Ceci redonne le résultat  de \cite{NST} évoqué ci-dessus.

\normalsize
\vspace{0.625cm}
\textsc{Sur le théorème 10.2 de \cite{AFP2}}

\bigskip
Avant d'énoncer une version édulcorée du théorème en question, nous rappelons la notion de $j$-syzygie (dans le cas particulier qui nous intéresse)~:

\bigskip
\begin{defi}\label{intro-syzygie}
Soient $M$ un $\mathrm{H}^{*}V$-module $\mathbb{N}$-gradué et $j\geq 0$ un entier. On dit que $M$ est une  {\em $\mathrm{H}^{*}V$-$j$-syzygie} s'il existe une suite exacte de $\mathrm{H}^{*}V$-modules $\mathbb{N}$-gradués
$$
0\to M\to L^{0}\to L^{1}\to\ldots\to L^{j-1}
$$
avec $L^{0},L^{1},\ldots,L^{j-1}$ libres (on convient que tout $M$ est une $\mathrm{H}^{*}V$-$0$-syzygie).
\end{defi}
\pagebreak

\medskip
\begin{theo}\label{genafp} {\em (Allday, Franz et Puppe \cite{AFP2})} Soient $V$ un $2$-groupe abélien élémentaire, $X$ un $V$-CW-complexe fini et $j\geq 0$ un entier. Les deux propriétés suivantes sont équivalentes~:
\begin{itemize}
\item[(i)] le $\mathrm{H}^{*}V$-module sous-jacent à $\mathrm{H}_{V}^{*}X$ est une $\mathrm{H}^{*}V$-$j$-syzygie~;
\item[(ii)] on a $\mathrm{H}^{p}\hspace{1pt}\widetilde{\mathrm{C}}_{\mathrm{top}}^{\bullet}X=0$ pour $p\leq j-2$.
\end{itemize}
\end{theo}

\medskip
Ce théorème est bien une généralisation du théorème \ref{introCtop}~: 

\medskip
Pour $j=n$ la propriété (i) ci-dessus coïncide avec la propriété (i) de \ref{introCtop}~; la propriété (ii) ci-dessus implique quant à elle la propriété (ii) de \ref{introCtop} grâce à un lemme facile \cite[Lemma 5.6]{AFP1} (lemme \ref{top-n-2} dans le cas qui nous occupe).

\bigskip
Allday, Franz et Puppe obtiennent  le théorème \ref{genafp} à l'aide de résultats d'algèbre commutative concernant la catégorie des $\mathrm{H}^{*}V$-modules $\mathbb{N}$-gradués. L'article \cite{AFP2} fait d'ailleurs suite à l'article \cite{AFP1} où le rôle de l'algèbre de polynômes $\mathbb{F}_{2}[U_{1},U_{2},\ldots,U_{n}]$ ($\simeq\mathrm{H}^{*}V$) était tenu par l'algèbre de polynômes $\mathbb{Q}[T_{1},T_{2},\ldots,T_{n}]$ ($\simeq\mathrm{H}^{*}(\mathrm{B}T;\mathbb{Q})$, $T$ tore de rang $n$).

\bigskip
Nous démontrons une version du théorème \ref{genafp} (Théorème \ref{introgentop} ci-dessous) en utilisant, comme pour le théorème \ref{introCtop}, la théorie des $\mathrm{H}^{*}V$-$\mathrm{A}$-modules instables (sous-produit des recherches sur la conjecture de Sullivan) et plus précisément la théorie des $\mathrm{H}^{*}V$-$\mathrm{A}$-modules instables qui sont de type fini comme $\mathrm{H}^{*}V$-modules. En fait la structure des $\mathrm{H}^{*}V$-$\mathrm{A}$-modules instables qui sont de type fini comme $\mathrm{H}^{*}V$-modules est beaucoup plus accessible que celle des $\mathrm{H}^{*}V$-modules $\mathbb{N}$-gradués de type fini généraux (voir la section 9 et en particulier la sous-section \ref{Serre1}).

\bigskip
Avant d'énoncer le théorème \ref{introgentop} nous devons introduire une définition~:

\bigskip
\begin{defi}\label{intro-syzygie-bis}
Soient $M$ un $\mathrm{H}^{*}V$-$\mathrm{A}$-module instable de type fini comme $\mathrm{H}^{*}V$-module et $j\geq 0$ un entier. Nous dirons que $M$ est une  {\em $(V_{\mathrm{tf}}\text{-}\mathcal{U})$-$j$-syzygie} s'il existe une suite exacte dans la catégorie $V_{\mathrm{tf}}\text{-}\mathcal{U}$
$$
0\to M\to L^{0}\to L^{1}\to\ldots\to L^{j-1}
$$
avec $L^{0},L^{1},\ldots,L^{j-1}$ libres comme $\mathrm{H}^{*}V$-modules (nous convenons que tout~$M$ est une $(V_{\mathrm{tf}}\text{-}\mathcal{U})$-$0$-syzygie).
\end{defi}

\pagebreak

\medskip
\begin{theo}\label{introgentop} Soient $V$ un $2$-groupe abélien élémentaire, $X$ un $V$-CW-complexe fini et $j\geq 0$ un entier. Les trois propriétés suivantes sont équi\-valentes~:
\begin{itemize}
\item[(i)] $\mathrm{H}_{V}^{*}X$ est une $(V_{\mathrm{tf}}\text{-}\mathcal{U})$-$j$-syzygie~;
\item[(ii)] le $\mathrm{H}^{*}V$-module sous-jacent à $\mathrm{H}_{V}^{*}X$ est une $\mathrm{H}^{*}V$-$j$-syzygie~;
\item[(iii)] on a $\mathrm{H}^{p}\hspace{1pt}\widetilde{\mathrm{C}}_{\mathrm{top}}^{\bullet}X=0$ pour $p\leq j-2$.
\end{itemize}
\end{theo}

Ce théorème est essentiellement conséquence de sa version algébrique~:

\begin{theo}\label{introgenalg} Soient $M$ un $\mathrm{H}^{*}V$-$\mathrm{A}$-module instable qui est de type fini comme $\mathrm{H}^{*}V$-module et $j\geq 0$ un entier. Les trois propriétés suivantes sont équivalentes~:
\begin{itemize}
\item[(i)] $M$ est une $(V_{\mathrm{tf}}\text{-}\mathcal{U})$-$j$-syzygie~;
\item[(ii)] le $\mathrm{H}^{*}V$-module sous-jacent à $M$ est une $\mathrm{H}^{*}V$-$j$-syzygie~;
\item[(iii)] on a $\mathrm{H}^{p}\hspace{1pt}\widetilde{\mathrm{C}}^{\bullet}M=0$ pour $p\leq j-2$.
\end{itemize}
\end{theo}

et de la généralisation suivante du point (b) de la proposition \ref{introCalgCtop}~:

\begin{pro}\label{introgenCalgCtop} Soient $V$ un $2$-groupe abélien élémentaire, $X$ un $V$-CW-complexe fini et $j\geq 0$ un entier. Si le $\mathrm{H}^{*}V$-module sous-jacent à $\mathrm{H}_{V}^{*}X$ est une $\mathrm{H}^{*}V$-$j$-syzygie, alors l'homomorphisme de $\mathrm{H}^{*}V$-$\mathrm{A}$-modules 
$$
\varkappa^{p}\hspace{2pt}:
\hspace{2pt}\widetilde{\mathrm{C}}_{\mathrm{alg}}^{p}\hspace{1pt}X
\longrightarrow
\widetilde{\mathrm{C}}_{\mathrm{top}}^{p}\hspace{1pt}X
$$
est un isomorphisme pour $p\leq j$.
\end{pro}

\bigskip
L'implication (ii)$\Rightarrow$(i) du théorème 0.9 est assez piquante. Nous montrons en fait un énoncé plus précis~:

\begin{pro}\label{harissa} Soient $M$ un $\mathrm{H}^{*}V$-$\mathrm{A}$-module instable qui est de type fini comme $\mathrm{H}^{*}V$-module et $j\geq 1$ un entier.

\smallskip
Il existe un complexe de cochaînes coaugmenté dans la catégorie $V_{\mathrm{tf}}\text{-}\mathcal{U}$
$$
\hspace{24pt}
\widetilde{\mathrm{L}}_{j}^{\bullet}M
\hspace{4pt}=\hspace{4pt}
(\hspace{2pt}M\to\mathrm{L}_{j}^{0}M\to\mathrm{L}_{j}^{1}M\to\ldots\to\mathrm{L}_{j}^{j-1}M\hspace{2pt})
\hspace{24pt},
$$
dépendant fonctoriellement de $M$, qui vérifie la propriété suivante~:

\smallskip
Si le $\mathrm{H}^{*}V$-module sous-jacent à $M$ est une $\mathrm{H}^{*}V$-$j$-syzygie alors $\widetilde{\mathrm{L}}_{j}^{\bullet}M$ est acyclique et les $\mathrm{H}^{*}V$-modules sous-jacents aux $\mathrm{L}_{j}^{k}M$, $0\leq k\leq j-1$, sont libres.
\end{pro}

(Pour une explicitation du complexe $\widetilde{\mathrm{L}}_{j}^{\bullet}M$, voir la proposition \ref{AFPfonctoriel}, la remarque \ref{j=1} et la sous-section \ref{Bpq}.)

\bigskip
Nous obtenons la proposition ci-dessus par des techniques d'algèbre homologique. Précisons un peu. Soit $\mathcal{W}$ l'ensemble des sous-espaces $W$ de $V$~; $\mathcal{W}$ est ordonné par inclusion et peut donc être vu comme une catégorie. Nous considérons la catégorie des foncteurs définis sur $\mathcal{W}$ et à valeurs dans la catégorie $\mathcal{E}$ des $\mathbb{F}_{2}$-espaces vectoriels. L'outil principal qui nous permet d'arriver à un énoncé du type \ref{harissa} est l'explicitation d'une résolution injective, dans la catégorie abélienne $\mathcal{E}^{\mathcal{W}}$, d'un objet arbitraire (sous-section \ref{resinj}).

\medskip
On observera incidemment (Remarque \ref{Bpq-5}) que les techniques évoquées ci-dessus apportent un nouvel éclairage sur le complexe $\mathrm{C}^{\bullet}M$.

\vspace{0.625cm}
Voici, en conclusion de cette introduction, le sommaire du mémoire~:

\medskip
\textbf{1.} Rappels sur la théorie des foncteurs $\mathrm{Fix}$\dotfill 11

\medskip
\textbf{2.} Rappels sur les injectifs des catégories $V\text{-}\mathcal{U}$ et  $V_{\mathrm{tf}}\text{-}\mathcal{U}$\dotfill 26

\medskip
\textbf{3.} Sur les dérivés du foncteur ``partie finie'' \dotfill 33

\medskip
\textbf{4.} Le complexe topologique\dotfill 45

\medskip
\textbf{5.} Le complexe algébrique\dotfill 56

\medskip
\textbf{6.} Comparaison entre les complexes algébrique et topologique\dotfill 75

\medskip
\textbf{7.} Illustrations\dotfill 89

\medskip
\textbf{8.} Modules de Steinberg et algèbre homologique\dotfill 118

\medskip
\textbf{9.} Filtration par la codimension du support et foncteurs $\mathrm{Fix}$\dotfill 149

\medskip
\textbf{Références}\dotfill 160

\vspace{1,25cm}
\textit{Remerciements}

\medskip
Dorra Bourguiba, Jean Lannes et Saïd Zarati tiennent à remercier la Société Mathématique de Tunisie dont les congrès annuels ont contribué à la maturation de ce travail.

\medskip
Jean Lannes et Saïd Zarati ont bénéficié du programme Recherches en Binôme du CIRM.

\medskip
Jean Lannes a bénéficié de l'hospitalité du CMLS.

\medskip
Lionel Schwartz a bénéficié du soutien du VIASM (Hanoï) lors de plusieurs visites ces dernières années, en particulier à l'automne 2019, visites durant lesquelles il a travaillé sur le contenu de ce mémoire.

\pagebreak

\sect{Rappels sur la théorie des foncteurs $\mathrm{Fix}$}

Soient $\ell$ un nombre premier (il n'y a aucune raison dans cette section de se limiter au cas $\ell=2$), $V$ un $\ell$--groupe abélien élémentaire et $W\subset V$ un sous-groupe~;  l'homomorphisme canonique, $q:V\to V/W$, induit un homomorphisme (injectif) de $\mathrm{A}$-algèbres instables $q^{*}:\mathrm{H}^{*}V/W\to\mathrm{H}^{*}V$. La notation $\mathrm{A}$ désigne ici l'algèbre de Steenrod modulo $\ell$, la cohomologie modulo $\ell$ d'un espace est le type même d'une $\mathrm{A}$-algèbre instable, pour une définition explicite de la catégorie $\mathcal{U}$ (resp. $\mathcal{K}$) des $\mathrm{A}$-modules instables (resp. $\mathrm{A}$-algèbres instables) voir par exemple \cite[1.7.1 (resp. 1.7.2)]{LaT}~;  rappelons que nous avons convenu d'abréger la notation $\mathrm{H}^{*}( -;\mathbb{F}_{\ell})$ en $\mathrm{H}^{*}$. On note
$$
\mathrm{Fix}_{(V,W)}:V\text{-}\mathcal{U}\to V/W\text{-}\mathcal{U}
$$
``le'' foncteur adjoint à gauche du foncteur ``extension des scalaires'' 
$$
V/W\text{-}\mathcal{U}\to V\text{-}\mathcal{U}
\hspace{12pt},\hspace{12pt}
N\mapsto\mathrm{H}^{*}V\otimes_{\mathrm{H}^{*}V/W}N
$$
(pour des commentaires concernant les guillemets qui entourent l'article défini voir \ref{remadj} et \ref{uniadj}, la notation $V\text{-}\mathcal{U}$ est la généralisation évidente de celle introduite pour $\ell=2$).

\medskip
Les foncteurs $\mathrm{Fix}_{(V,W)}$ sont introduits et étudiés dans \cite{LZtorext} (ce sont des avatars des foncteurs $\mathrm{Fix}_{V',V''}$, définis pour $V=V'\oplus V'' $, introduits et étudiés dans \cite{LZsmith}). Le foncteur $\mathrm{Fix}_{(V,V)}:V\text{-}\mathcal{U}\to\mathcal{U}$ est aussi noté $\mathrm{Fix}_{V}$ dans \cite{LZsmith}. Ce foncteur est introduit et étudié dans \cite{LaT} (où il est d'ailleurs simplement noté $\mathrm{Fix}$). On rappelle ci-dessous l'origine de cette notation.

\medskip
Soit $X$ un espace muni d'une action de $V$. Soient $\underline{X}$ l'espace fonctionnel $\mathbf{hom}(\mathrm{E}V,X)$ et $\mathrm{i}:X\to\underline{X}$ l'application induite par l'application $\mathrm{E}V\to\mathrm{pt}$~; on observera que l'inclusion $\mathrm{i}$ est une équivalence d'homotopie. L'action de $V$, au but et à la source, munit $\mathbf{hom}(\mathrm{E}V,X)$ d'une action de $V\times V^{\mathrm{op}}$ et induit \textit{via} l'homomorphisme $V\to V\times V^{\mathrm{op}} , v\mapsto (v,v^{-1})$, une action de $V$ sur $\underline{X}$ telle que $\mathrm{i}$ est $V$-équivariante. On observera également que $\mathrm{H}^{*}_{V}\mathrm{i}:\mathrm{H}^{*}_{V}\underline{X}\to\mathrm{H}^{*}_{V}X$ est un isomorphisme. L'espace $(\underline{X})^{V}(=\mathbf{hom}_{V}(\mathrm{E}V,X))$ est appelé l'espace des points fixes homotopiques de l'action de $V$ sur $X$ et noté $X^{\mathrm{h}V}$.

\medskip
L'inclusion $X^{V}\hookrightarrow X$ induit un homomorphisme de $\mathrm{H}^{*}V$-$\mathrm{A}$-modules instables $\mathrm{H}^{*}_{V}X\to\mathrm{H}^{*}_{V}X^{V}=\mathrm{H}^{*}V\otimes\mathrm{H}^{*}X^{V}$ et par adjonction un homomorphisme de $\mathrm{H}^{*}V$-$\mathrm{A}$-modules instables $\nu_{X}:\mathrm{Fix}_{V}\mathrm{H}^{*}_{V}X\to\mathrm{H}^{*}X^{V}$. On note $\lambda_{X}$ l'homomorphisme composé
$$
\begin{CD}
\hspace{24pt}
\mathrm{Fix}_{V}\mathrm{H}^{*}_{V}X
@>{(\mathrm{Fix}_{V}(\mathrm{H}^{*}_{V}\mathrm{i}))}^{-1}>\cong>\mathrm{Fix}_{V}\mathrm{H}^{*}_{V}\underline{X}
@>\nu_{\underline{X}}>>
\mathrm{H}^{*}(\underline{X})^{V}
=
\mathrm{H}^{*}X^{\mathrm{h}V}
\hspace{24pt}.\end{CD}
$$
On montre dans  \cite[Chap. 4]{LaT} que $\lambda_{X}$ est en fait sous-jacent à un homomorphisme de $\mathrm{A}$-algèbres instables qui est ``souvent'' un isomorphisme.

\medskip
Pareillement, on définit un homomorphisme de $\mathrm{H}^{*}V/W$-$\mathrm{A}$-modules instables
$$
\nu_{X}:\mathrm{Fix}_{(V,W)}\mathrm{H}^{*}_{V}X\to\mathrm{H}^{*}_{V/W}X^{W}
$$ 
comme adjoint de l'homomorphisme $\mathrm{H}^{*}_{V}X\to\mathrm{H}^{*}_{V}X^{W}\cong\mathrm{H}^{*}V\otimes_{\mathrm{H}^{*}V/W}\mathrm{H}^{*}_{V/W}X^{W}$ induit par l'inclusion $X^{W}\hookrightarrow X$. (\textit{Mutatis mutandis} on peut définir une transformation naturelle $\lambda_{X}:\mathrm{Fix}_{(V,W)}\mathrm{H}^{*}_{V}X\to\mathrm{H}^{*}_{V/W}X^{\mathrm{h}W}$ qui vérifie des propriétés analogues à celles que l'on a pour $W=V$.)

\begin{pro}\label{cohFix}
Si $X$ est un $V$-CW-complexe fini alors l'homomorphisme
$$
\nu_{X}:\mathrm{Fix}_{(V,W)}\mathrm{H}^{*}_{V}X\to\mathrm{H}^{*}_{V/W}X^{W}
$$
est un isomorphisme.
\end{pro}

\begin{scho}\label{etaFix} Soit $X$ un $V$-CW-complexe fini. L'unité d'adjonction
$$
\mathrm{H}_{V}^{*}X\to\mathrm{H}^{*}V\otimes_{\mathrm{H}^{*}V/W}\mathrm{Fix}_{(V,W)}\mathrm{H}_{V}^{*}X
$$
s'identifie à l'homomorphisme de $\mathrm{H}^{*}V$-$\mathrm{A}$-modules instables $\mathrm{H}_{V}^{*}X\to\mathrm{H}_{V}^{*}X^{W}$.
\end{scho}

\textit{Démonstration.} ``Abstract nonsense ''. Soient $M$ un $\mathrm{H}^{*}V$-$\mathrm{A}$-module instable,  $N$ un $\mathrm{H}^{*}V/W$-$\mathrm{A}$-module instable, $f:M\to\mathrm{H}^{*}V\otimes_{\mathrm{H}^{*}V/W}N$ un homomorphisme et $\widetilde{f}:\mathrm{Fix}_{(V,W)}M\to N$ son adjoint, alors $f$ est le composé de l'unité d'adjonction $M\to\mathrm{H}^{*}V\otimes_{\mathrm{H}^{*}V/W}\mathrm{Fix}_{(V,W)}M$ et de l'homomorphisme $\mathrm{H}^{*}V\otimes_{\mathrm{H}^{*}V/W}\widetilde{f}$.
\hfill$\square$

\bigskip
La démonstration de la proposition \ref{cohFix} est renvoyée à la fin de cette section. Avant cela nous rappelons, de manière assez détaillée et relativement ``self-contained'', la théorie (algébrique)  des foncteurs $\mathrm{Fix}$. C'est en fait un sous-produit de celle des foncteurs $\mathrm{T}$ pour laquelle nous renvoyons à \cite{LaT}. La présentation que nous en donnons ci-après est légèrement différente de celle de \cite{LZsmith}, \cite{LZtorext} et \cite[ Chap. 4]{LaT} (pour le cas $W=V$)~; le lecteur est invité à la comparer à \cite[\S 2]{DWsmith2}.

\bigskip
On commence par quelques observations très simples concernant le produit fibré $V\underset{V/W}{\times}V$ de deux copies de $V$ au-dessus de $V/W$ (en clair, le sous-groupe de $V\times V$ constitué des couples $(z,t)$ avec $z\equiv t \bmod{W}$)~; on observe que l'homomorphisme canonique de  $\mathrm{A}$-algèbres instables 
$$
\mathrm{H}^{*}V\underset{\mathrm{H}^{*}V/W}{\otimes}\mathrm{H}^{*}V\to\mathrm{H}^{*}(V\underset{V/W}{\times}V)
$$
est un isomorphisme. On observe également que l'homomorphisme de groupes $\iota:W\oplus V\to V\underset{V/W}{\times}V,(x,y)\mapsto(x+y,y)$ est un isomorphisme. On a donc un isomorphisme canonique de $\mathrm{A}$-algèbres instables 
$$
\hspace{24pt}
\mathrm{H}^{*}V\underset{\mathrm{H}^{*}V/W}{\otimes}\mathrm{H}^{*}V
\hspace{4pt}\cong\hspace{4pt}
\mathrm{H}^{*}W\otimes\mathrm{H}^{*}V
\hspace{24pt}.
$$

\begin{scho}\label{observation}
Soient $\varphi:\mathrm{H}^{*}V\to\mathrm{H}^{*}W\otimes\mathrm{H}^{*}V$ l'homomorphisme de $\mathrm{A}$-algèbres instables induit par l'application linéaire $W\oplus V\to V,(x,y)\mapsto x+y$ et $N$ un $\mathrm{H}^{*}V$-$\mathrm{A}$-module instable, alors on a un isomorphisme canonique de $\mathrm{H}^{*}V$-$\mathrm{A}$-modules instables, naturel en $N$,
$$
\hspace{24pt}
\mathrm{H}^{*}V\otimes_{\mathrm{H}^{*}V/W}N
\hspace{4pt}\cong\hspace{4pt}
\mathrm{H}^{*}W\otimes N
\hspace{24pt},
$$
$\mathrm{H}^{*}W\otimes N$ étant un $\mathrm{H}^{*}V$-module \textit{via} $\varphi$.
\end{scho}

\medskip
\textit{Démonstration.} Soient $\alpha$ et $\beta$ les deux homomorphismes $V\underset{V/W}{\times}V\to V$  respectivement induits par la première et seconde projection. On constate que l'on a des isomorphismes canoniques de $\mathrm{H}^{*}V$-$\mathrm{A}$-modules instables
$$
\hspace{24pt}
\mathrm{H}^{*}V\otimes_{\mathrm{H}^{*}V/W}N
\hspace{4pt}\cong\hspace{4pt}
(\mathrm{H}^{*}V\otimes_{\mathrm{H}^{*}V/W}\mathrm{H}^{*}V)\otimes_{\mathrm{H}^{*}V}N
\hspace{24pt},
$$
$\mathrm{H}^{*}V\otimes_{\mathrm{H}^{*}V/W}\mathrm{H}^{*}V$ étant un $\mathrm{H}^{*}V$-module à gauche \textit{via} $\alpha^{*}$ et un $\mathrm{H}^{*}V$-module à droite \textit{via} $\beta^{*}$. On a donc d'après ce qui précède
$$
\hspace{24pt}
\mathrm{H}^{*}V\otimes_{\mathrm{H}^{*}V/W}N
\hspace{4pt}\cong\hspace{4pt}
(\mathrm{H}^{*}W\otimes\mathrm{H}^{*}V)\otimes_{\mathrm{H}^{*}V}N
\hspace{24pt},
$$
$\mathrm{H}^{*}W\otimes\mathrm{H}^{*}V$ étant un $\mathrm{H}^{*}V$-module à gauche \textit{via} $(\alpha\circ\iota)^{*}$, c'est-à-dire\linebreak l'homomorphisme de $\mathrm{A}$-algèbres instables induit par l'application linéaire $W\oplus V\to \nolinebreak  V,(x,y)\mapsto x+y$, et un $\mathrm{H}^{*}V$-module à droite \textit{via} $(\beta\circ\iota)^{*}$, c'est-à-dire l'homomorphisme de $\mathrm{A}$-algèbres instables induit par la projection $W\oplus V\to V$.
\hfill$\square$

\bigskip
On rappelle que l'on note $\mathrm{T}_{W}:\mathcal{U}\to\mathcal{U}$ l'adjoint à gauche du foncteur $\mathcal{U}\to\mathcal{U},N\mapsto\mathrm{H}^{*}W\otimes N$ ($\mathrm{T}_{W}$ est ``canonisé'' dans \cite{LaT}, d'où l'article défini). Soit $M$ un $\mathrm{A}$-module instable~; comme $\mathrm{T}_{W}$ ``préserve les produits tensoriels'', $\mathrm{T}_{W}M$ est naturellement un $\mathrm{T}_{W}\mathrm{H}^{*}V$-$\mathrm{A}$-module instable. Le scholie ci-dessus entraîne~:

\begin{pro}\label{adjonction}
Soient $M$ et $N$ deux $\mathrm{H}^{*}V$-$\mathrm{A}$-modules instables. On a un isomorphisme, naturel en $M$ et $N$,
$$
\hspace{12pt}
\mathrm{Hom}_{V\text{-}\mathcal{U}}(M,\mathrm{H}^{*}V\otimes_{\mathrm{H}^{*}V/W}N)
\hspace{4pt}\cong\hspace{4pt}
\mathrm{Hom}_{V\text{-}\mathcal{U}}(\mathrm{H}^{*}V\otimes_{\mathrm{T}_{W}\mathrm{H}^{*}V}\mathrm{T}_{W}M, N)
\hspace{12pt},
$$
$\mathrm{H}^{*}V$ étant dans le membre de droite un $\mathrm{T}_{W}\mathrm{H}^{*}V$-module \textit{via} l'homomorphisme de $\mathrm{A}$-algèbres instables $\widetilde{\varphi}:\mathrm{T}_{W}\mathrm{H}^{*}V\to\mathrm{H}^{*}V$ adjoint de $\varphi$. En d'autres termes le foncteur $V\text{-}\mathcal{U}\to V\text{-}\mathcal{U},M\mapsto\mathrm{H}^{*}V\otimes_{\mathrm{T}_{W}\mathrm{H}^{*}V}\mathrm{T}_{W}M$ est un adjoint à gauche du foncteur $V\text{-}\mathcal{U}\to V\text{-}\mathcal{U},N\mapsto\mathrm{H}^{*}V\otimes_{\mathrm{H}^{*}V/W}N$.

\end{pro}

\medskip
\textit{Démonstration.} On a
$$
\mathrm{Hom}_{V\text{-}\mathcal{U}}(M,\mathrm{H}^{*}V\otimes_{\mathrm{H}^{*}V/W}N)\cong\mathrm{Hom}_{V\text{-}\mathcal{U}}(M,\mathrm{H}^{*}W\otimes N)
$$
d'après \ref{observation}. Or l'inclusion $\mathrm{Hom}_{V\text{-}\mathcal{U}}(M,\mathrm{H}^{*}W\otimes N)\subset\mathrm{Hom}_{\mathcal{U}}(M,\mathrm{H}^{*}W\otimes N)$ s'identifie par adjonction à l'inclusion
$$
\hspace{24pt}
\mathrm{Hom}_{V\text{-}\mathcal{U}}(\mathrm{H}^{*}V\otimes_{\mathrm{T}_{W}\mathrm{H}^{*}V}\mathrm{T}_{W}M, N)
\hspace{4pt}\subset\hspace{4pt}
\mathrm{Hom}_{\mathcal{U}}(\mathrm{T}_{W}M, N)
\hspace{24pt}.
$$
On s'en convainc grâce à \cite[Proposition 1.3.1]{LZsmith}~: faire $K=L=\mathrm{H}^{*}V$, $E=W$ et prendre pour $\varphi:K\to\mathrm{H}^{*}E\otimes L$ l'homomorphisme de $\mathrm{A}$-algèbres instables que l'on a intentionnellement noté $\varphi$ dans \ref{observation}.
\hfill$\square$

\bigskip
L'endofoncteur de la catégorie $V\text{-}\mathcal{U}$
$$
M
\hspace{4pt}\mapsto\hspace{4pt}
\mathrm{H}^{*}V\otimes_{\mathrm{T}_{W}\mathrm{H}^{*}V}\mathrm{T}_{W}M
$$
($\mathrm{H}^{*}V$ étant un $\mathrm{T}_{W}\mathrm{H}^{*}V$-module \textit{via} l'adjoint de $\varphi$), qui apparaît dans la proposition \ref{adjonction}, jouera un rôle important dans ce mémoire~; pour alléger  nous le noterons $\mathrm{EFix}_{(V,W)}$. Cette notation entend suggérer le lien avec le foncteur $\mathrm{Fix}_{(V,W)}$ que nous précisons ci-après (\cite[Proposition 3.1]{LZtorext}). (Dans \cite{DWsmith2}, $\mathrm{EFix}_{(V,W)}M$ est noté $\mathrm{T}^{W}_{i^{*}}M$, $i^{*}:\mathrm{H}^{*}V\to\mathrm{H}^{*}W$ désignant  l'homomorphisme de $\mathrm{A}$-algèbres instables induit par l'inclusion $i$ de $W$ dans $V$.)

\begin{pro-def}\label{EFix} Soit $M$ un $\mathrm{H}^{*}V$-$\mathrm{A}$-module instable~; on a un isomorphisme de $\mathrm{H}^{*}V$-$\mathrm{A}$-modules instables, naturel en $M$
$$
\hspace{24pt}
\mathrm{EFix}_{(V,W)}M
\hspace{4pt}\cong\hspace{4pt}
\mathrm{H}^{*}V\otimes_{\mathrm{H}^{*}V/W}\mathrm{Fix}_{(V,W)}M
\hspace{24pt}.
$$
En d'autres termes on a un isomorphisme naturel, entre endofoncteurs de~$\mathrm{V}\text{-}\mathcal{U}$,
$$
\hspace{24pt}
\mathrm{EFix}_{(V,W)}
\hspace{4pt}\cong\hspace{4pt}
\mathrm{e}_{(V,W)}\circ\mathrm{Fix}_{(V,W)}
\hspace{24pt},
$$
$\mathrm{e}_{(V,W)}:V/W\text{-}\mathcal{U}\to V\text{-}\mathcal{U}$ désignant le foncteur $N\mapsto\mathrm{H}^{*}V\otimes_{\mathrm{H}^{*}V/W}N$.
\end{pro-def}

\bigskip
\textit{Démonstration.} Soit $\mathrm{E}_{(V,W)}:V\text{-}\mathcal{U}\to V\text{-}\mathcal{U}$ le foncteur $N\mapsto\mathrm{H}^{*}V\otimes_{\mathrm{H}^{*}V/W}N$~; la proposition \ref{adjonction} dit que $\mathrm{EFix}_{(V,W)}$ est adjoint à gauche de $\mathrm{E}_{(V,W)}$. On a un isomorphisme naturel, entre endofoncteurs de $\mathrm{V}\text{-}\mathcal{U}$, 
$$
\hspace{24pt}
\mathrm{E}_{(V,W)}
\hspace{4pt}\cong\hspace{4pt}
\mathrm{e}_{(V,W)}\circ\mathcal{O}_{(V,W)}
\hspace{24pt},
$$
$\mathcal{O}_{(V,W)}: V\hspace{-2pt}\text{-}\mathcal{U}\to V/W\hspace{-2pt}\text{-}\mathcal{U}$ désignant le foncteur oubli évident. Comme les foncteurs $\mathrm{Fix}_{(V,W)}: V\hspace{-2pt}\text{-}\mathcal{U}\to V/W\hspace{-2pt}\text{-}\mathcal{U}$ et $\mathrm{e}_{(V,W)}: V/W\text{-}\mathcal{U}\to V\hspace{-2pt}\text{-}\mathcal{U}$ sont respectivement adjoints à gauche des foncteurs $\mathrm{e}_{(V,W)}$ et $\mathcal{O}_{(V,W)}$,  on a un isomorphisme naturel
$\mathrm{EFix}_{(V,W)}\cong\mathrm{e}_{(V,W)}\circ\mathrm{Fix}_{(V,W)}$ (spécialisation de \cite[Chap. IV, \S8, Theorem~1]{McL}, noter l'ordre des facteurs dans la composition ci-dessus~!).
\hfill$\square$

\bigskip
\begin{rem}\label{RemEfix}
Pour $W=V$ l'isomorphisme naturel de \ref{EFix} est celui de la proposition 4.5 de \cite{LaT}. La démonstration que nous en donnons ci-dessus est plus directe que celle de \cite{LaT} qui imite la démonstration de la proposition ``topologique'' 4.2 de cette référence.
\end{rem}

\bigskip
Nous dégageons incidemment une variante de la proposition \ref{cohFix} dont la démonstration sera évoquée à la fin de cette section. Soient $X$ un $V$-espace et $W$ un sous-groupe de $V$. On note $\mathrm{a}_{W}:W\oplus V\to V$ l'homomorphisme de groupes $(x,y)\mapsto x+y$~; on fait agir $W\oplus W$ sur $X$ \textit{via}~$\mathrm{a}_{W}$. On note $\delta_{W,X}:\mathrm{H}^{*}_{V}X\to\mathrm{H}^{*}V\otimes_{\mathrm{H}^{*}V/W}\mathrm{H}^{*}_{V}X^{W}$ l'homomorphisme composé
$$
\hspace{4pt}
\mathrm{H}^{*}_{V}X\to\mathrm{H}^{*}_{W\oplus V}X\to\mathrm{H}^{*}_{W\oplus V}X^{W}\cong\mathrm{H}^{*}W\otimes\mathrm{H}^{*}_{V}X^{W}
\cong\mathrm{H}^{*}V\otimes_{\mathrm{H}^{*}V/W}\mathrm{H}^{*}_{V}X^{W}
\hspace{4pt},
$$
la flèche de gauche étant induite par $\mathrm{a}_{W}$ et celle de droite par l'inclusion de $X^{W}$ dans $X$~;  $\delta_{W,X}$ est un morphisme dans la catégorie $V\hspace{-2pt}\text{-}\mathcal{U}$. On note
$$
\mathrm{E}\hspace{1pt}\nu_{X}:\mathrm{EFix}_{(V,W)}\mathrm{H}^{*}_{V}X\to\mathrm{H}^{*}_{V}X^{W}
$$
l'homomorphisme de $\mathrm{H}^{*}V$-$\mathrm{A}$-modules instables adjoint de $\delta_{W,X}$ .

\begin{pro}\label{cohEFix}
Si $X$ est un $V$-CW-complexe fini alors l'homomorphisme
$$
\mathrm{E}\hspace{1pt}\nu_{X}:\mathrm{EFix}_{(V,W)}\mathrm{H}^{*}_{V}X\to\mathrm{H}^{*}_{V}X^{W}
$$
est un isomorphisme.
\end{pro}

\bigskip
Les formules de \ref{EFix} ``s'inversent'' aisément~:

\begin{pro}\label{Fix} Soit $s:V/W\to V$ une section linéaire de $q:V\to V/W$. On considère $\mathrm{H}^{*}V/W$ comme un $(\mathrm{H}^{*}V/W,\mathrm{H}^{*}V)$-bimodule, la structure de $\mathrm{H}^{*}V/W$-module à gauche étant la structure évidente et la structure de $\mathrm{H}^{*}V$-module à droite étant induite par l'homomorphisme (surjectif) de $\mathrm{A}$-algèbres instables $s^{*}:\mathrm{H}^{*}V\to\mathrm{H}^{*}V/W$.  Soit $M$ un $\mathrm{H}^{*}V$-$\mathrm{A}$-module instable~; on a un isomorphisme de $\mathrm{H}^{*}V/W$-$\mathrm{A}$-modules instables, naturel en $M$~:
$$
\hspace{24pt}
\mathrm{Fix}_{(V,W)}M
\hspace{4pt}\cong\hspace{4pt}
\mathrm{H}^{*}V/W\otimes_{\mathrm{H}^{*}V}\mathrm{EFix}_{(V,W)}M
\hspace{24pt}.
$$
\end{pro}

\footnotesize
\begin{com}\label{remadj}
Il n'est pas difficile de se convaincre \textit{a priori} de l'existence d'un adjoint à gauche pour le foncteur $\mathrm{e}_{(V,W)}$ (existence que nous avons implicitement admise dans notre exposition). En fait, les propriétés d'adjonction des foncteurs $\mathrm{T}$ impliquent cette existence. Précisons un peu. Soit $\mathrm{e}_{s}:V/W\text{-}\mathcal{U}\to V\text{-}\mathcal{U}$ le foncteur défini \textit{via} $s^{*}$. On fait les deux observations suivantes~:

\smallskip
($\mathrm{O}_{1}$) L'endofoncteur  $\mathcal{O}_{(V,W)}\circ\mathrm{e}_{s}:V/W\text{-}\mathcal{U}\to V/W\text{-}\mathcal{U}$ est l'identité.

\smallskip
($\mathrm{O}_{2}$) Soit $\mathrm{f}_{s}:V\text{-}\mathcal{U}\to V/W\text{-}\mathcal{U}, M\mapsto\mathrm{H}^{*}V/W\otimes_{\mathrm{H}^{*}V}M$, le foncteur défini grâce à la structure de bimodule sur $\mathrm{H}^{*}V/W$ considérée en \ref{Fix}~;  $\mathrm{f}_{s}$ est adjoint à gauche de $\mathrm{e}_{s}$.

\smallskip
On a $\mathrm{e}_{(V,W)}=\mathrm{e}_{(V,W)}\circ\mathcal{O}_{(V,W)}\circ\mathrm{e}_{s}$ d'après $(\mathrm{O}_{1})$. Le foncteur $\mathrm{e}_{(V,W)}$ admet donc comme adjoint à gauche le composé $\mathrm{f}_{s}\circ\mathrm{EFix}_{(V,W)}$, d'après l'observation $(\mathrm{O}_{2})$ et la proposition \ref{adjonction}. On retrouve ainsi l'énoncé \ref{Fix}. On vérifie que le foncteur  $\mathrm{f}_{s}\circ\mathrm{EFix}_{(V,W)}$ s'identifie au foncteur $\mathrm{Fix}_{\hspace{1pt}W,s(V/W)}$ de \cite[1.3.4]{LZsmith}\footnote{Signalons incidemment un misprint dans la dernière formule de la page 14 de cette référence~: $\mathrm{Fix}'_{V}$ doit être remplacé par  $\mathrm{Fix}_{V'}$.} \textit{via} l'isomorphisme évident $\mathrm{H}^{*}V/W\cong\mathrm{H}^{*}s(V/W)$.
\end{com}

\medskip
{\begin{com}\label{uniadj}
L'expression que donne la proposition \ref{Fix} pour le foncteur $\mathrm{Fix}_{(V,W)}$ dépend du choix de $s$. Cela peut paraître inesthétique, mais ce n'est pas surprenant car un foncteur adjoint est défini à isomorphisme fonctoriel (canonique) près. Expliquons le phénomène dans le cas qui nous occupe. Notons $\mathrm{H}^{*}V/W\text{-}s$ le $(\mathrm{H}^{*}V/W,\mathrm{H}^{*}V)$-bimodule introduit en \ref{Fix}~; soient $s_{0}$ et $s_{1}$ deux sections linéaires de $q$, alors les deux bimodules  $\mathrm{H}^{*}V/W\text{-}s_{0}$ et  $\mathrm{H}^{*}V/W\text{-}s_{1}$ (avec toutes leurs structures) sont canoniquement isomorphes.
\end{com}
\normalsize

\bigskip
On fait  maintenant la liste des propriétés des foncteurs $\mathrm{Fix}$ (et $\mathrm{EFix}$) immé\-diatement impliquées par celles des foncteurs $\mathrm{T}$.

\begin{pro}\label{exactEFix} Le foncteur $\mathrm{EFix}_{(V,W)}$ est exact.
\end{pro}

\medskip
\textit{Démonstration.} Cette exactitude résulte de celle du foncteur $\mathrm{T}_{W}$ et du fait  que $\mathrm{H}^{*}V$ est un $\mathrm{T}_{W}\mathrm{H}^{*}V$-module plat \textit{via} $\widetilde{\varphi}$.
\hfill$\square$

\begin{pro}\label{exactFix} Le foncteur $\mathrm{Fix}_{(V,W)}$ est exact.
\end{pro}

\medskip
\textit{Démonstration.}  Compte tenu de la proposition \ref{EFix}, cette exactitude résulte de celle du foncteur $\mathrm{EFix}_{(V,W)}$ et du fait  que $\mathrm{H}^{*}V$ est  un $\mathrm{H}^{*}V/W$-module fidèlement plat.
\hfill$\square$

\begin{pro}\label{tensEFix} Soient $M_{1}$ et $M_{2}$ deux $\mathrm{H}^{*}V$-$\mathrm{A}$-modules instables. On a un isomorphisme de $\mathrm{H}^{*}V$-$\mathrm{A}$-modules instables, naturel en $M_{1}$ et $M_{2}$~:
$$
\hspace{24pt}
\mathrm{EFix}_{(V,W)}(M_{1}\otimes_{\mathrm{H}^{*}V}M_{2})
\hspace{4pt}\cong\hspace{4pt}
\mathrm{EFix}_{(V,W)}(M_{1})\otimes_{\mathrm{H}^{*}V}
\mathrm{EFix}_{(V,W)}(M_{2})
\hspace{24pt}.
$$
\end{pro}

\textit{Démonstration.} Résulte du fait que le foncteur $\mathrm{T}_{W}$ préserve les produits tensoriels.
\hfill$\square$

\bigskip
Compte tenu de \ref{Fix} (ou \ref{EFix}), la proposition \ref{tensEFix} fournit aisément l'énoncé suivant (qui généralise le théorème 4.6.2.1 de \cite{LaT})~:

\begin{cor}\label{tensFix} Soient $M_{1}$ et $M_{2}$ deux $\mathrm{H}^{*}V$-$\mathrm{A}$-modules instables. On a un isomorphisme de $\mathrm{H}^{*}V/W$-$\mathrm{A}$-modules instables, naturel en $M_{1}$ et $M_{2}$~:
$$
\hspace{24pt}
\mathrm{Fix}_{(V,W)}(M_{1}\otimes_{\mathrm{H}^{*}V}M_{2})
\hspace{4pt}\cong\hspace{4pt}
\mathrm{Fix}_{(V,W)}(M_{1})\otimes_{\mathrm{H}^{*}V/W}
\mathrm{Fix}_{(V,W)}(M_{2})
\hspace{24pt}.
$$
\end{cor}

Le slogan pour la proposition \ref{tensEFix} (resp. le corollaire \ref{tensFix}) est que les foncteurs $\mathrm{EFix}$ (resp. $\mathrm{Fix}$) ``préservent les produits tensoriels''. Voici quelques conséquences de cette préservation.

\medskip
Soit $S$ un $\mathrm{A}$-module instable~; les deux énoncés qui suivent font intervenir le foncteur de $V\text{-}\mathcal{U}\to V\text{-}\mathcal{U},M\mapsto S\otimes M$. On observera que pour $S=\Sigma\hspace{1pt}\mathbb{F}_{\ell}$ ce~foncteur est le {\em foncteur suspension}.

\begin{cor}\label{suspensionEFix} Soit $S$ un $\mathrm{A}$-module instable.

\smallskip
{\em (a)} On a un isomorphisme de $\mathrm{H}^{*}V$-$\mathrm{A}$-modules instables, naturel en le $\mathrm{H}^{*}V$-$\mathrm{A}$-module instable $M$~:
$$
\hspace{24pt}
\mathrm{EFix}_{(V,W)}(S\otimes M)
\hspace{4pt}\cong\hspace{4pt}
\mathrm{T}_{W}S\otimes\mathrm{EFix}_{(V,W)}M
\hspace{24pt}.
$$

\smallskip
{\em (b)} Si $S$ est fini, alors on a un isomorphisme de $\mathrm{H}^{*}V$-$\mathrm{A}$-modules instables, naturel en le $\mathrm{H}^{*}V$-$\mathrm{A}$-module instable $M$~:
$$
\hspace{24pt}
\mathrm{EFix}_{(V,W)}(S\otimes M)
\hspace{4pt}\cong\hspace{4pt}
S\otimes\mathrm{EFix}_{(V,W)}M
\hspace{24pt}.
$$
\end{cor}

\textit{Démonstration.} Le point (b) est conséquence du point (a) car si $S$ est fini alors on a un isomorphisme (canonique) de $\mathrm{A}$-modules instables $\mathrm{T}_{W}S\cong S$. Pour démontrer le point (a) on observe que l'on a $S\otimes M\cong(S\otimes\mathrm{H}^{*}V)\otimes_{\mathrm{H}^{*}V}M$, on vérifie l'isomorphisme $\mathrm{EFix}_{(V,W)}(S\otimes\mathrm{H}^{*}V)\cong\mathrm{T}_{W}S\otimes\mathrm{H}^{*}V$ et on invoque la proposition \ref{tensEFix}.
\hfill$\square$

\begin{cor}\label{suspensionFix} Soit $S$ un $\mathrm{A}$-module instable.

\smallskip
{\em (a)} On a un isomorphisme de $\mathrm{H}^{*}V$-$\mathrm{A}$-modules instables, naturel en le $\mathrm{H}^{*}V$-$\mathrm{A}$-module instable $M$~:
$$
\hspace{24pt}
\mathrm{Fix}_{(V,W)}(S\otimes M)
\hspace{4pt}\cong\hspace{4pt}
\mathrm{T}_{W}S\otimes\mathrm{Fix}_{(V,W)}M
\hspace{24pt}.
$$

\smallskip
{\em (b)} Si $S$ est fini, alors on a un isomorphisme de $\mathrm{H}^{*}V$-$\mathrm{A}$-modules instables, naturel en le $\mathrm{H}^{*}V$-$\mathrm{A}$-module instable $M$~:
$$
\hspace{24pt}
\mathrm{Fix}_{(V,W)}(S\otimes M)
\hspace{4pt}\cong\hspace{4pt}
S\otimes\mathrm{Fix}_{(V,W)}M
\hspace{24pt}.
$$
\end{cor}

\begin{pro-def}\label{torFix} Soient $V$ un $\ell$-groupe abélien élémentaire, $W\subset V$ un sous-groupe, $M_{1}$ et $M_{2}$ deux $\mathrm{H}^{*}V$-$\mathrm{A}$-modules instables et $p$ un entier naturel.

\medskip
{\em (a)} La structure de $\mathrm{H}^{*}V$-module gradué de $\mathrm{Tor}_{p}^{\mathrm{H}^{*}V}(M_{1},M_{2})$ peut être naturellement enrichie en une structure de $\mathrm{H}^{*}V$-$\mathrm{A}$-module instable. Le $\mathrm{H}^{*}V$-$\mathrm{A}$-module instable ainsi obtenu sera toujours noté  $\mathrm{Tor}_{p}^{\mathrm{H}^{*}V}(M_{1},M_{2})$.

\medskip
{\em (b)} On a un isomorphisme de $\mathrm{H}^{*}V/W$-$\mathrm{A}$-modules instables
$$
\hspace{24pt}
\mathrm{Fix}_{(V,W)}\hspace{1pt}\mathrm{Tor}_{p}^{\mathrm{H}^{*}V}(M_{1},M_{2})
\hspace{4pt}\cong\hspace{4pt}
\mathrm{Tor}_{p}^{\mathrm{H}^{*}V/W}(
\mathrm{Fix}_{(V,W)}M_{1},\mathrm{Fix}_{(V,W)}M_{2})
\hspace{24pt},
$$
naturel en $M_{1}$ et $M_{2}$.
\end{pro-def}

\medskip
\textit{Démonstration.} Rappelons pour commencer la théorie des projectifs de la catégorie abélienne $V\text{-}\mathcal{U}$. Cette théorie est banale~:

\smallskip
-- Le foncteur $M\mapsto M^{k}$, $k\in\mathbb{N}$, défini sur la catégorie $V\text{-}\mathcal{U}$ et à valeurs dans la catégorie des $\mathbb{F}_{\ell}$-espaces vectoriels, est représentable~:
$$
\hspace{24pt}
M^{k}
\hspace{4pt}\cong\hspace{4pt}
\mathrm{Hom}_{V\text{-}\mathcal{U}}(\mathrm{H}^{*}V \otimes\mathrm{F}(k),M)
\hspace{24pt},
$$
$\mathrm{F}(k)$ désignant le $\mathrm{A}$-module instable librement engendré par un élément de degré $k$~; $\mathrm{H}^{*}V \otimes\mathrm{F}(k)$ est donc ``tautologiquement'' un projectif de $V\text{-}\mathcal{U}$ (compte tenu de sa définition, on dit aussi qu'il est libre).

\smallskip
-- Tout objet de $V\text{-}\mathcal{U}$ est fonctoriellement quotient d'une somme directe de ces projectifs tautologiques. En particulier $V\text{-}\mathcal{U}$ a assez de projectifs.

\smallskip
-- Tout projectif de $V\text{-}\mathcal{U}$ est isomorphe à une telle somme directe.

\medskip
Après ce préalable, passons à la démonstration.

\medskip
Soit $M_{i}\leftarrow P_{i,\bullet}$, $i=1,2$, une résolution projective de $M_{i}$ dans la catégorie $V\text{-}\mathcal{U}$~; on considère le  $\mathrm{H}^{*}V$-$\mathrm{A}$-module instable
$$
\Theta_{p}
\hspace{4pt}:=\hspace{4pt}
\mathrm{H}_{p}\hspace{1pt}\mathrm{Tot}\hspace{1pt}(P_{1,\bullet}\otimes_{\mathrm{H}^{*}V}P_{2,\bullet})
$$
(la notation $\mathrm{Tot}$ désigne le totalisé d'un bicomplexe). Les arguments habituels montrent que $\Theta_{p}$ ``est indépendant'' du choix des résolutions et que l'on a des isomorphismes canoniques
$$
\hspace{24pt}
\Theta_{p}
\hspace{4pt}\cong\hspace{4pt}
\mathrm{H}_{p}\hspace{1pt}(P_{1,\bullet}\otimes_{\mathrm{H}^{*}V}M_{2})
\hspace{12pt}\text{et}\hspace{12pt}
\Theta_{p}
\hspace{4pt}\cong\hspace{4pt}
\mathrm{H}_{p}\hspace{1pt}(M_{1}\otimes_{\mathrm{H}^{*}V}P_{2,\bullet})
\hspace{24pt}.
$$
On fait les deux observations suivantes~:

\smallskip
-- Le $\mathrm{H}^{*}V$-module sous-jacent à un projectif de $V\text{-}\mathcal{U}$ est libre.

\smallskip
-- Le foncteur $\mathrm{Fix}_{(V,W)}:V\text{-}\mathcal{U}\to V/W\text{-}\mathcal{U}$ transforme projectif en projectif. Ceci résulte formellement du fait que $\mathrm{Fix}_{(V,W)}$ est adjoint à gauche d'un foncteur exact.

\medskip
La première observation montre que $M_{i}\leftarrow P_{i,\bullet}$ est une résolution libre du $\mathrm{H}^{*}V$-module sous-jacent à $M_{i}$ et conduit au point (a)~:  le $\mathrm{H}^{*}V$-module sous-jacent à $\Theta_{p}$ est bien $\mathrm{Tor}_{p}^{\mathrm{H}^{*}V}(M_{1},M_{2})$. La seconde observation et l'exactitude de $\mathrm{Fix}_{(V,W)}$ montrent que $\mathrm{Fix}_{(V,W)}M_{i}\leftarrow \mathrm{Fix}_{(V,W)}P_{i,\bullet}$ est une résolution projective dans la catégorie $V/W\text{-}\mathcal{U}$ et l'on obtient le point (b) en invoquant le fait que $\mathrm{Fix}_{(V,W)}$ ``préserve les produits tensoriels'' (Corollaire \ref{tensFix}) et le point (a).
\hfill$\square$

\medskip
\begin{rem}\label{derivetor} Soit $\mathrm{gtens}_{M_{1}}:V\text{-}\mathcal{U}\to V\text{-}\mathcal{U}$ (resp. $\mathrm{dtens}_{M_{2}}:V\text{-}\mathcal{U}\to V\text{-}\mathcal{U}$) le foncteur $M\mapsto M_{1}\otimes_{\mathrm{H}^{*}V}M$ (resp. $M\mapsto M\otimes_{\mathrm{H}^{*}V}M_{2}$)~; $\mathrm{gtens}_{M_{1}}$ (resp. $\mathrm{dtens}_{M_{2}}$) est exact à droite et l'on a
$$
\hspace{18pt}
\mathrm{Tor}_{p}^{\mathrm{H}^{*}V}(M_{1},M_{2})
\hspace{4pt}=\hspace{4pt}
(\mathrm{L}_{p}\hspace{1pt}\mathrm{gtens}_{M_{1}})(M_{2})
\hspace{4pt}=\hspace{4pt}
(\mathrm{L}_{p}\hspace{1pt}\mathrm{dtens}_{M_{2}})(M_{1})\hspace{4pt}=\hspace{4pt}
\Theta_{p}
\hspace{18pt}.
$$
\end{rem}

\medskip
Voici une spécialisation du point (b) de la proposition précédente que nous utiliserons dans les sections 4 et 5.

\begin{cor}\label{sptorFix} Soient $M$ un $\mathrm{H}^{*}V$-$\mathrm{A}$-module instable et $p$ un entier naturel~; on a un isomorphisme de $\mathrm{H}^{*}V/W$-$\mathrm{A}$-modules instables
$$
\hspace{24pt}
\mathrm{Fix}_{(V,W)}\hspace{1pt}\mathrm{Tor}_{p}^{\mathrm{H}^{*}V}(\mathrm{H}^{*}W,M)
\hspace{4pt}\cong\hspace{4pt}
\mathrm{Tor}_{p}^{\mathrm{H}^{*}V/W}(\mathbb{F}_{\ell},\mathrm{Fix}_{(V,W)}M)
\hspace{24pt},
$$
naturel en $M$.
\end{cor}

\bigskip
\textit{Démonstration.} On a $\mathrm{Fix}_{(V,W)}\mathrm{H}^{*}W\cong\mathbb{F}_{\ell}$ (pour s'en convaincre on peut voir $\mathrm{H}^{*}W$ comme le produit tensoriel $\mathrm{H}^{*}V\otimes_{\mathrm{H}^{*}V/W}\mathbb{F}_{\ell}$ et invoquer \ref{corproclef}).
\hfill$\square$

\bigskip
On clôt la liste déroulée ci-dessus par la proposition-définition \ref{Wfonctoriel} ci-après~; cette proposition traite de la ``fonctorialité en $W$'' des $\mathrm{H}^{*}V$-$\mathrm{A}$-modules instables $\mathrm{EFix}_{(V,W)}M$, sa démonstration est laissée au lecteur.

\begin{pro-def} \label{Wfonctoriel} Soient $W_{0}$ et $W_{1}$ deux sous-groupes de $V$ avec $W_{0}\subset W_{1}$. La transformation naturelle $\mathrm{T}_{W_{0}}\to\mathrm{T}_{W_{1}}$ induit une transformation naturelle
$$
\hspace{24pt}
\rho\hspace{1pt}(W_{0},W_{1}):\mathrm{EFix}_{(V,W_{0})}\to\mathrm{EFix}_{(V,W_{1})}
\hspace{24pt}.
$$
Soient $W_{0}$, $W_{1}$ et $W_{2}$ trois sous-groupes de $V$ avec $W_{0}\subset W_{1}\subset W_{2}$~; on~a~:
$$
\hspace{24pt}
\rho\hspace{1pt}(W_{0},W_{2})
\hspace{4pt}=\hspace{4pt}
\rho\hspace{1pt}(W_{1},W_{2})\circ\rho\hspace{1pt}(W_{0},W_{1})
\hspace{24pt}.
$$
Soit $W$ un sous-groupe de $V$ ; la transformation naturelle $\rho(0,W)$ s'identifie à une transformation naturelle $\mathrm{id}\to\mathrm{EFix}_{(V,W)}$ qui sera aussi notée $\rho_{(V,W)}$.
\end{pro-def}

\bigskip
\begin{rem}\label{Wfonctoriel-2} D'après \ref{adjonction}, les endofoncteurs $\mathrm{EFix}_{(V,W_{0})}$ et $\mathrm{EFix}_{(V,W_{1})}$ sont respectivement adjoints à gauche des endofoncteurs $\mathrm{E}_{(V,W_{0})}$ et $\mathrm{E}_{(V,W_{1})}$ (on rappelle que $\mathrm{E}_{(V,W)}$ est le foncteur $V\text{-}\mathcal{U}\to V\text{-}\mathcal{U},N\mapsto\mathrm{H}^{*}V\otimes_{\mathrm{H}^{*}V/W}N$). On vérifie que la transformation naturelle $\rho\hspace{1pt}(W_{0},W_{1})$ correspond par adjonction à la transformation naturelle $\sigma\hspace{1pt}(W_{1},W_{0}):\mathrm{E}_{(V,W_{1})}\to\mathrm{E}_{(V,W_{0})}$ telle que $\sigma\hspace{1pt}(W_{1},W_{0})_{N}$ est l'épimorphisme canonique $\mathrm{H}^{*}V\otimes_{\mathrm{H}^{*}V/W_{1}}N\to\mathrm{H}^{*}V\otimes_{\mathrm{H}^{*}V/W_{0}}N$ (spécialisation de \cite[Chap. IV, \S7, Theorem~2]{McL}).
\end{rem}

\bigskip
\begin{rem}\label{Wfonctoriel-3} Soient $W_{0}$, $W_{1}$ deux sous-groupes de $V$ avec $W_{0}\subset W_{1}$ et~$X$ un $V$-CW-complexe fini. On vérifie que le diagramme
$$
\begin{CD}
\mathrm{EFix}_{(V,W_{0})}\mathrm{H}^{*}_{V}X
@>\mathrm{E}\hspace{1pt}\nu_{W_{0},X}>\cong>
\mathrm{H}^{*}_{V}X^{W_{0}} \\
@V\rho\hspace{1pt}(W_{0},W_{1})VV @VVV \\
\mathrm{EFix}_{(V,W_{1})}\mathrm{H}^{*}_{V}X
@>\mathrm{E}\hspace{1pt}\nu_{W_{1},X}>\cong>
\mathrm{H}^{*}_{V}X^{W_{1}}
\end{CD}
$$
est commutatif (pour des raisons évidentes on a précisé ci-dessus la notation~$\mathrm{E}\hspace{1pt}\nu_{X}$ de \ref{cohEFix} en $\mathrm{E}\hspace{1pt}\nu_{W,X}$).
\end{rem}

\bigskip
\begin{rem}\label{Wfonctoriel-4} L'homomorphisme de restriction $\mathrm{H}^{*}_{V}X^{W_{0}}\to\mathrm{H}^{*}_{V}X^{W_{1}}$ qui apparaît dans la remarque précédente peut s'écrire $\mathrm{H}^{*}V\otimes_{\mathrm{H}^{*}V/W_{0}}r$, $r$ désignant l'homomorphisme de restriction $\mathrm{H}^{*}_{V/W_{0}}X^{W_{0}}\to\mathrm{H}^{*}_{V/W_{0}}X^{W_{1}}$. La version ``purement algébrique'' de cette observation est la suivante~:

\smallskip
Soit $M$ un $\mathrm{H}^{*}V$-$\mathrm{A}$-module instable. On considère l'homomorphisme de $\mathrm{H}^{*}V$-$\mathrm{A}$-modules instables
$$
{\rho(W_{0},W_{1})}_{M}
\hspace{2pt}:\hspace{2pt}
\mathrm{EFix}_{(V,W_{0})}\hspace{1pt}M
\longrightarrow
\mathrm{EFix}_{(V,W_{1})}\hspace{1pt}M
$$
et  l'homomorphisme de $\mathrm{H}^{*}V/W_{0}$-$\mathrm{A}$-modules instables
$$
{\rho_{(V/W_{0},W_{1}/W_{0})}}_{\mathrm{Fix}_{(V,W_{0})}M}
\hspace{2pt}:\hspace{2pt}
\mathrm{Fix}_{(V,W_{0})}M
\longrightarrow
\mathrm{EFix}_{(V/W_{0},W_{1}/W_{0})}\hspace{1pt}\mathrm{Fix}_{(V,W_{0})}M
$$
(la notation $\rho_{(-,-)}$ est introduite à la toute fin de \ref{Wfonctoriel})~; on a l'identification
$$
\hspace{24pt}
{\rho(W_{0},W_{1})}_{M}
\hspace{4pt}=\hspace{4pt}
\mathrm{H}^{*}V\otimes_{\mathrm{H}^{*}V/W_{0}}\hspace{1pt}{\rho_{(V/W_{0},W_{1}/W_{0})}}_{\mathrm{Fix}_{(V,W_{0})}M}
\hspace{24pt}.
$$

On s'en convainc en contemplant les identités~:
$$
\mathrm{E}_{(V,W_{1})}=\mathrm{e}_{(V,W_{0})}\circ\mathrm{E}_{(V/W_{0},W_{1}/W_{0})}\circ\mathcal{O}_{(V,W_{0})}
\hspace{8pt},\hspace{8pt}
\mathrm{E}_{(V,W_{0})}=\mathrm{e}_{(V,W_{0})}\circ\mathrm{id}\circ\mathcal{O}_{(V,W_{0})}
$$
($\mathrm{id}$ désigne ici le foncteur identique de la catégorie $V/W_{0}\text{-}\mathcal{U}$) et en observant que la transformation naturelle $\sigma\hspace{1pt}(W_{1},W_{0}):\mathrm{E}_{(V,W_{1})}\to\mathrm{E}_{(V,W_{0})}$ de~\ref{Wfonctoriel-2} est induite par  $\sigma\hspace{1pt}(W_{1}/W_{0},0):\mathrm{E}_{(V/W_{0},W_{1}/W_{0})}\to\mathrm{E}_{(V/W_{0\hspace{1pt}},\hspace{1pt}0)}=\mathrm{id}$.
\end{rem}

\bigskip
Compte tenu de \ref{EFix} la transformation naturelle $\rho_{(V,W)}:\mathrm{id}\to\mathrm{EFix}_{(V,W)}$ s'identifie à une transformation naturelle $\mathrm{id}\to\mathrm{e}_{(V,W)}\circ\mathrm{Fix}_{(V,W)}$~; on constate sans surprise que l'on a l'énoncé suivant~: 

\begin{pro}\label{commutativity} Soit $\kappa_{(V,W)}:\mathrm{EFix}_{(V,W)}\to\mathrm{e}_{(V,W)}\circ\mathrm{Fix}_{(V,W)}$ l'isomorphisme naturel de \ref{EFix}, alors la transformation naturelle composée
$$
\begin{CD}
\mathrm{id}@>\rho_{(V,W)}>>\mathrm{EFix}_{(V,W)}@>\kappa_{(V,W)}>>\mathrm{e}_{(V,W)}\circ\mathrm{Fix}_{(V,W)}
\end{CD}
$$
est l'unité de l'adjonction du couple de foncteurs adjoints $(\mathrm{Fix}_{(V,W)},\mathrm{e}_{(V,W)})$.

\end{pro}

\textit{Démonstration.} Elle fait intervenir les unités des adjonctions, pour les paires de foncteurs adjoints  $(\mathrm{Fix}_{(V,W)},\mathrm{e}_{(V,W)})$ et  $(\mathrm{EFix}_{(V,W)},\mathrm{e}_{(V,W)}\circ\mathcal{O}_{(V,W)})$ et la co-unité de l'adjonction de la paire de foncteurs adjoints $(\mathrm{e}_{(V,W)},\mathcal{O}_{(V,W)}$). Les détails (fastidieux) sont laissés au lecteur.
\hfill$\square$

\bigskip
On en vient maintenant  à des propriétés plus techniques des foncteurs $\mathrm{Fix}$ (\ref{proclef}, \ref{corproclef} et \ref{proclefbis}) qui auront un rôle à jouer dans notre mémoire.

\bigskip
\footnotesize
Compte tenu de \ref{cohFix}, \ref{EFix} et \ref{Fix}, les énoncés \ref{proclef}, \ref{corproclef} et \ref{proclefbis} ci-après constituent le pendant algébrique de la simple observation suivante~:

\medskip
Soit $X$ un $V/U$-CW-complexe fini tel que l'action de $V/U$ sur $X$ est libre. Si l'on considère $X$ comme un $V$-CW-complexe (fini) alors on a
$$
X^{W}
\hspace{4pt}=\hspace{4pt}
\begin{cases}
X
& \text{pour $W\subset U$,} \\
\emptyset
& \hspace{1pt}\text{pour $W\not\subset U$.}
\end{cases}
$$
\normalsize

\begin{pro}\label{proclef} Soient $V$ un $\ell$-groupe abélien élémentaire et $W$, $U$ des sous-groupes~; soit $N$ un $\mathrm{H}^{*}V/U$-$\mathrm{A}$-module instable fini.

\medskip
{\em (a)} Si l'on a $W\subset U$ alors l'homomorphisme de $\mathrm{H}^{*}V$-$\mathrm{A}$-modules instables
$$
{\rho_{(V,W)}}_{\mathrm{H}^{*}V\otimes_{\mathrm{H}^{*}V/U}N}
\hspace{4pt}:\hspace{4pt}
\mathrm{H}^{*}V\otimes_{\mathrm{H}^{*}V/U}N\to
\mathrm{EFix}_{(V,W)}(\mathrm{H}^{*}V\otimes_{\mathrm{H}^{*}V/U}N)
$$
(introduit en \ref{Wfonctoriel}) est un isomorphisme.

\medskip
{\em (b)} Si l'on a $W\not\subset U$ alors $\mathrm{EFix}_{(V,W)}(\mathrm{H}^{*}V\otimes_{\mathrm{H}^{*}V/U}N)$ est nul.
\end{pro}

\bigskip
\textit{Démonstration.} Il faut montrer en particulier que l'on a
$$
\mathrm{H}^{*}V\otimes_{\mathrm{T}_{W}\mathrm{H}^{*}V}\mathrm{T}_{W}(\mathrm{H}^{*}V\otimes_{\mathrm{H}^{*}V/U}N)
\hspace{4pt}\cong\hspace{4pt}
\begin{cases}
\mathrm{H}^{*}V\otimes_{\mathrm{H}^{*}V/U}N
& \text{pour $W\subset U$,} \\
0
& \hspace{1pt}\text{pour $W\not\subset U$,}
\end{cases}
\leqno{(*)}
$$
$\mathrm{H}^{*}V$ étant ci-dessus  un $\mathrm{T}_{W}\mathrm{H}^{*}V$-module à droite \textit{via} $\widetilde{\varphi}$ ou ce qui revient au même \textit{via} l'homomorphisme composé
$$
\begin{CD}
\hspace{24pt}\mathrm{T}_{W}\mathrm{H}^{*}V
\hspace{2pt}\cong\hspace{2pt}
(\mathrm{H}^{*}V)^{\hspace{1pt}\mathrm{Hom}(W,V)}
@>\pi_{i}>>\mathrm{H}^{*}V\hspace{24pt},
\end{CD}
$$
$\pi_{i}$ désignant la projection sur la composante indexée par l'inclusion $i$ de $W$ dans $V$. Comme $\mathrm{T}_{W}$ ``commute aux produits tensoriels'', on dispose d'un isomorphisme (naturel en $N$) de $\mathrm{T}_{W}\mathrm{H}^{*}V$-$\mathrm{A}$-modules instables
$$
\hspace{24pt}
\mathrm{T}_{W}(\mathrm{H}^{*}V\otimes_{\mathrm{H}^{*}V/U}N)
\hspace{4pt}\cong\hspace{4pt}
\mathrm{T}_{W}\mathrm{H}^{*}V\otimes_{\mathrm{T}_{W}\hspace{1pt}\mathrm{H}^{*}V/U}\mathrm{T}_{W}N
\hspace{24pt}.
$$
Si $N$ est fini alors on a $\mathrm{T}_{W}N\cong N$, ce que l'on peut préciser ainsi~: l'homorphisme canonique $\mathrm{T}_{0}N\to\mathrm{T}_{W}N$ induit par l'inclusion de $0$ dans $W$ est un isomorphisme. Il en résulte
$$
\hspace{24pt}
\mathrm{T}_{W}(\mathrm{H}^{*}V\otimes_{\mathrm{H}^{*}V/U}N)
\hspace{4pt}\cong\hspace{4pt}
\mathrm{T}_{W}\mathrm{H}^{*}V\otimes_{\mathrm{T}_{W}\hspace{1pt}\mathrm{H}^{*}V/U}N
\hspace{24pt},
$$
$N$ étant  $\mathrm{T}_{W}\hspace{1pt}\mathrm{H}^{*}V/U$-module à gauche \textit{via} l'homomorphisme composé
$$
\begin{CD}
\hspace{24pt}\mathrm{T}_{W}\hspace{1pt}\mathrm{H}^{*}V/U
\hspace{2pt}\cong\hspace{2pt}
(\mathrm{H}^{*}V/U)^{\hspace{1pt}\mathrm{Hom}(W,V/U)}
@>\pi_{0}>>\mathrm{H}^{*}V/U\hspace{24pt},
\end{CD}
$$
$\pi_{0}$ désignant la projection sur la composante indexée par l'homorphisme nul de $W$ dans $V/U$. On note $\kappa:W\to V/U$ l'homorphisme composé $q_{U}\circ i$, $q_{U}$ désignant la surjection canonique $V\to V/U$~; on obtient au bout du compte un isomorphisme (naturel en $N$) de $\mathrm{H}^{*}V$-$\mathrm{A}$-modules instables
$$
\mathrm{H}^{*}V\otimes_{\mathrm{T}_{W}\mathrm{H}^{*}V}\mathrm{T}_{W}(\mathrm{H}^{*}V\otimes_{\mathrm{H}^{*}V/U}N)
\hspace{4pt}\cong\hspace{4pt}
\mathrm{H}^{*}V\otimes_{\mathrm{H}^{*}V/U}K\otimes_{\mathrm{H}^{*}V/U}N
$$
avec
$$
\hspace{24pt}
K
\hspace{4pt}:=\hspace{4pt}
\mathrm{H}^{*}V/U
\otimes_{(\mathrm{H}^{*}V/U)^{\mathrm{Hom}(W,V/U)}}
\mathrm{H}^{*}V/U
\hspace{24pt},
$$
$\mathrm{H}^{*}V/U$ étant à gauche un $(\mathrm{H}^{*}V/U)^{\mathrm{Hom}(W,V/U)}$-module à droite \textit{via} $\pi_{\kappa}$ (la projection sur la composante indexée par $\kappa$) et $\mathrm{H}^{*}V/U$ étant à droite un $(\mathrm{H}^{*}V/U)^{\mathrm{Hom}(W,V/U)}$-module à gauche \textit{via} $\pi_{0}$. On constate que l'on a
$$
K
\hspace{4pt}\cong\hspace{4pt}
\begin{cases}
\mathrm{H}^{*}V/U
& \text{pour $\kappa=0$,} \\
0
& \hspace{1pt}\text{pour $\kappa\not=0$,}
\end{cases}
$$
ce qui conduit bien à $(*)$.

\bigskip
Pour achever la démonstration, il reste à vérifier que dans le cas $W\subset U$ l'inverse de l'isomorphisme (*) coïncide avec la transformation naturelle $\rho_{(V,W)}$ introduite en \ref{Wfonctoriel}, ou encore que le composé de $\rho_{(V,W)}$ et de l'isomorphisme (*) est l'identité. Cette dernière vérification est immédiate.
\hfill$\square$

\medskip
\begin{cor}\label{corproclef} Soient $V$ un $\ell$-groupe abélien élémentaire et $W$, $U$ des sous-groupes~; soit $N$ un $\mathrm{H}^{*}V/U$-$\mathrm{A}$-module instable. Si $N$ est fini, alors on a un isomorphisme de $\mathrm{H}^{*}V/W$-$\mathrm{A}$-modules instables, naturel en $N$
$$
\mathrm{Fix}_{(V,W)}(\mathrm{H}^{*}V\otimes_{\mathrm{H}^{*}V/U}N)
\hspace{4pt}\cong\hspace{4pt}
\begin{cases}
\mathrm{H}^{*}V/W\otimes_{\mathrm{H}^{*}V/U}N
& \text{pour $W\subset U$,} \\
0
& \hspace{1pt}\text{pour $W\not\subset U$.}
\end{cases}
$$
\end{cor}

\medskip
\textit{Démonstration.} Conséquence de \ref{Fix} et \ref{proclef}.
\hfill$\square$

\bigskip
Compte tenu de \ref{commutativity}, la proposition \ref{proclef} peut être reformulée ainsi~:

\medskip
\begin{pro}\label{proclefbis} Soient $V$ un $\ell$-groupe abélien élémentaire et $W$, $U$ des sous-groupes~; soit $N$ un $\mathrm{H}^{*}V/U$-$\mathrm{A}$-module instable fini. Alors l'unité d'adjonction
$$
\begin{CD}
\mathrm{H}^{*}V\otimes_{\mathrm{H}^{*}V/U}N
@>\eta_{(V,W)}>>
\mathrm{H}^{*}V\otimes_{\mathrm{H}^{*}V/W}\mathrm{Fix}_{(V,W)}(\mathrm{H}^{*}V\otimes_{\mathrm{H}^{*}V/U}N)
\end{CD}
$$
est un isomorphisme si l'on a $W\subset U$  et est nulle si l'on a $W\not\subset U$.
\end{pro}

\vspace{1,5cm}
\textit{Démonstration de la proposition \ref{cohFix}.}

\medskip
Elle est analogue à celle que l'on trouve dans \cite[paragraphe 4.7]{LaT}. On la divise en trois étapes.

\medskip
1) On étend la transformation naturelle $\nu_{X}$ en une transformation naturelle entre foncteurs définis sur la catégorie des paires de $V$-espaces
$$
\hspace{24pt}
\nu_{(X,Y)}:\mathrm{Fix}_{(V,W)}\mathrm{H}^{*}_{V}(X,Y)\to\mathrm{H}^{*}_{V/W}(X^{W},Y^{W})
\hspace{24pt}.
$$
On observe que l'énoncé suivant est vérifié~:

\begin{pro}\label{Fixconnectant} Soit $(X,Y)$ une paire de $V$-espaces, alors le diagramme d'homomorphismes de $\mathrm{H}^{*}V/W$-$\mathrm{A}$-modules instables
$$
\begin{CD}
\mathrm{Fix}_{(V,W)}\mathrm{H}^{*}_{V}(X,Y)
@>\mathrm{Fix}_{(V,W)}\hspace{1pt}\partial>>
\Sigma\hspace{1pt}\mathrm{Fix}_{(V,W)}\mathrm{H}^{*}_{V}Y \\
@V\nu_{(X,Y)}VV @V\Sigma\hspace{1pt}\nu_{Y}VV \\
\mathrm{H}^{*}_{V/W}(X^{W},Y^{W})
@>\partial>>
\Sigma\hspace{1pt}\mathrm{H}^{*}_{V/W}Y^{W}
\end{CD}
$$
est commutatif.

\smallskip
(Au-dessus de la flèche horizontale du haut $\partial$ désigne le connectant, en cohomologie $V$-équivariante de la paire $(X,Y)$, au-dessus de la flèche horizontale du bas $\partial$ désigne le connectant, en cohomologie $V/W$-équivariante de la paire $(X^{W},Y^{W})$.)
\end{pro}

\medskip
On en déduit que si $\nu_{Y}$ et $\nu_{(X,Y)}$ sont des isomorphisme alors il en est de même pour $\nu_{X}$.

\medskip
2) On vérifie que $\nu_{(X,Y)}$ est un isomorphisme pour $(X,Y)=(\mathrm{D}^{m},\mathrm{S}^{m-1})\times V/U$ avec $m$ un entier naturel et $U$ un sous-groupe de $V$. Compte tenu du point (b) de \ref{suspensionFix} (avec $S=\Sigma^{m}\mathbb{F}_{\ell}$), il suffit d'effectuer cette vérification pour $m=0$, c'est-à-dire de se convaincre que $\nu_{X}$ est un isomorphisme pour $X=V/U$.

\medskip
On pose donc $X=V/U$. On a $\mathrm{H}^{*}_{V}X\cong\mathrm{H}^{*}V\otimes_{\mathrm{H}^{*}V/U}\mathrm{H}^{*}_{V/U}X\cong\mathrm{H}^{*}V\otimes_{\mathrm{H}^{*}V/U}\mathbb{F}_{\ell}$, si bien que l'on peut appliquer \ref{corproclef}~; on obtient~:
$$
\mathrm{Fix}_{(V,W)}\mathrm{H}^{*}_{V}X
\hspace{4pt}\cong\hspace{4pt}
\begin{cases}
\mathrm{H}^{*}_{V/W}X & \text{pour $W\subset U$,} 
\\
0 & \hspace{1pt}\text{pour $W\not\subset U$.}
\end{cases}
$$
Dans les deux cas le second membre est égal à $\mathrm{H}^{*}_{V/W}X^{W}$, en effet~:
$$
X^{W}
\hspace{4pt}=\hspace{4pt}
\begin{cases}
X
& \text{pour $W\subset U$,} \\
\emptyset
& \hspace{1pt}\text{pour $W\not\subset U$.}
\end{cases}
$$
Dans le cas $W\not\subset U$ il n'y a plus rien à démontrer~; passons au cas $W\subset U$. Par ``abstract nonsense'' (même argument que pour \ref{etaFix}), le diagramme
$$
\begin{CD}
\mathrm{H}^{*}_{V}X @>\eta_{(V,W)}>>
\mathrm{H}^{*}V\otimes_{\mathrm{H}^{*}V/W}\mathrm{Fix}_{(V,W)}\mathrm{H}^{*}_{V}X
\\
@VVV @VV\mathrm{H}^{*}{V}\otimes_{\mathrm{H}^{*}{V}/W}\hspace{1pt}\nu_{X}V
\\
\mathrm{H}^{*}_{V}X^{W} @<<<
\mathrm{H}^{*}V\otimes_{\mathrm{H}^{*}V/W}\mathrm{H}^{*}_{V/W}X^{W}
\end{CD}
$$
est commutatif. Il est clair que la flèche verticale de gauche et la flèche horizontale du bas sont des isomorphismes~;  l'unité d'adjonction $\eta_{(V,W)}$ est un isomorphisme d'après \ref{proclefbis}. Il en résulte que $\mathrm{H}^{*}{V}\otimes_{\mathrm{H}^{*}V/W}\hspace{1pt}\nu_{X}$ est un isomorphisme et donc que $\nu_{X}$ en est un aussi puisque $\mathrm{H}^{*}{V}$ est un $\mathrm{H}^{*}{V/W}$-module fidèlement plat.

\medskip
3) Soit maintenant $X$ un $V$-CW-complexe fini  arbitraire. Soit $\mathrm{Sk}_{m} X$ son $m$-ième squelette~; on montre que $\nu_{\mathrm{Sk}_{m}X}$ est un isomorphisme par récurrence sur l'entier $m$ grâce à la deuxième étape. Comme l'on a par hypothèse $X=\mathrm{Sk}_{m} X$ pour $m$ assez grand la démonstration de la proposition \ref{cohFix} est achevée.

\medskip
\begin{cor}\label{cohFixpaire} Si $(X,Y)$ est une paire de $V$-CW-complexes finis alors $\nu_{(X,Y)}$ est un isomorphisme.
\end{cor}

\vspace{1,5cm}
\textit{Démonstration de la proposition \ref{cohEFix}.}

\medskip
On étend la transformation naturelle $\mathrm{E}\nu_{X}$ en une transformation naturelle entre foncteurs définis sur la catégorie des paires de $V$-espaces
$$
\mathrm{E}\nu_{(X,Y)}:\mathrm{EFix}_{(V,W)}\mathrm{H}^{*}_{V}(X,Y)\to\mathrm{H}^{*}_{V}(X^{W},Y^{W})
$$
et on procède \textit{mutatis mutandis} comme précédemment. On obtient du même coup l'énoncé suivant~:

\medskip
\begin{cor}\label{cohEFixpaire} Si $(X,Y)$ est une paire de $V$-CW-complexes finis alors $\mathrm{E}\nu_{(X,Y)}$ est un isomorphisme.
\end{cor}

\pagebreak

\sect{Rappels sur les injectifs des catégories $V\text{-}\mathcal{U}$ et  $V_{\mathrm{tf}}\text{-}\mathcal{U}$}

Rappelons tout d'abord que les notations $V\text{-}\mathcal{U}$ et  $V_{\mathrm{tf}}\text{-}\mathcal{U}$ désignent respectivement la catégorie des $\mathrm{H}^{*}V$-$\mathrm{A}$-modules instables et sa sous-catégorie pleine dont les objets sont les $\mathrm{H}^{*}V$-$\mathrm{A}$-modules instables qui sont de type fini comme $\mathrm{H}^{*}V$-module. Cette dernière hypothèse intervient très souvent dans ce mé\-moire, aussi nous abrégerons souvent ``$\hspace{1pt}\mathrm{H}^{*}V$-$\mathrm{A}$-module instable de type fini comme $\mathrm{H}^{*}V$-module'' en {\em ``$\hspace{2pt}\mathrm{H}^{*}V_{\mathrm{tf}}$-$\mathrm{A}$-module instable''}. Les catégories $V\text{-}\mathcal{U}$ et  $V_{\mathrm{tf}}\text{-}\mathcal{U}$ sont toutes deux des catégories abéliennes qui ont assez d'injectifs. C'est formel dans le cas de $V\text{-}\mathcal{U}$ (voir ci-dessous)~; cela l'est moins dans le cas de $V_{\mathrm{tf}}\text{-}\mathcal{U}$ (pour une généralisation voir \cite[Theorem 0.1]{Hennalg}). La référence principale pour cette section est \cite{LZsmith}~; nous pourrions comme dans la section précédente supposer que $V$ est un $\ell$-groupe abélien élémentaire avec $\ell$ un nombre premier arbitraire mais par paresse nous revenons comme dans l'introduction à $\ell=2$.

\bigskip
On commence par exhiber quelques injectifs de ces catégories.

\medskip
1) On note $\mathrm{J}_{V}(k)$, $k\in\mathbb{N}$,  le $\mathrm{H}^{*}V$-$\mathrm{A}$-module instable caractérisé, à isomorphisme près, par l'isomorphisme fonctoriel en le $\mathrm{H}^{*}V$-$\mathrm{A}$-module instable $M$~:
$$
\mathrm{Hom}_{V\text{-}\mathcal{U}}(M,\mathrm{J}_{V}(k))
\hspace{4pt}\cong\hspace{4pt}
\mathrm{Hom}_{\mathbb{F}_{2}}(M^{k},\mathbb{F}_{2})\hspace{4pt}:=\hspace{4pt}
\mathrm{H}_{k}M
$$
($M^{k}$ désigne ci-dessus le $\mathrm{F}_{2}$-espace vectoriel constitué des éléments de degré $k$ de $M$). Il est clair que $\mathrm{J}_{V}(k)$ est un injectif de $V\text{-}\mathcal{U}$. On constate qu'il est fini et donc \textit{a fortiori} de type fini comme $\mathrm{H}^{*}V$-module. L'homomorphisme canonique $M\to\prod_{k\in\mathbb{N}}\prod_{u\in\mathrm{H}_{k}M}\mathrm{J}_{V}(k)$ est par construction injectif~; il en résulte que la catégorie $V\text{-}\mathcal{U}$ a assez d'injectifs.

\medskip
2) Soient $W$ un sous-groupe de $V$ et $k\geq 0$ un entier. L'isomorphisme, 
$$
\mathrm{Hom}_{V\text{-}\mathcal{U}}(M,\mathrm{H}^{*}V\otimes_{\mathrm{H}^{*}V/W}\mathrm{J}_{V/W}(k))
\hspace{4pt}\cong\hspace{4pt}
\mathrm{Hom}_{V/W\text{-}\mathcal{U}}(\mathrm{Fix}_{(V,W)}M,\mathrm{J}_{V/W}(k))
$$
et l'exactitude du foncteur $\mathrm{Fix}_{(V,W)}$ montrent que $\mathrm{H}^{*}V\otimes_{\mathrm{H}^{*}V/W}\mathrm{J}_{V/W}(k)$ est un injectif de $V\text{-}\mathcal{U}$. On constate à nouveau qu'il est de type fini comme $\mathrm{H}^{*}V$-module.

\medskip
Ce qui précède se généralise~; l'exactitude du foncteur $\mathrm{Fix}_{(V,W)}$ est équivalente à l'énoncé suivant~:

\begin{pro}\label{exactFixbis} Si $I$ est un $V/W\text{-}\mathcal{U}$-injectif alors $\mathrm{H}^{*}V\otimes_{\mathrm{H}^{*}V/W}I$ est un $V\text{-}\mathcal{U}$-injectif.
\end{pro}

\medskip
3) Soient $E$ un $2$-groupe abélien élémentaire et $I$ un injectif de $V\text{-}\mathcal{U}$. On montre que $\mathrm{H}^{*}E\otimes I$ est un injectif de $V\text{-}\mathcal{U}$ (comme dans les énoncés \ref{suspensionEFix} et \ref{suspensionFix}, la structure de $\mathrm{H}^{*}V$-$\mathrm{A}$-module instable de $\mathrm{H}^{*}E\otimes I$ provient ``exclusivement'' de celle de $I$). Evidemment, si $L$ est un facteur direct, comme $\mathrm{A}$-module (instable), de $\mathrm{H}^{*}E$, alors $L\otimes I$ est encore un injectif de $V\text{-}\mathcal{U}$~; $L\otimes I$ est de type fini comme $\mathrm{H}^{*}V$-module  si et seulement si $L$ est isomorphe à $(\mathbb{F}_{2})^{m}$ avec $m\in\mathbb{N}$ et  si $I$ est de type fini comme $\mathrm{H}^{*}V$-module.

\bigskip
On énonce maintenant le théorème de classification des $\mathrm{H}^{*}V$-$\mathrm{A}$-modules instables injectifs  \cite[Théorème 0.12]{LZsmith}. Cet énoncé nécessite l'introduction de deux notations~:

\smallskip
--\hspace{8pt}$\mathcal{L}$ désigne un système de représentants pour les classes de $\mathcal{U}$-isomorphismes des facteurs directs indécomposables de $\mathrm{H}^{*}(\mathbb{Z}/2\mathbb{Z})^{d}$, $d$ parcourant $\mathbb{N}$ (chacune de ces classes est donc représentée dans $\mathcal{L}$ une et une seule fois)~;

\smallskip
--\hspace{8pt}$\mathcal{W}$ désigne l'ensemble des sous-groupes de $V$.

\begin{theo}\label{injectifs1} Soit $I$ un $V\text{-}\mathcal{U}$-injectif. Alors il existe une unique famille de cardinaux $(a_{L,W,k})_{(L,W,k)\in\mathcal{L}\times\mathcal{W}\times\mathbb{N}}$ telle que $I$ est isomorphe à la somme directe
$$
\hspace{24pt}
\bigoplus_{(L,W,k)\in\mathcal{L}\times\mathcal{W}\times\mathbb{N}}
{(\hspace{1pt}L\otimes (\mathrm{H}^{*}V\otimes_{\mathrm{H}^{*}V/W}\mathrm{J}_{V/W}(k))\hspace{1pt})}^{\oplus\hspace{1pt}a_{L,W,k}}
\hspace{24pt}.
$$
(Réciproquement tout $\mathrm{H}^{*}V$-$\mathrm{A}$-module instable de cette forme est un $V\text{-}\mathcal{U}$-injectif.)

\end{theo}
 
On montre dans \cite{LZsmith} que le théorème ci-dessus a pour sous-produit l'énoncé suivant (voir \cite[Théorème 0.13]{LZsmith})~:

\medskip
\begin{theo}\label{LZsmith0.13} Soient $M$ un $\mathrm{H}^{*}V$-$\mathrm{A}$-module instable et $i:M\to E$ une enveloppe injective de $M$ dans la catégorie $V\text{-}\mathcal{U}$. Si $M$ est de type fini comme $\mathrm{H}^{*}V$-module alors $E$ est isomorphe à une somme directe finie de $V\text{-}\mathcal{U}$-injectifs de la forme $\mathrm{H}^{*}V\otimes_{\mathrm{H}^{*}V/W}\mathrm{J}_{V/W}(k)$. En particulier $E$ est aussi de type fini comme $\mathrm{H}^{*}V$-module.
\end{theo}

\bigskip
En observant qu'un $V\text{-}\mathcal{U}$-injectif qui est de type fini comme $\mathrm{H}^{*}V$-module est un $V_{\mathrm{tf}}\text{-}\mathcal{U}$-injectif ($V_{\mathrm{tf}}\text{-}\mathcal{U}$-injectif est une abréviation pour injectif de  la catégorie $V_{\mathrm{tf}}\text{-}\mathcal{U}$), on peut dédoubler (voire détripler) le théorème \ref{LZsmith0.13}~:

\medskip
\begin{theo}\label{injectifs2} La catégorie $V_{\mathrm{tf}}\text{-}\mathcal{U}$ a assez d'injectifs.
\end{theo}

\pagebreak

\begin{theo}\label{injectifs3} Soit $I$ un $V_{\mathrm{tf}}\text{-}\mathcal{U}$-injectif. Alors il existe une unique application $a:\mathcal{W}\times\mathbb{N}\to\mathbb{N}$, à support fini, en clair avec $a(W,k)=0$ en dehors d'un sous-ensemble fini de $\mathcal{W}\times\mathbb{N}$, telle que $I$ est isomorphe à la somme directe
$$
\hspace{24pt}
\bigoplus_{(W,k)\in\mathcal{W}\times\mathbb{N}}
{(\hspace{1pt}\mathrm{H}^{*}V\otimes_{\mathrm{H}^{*}V/W}\mathrm{J}_{V/W}(k)\hspace{1pt})}^{\hspace{1pt}a(W,k)}
\hspace{24pt}.
$$
(Réciproquement tout $\mathrm{H}^{*}V$-$\mathrm{A}$-module instable de cette forme est un $V_{\mathrm{tf}}\text{-}\mathcal{U}$-injectif.)
\end{theo}

\medskip
\textit{Démonstration.} La partie ``unicité'' est conséquence de celle de \ref{injectifs1}.
\hfill$\square$

\begin{scho}\label{smith0.13bis} Un $V_{\mathrm{tf}}\text{-}\mathcal{U}$-injectif est aussi un $V\text{-}\mathcal{U}$-injectif. Soient $M$ un objet de $V_{\mathrm{tf}}\text{-}\mathcal{U}$ et $i:M\to E$ un morphisme de $V\text{-}\mathcal{U}$ alors les deux propriétés suivantes sont équivalentes~:
\begin{itemize}
\item[(i)] $i$ est une enveloppe injective dans la catégorie $V\text{-}\mathcal{U}$~;
\item[(ii)] $i$ est une enveloppe injective dans la catégorie $V_{\mathrm{tf}}\text{-}\mathcal{U}$.
\end{itemize}
\end{scho}

\begin{pro-def}\label{definition-de-a} L'application à support fini $a:\mathcal{W}\times\mathbb{N}\to\mathbb{N}$ qui apparaît dans \ref{injectifs3} ne dépend que de la classe d'isomorphisme de $I$~; on la note $\mathbf{a}_{I}$. Soient $M$ un objet de $V_{\mathrm{tf}}\text{-}\mathcal{U}$ et $i:M\to E$ une enveloppe injective dans cette catégorie, alors l'application $\mathbf{a}_{E}$ ne dépend que de la classe d'isomorphisme de $M$~; on la note $\mathbf{a}_{M}$. Les conditions $M=0$ et $\mathbf{a}_{M}=0$ sont équivalentes.

\end{pro-def}

\textit{Démonstration.} La première partie de la proposition résulte de l'unicité qui figure dans l'énoncé \ref{injectifs3}, la deuxième de ``l'unicité'' de l'enveloppe injective, la troisième est triviale.
\hfill$\square$

\medskip
\begin{pro}\label{Fixinjectifs} Soit $W$un sous-groupe de $V$.

\medskip
{\em (a)} Soit $I$ un $V\text{-}\mathcal{U}$-injectif. Alors $\mathrm{EFix}_{(V,W)}\hspace{1pt}I$ (resp. $\mathrm{Fix}_{(V,W)}\hspace{1pt}I$) est un $V\text{-}\mathcal{U}$-injectif (resp. $V/W\text{-}\mathcal{U}$-injectif).

\medskip
{\em (b)} Soit $I$ un $V_{\mathrm{tf}}\text{-}\mathcal{U}$-injectif. Alors $\mathrm{EFix}_{(V,W)}\hspace{1pt}I$ (resp. $\mathrm{Fix}_{(V,W)}\hspace{1pt}I$) est un $V_{\mathrm{tf}}\text{-}\mathcal{U}$-injectif (resp. $V/W_{\mathrm{tf}}\text{-}\mathcal{U}$-injectif).

\end{pro}

\textit{Démonstration.} Pour le (a) utiliser \ref{injectifs1},  \ref{suspensionEFix} (resp. \ref{suspensionFix}),  l'isomorphisme $\mathrm{T}_{W}\mathrm{H}^{*}E\cong(\mathrm{H}^{*}E)^{\hspace{1pt}\mathrm{Hom}(W,E)}$, et \ref{proclef} (resp. \ref{corproclef}). Pour le (b) utiliser \ref{injectifs3} et \ref{proclef} (resp. \ref{corproclef}).
\hfill$\square$

\bigskip
Comme la catégorie $V_{\mathrm{tf}}\text{-}\mathcal{U}$ a assez d'injectifs et que les algèbres $\mathrm{H}^{*}V$ et  $\mathrm{H}^{*}V/W$ sont noethériennes,  le point (b) de \ref{Fixinjectifs} implique l'énoncé suivant (pour une preuve alternative dans le cas de $\mathrm{Fix}_{V}=\mathrm{Fix}_{(V,V)}$ voir \cite[Lemme 2.4.2]{LZsmith}, pour une vaste généralisation, avec une approche très différente, voir \cite[Theorem 1.4]{DWsmith2})~:

\begin{pro}\label{typefiniFix} Soit $M$ un $\mathrm{H}^{*}V$-$\mathrm{A}$-module instable. Si $M$ est de type fini comme $\mathrm{H}^{*}V$-module alors $\mathrm{EFix}_{(V,W)}M$ (resp. $\mathrm{Fix}_{(V,W)}M$) est de type fini comme $\mathrm{H}^{*}V$-module (resp.  $\mathrm{H}^{*}V/W$-module).
\end{pro}

\medskip
Dans le cas particulier $W=V$ on obtient (voir à nouveau  \cite[Lemme 2.4.2]{LZsmith})~:

\begin{cor} Soit $M$ un $\mathrm{H}^{*}V$-$\mathrm{A}$-module instable. Si $M$ est de type fini comme $\mathrm{H}^{*}V$-module alors le $\mathrm{A}$-module instable $\mathrm{Fix}_{V}M$ est fini.
\end{cor}

\bigskip
La proposition \ref{typefiniFix} montre que le foncteur $\mathrm{EFix}_{(V,W)}:V\text{-}\mathcal{U}\to V\text{-}\mathcal{U}$ (resp. $\mathrm{Fix}_{(V,W)}:V\text{-}\mathcal{U}\to V/W\text{-}\mathcal{U}$) induit un foncteur $\mathrm{EFix}_{(V,W)}:V_{\mathrm{tf}}\text{-}\mathcal{U}\to V_{\mathrm{tf}}\text{-}\mathcal{U}$ (resp. $\mathrm{Fix}_{(V,W)}:V_{\mathrm{tf}}\text{-}\mathcal{U}\to V/W_{\mathrm{tf}}\text{-}\mathcal{U}$).  Compte tenu de \ref{cohEFix}, le point (b) de la proposition-définition ci-après est le pendant algébrique de l'énoncé trivial suivant~: soient $X$ un $V$-CW-complexe fini et $W_{1}$, $W_{2}$ deux sous-groupes de~$V$, on a ${(X^{W_{2}})}^{W_{1}}=X^{W_{1}+W_{2}}$.

\bigskip
\begin{pro-def}\label{compoEFix} Soient $W_{1}$ et $W_{2}$ deux sous-groupes de $V$.

\medskip
{\em (a)} L'isomorphisme fonctoriel $\mathrm{T}_{W_{1}}\circ\mathrm{T}_{W_{2}}\cong\mathrm{T}_{W_{1}\oplus W_{2}}$ et le morphisme fonctoriel $\mathrm{T}_{W_{1}\oplus W_{2}}\to\mathrm{T}_{W_{1}+W_{2}}$ induisent un morphisme fonctoriel
$$
\mu:\mathrm{EFix}_{(V,W_{1})}\circ\mathrm{EFix}_{(V,W_{2})}\to\mathrm{EFix}_{(V,W_{1}+W_{2})}
$$
tel que $\mu_{M}$ est un épimorphisme pour tout $\mathrm{H}^{*}V$-$\mathrm{A}$-module instable.

\medskip
{\em (b)} Si $M$ est un $\mathrm{H}^{*}V_{\mathrm{tf}}$-$\mathrm{A}$-module instable alors $\mu_{M}$ est un isomorphisme~; en d'autres termes $\mu$ induit  isomorphisme fonctoriel
$$
\mathrm{EFix}_{(V,W_{1})}\circ\mathrm{EFix}_{(V,W_{2})}
\hspace{4pt}\cong\hspace{4pt}
\mathrm{EFix}_{(V,W_{1}+W_{2})}
$$
entre endofoncteurs de $V_{\mathrm{tf}}\text{-}\mathcal{U}$.

\end{pro-def}

\textit{Démonstration du (a).} On montre l'existence d'un morphisme fonctoriel canonique $\mu:\mathrm{EFix}_{(V,W_{1})}\circ\mathrm{EFix}_{(V,W_{2})}\to\mathrm{EFix}_{(V,W_{1}+W_{2})}$ tel que $\mu_{M}$ est un épimorphisme pour tout $\mathrm{H}^{*}V$-$\mathrm{A}$-module instable~; on laisse au lecteur le soin de vérifier qu'il est bien induit par l'isomorphisme fonctoriel $\mathrm{T}_{W_{1}}\circ\mathrm{T}_{W_{2}}\cong\mathrm{T}_{W_{1}\oplus W_{2}}$ et le morphisme fonctoriel $\mathrm{T}_{W_{1}\oplus W_{2}}\to\mathrm{T}_{W_{1}+W_{2}}$.

\medskip
La proposition \ref{adjonction} dit que le foncteur $\mathrm{EFix}_{(V,W)}:V\text{-}\mathcal{U}\to V\text{-}\mathcal{U}$ est l'adjoint à gauche du foncteur $\mathrm{E}_{(V,W)}:V\text{-}\mathcal{U}\to V\text{-}\mathcal{U},N\mapsto\mathrm{H}^{*}V\otimes_{\mathrm{H}^{*}V/W}N$. Le scholie \ref{observation} dit que l'on a un isomorphisme naturel de $\mathrm{H}^{*}V$-$\mathrm{A}$-modules instables $\mathrm{E}_{(V,W)}N\cong\mathrm{H}^{*}W\otimes N$ le membre de gauche étant muni de la structure de $\mathrm{H}^{*}V$-module induite par l'homomorphisme $W\oplus V\to V,(x,y)\mapsto x+y$. Il~en résulte que le $\mathrm{H}^{*}V$-$\mathrm{A}$-module instable $(\mathrm{E}_{(V,W_{2})}\circ\mathrm{E}_{(V,W_{1})})(N)$  est naturellement isomorphe à $\mathrm{H}^{*}(W_{2}\oplus W_{1})\otimes N$, cet $\mathrm{A}$-module instable étant muni de la structure de $\mathrm{H}^{*}V$-module induite par l'homomorphisme $W_{2}\oplus W_{1}\oplus V\to V,\linebreak (x_{2},x_{1},y)\mapsto x_{2}+x_{1}+y$. En considérant le monomorphisme évident\linebreak $\mathrm{H}^{*}(W_{2}+W_{1})\to\mathrm{H}^{*}(W_{2}\oplus W_{1})$ on constate au bout du compte que l'on dispose d'un monomorphisme naturel $\widetilde{\mu}_{N}:\mathrm{E}_{(V,W_{1}+W_{2})}N\hookrightarrow(\mathrm{E}_{(V,W_{2})}\circ\mathrm{E}_{(V,W_{1})})(N)$. Par~adjonction on obtient un épimorphisme naturel
$$
\hspace{24pt}
\mu_{M}:(\mathrm{EFix}_{(V,W_{1})}\circ\mathrm{EFix}_{(V,W_{2})})(M)\twoheadrightarrow\mathrm{EFix}_{(V,W_{1}+W_{2})}M
\hspace{24pt}.
$$

\medskip
\textit{Démonstration du (b).} On prend tout d'abord $M=\mathrm{H}^{*}V\otimes_{\mathrm{H}^{*}V/W}N$ avec $N$ un $\mathrm{H}^{*}V/W$-$\mathrm{A}$-module instable fini. La proposition \ref{proclef} entraîne que l'on a dans ce cas un $(V_{\mathrm{tf}}\text{-}\mathcal{U})$-isomorphisme $(\mathrm{EFix}_{(V,W_{1})}\circ\mathrm{EFix}_{(V,W_{2})})(M)\cong\mathrm{E}_{(V,W_{1}+W_{2})}M$~; $\mu_{M}$ est nécessairement un $(V_{\mathrm{tf}}\text{-}\mathcal{U})$-isomorphisme, en effet, en un degré donné, $\mu_{M}$ est une surjection entre deux $\mathbb{F}_{2}$-espaces vectoriels de même dimension finie. Il en résulte compte tenu de \ref{injectifs3} que $\mu_{I}$ est un isomorphisme pour $I$ un $(V_{\mathrm{tf}}\text{-}\mathcal{U})$-injectif. On obtient le cas général en considérant le début $0\to M\to I^{0}\to I^{1}$ d'une résolution injective de $M$ dans la catégorie $V_{\mathrm{tf}}\text{-}\mathcal{U}$.
\hfill$\square$

\bigskip
L'énoncé ci-dessous précise le point (b) de \ref{Fixinjectifs} \cite[Proposition 3.3.3]{LZtorext}~:

\begin{pro}\label{enveloppeFix}  Les deux foncteurs $\mathrm{EFix}_{(V,W)}:V_{\mathrm{tf}}\text{-}\mathcal{U}\to V_{\mathrm{tf}}\text{-}\mathcal{U}$ et  $\mathrm{Fix}_{(V,W)}:V_{\mathrm{tf}}\text{-}\mathcal{U}\to V/W_{\mathrm{tf}}\text{-}\mathcal{U}$ préservent les enveloppes injectives.

\end{pro}

\textit{Démonstration.} On suit celle de \cite{LZtorext}\footnote{
Signalons incidemment un misprint dans la version publiée de cette démonstration~: dans le diagramme qui y apparaît, les flèches verticales pointent dans la mauvaise direction.}. Soient $M$ un $\mathrm{H}^{*}V_{\mathrm{tf}}$-$\mathrm{A}$-module instable et $i:M\to E$ une enveloppe injective de $M$ dans la catégorie $V_{\mathrm{tf}}\text{-}\mathcal{U}$.

\medskip
1) Le cas du foncteur $\mathrm{EFix}_{(V,W)}$. On note $\eta_{M}$ la transformation naturelle $M\to\mathrm{H}^{*}V\otimes_{\mathrm{H}^{*}V/W}\mathrm{Fix}_{(V,W)}M=\mathrm{EFix}_{(V,W)}M$ et on identifie $i$ et $\mathrm{EFix}_{(V,W)}(i)$ avec des inclusions $M\subset E$ et $\mathrm{EFix}_{(V,W)}M\subset\mathrm{EFix}_{(V,W)}E$. La proposition \ref{proclefbis} montre que le $\mathrm{H}^{*}V$-$\mathrm{A}$-module instable $E$ se décompose naturellement en une somme directe $\mathrm{EFix}_{(V,W)}E\oplus\ker\eta_{E}$ et l'on a $\ker\eta_{M}=M\cap\ker\eta_{E}$. Soit maintenant $P$ un sous-module de $\mathrm{EFix}_{(V,W)}E$ avec $P\cap\mathrm{EFix}_{(V,W)}M=0$~; on~constate que $P$ s'identifie avec un sous-module de $E$ vérifiant $P\cap M=0$. On a donc $P=0$.
\hfill$\square$

\medskip
2) Le cas du foncteur $\mathrm{Fix}_{(V,W)}$. Soit $Q$ un sous-module de $\mathrm{Fix}_{(V,W)}E$ avec $Q\cap\mathrm{Fix}_{(V,W)}M=0$, alors le produit tensoriel $\mathrm{H}^{*}V\otimes_{\mathrm{H}^{*}V/W}Q$ s'identifie à un sous-module de  $\mathrm{EFix}_{(V,W)}E$ dont l'intersection avec $\mathrm{EFix}_{(V,W)}E$ est triviale, on a donc $\mathrm{H}^{*}V\otimes_{\mathrm{H}^{*}V/W}Q=0$ d'après le 1). Ceci implique $Q=0$ puisque $\mathrm{H}^{*}V$ est un $\mathrm{H}^{*}V/W$-module fidèlement plat.
\hfill$\square$

\bigskip
Soit $M$ un $\mathrm{H}^{*}V_{\mathrm{tf}}$-$\mathrm{A}$-module instable~; les propositions \ref{proclef}, \ref{corproclef} et \ref{enveloppeFix} permettent d'exprimer l'invariant $\mathbf{a}$ (défini\-tion \ref{definition-de-a}) de $\mathrm{EFix}_{(V,W)}M$ et $\mathrm{Fix}_{(V,W)}M$ en fonction de celui de $M$~:

\begin{scho}\label{aFix}  Soit $M$ un $\mathrm{H}^{*}V_{\mathrm{tf}}$-$\mathrm{A}$-module instable.

\medskip
{\em (a)} Soient $U$ un sous-groupe de $V$ et $k$ un entier naturel~; on a
$$
\mathbf{a}_{\mathrm{EFix}_{(V,W)}M}(U,k)
\hspace{4pt}=\hspace{4pt}
\begin{cases}
\mathbf{a}_{M}(U,k)
& \text{pour $W\subset U$,} \\
0
& \hspace{1pt}\text{pour $W\not\subset U$.}
\end{cases}
$$

\medskip
{\em (b)} Soient $U$ un sous-groupe de $V/W$ et $k$ un entier naturel~; on a
$$
\hspace{24pt}
\mathbf{a}_{\mathrm{Fix}_{(V,W)}M}(U,k)
\hspace{4pt}=\hspace{4pt}
\mathbf{a}_{M}(q^{-1}(U),k)
\hspace{24pt},
$$
$q$ désignant la surjection  canonique de $V$ dans $V/W$.
\end{scho}

\medskip
Le scholie ci-dessus conduit à l'énoncé suivant (qui est implicite dans la démonstration de \cite[Lemme 3.2.3]{LZtorext}, pour une généralisation voir \cite[Corollary 2.13]{Hennalg})~:

\begin{pro}\label{caracterisationfini}  Soit $M$ un $\mathrm{H}^{*}V_{\mathrm{tf}}$-$\mathrm{A}$-module instable. Les conditions suivantes sont équivalentes~:
\begin{itemize}
\item[(i)] $\mathbf{a}_{M}(U,k)=0$ pour tout $(U,k)$ avec $U\not=0$~;
\item[(ii)] $\mathrm{EFix}_{(V,W)}M=0$ pour tout $W\not=0$~;
\item[(ii-bis)] $\mathrm{EFix}_{(V,W)}M=0$ pour tout $W$ avec $\dim W=1$~;
\item[(iii)] $\mathrm{Fix}_{(V,W)}M=0$ pour tout $W\not=0$~;
\item[(iii-bis)] $\mathrm{Fix}_{(V,W)}M=0$ pour tout $W$ avec $\dim W=1$~;
\item[(iv)] $M$ est fini.
\end{itemize}
\end{pro}

\medskip
\textit{Démonstration.} D'après \ref{aFix} la condition (i) est équivalente à chacune des conditions (ii), (ii-bis), (iii) et (iii-bis). Si (i) est satisfaite alors $M$ s'injecte dans $\bigoplus_{k\in\mathbb{N}}(\mathrm{J}_{V}(k))^{\mathbf{a}_{M}(0,k)}$ qui est fini~; on a donc $(i)\Rightarrow (iv)$. On montre facilement  qu'un $\mathrm{H}^{*}V$-$\mathrm{A}$-module instable fini s'injecte dans une somme directe finie de $\mathrm{J}_{V}(k)$ (voir la proposition \ref{resolutionfini} ci-après)~; on en déduit $(iv)\Rightarrow(i)$.
\hfill$\square$

\bigskip
\begin{rem}\label{implicationelementaire} L'implication (iv)$\Rightarrow$(iii) de la proposition \ref{caracterisationfini} est en fait très élémentaire~:

\smallskip
Soit $M$ un $\mathrm{H}^{*}V$-$\mathrm{A}$-module instable fini et $N$ un $\mathrm{H}^{*}V/W$-$\mathrm{A}$-module instable. On a $\mathrm{Hom}_{\mathrm{H}^{*}V}(M,\mathrm{H}^{*}V\otimes_{\mathrm{H}^{*}V/W}N)=0$ (homomorphismes dans la catégorie des $\mathrm{H}^{*}V$-modules gradués). En effet il est facile de se convaincre que si $W$ est non nul alors tout sous-$\mathrm{H}^{*}V$-module fini de $\mathrm{H}^{*}V\otimes_{\mathrm{H}^{*}V/W}N$ est trivial (voir Lemme \ref{pre-annPf}). On a \textit{a fortiori}  $\mathrm{Hom}_{\mathrm{V}\text{-}\mathcal{U}}(M,\mathrm{H}^{*}V\otimes_{\mathrm{H}^{*}V/W}N)=0$ et donc $\mathrm{Fix}_{(V,W)}M=0$.

\medskip
Compte tenu de \ref{EFix} la même remarque vaut pour l'implication (iv)$\Rightarrow$(ii).
\end{rem}

\bigskip
Soit $M$ un $\mathrm{H}^{*}V$-$\mathrm{A}$-module instable, on note $\Vert M\Vert$ l'élément de $\mathbb{N}\cup\{+\infty\}$ défini par $\Vert M\Vert:=\sup\{k;M^{k}\not=0\}$ (on a donc  $\Vert 0\Vert=0$).

\begin{pro}\label{resolutionfini} Soit $M$ un $\mathrm{H}^{*}V$-$\mathrm{A}$-module instable fini~; $M$ admet une résolution injective dans la catégorie $V_{\mathrm{tf}}\text{-}\mathcal{U}$ dont chaque terme est une somme finie de $\mathrm{J}_{V}(k)$ et dont la longueur est inférieure ou égale à $\Vert M\Vert$.
\end{pro}

\medskip
\textit{Démonstration.} Celle-ci généralise celle de \cite[Proposition 6.1.3]{LZens}. Le cas  $\Vert M\Vert=0$ est trivial. On suppose $\Vert M\Vert\geq 1$ et on considère l'homomorphisme canonique
$$
\hspace{24pt}
i\hspace{2pt}:\hspace{2pt}M\longrightarrow\bigoplus_{k=0}^{\Vert M\Vert}\hspace{4pt}M^{k}\otimes\mathrm{J}_{V}(k)
\hspace{24pt}.
$$
Par construction $i$ est injectif. De plus $i$ est un isomorphisme en degré  $\Vert M\Vert$ si bien que l'on a $\Vert\mathop{\mathrm{coker}}i\Vert\leq\Vert M\Vert-1$. On conclut par récurrence.
\hfill$\square$

\bigskip
On achève cette section par deux commentaires~:

\medskip
\begin{com}\label{DoSa} La théorie des  $V_{\mathrm{tf}}\text{-}\mathcal{U}$-injectifs est  un outil important dans la preuve, par le premier et quatrième auteur \cite{BZ}, de la conjecture de Landweber et Stong \cite{LSt} sur la profondeur d'un anneau d'invariants $\mathbb{F}_{p}[X_{1},X_{2},\ldots,X_{n}]^{G}$, $G$ désignant un sous-groupe de $\mathrm{GL}_{n}(\mathbb{F}_{p})$. Nous verrons également en \ref{Serre1} que cette théorie fournit une démonstration (assez détournée~!) du résultat de Serre concernant les idéaux homogènes de $\mathrm{H}^{*}V$ stables sous l'action de l'al\-gèbre de Steenrod.
\end{com}

\medskip
\begin{com}\label{proj} Comme nous l'avons rappelé en préambule à la démon\-stration de \ref{tensFix} le fait que $V\text{-}\mathcal{U}$ a assez de projectifs est formel. La situation est différente dans le cas de $V_{\mathrm{tf}}\text{-}\mathcal{U}$~: on peut se convaincre de ce qu'un objet $P$ de $V_{\mathrm{tf}}\text{-}\mathcal{U}$ est projectif si et seulement l'on a $P\simeq(\mathrm{H}^{*}V)^{\oplus m}$ avec $m\in\mathbb{N}$. Il en résulte que  $V_{\mathrm{tf}}\text{-}\mathcal{U}$ n'a pas assez de projectifs, en effet on a $\mathrm{Hom}_{V_{\mathrm{tf}}\text{-}\mathcal{U}}(\mathrm{H}^{*}V,M)\cong M^{0}$ si bien qu'un objet $M$ de $V_{\mathrm{tf}}\text{-}\mathcal{U}$, avec $M^{0}=0$, qui est quotient d'un projectif, est nul.
\end{com}

\sect{Sur les dérivés du foncteur ``partie finie''}

Soit $M$ un $\mathrm{H}^{*}V_{\mathrm{tf}}$-$\mathrm{A}$-module instable. Il n'est pas difficile de se convaincre que l'ensemble des sous-$\mathrm{H}^{*}V$-$\mathrm{A}$-modules finis de $M$ (ordonné par inclusion) possède un plus grand élément (voir \ref{defPf})~; nous l'appelons la {\em partie finie} de~$M$. L'endofoncteur de $V_{\mathrm{tf}}\text{-}\mathcal{U}$ qui associe à $M$ sa partie finie est un foncteur (additif) exact à gauche~; l'étude de ses foncteurs dérivés à droite est le principal objet de cette section.

\medskip
Cette théorie est réminiscente de celle de la cohomologie locale dont nous rappelons la définition ci-après. 

\medskip
Soient $R$ un anneau (commutatif unitaire) noethérien, $\mathfrak{a}\subset R$ un idéal et $M$ un $R$-module. On pose
$$
\hspace{24pt}
\Gamma_{\mathfrak{a}}(M)
\hspace{4pt}:=\hspace{4pt}\bigcup_{t\in\mathbb{N}}\hspace{4pt}
\mathrm{Ann}_{M}(\mathfrak{a}^{t})
\hspace{24pt},
$$
$\mathrm{Ann}_{M}(\mathfrak{a}^{t})$ désignant le sous-module de $M$ constitué des éléments annulés par la puissance $t$-ième de l'idéal $\mathfrak{a}$~;  $\Gamma_{\mathfrak{a}}(M)$ est souvent appelé la {\em $\mathfrak{a}$-torsion} de $M$. Deux observations~:

\smallskip
-- Si $M$ est de type fini alors la suite croissante de sous-modules $(\mathrm{Ann}_{M}(\mathfrak{a}^{t}))_{t\in\mathbb{N}}$ est stationnaire si bien que l'on a $\Gamma_{\mathfrak{a}}(M)=\mathrm{Ann}_{M}(\mathfrak{a}^{t_{0}})$ pour un certain $t_{0}$.

\smallskip
-- On a $\Gamma_{\mathfrak{a}}(M)=\Gamma_{\sqrt{\mathfrak{a}}}\hspace{2pt}(M)$, $\sqrt{\mathfrak{a}}$ désignant le radical de $\mathfrak{a}$ (voir \ref{Serre1}).

\smallskip
Il est clair que le foncteur $\Gamma_{\mathfrak{a}}$ est exact à gauche~; l'étude de ses foncteurs dérivés à droite est la théorie de la cohomologie locale initiée par Alexander Grothendieck \cite{Gr}. On pourra trouver une exposition approfondie de cette théorie, de sa version $\mathbb{Z}$-graduée et de ses liens avec la géométrie algébrique dans \cite{BS}.

\bigskip
On prend maintenant pour $R$ l'anneau $\mathrm{H}^{*}V$.

\medskip
\begin{pro}\label{GammaSteenrod} Soit $M$ un $\mathrm{H}^{*}V$-$\mathrm{A}$-module instable et $\mathfrak{a}\subset\mathrm{H}^{*}V$ un idéal homogène et stable sous l'action de $\mathrm{A}$.

\medskip
{\em (a)} Le sous-$\mathrm{H}^{*}V$-module $\mathrm{Ann}_{M}(\mathfrak{a})$ de $M$ est stable sous l'action de $\mathrm{A}$.

\medskip
{\em (b)} Le sous-$\mathrm{H}^{*}V$-module $\Gamma_{\mathfrak{a}}(M)$ de $M$ est stable sous l'action de $\mathrm{A}$.
\end{pro}

\pagebreak

\textit{Démonstration.} On peut montrer $\mathrm{Sq}^{i}\hspace{1pt}\mathrm{Ann}_{M}(\mathfrak{a})\subset\mathrm{Ann}_{M}(\mathfrak{a})$ par récurrence sur l'entier $i$ en contemplant l'égalité $\mathrm{Sq}^{i}ax=\sum_{j+k=i}(\mathrm{Sq}^{j}a)(\mathrm{Sq}^{k}x)$, $a$ dans $\mathfrak{a}$ et $x$ dans $M$. Le (b) est conséquence du (a) puisque si $\mathfrak{a}$ est homogène et stable sous l'action de $\mathrm{A}$ alors il en est de même pour $\mathfrak{a}^{t}$.
\hfill$\square$

\bigskip
\begin{rem}\label{listeSerre} Si $\mathfrak{a}$ est un idéal de $\mathrm{H}^{*}V$ homogène et stable sous l'action de~$\mathrm{A}$ alors il en est de même pour $\sqrt{\mathfrak{a}}$ (voir \ref{racine}). Comme l'on a $\Gamma_{\mathfrak{a}}=\Gamma_{\sqrt{\mathfrak{a}}}$, on peut supposer lorsque l'on considère les foncteurs $\Gamma_{\mathfrak{a}}$ qui apparaissent dans \ref{GammaSteenrod}, que les idéaux $\mathfrak{a}$ sont radiciels ($\mathfrak{a}=\sqrt{\mathfrak{a}}$). La liste de ces idéaux a été déterminée par Serre (voir \ref{Serre2})~; elle est indexée par les parties $P$ de $\mathcal{W}$ telles que la relation d'ordre sur $P$, définie par l'inclusion, est l'égalité.
\end{rem}

\bigskip
On prend enfin pour $\mathfrak{a}$ l'idéal $\widetilde{\mathrm{H}}^{*}V$ (constitué des éléments de $\mathrm{H}^{*}V$ de degré strictement positifs).

\begin{cor-def}\label{defPf} Soit $M$ un $\mathrm{H}^{*}V$-$\mathrm{A}$-module instable. Si $M$ est de type fini comme $\mathrm{H}^{*}V$-module alors $\Gamma_{\widetilde{\mathrm{H}}^{*}V}(M)$ est le plus grand sous-$\mathrm{H}^{*}V$-$\mathrm{A}$-module fini de $M$~; nous l'appelons la {\em partie finie} de $M$. Nous la notons $\mathrm{Pf}(M)$ (ou $\mathrm{Pf}_{V}(M)$ quand la mention de $V$ nous semble utile voire nécessaire).
\end{cor-def}

\medskip
\textit{Démonstration.} Nous savons déjà que le sous-$\mathrm{H}^{*}V$-module $\Gamma_{\widetilde{\mathrm{H}}^{*}V}(M)$ de $M$ est stable sous l'action de $\mathrm{A}$. Si $M$ est de type fini comme $\mathrm{H}^{*}V$-module alors $\Gamma_{\widetilde{\mathrm{H}}^{*}V}(M)$ est fini. En effet il existe dans ce cas un entier naturel $t_{0}$ tel que $\Gamma_{\widetilde{\mathrm{H}}^{*}V}(M)$ est un $\mathrm{H}^{*}V/(\widetilde{\mathrm{H}}^{*}V)^{t_{0}}$-module de type fini et l'algèbre $\mathrm{H}^{*}V/(\widetilde{\mathrm{H}}^{*}V)^{t_{0}}$ est finie. De plus, si $N$ est un sous-$\mathrm{H}^{*}V$-$\mathrm{A}$-module fini de $M$ alors il existe un entier $t_{N}$ tel que $N$ est annulé par $(\widetilde{\mathrm{H}}^{*}V)^{t_{N}}$.
\hfill$\square$

\bigskip
Nous considérons dans ce mémoire les dérivés à droite du foncteur $\mathrm{Pf}$ dans la catégorie abélienne $V_{\mathrm{tf}}\text{-}\mathcal{U}$, et non dans la catégorie des $\mathrm{H}^{*}V$-modules, la graduation de $\mathrm{H}^{*}V$ étant oubliée, ou dans celle des $\mathrm{H}^{*}V$-modules $\mathbb{Z}$-gradués. Il existe cependant, comme le souligne Henn dans \cite{Hennalg}, une analogie formelle entre les énoncés de la présente section et certains des énoncés de la théorie de la cohomologie locale. Le contexte dans lequel travaille Henn est bien plus général que le notre~: $\mathrm{H}^{*}V$ est remplacée par une $\mathrm{A}$-algèbre instable $K$ qui est de type fini comme algèbre graduée, et la catégorie $V_{\mathrm{tf}}\text{-}\mathcal{U}$ par celle des $K$-$\mathrm{A}$-modules instables qui sont de type fini comme $K$-modules gradués. Nous recommandons au lecteur désireux d'approfondir la discussion précédente de consulter l'introduction de \cite{Hennalg}.

\bigskip
Nous allons relier les $\mathrm{R}^{k}\mathrm{Pf}M$ ($M$ objet de $V_{\mathrm{tf}}\text{-}\mathcal{U}$) au foncteur 
$W\mapsto\mathrm{EFix}_{(V,W)}M$ implicitement introduit dans la proposition-définition \ref{Wfonctoriel}.

\pagebreak

\begin{pro-def}\label{defPsi} On pose $\mathcal{W}_{0}:=\mathcal{W}-\{0\}$ (on rappelle que la notation $\mathcal{W}$ désigne l'ensemble des sous-groupes $W\subset V$)~; $\mathcal{W}$ (resp. $\mathcal{W}_{0}$) est ordonné par inclusion et peut donc être vu comme l'ensemble des objets d'une catégorie que l'on note encore $\mathcal{W}$ (resp. $\mathcal{W}_{0}$).

\medskip
On considère les deux catégories $\mathcal{W}$ et $V_{\mathrm{tf}}\text{-}\mathcal{U}$ et on se donne un $\mathrm{H}^{*}V_{\mathrm{tf}}$-$\mathrm{A}$-module instable $M$. Alors l'application entre objets
$$
\hspace{24pt}
W
\hspace{4pt}
\mapsto
\hspace{4pt}
\mathrm{EFix}_{(V,W)}M
\cong
\mathrm{H}^{*}V\otimes_{\mathrm{H}^{*}V/W}\mathrm{Fix}_{(V,W)}M
\hspace{24pt}
$$
se prolonge canoniquement en un foncteur de $\mathcal{W}$ dans $V_{\mathrm{tf}}\text{-}\mathcal{U}$. On note $\Psi_{M}$ la restriction de ce foncteur à $\mathcal{W}_{0}$. On note $\Psi:V_{\mathrm{tf}}\text{-}\mathcal{U}\to (V_{\mathrm{tf}}\text{-}\mathcal{U})^{\mathcal{W}_{0}}$ le foncteur $M\mapsto\Psi_{M}$ ($\Psi_{M}$ pourra donc aussi être noté $\Psi(M)$).
\end{pro-def}

\bigskip
\textit{Démonstration.} La ``partie proposition" de l'énoncé ci-dessus résulte de la proposition-définition \ref{Wfonctoriel}.
\hfill$\square$

\bigskip
En préalables respectifs aux énoncés \ref{derivePf} et \ref{annderivePf1}, nous rappelons ci-dessous quelques points qui nous seront utiles de la théorie des ``dérivés de la limite'' d'un foncteur à valeurs dans une catégorie abélienne.

\vspace{0.75cm}
\begin{raps}\label{rappelslim}

\medskip
(R.1) Soit $\mathcal{C}$ une petite catégorie finie (en clair nous supposons que l'ensemble des objets et l'ensemble des morphismes de $\mathcal{C}$ sont finis, hypothèse qui simplifiera notre exposition)~; soit $\mathcal{A}$ une catégorie abélienne avec assez d'injectifs. On note $\mathcal{A}^{\mathcal{C}}$ la catégorie des foncteurs définis sur $\mathcal{C}$ et à valeurs dans $\mathcal{A}$~; la catégorie $\mathcal{A}^{\mathcal{C}}$ est une catégorie abélienne avec assez d'injectifs. Le foncteur $\lim_{\mathcal{C}}:\mathcal{A}^{\mathcal{C}}\to\mathcal{A}$ est l'adjoint à droite du ``foncteur diagonal'' $\Delta:\mathcal{A}\to\mathcal{A}^{\mathcal{C}}$,  qui associe à un objet $a$ de $\mathcal{A}$ le foncteur $\Delta_{a}:\mathcal{C}\to\mathcal{A}$ envoyant tout $\mathcal{C}$-objet sur~$a$ et tout $\mathcal{C}$-morphisme sur $\mathrm{id}_{a}$. Le foncteur $\lim_{\mathcal{C}}$ est un foncteur additif exact à gauche~; son $k$-ième dérivé à droite est noté $\lim_{\mathcal{C}}^{k}$ (on a donc $\lim_{\mathcal{C}}^{0}=\lim_{\mathcal{C}}$).

\medskip
(R.2) Soit maintenant $\mathcal{C}$ un ensemble ordonné fini (que l'on peut voir comme une petite catégorie finie), nous rappelons l'expression très concrète que l'on a dans ce cas pour les $\lim_{\mathcal{C}}^{k}$.

\smallskip
Soit $F$ un foncteur de $\mathcal{C}$ dans $\mathcal{A}$~; on note $\mathrm{L}^{\bullet}(F)$ le $\mathcal{A}$-complexe de cochaînes défini de la façon suivante~:

\smallskip
-- Soit $k\geq 0$ un entier~; on note $\mathrm{S}_{k}\mathcal{C}$ l'ensemble des $k$-simplexes de $\mathcal{C}$ (un $k$-simplexe de $\mathcal{C}$ est un sous-ensemble $\sigma\subset\mathcal{C}$, totalement ordonné à $k+1$ éléments) et on pose
$$
\hspace{24pt}
\mathrm{L}^{k}(F)
\hspace{4pt}:=\hspace{4pt}
\prod_{\sigma\in\mathrm{S}_{k}\mathcal{C}}
F(\sup\sigma)
\hspace{24pt}.
$$
\footnotesize
On observera que l'ensemble $\mathrm{S}_{k}\mathcal{C}$ est fini si bien que le produit ci-dessus existe sans hypothèse supplémentaire sur la catégorie abélienne $\mathcal{A}$~; la raison d'être de cette observation est que nous allons prendre $\mathcal{A}=V_{\mathrm{tf}}\text{-}\mathcal{U}$ et que dans cette catégorie seuls les produits finis existent.
\normalsize

\smallskip
-- Le cobord $\mathrm{d}:\mathrm{L}^{k}(F)\to\mathrm{L}^{k+1}(F)$ est défini de la manière habituelle~:

\smallskip
On identifie le groupe abélien $\mathrm{Hom}_{\mathcal{A}}(\mathrm{L}^{k}(F),\mathrm{L}^{k+1}(F))$ au produit
$$
\hspace{24pt}
\prod_{(\sigma,\tau)\in\mathrm{S}_{k}\mathcal{C}\times\mathrm{S}_{k+1}\mathcal{C}}
\mathrm{Hom}_{\mathcal{A}}(F(\sup\sigma),F(\sup\tau))
\hspace{24pt};
$$
on a $\mathrm{d}_{(\sigma,\tau)}=0$ pour $\sigma\not\subset\tau$ et $\mathrm{d}_{(\sigma,\tau)}=(-1)^{\nu(\sigma,\tau)}F((\sup\sigma,\sup\tau))$ pour $\sigma\subset\tau$. Précisons la notation : $(\sup\sigma,\sup\tau)$ désigne l'unique $\mathcal{C}$-morphisme de $\sup\sigma$ dans $\sup\tau$ et $\nu(\sigma,\tau)$ le cardinal du sous-ensemble de $\sigma$ constitué des éléments $c$ vérifiant $c<\gamma(\sigma,\tau)$ avec $\{\gamma(\sigma,\tau)\}=\tau-\sigma$.

\medskip
Pour tout $k\geq 0$, on a un isomorphisme $\lim_{\mathcal{C}}^{k}F\cong\mathrm{H}^{k}\hspace{1pt}\mathrm{L}^{\bullet}(F)$, canonique et naturel en $F$. Pour $k=0$ ceci résulte de la définition même de $\lim_{\mathcal{C}}F$. Pour $k>0$ on peut s'en convaincre de la façon suivante~:

\smallskip
1) Soient $J$ un objet de $\mathcal{A}$ et $\gamma$ un élément de $\mathcal{C}$. On considère le foncteur  $J^{\mathrm{Hom}_{\mathcal{C}}(-,\gamma)}$ de $\mathcal{C}$ dans $\mathcal{A}$~; ce foncteur est déterminé, à isomorphisme près, par l'isomorphisme $\mathrm{Hom}_{\mathcal{A}^{\mathcal{C}}}(F,J^{\mathrm{Hom}_{\mathcal{C}}(-,\gamma)})\cong\mathrm{Hom}_{\mathcal{A}}(F(\gamma),J)$ naturel en $F$. On vérifie que l'on a $\mathrm{H}^{k}\hspace{1pt}\mathrm{L}^{\bullet}(J^{\mathrm{Hom}_{\mathcal{C}}(-,\gamma)})=0$ pour $k>0$.

\smallskip
\footnotesize
Pour ce faire on peut, par exemple, invoquer les arguments ci-après. Soient $\mathcal{C}_{\leq\gamma}$ (resp.~$\mathcal{C}_{<\gamma}$) le sous-ensemble ordonné de $\mathcal{C}$ constitué des éléments $c\leq\gamma$ (resp.  $c<\gamma$) et $\Vert\mathcal{C}_{\leq\gamma}\Vert$ (resp. $\Vert\mathcal{C}_{<\gamma}\Vert$) le polyèdre associé~; soit $K$ un injectif (arbitraire) de $\mathcal{A}$. On constate que le complexe $\mathrm{Hom}_{\mathcal{A}}(\mathrm{L}^{\bullet}(J^{\mathrm{Hom}_{\mathcal{C}}(-,\gamma)}),K)$ s'identifie au complexe de chaînes polyédrales de $\Vert\mathcal{C}_{\leq\gamma}\Vert$ à coefficients $\mathrm{Hom}_{\mathcal{A}}(J,K)$. On obtient $\mathrm{Hom}_{\mathcal{A}}(\mathrm{H}^{k}\hspace{1pt}\mathrm{L}^{\bullet}(J^{\mathrm{Hom}_{\mathcal{C}}(-,\gamma)}),K)=0$ pour $k>0$ en observant que $\Vert\mathcal{C}_{\leq\gamma}\Vert$ est isomorphe au cône de $\Vert\mathcal{C}_{<\gamma}\Vert$ (pour $k=0$ on obtient $\mathrm{Hom}_{\mathcal{A}}(\mathrm{H}^{k}\hspace{1pt}\mathrm{L}^{\bullet}(J^{\mathrm{Hom}_{\mathcal{C}}(-,\gamma)}),K)=\mathrm{Hom}_{\mathcal{A}}(J,K)$, ce qui est bien compatible avec l'isomorphisme canonique $\lim_{\mathcal{C}}J^{\mathrm{Hom}_{\mathcal{C}}(-,\gamma)})\cong J$, voir le point (b) de \ref{co-induit}).
\normalsize

\smallskip
On observera que l'isomorphisme $\mathrm{Hom}_{\mathcal{A}^{\mathcal{C}}}(F,J^{\mathrm{Hom}_{\mathcal{C}}(-,\gamma)})\cong\mathrm{Hom}_{\mathcal{A}}(F(\gamma),J)$ montre que $J^{\mathrm{Hom}_{\mathcal{C}}(-,\gamma)}$ est un injectif de $\mathcal{A}^{\mathcal{C}}$ si $J$ est un injectif de $\mathcal{A}$.

\smallskip
2) Soit $F$ un objet de $\mathcal{A}^{\mathcal{C}}$~; soit $(i_{c}:F(c)\to J_{c})_{c\in\mathcal{C}}$ une famille d'homomorphismes injectifs avec $J_{c}$ un injectif de $\mathcal{A}$. Cette famille induit tautologiquement  un homomorphisme injectif $i:F\to\prod_{c\in\mathcal{C}}J_{c}^{\mathrm{Hom}_{\mathcal{C}}(-,c)}$ (c'est là l'argument que  l'on invoque pour montrer que si $\mathcal{A}$ a assez d'injectifs alors il en est de même pour  $\mathcal{A}^{\mathcal{C}}$). En particulier si $F$ est un injectif de $\mathcal{A}^{\mathcal{C}}$ alors $F$ est facteur direct dans un produit (fini) d'injectifs de $\mathcal{A}^{\mathcal{C}}$ du type considéré précédemment~; on a donc encore dans ce cas $\mathrm{H}^{k}\hspace{1pt}\mathrm{L}^{\bullet}(F)=0$ pour $k>0$.

\smallskip
3) Soient $F$ un objet de $\mathcal{A}^{\mathcal{C}}$ et $F\to I^{\bullet}$ une résolution injective dans cette catégorie. On considère le bicomplexe $\mathrm{L}^{\bullet}(I^{\bullet})$ (on utilise ici ``the usual sign trick'' pour transformer un complexe de complexes en bicomplexe, voir par exemple \cite[1.2.5]{We}, la différentielle horizontale $\mathrm{d}_{\mathrm{h}}^{p,q}:\mathrm{L}^{p}(I^{q})\to \mathrm{L}^{p+1}(I^{q})$ est donnée par la différentielle de $\mathrm{L}^{\bullet}(I^{q})$ et la différentielle verticale $\mathrm{d}_{\mathrm{v}}^{p,q}:\mathrm{L}^{p}(I^{q})\to \mathrm{L}^{p}(I^{q+1})$ est $(-1)^{p}\mathrm{L}^{p}(I^{q}\to I^{q+1})$). Ce qui précède montre que pour tout $k\geq 0$, le $\mathcal{A}$-objet $\mathrm{H}^{k}\mathrm{Tot}\hspace{2pt}\mathrm{L}^{\bullet}(I^{\bullet})$ est à la fois canoniquement isomorphe au $\mathcal{A}$-objet $\mathrm{H}^{k}\lim_{\mathcal{C}}I^{\bullet}$ (``commencer par dériver horizontalement'') et au $\mathcal{A}$-objet $\mathrm{H}^{k}\mathrm{L}^{\bullet}(F)$ (``commencer par dériver verticalement''). Précisons. Par construction le complexe $\mathrm{L}^{\bullet}(F)$ est muni d'une coaugmentation $\lim_{\mathcal{C}}F\to\mathrm{L}^{\bullet}(F)$~; en d'autres termes on dispose d'un homomorphisme de $\mathcal{A}$-complexes $\mathrm{c}^{\bullet}\lim_{\mathcal{C}}F\to\mathrm{L}^{\bullet}(F)$, $\mathrm{c}^{\bullet}a$, $a$ objet de $\mathcal{A}$, désignant le $\mathcal{A}$-complexe $(a\to 0\to 0\ldots)$. En appliquant le foncteur $\mathrm{Tot}$ à l'homomorphisme de bicomplexes $\mathrm{c}^{\bullet}\lim_{\mathcal{C}}I^{\bullet}\to\mathrm{L}^{\bullet}(I^{\bullet})$ on obtient un homomorphisme de $\mathcal{A}$-complexes $\eta_{\mathrm{v}}:\lim_{\mathcal{C}}I^{\bullet}\to\mathrm{Tot}\hspace{2pt}\mathrm{L}^{\bullet}(I^{\bullet})$. Pareillement, en appliquant le foncteur $\mathrm{Tot}$ à l'homomorphisme de bicomplexes $\mathrm{L}^{\bullet}(\mathrm{c}^{\bullet}F)\to\mathrm{L}^{\bullet}(I^{\bullet})$, $\mathrm{c}^{\bullet}F$ désignant le $\mathcal{A}^{\mathcal{C}}$-complexe $(F\to 0\to 0\ldots)$, on obtient un homomorphisme de $\mathcal{A}$-complexes $\eta_{\mathrm{h}}:\mathrm{L}^{\bullet}(F)\to\mathrm{Tot}\hspace{2pt}\mathrm{L}^{\bullet}(I^{\bullet})$. La précision promise est la suivante~: $\mathrm{H}^{k}(\eta_{\mathrm{v}})$ et $\mathrm{H}^{k}(\eta_{\mathrm{h}})$ sont des isomorphismes pour tout $k\geq 0$.

\medskip
(R.3) Soit $0\to F'\to F\to F''\to 0$ une suite exacte dans la catégorie~$\mathcal{A}^{\mathcal{C}}$. Comme le foncteur $\mathrm{L}^{\bullet}(-)$ est exact, on a une suite exacte de $\mathcal{A}$-complexes $0\to\mathrm{L}^{\bullet}(F')\to\mathrm{L}^{\bullet}(F')\to\mathrm{L}^{\bullet}(F')\to 0$ et donc un connectant $\partial:\mathrm{H}^{k}\mathrm{L}^{\bullet}(F'')\to\mathrm{H}^{k}\mathrm{L}^{\bullet}(F')$. En reprenant l'argument de bicomplexe détaillé ci-dessus, on se convainc que ce connectant s'identifie au connectant $\partial:\lim_{\mathcal{C}}^{k}F''\to\lim_{\mathcal{C}}^{k+1}\hspace{-2pt}F'$.

\end{raps}

\bigskip
On rassemble dans l'énoncé ci-dessous les points évoqués dans le 1) du (R.2)~:

\begin{scho-def}\label{co-induit} Soient $\mathcal{A}$ une catégorie abélienne avec assez d'injectifs et $\mathcal{C}$ un ensemble ordonné fini~; soient $J$ un objet de $\mathcal{A}$ et $\gamma$ un objet de $\mathcal{C}$.

\medskip
{\em (a)} On a $\lim_{C}^{k}J^{\mathrm{Hom}_{\mathcal{C}}(-,\gamma)}=0$ pour $k>0$.

\medskip
{\em (b)}  Soit $a$ un objet de $\mathcal{A}$~; on note $\mathrm{Y}_{a}:\mathrm{Hom}_{\mathcal{A}^{\mathcal{C}}}(\Delta_{a},J^{\mathrm{Hom}_{\mathcal{C}}(-,\gamma)})\overset{\cong}{\to}\mathrm{Hom}_{\mathcal{A}}(a,J)$ et $\mathrm{D}_{a}:\mathrm{Hom}_{\mathcal{A}^{\mathcal{C}}}(\Delta_{a},J^{\mathrm{Hom}_{\mathcal{C}}(-,\gamma)})\overset{\cong}{\to}\mathrm{Hom}_{\mathcal{A}}(a,\lim_{\mathcal{C}}J)$ les isomorphismes respectivement donnés par Yoneda et la définition même de $\lim_{\mathcal{C}}$. Alors le $\mathcal{A}$-morphisme $(\mathrm{D}_{J}\circ\mathrm{Y}_{J}^{-1})(\mathrm{id}_{J}):J\to\lim_{C}J^{\mathrm{Hom}_{\mathcal{C}}(-,\gamma)}$ est un isomorphisme.

\medskip
{\em (c)} Si $J$ est un injectif de $\mathcal{A}$ alors $J^{\mathrm{Hom}_{\mathcal{C}}(-,\gamma)}$ est un injectif de $\mathcal{A}^{\mathcal{C}}$ (ce type d'injectif est dit {\em tautologique}).

\medskip
{\em (d)} Si $F$ est un  injectif de $\mathcal{A}^{\mathcal{C}}$ alors $F(\gamma)$ est un injectif de $\mathcal{A}$ pour tout objet $\gamma$ de $\mathcal{C}$.

\bigskip
Un foncteur de $\mathcal{C}$ dans $\mathcal{A}$ isomorphe à un foncteur de la forme $J^{\mathrm{Hom}_{\mathcal{C}}(-,\gamma)}$ ($J$~arbitraire) est dit {\em co-induit}.
\end{scho-def}

\medskip
\begin{rem}\label{co-induit-bis} Les points (b) et (c) de \ref{co-induit} (qui sont triviaux) impliquent facilement le point (a) si l'on adopte le point de vue du (R.1) sur les dérivés de~$\lim$. En effet, si $J\to I^{\bullet}$ est une résolution injective dans la catégorie $\mathcal{A}$ alors $J^{\mathrm{Hom}_{\mathcal{C}}(-,\gamma)}\to {I^{\bullet}}\hspace{1pt}^{\mathrm{Hom}_{\mathcal{C}}(-,\gamma)}$ est une résolution injective dans la catégorie~$\mathcal{A}^{\mathcal{C}}$ et le complexe $\lim_{\mathcal{C}}{I^{\bullet}}\hspace{1pt}^{\mathrm{Hom}_{\mathcal{C}}(-,\gamma)}$ s'identifie à $I^{\bullet}$.
\end{rem}

\vspace{0.75cm}
Revenons maintenant au foncteur $\Psi_{M}$. Par définition $\Psi_{M}$ est la restriction à $\mathcal{W}_{0}$ du foncteur $\mathcal{W}\to V_{\mathrm{tf}}\text{-}\mathcal{U},W\mapsto\mathrm{EFix}_{(V,W)}M$, disons $\widehat{\Psi}_{M}$. Les homomorphismes $M=\widehat{\Psi}_{M}(0)\to\widehat{\Psi}_{M}(W)=\Psi_{M}(W)$, $W$ décrivant $\mathcal{W}_{0}$, fournissent un $(V_{\mathrm{tf}}\text{-}\mathcal{U})^{\mathcal{W}_{0}}$-morphisme $\Delta_{M}\to\Psi_{M}$ (la notation $\Delta_{M}$ est introduite dans le rappel (R.1) de \ref{rappelslim}), soit encore un $V_{\mathrm{tf}}\text{-}\mathcal{U}$-morphisme $M\to\lim_{\mathcal{W}_{0}}\Psi_{M}$, naturel en le $\mathrm{H}^{*}V_{\mathrm{tf}}$-$\mathrm{A}$-module instable $M$, que nous notons~$\rho_{M}$. Cette définition peut être paraphrasée ainsi~: puisque $0$ est un objet initial de $\mathcal{W}$, on dispose d'un homomorphisme canonique $\widehat{\Psi}_{M}(0)\to\lim_{\mathcal{W}}\widehat{\Psi}_{M}$ et cet homomorphisme est un isomorphisme~; $\rho_{M}$ est l'homomorphisme
$$
\hspace{24pt}
M=\widehat{\Psi}_{M}(0)\cong{\lim}_{\mathcal{W}}\widehat{\Psi}_{M}
\to{\lim}_{\mathcal{W}_{0}}(\widehat{\Psi}_{M})_{\vert\mathcal{W}_{0}}
={\lim}_{\mathcal{W}_{0}}\Psi_{M}
\hspace{24pt}.
$$

\medskip
\begin{rem}
La proposition \ref{commutativity} montre que $\rho_{M}$ est aussi induite par le produit des unités d'adjonction $M\to\mathrm{H}^{*}V\otimes_{\mathrm{H}^{*}V/W}\mathrm{Fix}_{(V,W)}M$.
\end{rem}

\begin{pro}\label{derivePf} Soit $M$ un $\mathrm{H}^{*}V_{\mathrm{tf}}$-$\mathrm{A}$-module instable~; on a dans $V_{\mathrm{tf}}\text{-}\mathcal{U}$ une suite exacte
$$
\begin{CD}
0@>>>\mathrm{Pf}M@>\subset>>M@>\rho_{M}>>\lim_{\mathcal{W}_{0}}\hspace{1pt}\Psi_{M}@>>>\mathrm{R}^{1}\mathrm{Pf}\hspace{1pt}M@>>>0
\end{CD}
$$
et pour $k>0$ un isomorphisme
$$
\hspace{24pt}
{\lim}_{\mathcal{W}_{0}}^{k}\hspace{1pt}\Psi_{M}
\hspace{4pt}\cong\hspace{4pt}
\mathrm{R}^{k+1}\mathrm{Pf}\hspace{1pt}M
\hspace{24pt}.
$$
\end{pro}

\textit{Démonstration.} Elle repose sur certaines  propriétés du foncteur $\Psi$ (introduit en \ref{defPsi}) que nous rassemblons dans l'énoncé ci-après.

\begin{pro}\label{proPsi} Le foncteur additif  $\Psi:V_{\mathrm{tf}}\text{-}\mathcal{U}\to(V_{\mathrm{tf}}\text{-}\mathcal{U})^{\mathcal{W}_{0}}$ possède les propriétés suivantes~:

\medskip
{\em (a)} $\Psi$ est exact~;

\medskip
{\em (b)} $\Psi$ envoie injectifs sur injectifs~;

\medskip
{\em (c)} pour tout injectif $I$ de $V_{\mathrm{tf}}\text{-}\mathcal{U}$, la suite de $V_{\mathrm{tf}}\text{-}\mathcal{U}$-morphismes
$$
\begin{CD}
0@>>>\mathrm{Pf}\hspace{1pt}I@>\subset>>I@>\rho_{I}>>{\lim}_{\mathcal{W}_{0}}\hspace{1pt}\Psi(I)@>>>0
\end{CD}
$$
est exacte.
\end{pro}

\pagebreak

\bigskip
\textit{Démonstration du (a).} L'exactitude du foncteur $\Psi$ résulte immédiatement de celle des foncteurs $\mathrm{EFix}_{(V,W)}$.
\hfill$\square$

\medskip
\textit{Démonstration du (b).} Compte tenu de \ref{injectifs3} il suffit de vérifier que $\Psi(I)$ est injectif pour $I=\mathrm{H}^{*}V\otimes_{\mathrm{H}^{*}V/W}\mathrm{J}_{V/W}(k)$ avec $(W,k)\in\mathcal{W}\times\mathbb{N}$. Soit $U$ un objet de $\mathcal{W}_{0}$, la proposition \ref{proclef} montre que l'on a
$$
\Psi(I)(U)
\hspace{4pt}=\hspace{4pt}
\begin{cases}
I^{\hspace{1pt}\mathrm{Hom}_{\mathcal{W}_{0}}(U,W)}
& \text{pour $W\not=0$,} \\
0
& \hspace{1pt}\text{pour $W=0$.}
\end{cases}
$$
On peut donc supposer $W\not=0$. Dans ce cas  la proposition \ref{proclef} implique plus précisément que l'on a un isomorphisme $\Psi(I)\cong I^{\hspace{1pt}\mathrm{Hom}_{\mathcal{W}_{0}}(-,W)}$ et donc que  $\Psi(I)$ est bien injectif (point (c) de \ref{co-induit}).
\hfill$\square$

\medskip
\textit{Démonstration du (c).} Là encore on peut supposer $I=\mathrm{H}^{*}V\otimes_{\mathrm{H}^{*}V/W}\mathrm{J}_{V/W}(k)$. Le cas $W=0$ est trivial~: on a $\mathrm{Pf}\hspace{1pt}I=I$ et $\Psi_{I}=0$ (voir la remarque \ref{implicationelementaire}, on rappelle que les notations $\Psi_{-}$ et $\Psi(-)$ sont interchangeables). Passons au cas $W\not=0$~:

\smallskip
-- D'une part on a $\mathrm{Pf}\hspace{1pt}I=0$~; cette égalité découle de l'énoncé \ref{annPf} (que nous vérifierons une fois achevée la démonstration en cours).

\smallskip
-- D'autre part le fait que $\rho_{I}$ est un isomorphisme est à nouveau conséquence de l'isomorphisme $\Psi(I)\cong I^{\hspace{1pt}\mathrm{Hom}_{\mathcal{W}_{0}}(-,W)}$. En effet on constate que l'isomorphisme $(\mathrm{D}_{I}\circ\mathrm{Y}_{I}^{-1})(\mathrm{id}_{I})$ évoqué dans le point (b) de \ref{co-induit} coïncide dans le cas présent avec $\rho_{I}$.
\hfill$\square$

\bigskip
\begin{lem}\label{annPf} Soit $N$ un $\mathrm{H}^{*}V/W_{\mathrm{tf}}$-$\mathrm{A}$-module instable. Si $W$ est non nul alors on a $\mathrm{Pf}_{V}(\mathrm{H}^{*}V\otimes_{\mathrm{H}^{*}V/W}N)=0$.
\end{lem}

\medskip
\textit{Démonstration.} Ce lemme ne concerne en fait que la structure de $\mathrm{H}^{*}V/W$-module de $N$~; il est conséquence du suivant~:

\medskip
\begin{lem}\label{pre-annPf} Soient $W$ un sous-groupe de $V$ et $r:V\to W$ une rétraction linéaire~; soit $N$ un $\mathrm{H}^{*}V/W$-module. Alors le $\mathrm{H}^{*}V$-module $\mathrm{H}^{*}V\otimes_{\mathrm{H}^{*}V/W}N$, vu comme un  $\mathrm{H}^{*}W$-module  {\em via} l'homomorphisme $r^{*}:\mathrm{H}^{*}W\to\mathrm{H}^{*}V$, est un $\mathrm{H}^{*}W$-module libre.
\end{lem}

\medskip
\textit{Démonstration.} On note $\phi:W\oplus V/W\to V$ l'isomorphisme induit par $r$. Comme $q\circ\phi$ (on rappelle que $q$ est l'homomorphisme canonique $V\to V/W$) et $r\circ\phi$ sont respectivement les projections sur $V/W$ et $W$, l'isomorphisme $\phi^{*}:\mathrm{H}^{*}V\to\mathrm{H}^{*}(W\oplus V/W)=\mathrm{H}^{*}W\otimes\mathrm{H}^{*}V/W$ est un isomorphisme de $\mathrm{H}^{*}V/W$-modules (à droite) et un isomorphisme de $\mathrm{H}^{*}W$-modules (à gauche), $\mathrm{H}^{*}V$ étant un  $\mathrm{H}^{*}W$-module (à gauche) \textit{via} $r^{*}$. On en déduit que $\phi^{*}$ induit~un isomorphisme de $\mathrm{H}^{*}W$-modules (à gauche) $\mathrm{H}^{*}V\otimes_{\mathrm{H}^{*}V/W}N\cong\mathrm{H}^{*}W\otimes N$.
\hfill$\square$

\pagebreak\bigskip
\textit{Fin de la démonstration de la proposition \ref{derivePf}.}

\medskip
On met en oeuvre la définition des $\lim_{\mathcal{W}_{0}}^{k}$ en termes d'injectifs de la catégorie $(V_{\mathrm{tf}}\text{-}\mathcal{U})^{\mathcal{W}_{0}}$ (voir le rappel (R.1) de \ref{rappelslim}).

\medskip
Soit $M\to I^{\bullet}$ une résolution injective dans la catégorie $V_{\mathrm{tf}}\text{-}\mathcal{U}$. D'après les points (a) et (b) de la proposition \ref{proPsi}, $\Psi(M)\to\Psi(I^{\bullet})$ est une résolution injective de $\Psi(M)$ dans la catégorie $(V_{\mathrm{tf}}\text{-}\mathcal{U})^{\mathcal{W}_{0}}$.  D'après le point (c) de cette proposition la suite de complexes de cochaînes dans la catégorie $V_{\mathrm{tf}}\text{-}\mathcal{U}$
$$
0\to\mathrm{Pf}(I^{\bullet})\to I^{\bullet}
\overset{\rho_{I^{\bullet}}}{\longrightarrow}{\lim}_{\mathcal{W}_{0}}\Psi(I^{\bullet})\to 0
$$
est exacte. On obtient le résultat de la proposition \ref{derivePf} en considérant la longue suite exacte associée des groupes de cohomologie et en observant que l'on a $\mathrm{H}^{0}I^{\bullet}=M$ et $\mathrm{H}^{k}I^{\bullet}=0$ pour $k>0$.
\hfill$\square\square$

\bigskip
Nous complétons notre étude des foncteurs dérivés de la partie finie par divers énoncés d'annulation de ces foncteurs que nous utiliserons par la suite.

\bigskip
\begin{pro}\label{annderivePf1} Les foncteurs dérivés $\mathrm{R}^{k}\mathrm{Pf}_{V}:V_{\mathrm{tf}}\text{-}\mathcal{U}\to V_{\mathrm{tf}}\text{-}\mathcal{U}$ sont nuls pour $k>\dim V$.
\end{pro}

\bigskip
\textit{Démonstration.} Le cas $V=0$ est trivial ($0_{\mathrm{tf}}\text{-}\mathcal{U}$ est la sous-catégorie pleine de $\mathcal{U}$ dont les objets sont les $\mathrm{A}$-modules instables finis et $\mathrm{Pf}_{0}$ est le foncteur identique)~; on suppose $V\not=0$. La proposition \ref{derivePf} dit que l'on a un isomorphisme $\mathrm{R}^{k}\mathrm{Pf}_{V}M\cong\lim_{\mathcal{W}_{0}}^{k-1}\Psi_{M}$ pour $k>\dim V$. Les formules pour les dérivés de $\lim$ que l'on a données dans le rappel (R.2) de \ref{rappelslim} montrent que les $\lim_{\mathcal{W}_{0}}^{k-1}$ sont nuls pour $k>\dim V$. En effet cette inégalité implique que l'ensemble $\mathrm{S}_{k-1}\mathcal{W}_{0}$ est vide.
\hfill$\square$

\bigskip
\begin{pro}\label{annderivePf2} Soit $M$ un $\mathrm{H}^{*}V$-$\mathrm{A}$-module instable fini. On a $\mathrm{R}^{k}\mathrm{Pf}M=\nolinebreak 0$ pour tout $k>0$.
\end{pro}

\medskip
\textit{Démonstration.} La proposition \ref{resolutionfini} dit en particulier que l'on dispose d'une résolution injective $M\to I^{\bullet}$ dans la catégorie $V_{\mathrm{tf}}\text{-}\mathcal{U}$ avec $I^{k}$ fini pour tout $k$. On a donc $\mathrm{Pf}\hspace{1pt}I^{\bullet}=I^{\bullet}$, d'où le résultat.
\hfill$\square$

\bigskip
\begin{pro}\label{annderivePf3} Soient $W$ un sous-groupe de $V$ et $M$ un $\mathrm{H}^{*}V/W_{\mathrm{tf}}$-$\mathrm{A}$-module instable. Si $W$ est non nul alors on a $\mathrm{R}^{k}\mathrm{Pf}_{V}(\mathrm{H}^{*}V\otimes_{\mathrm{H}^{*}V/W}M)=0$ pour tout $k\geq 0$.
\end{pro}

\pagebreak

\medskip
\textit{Démonstration.} Soit $M\to I^{\bullet}$ une résolution injective de $M$ dans la catégorie $V/W_{\mathrm{tf}}\text{-}\mathcal{U}$. Comme $\mathrm{H}^{*}V$ est plat sur $\mathrm{H}^{*}V/W$, la proposition \ref{exactFixbis} implique que $\mathrm{H}^{*}V\otimes_{\mathrm{H}^{*}V/W}M\to\mathrm{H}^{*}V\otimes_{\mathrm{H}^{*}V/W}I^{\bullet}$ est une résolution injective dans la catégorie $V_{\mathrm{tf}}\text{-}\mathcal{U}$. Or on a $\mathrm{Pf}_{V}(\mathrm{H}^{*}V\otimes_{\mathrm{H}^{*}V/W}I^{\bullet})=0$ d'après \ref{annPf}.
\hfill$\square$

\bigskip
En vue d'une future référence, nous énonçons une variante de la proposition ci-dessus, dont la démonstration utilise le même argument.

\bigskip
\begin{pro}\label{annderivePf3var} Soient $N$ un $\mathrm{H}^{*}V$-$\mathrm{A}$-module instable fini, $W$ un sous-groupe non nul de $V$ et $M$ un $\mathrm{H}^{*}V/W_{\mathrm{tf}}$-$\mathrm{A}$-module instable. Alors on a
$$
\mathrm{Ext}^{k}_{V_{\mathrm{tf}}\text{-}\mathcal{U}}(N,\mathrm{H}^{*}V\otimes_{\mathrm{H}^{*}V/W}M)
\hspace{4pt}=\hspace{4pt}
0
$$
pour tout $k\geq 0$.
\end{pro}

\bigskip
\textit{Démonstration.} On a $\mathrm{Hom}_{V_{\mathrm{tf}}\text{-}\mathcal{U}}(N,J^{\bullet})=\mathrm{Hom}_{V_{\mathrm{tf}}\text{-}\mathcal{U}}(N,\mathrm{Pf}_{V}J^{\bullet})=0$, $J^{\bullet}$ dési\-gnant la résolution injective de la démonstration précédente.
\hfill$\square$

\bigskip
\begin{pro}\label{annderivePf4} Soit $M$ un $\mathrm{H}^{*}V_{\mathrm{tf}}$-$\mathrm{A}$-module instable~; soit $\mathrm{dp}_{V}M$ la dimension projective du $\mathrm{H}^{*}V$-module sous-jacent à $M$. Alors on a $\mathrm{R}^{k}\mathrm{Pf}_{V}M=0$ pour $0\leq k<\dim V-\mathrm{dp}_{V}M$.
\end{pro}

\medskip
La démonstration de la proposition ci-dessus utilise essentiellement la proposition \ref{DS} ci-dessous dont la démonstration est renvoyée à la fin de cette section (la proposition \ref{DS} est équivalente à  la proposition 3.2.1 de \cite{BZ} qui est un point clef de cette référence).

\begin{pro}\label{DS} Soit $M$ un $\mathrm{H}^{*}V_{\mathrm{tf}}$-$\mathrm{A}$-module instable~; soit $i:M\to E$ une enveloppe injective de $M$ dans la catégorie $V_{\mathrm{tf}}\text{-}\mathcal{U}$. Alors on a l'inégalité
$$
\hspace{24pt}
\mathrm{dp}_{V}E
\hspace{4pt}\leq\hspace{4pt}
\mathrm{dp}_{V}M
\hspace{24pt}.
$$
\end{pro}

\textit{Démonstration de \ref{annderivePf4}.} Soit $i:M\to E$ une enveloppe injective de $M$ dans la catégorie $V_{\mathrm{tf}}\text{-}\mathcal{U}$. La proposition \ref{DS} implique l'inégalité
$$
\hspace{24pt}
\mathrm{dp}_{V}(\mathop{\mathrm{coker}}i)
\hspace{4pt}\leq\hspace{4pt}\mathrm{dp}_{V}M+1
\hspace{24pt};
$$
on en déduit par récurrence que $M$ admet une résolution injective dans la catégorie $V_{\mathrm{tf}}\text{-}\mathcal{U}$ 
$$
0\to M\to I^{0}\to I^{1}\to\ldots\to I^{k}\to\ldots
$$
avec $\mathrm{dp}_{V}I^{k}\hspace{2pt}\leq\hspace{2pt}\mathrm{dp}_{V}M+k$. On conclut à l'aide de l'énoncé ci-après (équivalent à un énoncé de \cite{BHZ}), concernant les $\mathrm{H}^{*}V$-modules gradués, qui implique $\mathrm{Pf}_{V}I^{k}=0$ pour $k<\dim V-\mathrm{dp}_{V}M$.

\begin{pro}\label{DSS} Soit $M$ un $\mathrm{H}^{*}V$-module gradué de type fini. Les deux conditions suivantes sont équivalentes~:
\begin{itemize}
\item[(i)] $\mathrm{dp}_{V}M=\dim V$~;
\item[(ii)] $\mathrm{Pf}_{V}M\not= 0$.
\end{itemize}

\smallskip
{\em (Evidemment $\mathrm{dp}_{V}(-)$ désigne ci-dessus la dimension projective d'un  $\mathrm{H}^{*}V$-module et  $\mathrm{Pf}_{V}(-)$ la ``partie finie'' d'un $\mathrm{H}^{*}V$-module de type fini.)}
\end{pro}

\medskip
\textit{Démonstration.} On pose, comme on l'a fait jusqu'à présent, $n=\dim V$. Soit $M$ un $\mathrm{H}^{*}V$-module gradué. En considérant une résolution libre ``minimale'' $M\leftarrow L_{\bullet}$ de $M$, c'est-à-dire telle que le bord du complexe  $\mathbb{F}_{2}\otimes_{\mathrm{H}^{*}V}L_{\bullet}$ est nul on constate que la condition (i) est équivalente à la suivante (on utilise ici le ``lemme de Nakayama gradué'')~:

\smallskip
(i-bis) $\mathrm{Tor}^{\mathrm{H}^{*}V}_{n}(\mathbb{F}_{2},M)\not=0$.

\smallskip
Or on a un isomorphisme de $\mathrm{H}^{*}V$-modules gradués~:
$$
\hspace{24pt}
\mathrm{Tor}^{\mathrm{H}^{*}V}_{n}(\mathbb{F}_{2},M)
\hspace{4pt}\cong\hspace{4pt}
\Sigma^{n}\hspace{2pt}\tau_{V}(M)
\hspace{24pt},
\leqno(*)
$$
$\tau_{V}(M)$ désignant le sous-module de $M$ constitué des éléments annulés par l'idéal d'augmentation de $\mathrm{H}^{*}V$. Pour s'en convaincre, se rappeler que $\mathrm{H}^{*}V$ est isomorphe à une algèbre de polynômes en $n$ indéterminées de degré $1$ et considérer une résolution libre ``à~la Koszul'' de $\mathbb{F}_{2}$. La condition (i) est donc encore équivalente à la suivante~:

\smallskip
(i-ter) $\tau_{V}(M)\not=0$.

\smallskip
On suppose maintenant $M$ de type fini. Il est clair que l'on a (i-ter)$\Rightarrow$(ii). On a aussi  (ii)$\Rightarrow$(i-ter), en effet si $x$ est un élément non nul de degré maximal de $\mathrm{Pf}_{V}M$ alors $x$ appartient à $\tau_{V}(M)$.
\hfill$\square\square$

\bigskip
Les propositions \ref{annderivePf1} et \ref{annderivePf4} impliquent~:

\medskip
\begin{cor}\label{DS3} Soit $M$ un $\mathrm{H}^{*}V$-$\mathrm{A}$-module instable. Si le $\mathrm{H}^{*}V$-module sous-jacent à $M$ est libre de dimension finie alors on a $\mathrm{R}^{k}\mathrm{Pf}_{V}M=0$ pour $k\not=\dim V$.

\end{cor}

\pagebreak

\medskip
\textit{Démonstration de la proposition \ref{DS}}

\medskip
La démonstration ci-après est une variante de celle de \cite[Proposition 3.2.1]{BZ}.

\medskip
Soit $i:M\to E$ une enveloppe injective de $M$ dans la catégorie $V_{\mathrm{tf}}\text{-}\mathcal{U}$. D'après \ref{injectifs3}, on peut supposer
$$
E
\hspace{4pt}=\hspace{4pt}
\bigoplus_{(U,m)\in\mathcal{W}\times\mathbb{N}}
\hspace{4pt}
(\mathrm{H}^{*}V\otimes_{\mathrm{H}^{*}V/U}\mathrm{J}_{V/U}(m))^{\oplus\hspace{1pt}a_{U,m}}
$$
avec $a_{U,m}\in\mathbb{N}$ et $a_{U,m}=0$ pour $m$ assez grand. Soit $\mathcal{W}_{E}$ le sous-ensemble de $\mathcal{W}$ constitué des $U$ tels qu'il existe $m$ avec $a_{U,m}\not=0$~; on constate que l'on a
$$
\hspace{24pt}
\mathrm{dp}_{V}(E)
\hspace{4pt}=\hspace{4pt}
\sup_{U\in\mathcal{W}_{E}}
\mathop{\mathrm{codim}}U
\hspace{24pt}.
$$
Soit $W\subset V$ pour lequel ce $\sup$ est atteint. On pose $p=\mathrm{dp}_{V}(E)$~; on a  donc $\mathop{\mathrm{codim}}W=p$. D'après \ref{enveloppeFix},
$\mathrm{Fix}_{(V,W)}(i):\mathrm{Fix}_{(V,W)}M\to\mathrm{Fix}_{(V,W)}E$ est une enveloppe injective de $\mathrm{Fix}_{(V,W)}M$ dans la catégorie $(V/W)_{\mathrm{tf}}\text{-}\mathcal{U}$. Le choix de $W$ fait qu'il existe un entier $m$ tel que $\mathrm{H}^{*}V\otimes_{\mathrm{H}^{*}V/W}\mathrm{J}_{V/W}(m)$ est un facteur direct de $E$ et donc, d'après \ref{corproclef}, que $\mathrm{J}_{V/W}(m)$ est un facteur direct de $\mathrm{Fix}_{(V,W)}E$. Comme $\mathrm{J}_{V/W}(m)$ est enveloppe injective de $\Sigma^{m}\mathbb{F}_{2}$ (dans la catégorie $V\text{-}\mathcal{U}$),  $\mathrm{Fix}_{(V,W)}M$ contient un sous-module isomorphe à un sous-module non nul de $\Sigma^{m}\mathbb{F}_{2}$ donc à $\Sigma^{m}\mathbb{F}_{2}$ lui-même. On voit donc que le $\mathrm{H}^{*}V/W$-module gradué sous-jacent à $\mathrm{Fix}_{(V,W)}M$ contient un sous-module isomorphe à $\Sigma^{m}\mathbb{F}_{2}$ pour un certain entier $m$, c'est-à-dire que $\tau_{V/W}(\mathrm{Fix}_{(V,W)}M)$ est non nul (la notation $\tau$ a été introduite au cours de la démonstration de \ref{DSS}) ou encore que l'on a $\mathrm{Tor}^{\mathrm{H}^{*}V/W}_{p}(\mathbb{F}_{2},\mathrm{Fix}_{(V,W)}M)\not=0$ (isomorphisme (*) de la démonstration précitée). D'après \ref{sptorFix} on a aussi $\mathrm{Tor}_{p}^{\mathrm{H}^{*}V} (\mathrm{H}^{*}W,M)\not=0$ ce qui implique bien l'inégalité $\mathrm{dp}_{V}(M)\geq p$.
\hfill$\square$

\vspace{0,75cm}
\textsc{Complément~: Localisation de $V_{\mathrm{tf}}\text{-}\mathcal{U}$ modulo les objets finis}

\medskip
Notre référence pour la théorie de la localisation dans les catégories abéliennes est \cite[Chap.  III]{Ga}.

\bigskip
Soit $\mathcal{F}$ la sous-catégorie pleine de $V_{\mathrm{tf}}\text{-}\mathcal{U}$ dont les objets sont les $\mathrm{H}^{*}V$-$\mathrm{A}$-modules instables finis. Cette sous-catégorie est épaisse (en clair~: si $0\to M'\to M\to M''\to 0$ est une suite exacte de $V_{\mathrm{tf}}\text{-}\mathcal{U}$ alors $M$ est fini si et seulement si $M'$ et $M''$ le sont) si bien que l'on peut considérer la catégorie quotient $(V_{\mathrm{tf}}\text{-}\mathcal{U})/\mathcal{F}$. Nous nous proposons, en suivant à nouveau \cite{Hennalg}, de reformuler dans le  langage de \cite[Chap.  III]{Ga} certains des résultats de cette section.

\medskip
Soit $M$ un $\mathrm{H}^{*}V_{\mathrm{tf}}$-$\mathrm{A}$-module instable. On dit que $M$ est {\em $\mathcal{F}$-fermé} si l'on a $\mathrm{Ext}_{V_{\mathrm{tf}}\text{-}\mathcal{U}}^{k}(N,M)=0$ pour $k=0$ et $k=1$, quand $N$ est fini. On observera que la condition pour $k=0$ équivaut à $\mathrm{Pf}\hspace{1pt}M=0$~; si cette condition est vérifiée on dit que $M$ est {\em $\mathcal{F}$-réduit}.

\medskip
On note $\mathrm{Q}:V_{\mathrm{tf}}\text{-}\mathcal{U}\to(V_{\mathrm{tf}}\text{-}\mathcal{U})/\mathcal{F}$ le foncteur canonique.

\medskip
-- La proposition \ref{derivePf} dit en particulier que le $(V_{\mathrm{tf}}\text{-}\mathcal{U})/\mathcal{F}$-morphisme
$$
\mathrm{Q}(\rho_{M})
\hspace{2pt}:\hspace{2pt}M\longrightarrow{\lim}_{\mathcal{W}_{0}}\Psi_{M}
$$
est un isomorphisme. En effet, le noyau et le conoyau de $\rho_{M}$ sont finis \cite[Chap. III, Lemme 2]{Ga}.

\medskip
-- Si $M$ est fini alors $\lim_{\mathcal{W}_{0}}\Psi_{M}$ est nul, par exemple parce que l'on a $\mathrm{Pf}M=M$ et $\mathrm{R}^{1}\mathrm{Pf}M=0$ (voir \ref{annderivePf2}). (En fait le foncteur $\Psi_{M}$ est lui-même nul, voir la remarque \ref{implicationelementaire}.) Le foncteur
$$
V_{\mathrm{tf}}\text{-}\mathcal{U}\to V_{\mathrm{tf}}\text{-}\mathcal{U}\hspace{12pt},\hspace{12pt}
M\mapsto{\lim}_{\mathcal{W}_{0}}\Psi_{M}
$$
induit donc \cite[Chap. III, Cor. 2]{Ga} un foncteur, disons
$\mathrm{S}:(V_{\mathrm{tf}}\text{-}\mathcal{U})/\mathcal{F}\to V_{\mathrm{tf}}\text{-}\mathcal{U}$ tel que le foncteur composé $\mathrm{Q}\circ\mathrm{S}$ est naturellement isomorphe au foncteur identique de $(V_{\mathrm{tf}}\text{-}\mathcal{U})/\mathcal{F}$ ($\mathrm{S}$ est pour ``section'').

\medskip
-- Pour tout $M$, le $\mathrm{H}^{*}V_{\mathrm{tf}}$-$\mathrm{A}$-module instable $\lim_{\mathcal{W}_{0}}\Psi_{M}$ est $\mathcal{F}$-fermé.

\smallskip
Vérifions cette affirmation. Par définition de $\lim$ on a une suite exacte dans $V_{\mathrm{tf}}\text{-}\mathcal{U}$~:
$$
0\to{\lim}_{\mathcal{W}_{0}}\Psi_{M}\to\mathrm{L}^{0}(\Psi_{M})
\overset{\mathrm{d}}{\longrightarrow}\mathrm{L}^{1}(\Psi_{M})
$$
(notations introduites dans le rappel (R.2) de \ref{rappelslim}). Comme d'après \ref{annderivePf3var} $\Psi_{M}(W)$ est $\mathcal{F}$-fermé pour tout $W$ dans $\mathcal{W}_{0}$ , il est de même pour $\mathrm{L}^{0}(\Psi_{M})$ et $\mathrm{L}^{1}(\Psi_{M})$. On a donc encore une suite exacte courte
$$
0\to{\lim}_{\mathcal{W}_{0}}\Psi_{M}\to\mathrm{L}^{0}(\Psi_{M})
\to\mathop{\mathrm{im}}\hspace{-1pt}\mathrm{d}\to 0
$$
avec $\mathrm{L}^{0}(\Psi_{M})$ $\mathcal{F}$-fermé et $\mathop{\mathrm{im}}\hspace{-1pt}\mathrm{d}$ $\mathcal{F}$-réduit~; il en résulte bien que $\lim_{\mathcal{W}_{0}}\Psi_{M}$ est $\mathcal{F}$-fermé.

\smallskip
Ceci implique formellement que le foncteur $\mathrm{S}$ est adjoint à droite du foncteur $\mathrm{Q}$ et plus précisément que la transformation naturelle
$$
\rho_{M}:M\to(\mathrm{S}\circ\mathrm{Q})(M)={\lim}_{\mathcal{W}_{0}}\Psi_{M}
$$
est l'unité de cette adjonction. En conclusion~:

\medskip
\begin{pro}
La sous-catégorie $\mathcal{F}$ de $V_{\mathrm{tf}}\text{-}\mathcal{U}$ est localisante et l'endofoncteur $V_{\mathrm{tf}}\text{-}\mathcal{U}\to V_{\mathrm{tf}}\text{-}\mathcal{U}\hspace{2pt},\hspace{2pt}M\mapsto\lim_{\mathcal{W}_{0}}\Psi_{M}$, muni de la transformation naturelle $\rho$, est un endofoncteur de localisation de $V_{\mathrm{tf}}\text{-}\mathcal{U}$ modulo $\mathcal{F}$.
\end{pro}

\sect{Le complexe topologique}
 
Soit $X$ un $V$-CW-complexe fini. Rappelons la définition, donnée dans l'introduction, des complexes de chaînes $\mathrm{C}^{\bullet}_{\mathrm{top}}X$ et $\widetilde{\mathrm{C}}^{\bullet}_{\mathrm{top}}X$. On définit une filtration croissante de $X$ par des sous-$V$-CW-complexes~:
$$
\emptyset=\mathrm{F}_{-1}\hspace{1pt}X\subset\mathrm{F}_{0}\hspace{1pt}X\subset\mathrm{F}_{1}\hspace{1pt}X\subset\ldots\subset\mathrm{F}_{p}\hspace{1pt}X\subset\ldots\subset\mathrm{F}_{n-1}\hspace{1pt}X\subset\mathrm{F}_{n}\hspace{1pt}X=X
$$
en posant
$$
\hspace{24pt}
\mathrm{F}_{p}\hspace{1pt}X
\hspace{4pt}:=\hspace{4pt}
\bigcup_{\mathop{\mathrm{codim}}W\hspace{1pt}\leq\hspace{1pt}p} X^{W}
\hspace{24pt},
$$
pour  $-1\leq p\leq n$~; cette filtration induit une filtration de la construction de Borel $\mathrm{EV}\times_{V}X$. Le complexe $\mathrm{C}^{\bullet}_{\mathrm{top}}X$ est le terme $\mathrm{E}_{1}^{\bullet,*}$ de la suite spectrale, associée à cette dernière filtration, convergeant vers la cohomologie modulo~$2$ de $\mathrm{EV}\times_{V}X$. En clair on pose $\mathrm{C}^{p}_{\mathrm{top}}X:=\mathop{\Sigma^{-p}}\mathrm{H}^{*}_{V}(\mathrm{F}_{p}\hspace{1pt}X,\mathrm{F}_{p-1}\hspace{1pt}X)$ et on prend pour cobord le connectant de la triade $\mathrm{E}V\times_{V}(\mathrm{F}_{p+1}\hspace{1pt}X,\mathrm{F}_{p}\hspace{1pt}X,\mathrm{F}_{p-1}\hspace{1pt}X)$. Les deux observations suivantes permettent de préciser la structure de $\mathrm{C}^{p}_{\mathrm{top}}X$~:

\medskip
-- Le théorème d'excision pour les CW-complexes montre que l'on a un isomorphisme (de $\mathrm{H}^{*}V$-$\mathrm{A}$-modules instables)
$$
\hspace{4pt}
\mathrm{H}_{V}^{*}(\mathrm{F}_{p}\hspace{1pt}X,\mathrm{F}_{p-1}\hspace{1pt}X)
\hspace{4pt}\cong\hspace{4pt}
\bigoplus_{\mathop{\mathrm{codim}}W=p}
\mathrm{H}^{*}V
\otimes_{\mathrm{H}^{*}V/W}
\mathrm{H}_{V/W}^{*}(X^{W},\mathrm{Sing}_{V/W}X^{W})
\hspace{4pt}.
$$
Décryptons le membre de droite. Le sous-$V$-CW-complexe $X^{W}\subset X$ est considéré ci-dessus comme un $V/W$-espace et $\mathrm{Sing}_{V/W}X^{W}\subset X^{W}$ désigne le sous-espace (stable sous l'action de $V/W$) constitué des points dont l'isotropie (pour l'action de $V/W$) est non triviale.

\bigskip
-- Les $\mathrm{H}^{*}V/W$-$\mathrm{A}$-modules instables $\mathrm{H}_{V/W}^{*}(X^{W},\mathrm{Sing}_{V/W}X^{W})$ sont finis.

\bigskip
Comme la propriété ci-dessus joue un rôle important dans notre mémoire nous en proposons une démonstration \textit{ab initio} dans l'intermède ci-après où nous analysons $\mathrm{H}^{*}_{V}(X,Y)$ pour une paire de $V$-CW-complexes finis telle que l'action de $V$ sur $X-Y$ est libre.

\pagebreak

\bigskip
\textsc{Intermède de Topologie Générale}

\begin{lem}\label{finitude1} Soient $G$ un groupe discret et $(X,Y)$ une paire de $G$-espaces~; on suppose que $Y$ est fermé dans $X$ et qu'il existe un voisinage ouvert de $Y$ dans $X$, stable sous l'action de $G$, qui se rétracte par déformation $G$-équivariante sur $Y$.

\smallskip
Si l'action de $G$ sur $X-Y$ est topologiquement libre alors l'homomorphisme canonique
$$
\mathrm{H}^{*}(G\backslash X,G\backslash Y)\longrightarrow
\mathrm{H}^{*}_{G}(X,Y)
$$
est un isomorphisme.
\end{lem}

\medskip
\footnotesize
(Une action continue d'un groupe topologique $G$ sur un espace $Z$ est dite {\em topologiquement libre} si pour tout point de $Z$ il existe un voisinage ouvert $U$ de ce point, stable sous l'action de~$G$, et une application continue $G$-équivariante de $U$ dans $G$~; si $G$ est discret fini et $Z$ séparé alors une action libre est topologiquement libre.)
\normalsize

\bigskip
\textit{Démonstration.} Soit $N$ le voisinage en question. On a $\mathrm{H}^{*}_{G}(X,Y)\cong\mathrm{H}^{*}_{G}(X,N)$ puis $\mathrm{H}^{*}_{G}(X,N)\cong\mathrm{H}^{*}_{G}(X-Y,N-Y)$ par excision. Comme l'action de $G$ sur $X-Y$ (et $N-Y$) est topologiquement libre on a $\mathrm{H}^{*}_{G}(X-Y,N-\nolinebreak Y)\cong\mathrm{H}^{*}(G\backslash (X-Y),G\backslash (N-Y))=\mathrm{H}^{*}(G\backslash X-G\backslash Y,G\backslash N-G\backslash Y)$. Puis on a $\mathrm{H}^{*}(G\backslash X-G\backslash Y,G\backslash N-G\backslash Y)\cong
\mathrm{H}^{*}(G\backslash X,G\backslash N)$, à nouveau par excision. Enfin, comme la rétraction par déformation $G$-équivariante de $N$ sur $Y$ induit une rétraction par déformation de $G\backslash N$ sur $G\backslash Y$, on a $\mathrm{H}^{*}(G\backslash X,G\backslash N)\cong\mathrm{H}^{*}(G\backslash X,G\backslash Y)$.
\hfill$\square$

\medskip
\begin{scho}\label{finitude2} Soient $X$ un $V$-CW-complexe fini et $Y$ un sous-$V$-CW-complexe. Si l'action de $V$ sur $X-Y$ est libre alors~:

\smallskip {\em (a)} l'homomorphisme canonique $\mathrm{H}^{*}(V\backslash X,V\backslash Y)\to\mathrm{H}^{*}_{V}(X,Y)$ est un isomorphisme~;

\smallskip {\em (b)} $\mathrm{H}^{*}_{V}(X,Y)$ est fini.
\end{scho}

\bigskip
\textit{Démonstration.} Le (b) résulte du (a) puisque $(V\backslash X,V\backslash Y)$ est une paire de CW-complexes finis. L'hypothèse du lemme \ref{finitude1} est satisfaite dans notre contexte (voir \cite[Chap. II, (1.17), Exc. 3]{tD}), d'où le point (b).
\hfill$\square$

\medskip
\begin{scho}\label{finitude3} Soient $X$ un $V$-CW-complexe fini et $Y$ un sous-$V$-CW-complexe. Si l'action de $V$ sur $X-Y$ est libre alors on a un isomorphisme canonique de $\mathrm{H}^{*}V$-$\mathrm{A}$-modules instables
$$
\hspace{24pt}
\mathrm{H}^{*}_{V}(X,Y)
\hspace{4pt}\cong\hspace{4pt}
\mathrm{H}^{*}_{\mathrm{c}}(V\backslash(X-Y))
\hspace{24pt}.
$$

\pagebreak

\smallskip
{\em (La notation $\mathrm{H}^{*}_{\mathrm{c}}$ désigne ci-dessus la cohomologie modulo $2$ à support compact. Posons $Z=X-Y$ et notons $\pi:Z\to V\backslash Z$ le passage au quotient~; la structure de  $\mathrm{H}^{*}V$-$\mathrm{A}$-module instable de $\mathrm{H}^{*}_{\mathrm{c}}(V\backslash Z)$ est donnée par l'isomorphisme
\begin{multline*}
\hspace{12pt}
\mathrm{H}^{*}_{\mathrm{c}}(V\backslash Z)
:=\mathop{\mathrm{colim}}_{K}\mathrm{H}^{*}(V\backslash Z,V\backslash Z-K)= \\
\mathop{\mathrm{colim}}_{K}\mathrm{H}^{*}(V\backslash Z,V\backslash (Z-\pi^{-1}(K)))
\cong
\mathop{\mathrm{colim}}_{K}\mathrm{H}^{*}_{V}(Z,Z-\pi^{-1}(K))\hspace{12pt},
\end{multline*}
$K$ décrivant l'ensemble des compacts de $V\backslash Z$ ordonné par inclusion.)}
\end{scho}

\bigskip
\textit{Démonstration.} On dispose du lemme suivant~:

\begin{lem}\label{finitude4} Soit $(A,B)$ une paire de CW-complexes finis. L'homomorphisme $\mathrm{H}^{*}(A,B)\to\mathrm{H}^{*}(A-B)$ induit par l'inclusion de paires $(A-B,\emptyset)\hookrightarrow(A,B)$ se factorise par un isomorphisme canonique $\mathrm{H}^{*}(A,B)\cong\mathrm{H}^{*}_{\mathrm{c}}(A-B)$.
\end{lem}

\medskip
\textit{Démonstration.} Soient $\mathcal{C}$ l'ensemble des compacts $K$ de $A-B$ ordonné par inclusion et $\mathcal{D}$ le sous-ensemble de $\mathcal{C}$ constitué des $K$ tels que $A-K$ se rétracte par déformation sur $B$. Comme $\mathcal{D}$ est cofinal dans $\mathcal{C}$ (voir \cite[Theorem 6.1]{LW}\cite[Proposition A.5]{Ha}), on a
$$
\hspace{24pt}
\mathrm{H}^{*}_{\mathrm{c}}(A-B)
\hspace{4pt}\cong\hspace{4pt}
\mathop{\mathrm{colim}}_{K\in\mathcal{D}}\mathrm{H}^{*}(A-B,(A-B)-K)
\hspace{24pt};
$$
or les inclusions de paires $(A,B)\hookrightarrow(A, A-K)\hookleftarrow (A-B,(A-B)-K)$ induisent des isomorphismes en cohomologie, la première par l'hypothèse faite sur $K$ et la seconde par excision.
\hfill$\square$

\medskip
Le lemme ci-dessus et le point (a) du scholie \ref{finitude2} montrent que l'on a un isomorphisme canonique de $\mathrm{A}$-modules instables
$$
\hspace{24pt}
\mathrm{H}^{*}_{V}(X,Y)
\hspace{4pt}\cong\hspace{4pt}
\mathrm{H}^{*}_{\mathrm{c}}(V\backslash X-V\backslash Y))
=
\mathrm{H}^{*}_{\mathrm{c}}(V\backslash(X-Y))
\hspace{24pt}.
$$
Pour en finir avec le scholie \ref{finitude3} il reste à se convaincre que c'est bien un isomorphisme de $\mathrm{H}^{*}V$-$\mathrm{A}$-modules instables. On pose $A=V\backslash X$ et $B=V\backslash Y$. Soit $N$ un voisinage ouvert de $Y$ dans $X$, stable sous l'action de $V$, qui se rétracte par déformation $V$-équivariante sur $Y$~; on note $L$ le compact $X-N$ (stable sous l'action de $V$) et $K$ son image dans $A$ (qui se rétracte par déformation sur~$B$). On considère le diagramme commutatif
$$
\begin{CD}
\mathrm{H}^{*}(A,B)@<\bar{\alpha}<<\mathrm{H}^{*}(A,A-K)
@>\bar{\beta}>>\mathrm{H}^{*}(A-B,(A-B)-K) \\
@V\nu_{\mathrm{g}}VV
@V\nu_{\mathrm{m}}VV
@V\nu_{\mathrm{d}}VV \\
\mathrm{H}^{*}_{V}(X,Y)@<\alpha<<\mathrm{H}^{*}_{V}(X,X-L)
@>\beta>>
\mathrm{H}^{*}_{V}(X-Y,(X-Y)-L)
\end{CD}
$$
(les indices $\mathrm{g}$, $\mathrm{m}$ et $\mathrm{d}$ sont pour gauche, milieu et droite~!). On observe que toutes les flèches de ce diagramme sont des isomorphismes. En effet $\bar{\alpha}$, $\bar{\beta}$, $\nu_{\mathrm{g}}$ et  $\alpha$ sont des isomorphismes d'après ce qui précède et $\nu_{\mathrm{d}}$ en est un parce que l'action de $V$ sur $X-Y$ est libre. L'isomorphisme du scholie \ref{finitude3} est le composé $\gamma\circ\nu_{\mathrm{d}}^{-1}\circ\beta\circ\alpha^{-1}:\mathrm{H}^{*}_{V}(X,Y)\to\mathrm{H}^{*}_{\mathrm{c}}(V\backslash (X,Y))$, $\gamma$ désignant l'homomorphisme canonique $\mathrm{H}^{*}(A-B,(A-B)-K) \to\mathrm{H}^{*}_{\mathrm{c}}(V\backslash (X,Y))$ qui comme on l'a vu plus haut est un isomorphisme. Or par définition même de la structure de $\mathrm{H}^{*}V$-$\mathrm{A}$-module instable de $\mathrm{H}^{*}_{\mathrm{c}}(V\backslash (X,Y))$, $\gamma\circ\nu_{\mathrm{d}}^{-1}$ est un homomorphisme de $\mathrm{H}^{*}V$-$\mathrm{A}$-modules instables.
\hfill$\square\square$

\bigskip
\begin{rem}\label{finitude5} On peut démontrer (de manière assez détournée~!) \ref{finitude2} en invoquant le corollaire \ref{cohFixpaire} et la proposition \ref{caracterisationfini} puisque l'hypothèse de \ref{finitude2} est équivalente à $X^{W}=Y^{W}$ pour tout sous-groupe $W$ de $V$ de dimension $1$.
\end{rem}

\bigskip
\textsc{Fin de l'intermède} 

\bigskip
Nous archivons les deux observations qui précèdent \ref{finitude1}~:

\begin{pro}\label{Ctopefini} Soient $V$ un $2$-groupe abélien élémentaire et $X$ un $V$-CW-complexe fini.

\medskip
{\em (a)} On a pour tout $p$ un isomorphisme de $\mathrm{H}^{*}V$-$\mathrm{A}$-modules instables,  naturel en $X$~:
$$
\hspace{24pt}
\mathop{\Sigma^{p}\hspace{2pt}\mathrm{C}^{p}_{\mathrm{top}}X
\hspace{4pt}\cong\hspace{4pt}
\bigoplus_{\mathop{\mathrm{codim}}W=p}
}\mathrm{H}^{*}V
\otimes_{\mathrm{H}^{*}V/W}
\mathrm{H}_{V/W}^{*}(X^{W},\mathrm{Sing}_{V/W}X^{W})
\hspace{24pt}.
$$

\medskip
{\em (b)} Le  $\mathrm{H}^{*}V/W$-$\mathrm{A}$-module instable $\mathrm{H}_{V/W}^{*}(X^{W},\mathrm{Sing}_{V/W}X^{W})$ est fini pour\linebreak tout $W$.
\end{pro}

\medskip
Le complexe $\widetilde{\mathrm{C}}_{\mathrm{top}}^{\bullet}\hspace{1pt}X$ est associé à la coaugmentation naturelle $\mathrm{H}_{V}^{*}X=:\linebreak\mathrm{C}_{\mathrm{top}}^{-1}\hspace{1pt}X\to\mathrm{C}_{\mathrm{top}}^{0}\hspace{1pt}X=\mathrm{H}^{*}_{V}X^{V}$ dont $\mathrm{C}_{\mathrm{top}}^{\bullet}\hspace{1pt}X$ est muni. Formalisons (lourdement) la définition de celle-ci. Soit $(\mathrm{F}_{p}^{\mathrm{tr}}X)_{-1\leq p\leq n}$ la filtration croissante ``triviale''  de $X$ définie par
$$
\mathrm{F}_{p}^{\mathrm{tr}}X
\hspace{4pt}=\hspace{4pt}
\begin{cases}
\emptyset & \text{pour $p=-1$}, \\
X & \text{pour $0\leq p\leq n$}.
\end{cases}
$$
Le complexe obtenu en remplaçant dans la définition de $\mathrm{C}_{\mathrm{top}}^{\bullet}\hspace{1pt}X$ la filtration  $(\mathrm{F}_{p}X)$ par la filtration $(\mathrm{F}_{p}^{\mathrm{tr}}X)$ est $\mathrm{H}^{*}_{V}X\to 0\to 0\ldots$~;  l'homomorphisme de ce dernier complexe dans $\mathrm{C}_{\mathrm{top}}^{\bullet}\hspace{1pt}X$,  induit par les inclusions $\mathrm{F}_{p}X\subset\mathrm{F}_{p}^{\mathrm{tr}}X$, est la coaugmentation évoquée plus haut.

\bigskip
Venons-en maintenant à la démonstration du théorème \ref{introCtop} dont nous reproduisons l'énoncé~:

\begin{theo}\label{Ctop}
Soient $V$ un $2$-groupe abélien élémentaire et $X$ un $V$-CW-complexe fini. Les deux propriétés suivantes sont équivalentes~:
\begin{itemize}
\item[(i)] Le $\mathrm{H}^{*}V$-module sous-jacent à $\mathrm{H}_{V}^{*}X$ est libre.
\item[(ii)] Le complexe coaugmenté $\widetilde{\mathrm{C}}_{\mathrm{top}}^{\bullet}\hspace{1pt}X$ est acyclique.
\end{itemize}
\end{theo}

\medskip
\textit{Démonstration de l'implication $(ii)\Rightarrow(i)$.} Soit $M$ un  $\mathrm{H}^{*}V$-module $\mathbb{N}$-gradué. Les deux conditions suivantes sont équivalentes (conséquence du ``lemme de Nakayama gradué'')~:

\smallskip
--\hspace{8pt}$M$ est libre comme $\mathrm{H}^{*}V$-module~;

\smallskip
--\hspace{8pt}$\mathrm{Tor}_{1}^{\mathrm{H}^{*}V}(\mathbb{F}_{2},M)=0$.

\smallskip
Il suffit donc de se convaincre de l'égalité $\mathrm{Tor}_{1}^{\mathrm{H}^{*}V}(\mathbb{F}_{2},\mathrm{H}_{V}^{*}X)=0$. Pour cela on observe que le point (a) de \ref{Ctopefini} entraîne que la dimension projective du  $\mathrm{H}^{*}V$-module sous-jacent à $\mathrm{C}_{\mathrm{top}}^{p}\hspace{1pt}X$ est inférieure ou égale à $p$ et on invoque le lemme \textit{ad hoc} suivant~:

\begin{lem}\label{base1} On considère une suite exacte de $\mathrm{H}^{*}V$-modules de la forme suivante~:
$$
\hspace{24pt}
0\to M\to C^{0}\to C^{1}\to\ldots\to C^{p}\to\ldots\to C^{r-1}\overset{\mathrm{d}^{r-1}}{\longrightarrow}C^{r}\to 0
\hspace{24pt}.
$$
Si l'on suppose $\mathrm{Tor}_{k}^{\mathrm{H}^{*}V}(\mathbb{F}_{2},C^{p})=0$ pour $0\leq p\leq r$ et $k>p$, alors on a $\mathrm{Tor}_{k}^{\mathrm{H}^{*}V}(\mathbb{F}_{2},M)=0$ pour $k>0$.
\end{lem}

\medskip
\textit{Démonstration.} On procède par récurrence sur l'entier $r$. Le cas $r=0$ est trivial. On franchit le pas de récurrence en considérant la suite exacte
$$
0\to M\to C^{0}\to C^{1}\to\ldots\to C^{p}\to\ldots\to \ker\mathrm{d}^{r-1}\to 0
$$
et en observant que l'on a $\mathrm{Tor}_{k}^{\mathrm{H}^{*}V}(\mathbb{F}_{2}, \ker\mathrm{d}^{r-1})=0$ pour $k>r-1$ (utiliser la suite exacte courte $0\to\ker\mathrm{d}^{r-1}\to C^{r-1}\to C^{r}\to 0$).
\hfill$\square$

\bigskip
\textit{Démonstration de l'implication $(i)\Rightarrow(ii)$} 

\medskip
Pour $V=0$ il n'y a rien à démontrer~; on suppose donc $V\not=0$.

\medskip
On constate tout d'abord que l'on a bien $\mathrm{H}^{-1}\hspace{1pt}\widetilde{\mathrm{C}}_{\mathrm{top}}^{\bullet}\hspace{1pt}X=0$, c'est-à-dire que la coaugmentation $\eta:\mathrm{H}^{*}_{V}X\to\mathrm{H}^{*}_{V}X^{V}$ est injective. En effet la théorie de Smith dit en particulier que $\ker\eta$ est un $\mathrm{H}^{*}V$-module de torsion (voir par exemple \cite[Theorem 4.2]{Qui}, théorème  dont la démonstration est très facile dans le cas où $X$ est un $V$-CW-complexe fini, et invoquer le résultat ``folklorique'' \ref{Serre5})~; il est donc nul si $\mathrm{H}^{*}_{V}X$ est un $\mathrm{H}^{*}V$-module libre.

\bigskip
\textit{Le cas $\dim V=1$.} D'après ce qui précède la longue suite exacte en cohomologie équivariante modulo $2$ de la paire $(X,X^{V})$ donne en fait  une suite exacte courte $0\to\mathrm{H}^{*}_{V}X\to\mathrm{H}^{*}_{V}X^{V}\to\mathop{\Sigma}^{-1}\mathrm{H}^{*}_{V}(X,X^{V})\to 0$ (sans hypothèse sur $\dim V$). Or pour $\dim V=1$, le complexe $\widetilde{\mathrm{C}}_{\mathrm{top}}^{\bullet}\hspace{1pt}X$ est le complexe de cochaînes $\mathrm{H}^{*}_{V}X\to\mathrm{H}^{*}_{V}X^{V}\to\mathop{\Sigma}^{-1}\mathrm{H}^{*}_{V}(X,X^{V})$ ; il est bien acyclique.

\bigskip
\textit{Le cas $\dim V\geq 2$.} On procède par récurrence sur l'entier $\dim V$.

\smallskip
Comme nous allons être amené à faire varier le $2$-groupe abélien élémentaire~$V$ nous préciserons, lorsque cela nous semblera nécessaire, les  notations $\mathrm{C}_{\mathrm{top}}^{\bullet}\hspace{1pt}X$, $\widetilde{\mathrm{C}}_{\mathrm{top}}^{\bullet}\hspace{1pt}X$ et $\mathrm{F}_{p}\hspace{1pt}X$ en $\mathrm{C}_{\mathrm{top},V}^{\bullet}\hspace{1pt}X$, $\widetilde{\mathrm{C}}_{\mathrm{top},V}^{\bullet}\hspace{1pt}X$ et $\mathrm{F}_{p,V}\hspace{1pt}X$.

\medskip
Explicitons l'hypothèse de récurrence~:

\smallskip
(HR) Soient $V'$ un $2$-groupe abélien élémentaire et $X$ un $V'$-CW-complexe fini avec $\mathrm{H}_{V'}^{*}X'$ libre comme $\mathrm{H}^{*}V'$-module. Si l'on a $\dim V'<\dim V$ alors le complexe $\widetilde{\mathrm{C}}_{\mathrm{top},V'}^{\bullet}\hspace{1pt}X'$ est acyclique.

\medskip
Nous divisons notre démonstration en deux étapes.

\bigskip
\textit{Première étape.} On montre que si (HR) est vérifiée alors on a $\mathrm{H}^{p}\hspace{1pt}\widetilde{\mathrm{C}}_{\mathrm{top}}^{\bullet}\hspace{1pt}X=0$ pour $p<\dim V-1$.

\bigskip
En préalable nous démontrons deux propositions, \ref{FixCtop} et \ref{PsiCtop}, concernant respectivement $\widetilde{\mathrm{C}}_{\mathrm{top}}^{\bullet}\hspace{1pt}X$ et $\mathrm{C}_{\mathrm{top}}^{\bullet}\hspace{1pt}X$, dans lesquelles l'hypothèse ``$\mathrm{H}_{V}^{*}X$ est un $\mathrm{H}^{*}V$-module libre'' n'intervient pas, ni d'ailleurs l'hypothèse $\dim V\geq 2$. Nous posons à nouveau $\dim V=:n$.

\bigskip
Soit $W\subset V$ un sous-groupe~; $X^{W}$ est un sous-$V$-CW-complexe de $X$ que l'on peut considérer comme un $V/W$-CW-complexe (fini), ce que nous faisons ci-après.

\medskip
On note $\mathop{\Sigma^{n}}\widetilde{\mathrm{C}}_{\mathrm{top}}^{\bullet}\hspace{1pt}X$ la $n$-ième suspension (terme à terme) du complexe $\widetilde{\mathrm{C}}_{\mathrm{top}}^{\bullet}\hspace{1pt}X$. On a  $(\mathop{\Sigma^{n}}\widetilde{\mathrm{C}}_{\mathrm{top}}^{\bullet}\hspace{1pt}X)^{-1}=\mathop{\Sigma^{n}}\mathrm{H}_{V}^{*}X $ et $(\mathop{\Sigma^{n}}\widetilde{\mathrm{C}}_{\mathrm{top}}^{\bullet}\hspace{1pt}X)^{p}=\mathop{\Sigma^{n-p}}\mathrm{H}_{V}^{*}(\mathrm{F}_{p}\hspace{1pt}X,\mathrm{F}_{p-1}\hspace{1pt}X)$ pour $0\leq p\leq n$, ce qui montre que $\mathop{\Sigma^{n}}\widetilde{\mathrm{C}}_{\mathrm{top}}^{\bullet}\hspace{1pt}X$ est un complexe de cochaînes coaugmenté, dans  la catégorie des $\mathrm{H}^{*}V$-$\mathrm{A}$-modules \underline{instables}.

\begin{pro}\label{FixCtop}
On a un isomorphisme canonique de complexes de $\mathrm{H}^{*}V/W$-$\mathrm{A}$-modules instables~:
$$
\hspace{24pt}
\mathrm{Fix}_{(V,W)}\hspace{2pt}(\mathop{\Sigma^{n}}\widetilde{\mathrm{C}}_{\mathrm{top}}^{\bullet}\hspace{1pt}X)
\hspace{4pt}\cong\hspace{4pt}
\mathop{\Sigma^{n}}\widetilde{\mathrm{C}}_{\mathrm{top},V/W}^{\bullet}\hspace{1pt}X^{W}
\hspace{24pt}.
$$
\end{pro}

\textit{Démonstration.} On vérifie tout d'abord que $(\mathrm{F}_{p}\hspace{1pt}X)^{W}$ s'identifie à $\mathrm{F}_{p,V/W}\hspace{1pt}X^{W}$ (il est opportun ici de noter que notre définition des $\mathrm{F}_{p}X$ fait sens pour tout $p\geq -1$ et donne $\mathrm{F}_{p}X=X$ pour $p\geq\dim V$). On invoque ensuite le fait que le foncteur $\mathrm{Fix}_{(V,W)}$ ``commute à la suspension'' (point (b) de \ref{suspensionFix}), la proposition \ref{Fixconnectant} et le corollaire \ref{cohFixpaire}.
\hfill$\square$

\bigskip
La proposition ci-dessous concerne le foncteur $\Psi:V_{\mathrm{tf}}\text{-}\mathcal{U}\to(V_{\mathrm{tf}}\text{-}\mathcal{U})^{\mathcal{W}_{0}}$ introduit en section 3.

\begin{pro}\label{PsiCtop} Le $(V_{\mathrm{tf}}\text{-}\mathcal{U})^{\mathcal{W}_{0}}$-complexe $\Psi(\mathop{\Sigma^{n}}\mathrm{C}_{\mathrm{top}}^{\bullet}\hspace{1pt}X)$ possède les propriétés suivantes~:

\medskip
{\em (a)} Pour $0\leq p<n$ le foncteur $\Psi(\mathop{\Sigma^{n}}\mathrm{C}_{\mathrm{top}}^{p}\hspace{1pt}X):\mathcal{W}_{0}\to V_{\mathrm{tf}}\text{-}\mathcal{U}$ est une somme directe finie de foncteurs co-induits (terminologie introduite en \ref{co-induit}).

\medskip
{\em (b)} Le foncteur $\Psi(\mathop{\Sigma^{n}}\mathrm{C}_{\mathrm{top}}^{n}\hspace{1pt}X)$ est nul.

\medskip
{\em (c)} Pour $0\leq p\leq n$ et $k>0$ on a $\lim_{\mathcal{W}_{0}}^{k}\Psi(\mathop{\Sigma^{n}}\mathrm{C}_{\mathrm{top}}^{p}\hspace{1pt}X)=0$.

\medskip
{\em (d)} Pour $0\leq p<n$ le $(V_{\mathrm{tf}}\text{-}\mathcal{U})$-morphisme
$$
\rho_{\hspace{2pt}\mathop{\Sigma^{n}\hspace{-1pt}}\mathrm{C}_{\mathrm{top}}^{p}\hspace{1pt}X}:\mathop{\Sigma^{n}}\mathrm{C}_{\mathrm{top}}^{p}\hspace{1pt}X\to
{\lim}_{\mathcal{W}_{0}}\Psi(\mathop{\Sigma^{n}}\mathrm{C}_{\mathrm{top}}^{p}\hspace{1pt}X)
$$
est un isomorphisme.
\end{pro}

\bigskip
\textit{Démonstration du (a).} Le $\mathrm{H}^{*}V_{\mathrm{tf}}$-$\mathrm{A}$-module instable $\mathop{\Sigma^{n}}\mathrm{C}_{\mathrm{top}}^{p}\hspace{1pt}X$ est isomorphe à une somme directe de la forme
$
\bigoplus_{\mathop{\mathrm{codim}}W=p}
\mathrm{H}^{*}V
\otimes_{\mathrm{H}^{*}V/W}
N_{W}
$
avec $N_{W}$ un  $\mathrm{H}^{*}V/W$-$\mathrm{A}$-module instable fini (voir \ref{Ctopefini}).

\smallskip
On pose $M_{W}=\mathrm{H}^{*}V\otimes_{\mathrm{H}^{*}V/W}N_{W}$ et on suppose $p\not=n$. La proposition \ref{proclef} implique $\Psi(M_{W})\cong M_{W}^{\hspace{1pt}\mathrm{Hom}_{\mathcal{W}_{0}}(-,W)}$ (même argument que dans la démonstration du point (b) de \ref{proPsi}).

\medskip
\textit{Démonstration du (b).}  On a $\mathrm{EFix}_{(V,W)}\hspace{2pt}(\mathop{\Sigma^{n}}\mathrm{C}_{\mathrm{top}}^{n}\hspace{1pt}X)=0$ pour $W\not=0$ puisque $\mathop{\Sigma^{n}}\mathrm{C}_{\mathrm{top}}^{n}\hspace{1pt}X$ est fini (voir \ref{implicationelementaire}).

\medskip
\textit{Démonstration du (c).} Conséquence du point (a) ci-dessus et du point (a) de~\ref{co-induit}.

\medskip
\textit{Démonstration du (d).} La suite exacte de \ref{derivePf} et le résultat d'annulation \ref{annderivePf3} montrent que $\rho_{M_{W}}$ est un isomorphisme.
\hfill$\square$

\bigskip
On suppose maintenant que $\mathrm{H}_{V}^{*}X$ est libre comme $\mathrm{H}^{*}V$-module.

\medskip
On observe que si $\mathrm{H}_{V}^{*}X$ est libre comme $\mathrm{H}^{*}V$-module alors  $\mathrm{H}_{V/W}^{*}X^{W}$ est libre comme $\mathrm{H}^{*}V/W$-module. En effet, on a un isomorphisme de  $\mathrm{H}^{*}V/W$-$\mathrm{A}$-modules instables $\mathrm{H}_{V/W}^{*}X^{W}\cong\mathrm{Fix}_{(V,W)}\mathrm{H}_{V}^{*}X$ (Proposition \ref{cohFix}) et le corollaire \ref{sptorFix} conduit à l'énoncé suivant~: 

\begin{pro}\label{libre}
Soit $M$ un $\mathrm{H}^{*}V$-$\mathrm{A}$-module instable~; si $M$ est libre comme $\mathrm{H}^{*}V$-module, alors le $\mathrm{H}^{*}V/W$-$\mathrm{A}$-module instable $\mathrm{Fix}_{(V,W)}M$ est libre comme $\mathrm{H}^{*}V/W$-module.
\end{pro}

\medskip
\textit{Démonstration.} On a $\mathrm{Tor}_{1}^{\mathrm{H}^{*}V/W}(\mathbb{F}_{2},\mathrm{Fix}_{(V,W)}M)=0$ d'après \ref{sptorFix}.
\hfill$\square$

\bigskip
La proposition \ref{FixCtop} et l'hypothèse de récurrence (HR) montrent que le complexe $\mathrm{Fix}_{(V,W)}\hspace{2pt}(\mathop{\Sigma^{n}}\widetilde{\mathrm{C}}_{\mathrm{top}}^{\bullet}\hspace{1pt}X)$ est acyclique pour $W\not=0$. On en déduit que le complexe $\Psi(\mathop{\Sigma^{n}}\widetilde{\mathrm{C}}_{\mathrm{top}}^{\bullet}\hspace{1pt}X)$ est également acyclique. En effet on a
$$
\Psi(\mathop{\Sigma^{n}}\widetilde{\mathrm{C}}_{\mathrm{top}}^{\bullet}X)(W)
\hspace{4pt}\cong\hspace{4pt}
\mathrm{H}^{*}V\otimes_{\mathrm{H}^{*}V/W}\hspace{1pt}
\mathrm{Fix}_{(V,W)}(\mathop{\Sigma^{n}}\widetilde{\mathrm{C}}_{\mathrm{top}}^{\bullet}X)
$$
et $\mathrm{H}^{*}V$ est un $\mathrm{H}^{*}V/W$-module plat.  La coaugmentation
$$
\Psi(\mathop{\Sigma^{n}\hspace{-1pt}}\mathrm{H}_{V}^{*}X)
\to
\Psi(\mathop{\Sigma^{n}}\mathrm{C}_{\mathrm{top}}^{\bullet}\hspace{1pt}X)
$$
peut donc être vue comme une résolution de l'objet $\Psi(\mathop{\Sigma^{n}\hspace{-1pt}}\mathrm{H}_{V}^{*}X)$ dans la catégorie $(V_{\mathrm{tf}}\text{-}\mathcal{U})^{\mathcal{W}_{0}}$. Le point (c) de la proposition \ref{PsiCtop} dit que cette résolution est une résolution $\lim_{\mathcal{W}_{0}}$-acyclique ($\lim_{\mathcal{W}_{0}}$ désigne ici le foncteur exact à gauche $F\mapsto\lim_{\mathcal{W}_{0}}F$ de $(V_{\mathrm{tf}}\text{-}\mathcal{U})^{\mathcal{W}_{0}}$ dans $V_{\mathrm{tf}}\text{-}\mathcal{U}$), résolution tout aussi adaptée pour étudier les foncteurs dérivés de $\lim_{\mathcal{W}_{0}}$ qu'une résolution injective, si bien que l'on a~:

\smallskip
--\hspace{8pt}$\mathrm{H}^{p}\hspace{1pt}\lim_{\mathcal{W}_{0}}\Psi(\mathop{\Sigma^{n}}\widetilde{\mathrm{C}}_{\mathrm{top}}^{\bullet}\hspace{1pt}X)=0$ pour $p=-1,0$~;

\smallskip
--\hspace{8pt}$\mathrm{H}^{p}\hspace{1pt}\lim_{\mathcal{W}_{0}}\Psi(\mathop{\Sigma^{n}}\widetilde{\mathrm{C}}_{\mathrm{top}}^{\bullet}\hspace{1pt}X)
\cong\lim_{\mathcal{W}_{0}}^{p}\Psi(\mathop{\Sigma^{n}}\mathrm{H}_{V}^{*}X)$ pour $p>0$.

\smallskip
Comme le point (b) de \ref{suspensionEFix} montre que le foncteur $\Psi:V_{\mathrm{tf}}\text{-}\mathcal{U}\to(V_{\mathrm{tf}}\text{-}\mathcal{U})^{\mathcal{W}_{0}}$ ``commute à la suspension'', l'isomorphisme ci-dessus peut être remplacé par~:

\smallskip
--\hspace{8pt}$\mathrm{H}^{p}\hspace{1pt}\lim_{\mathcal{W}_{0}}\Psi(\mathop{\Sigma^{n}}\widetilde{\mathrm{C}}_{\mathrm{top}}^{\bullet}\hspace{1pt}X)
\cong\mathop{\Sigma^{n}}\lim_{\mathcal{W}_{0}}^{p}\Psi(\mathrm{H}_{V}^{*}X)$ pour $p>0$.

\medskip
La proposition \ref{derivePf} et le corollaire \ref{DS3} entraînent que l'homomorphisme $\rho_{\hspace{1pt}\mathrm{H}_{V}^{*}X}:\mathrm{H}_{V}^{*}X\to\lim_{\mathcal{W}_{0}}\Psi(\mathrm{H}_{V}^{*}X)$ est un isomorphisme ou ce qui revient au même que l'homomorphisme $\rho_{\hspace{1pt}\mathop{\Sigma^{n}}\mathrm{H}_{V}^{*}X}:\mathop{\Sigma^{n}}\mathrm{H}_{V}^{*}X\to\lim_{\mathcal{W}_{0}}\Psi(\mathop{\Sigma^{n}}\mathrm{H}_{V}^{*}X)$ est un isomorphisme. Ce dernier point et les points (d) et (b) de \ref{PsiCtop} montrent que l'homomorphisme de $(V_{\mathrm{tf}}\text{-}\mathcal{U})$-complexes
$$
\rho_{\hspace{2pt}\mathop{\Sigma^{n}}\hspace{-1pt}\widetilde{\mathrm{C}}_{\mathrm{top}}^{\bullet}\hspace{1pt}X}:
\mathop{\Sigma^{n}}\hspace{-1pt}\widetilde{\mathrm{C}}_{\mathrm{top}}^{\bullet}\hspace{1pt}X
\longrightarrow
{\lim}_{\mathcal{W}_{0}}\Psi(\mathop{\Sigma^{n}}\hspace{-1pt}\widetilde{\mathrm{C}}_{\mathrm{top}}^{\bullet}\hspace{1pt}X)
$$
est un épimorphisme dont le noyau est le $(V_{\mathrm{tf}}\text{-}\mathcal{U})$-complexe, disons $\mathrm{N}^{\bullet}$, défini par $\mathrm{N}^{n}=\mathop{\Sigma^{n}}\hspace{-1pt}\mathrm{C}_{\mathrm{top}}^{n}\hspace{1pt}X$ et $\mathrm{N}^{p}=0$ pour $p\not=n$.

\medskip
Les énoncés \ref{derivePf} et  \ref{DS3} montrent aussi que l'on a $\lim_{\mathcal{W}_{0}}^{p}\Psi(\mathrm{H}_{V}^{*}X)=0$ pour $0<p<n-1$. On obtient donc $\mathrm{H}^{p}\hspace{1pt}\lim_{\mathcal{W}_{0}}\Psi(\mathop{\Sigma^{n}}\widetilde{\mathrm{C}}_{\mathrm{top}}^{\bullet}\hspace{1pt}X)=0$ pour $p<n-1$ ce qui compte tenu de ce qui précède implique bien $\mathrm{H}^{p}\hspace{1pt}\mathop{\Sigma^{n}}\widetilde{\mathrm{C}}_{\mathrm{top}}^{\bullet}\hspace{1pt}X=0$, soit encore $\mathrm{H}^{p}\hspace{1pt}\widetilde{\mathrm{C}}_{\mathrm{top}}^{\bullet}\hspace{1pt}X=0$, pour $p<n-1$.
\hfill$\square$

\pagebreak

\bigskip
\textit{Seconde étape.} On montre que l'on a aussi $\mathrm{H}^{p}\hspace{1pt}\widetilde{\mathrm{C}}_{\mathrm{top}}^{\bullet}\hspace{1pt}X=0$ pour $p=n-1$ et~$p=n$.  En fait, cette étape n'utilise pas l'hypothèse ``$\mathrm{H}_{V}^{*}X$ libre comme $\mathrm{H}^{*}V$-module''~:

\medskip
\begin{lem}\label{top-n-2} Soient $V$ un $2$-groupe abélien élémentaire et $X$ un $V$-CW-complexe fini. Si l'on a $\mathrm{H}^{p}\hspace{1pt}\widetilde{\mathrm{C}}_{\mathrm{top}}^{\bullet}\hspace{1pt}X=0$ pour $p\leq\dim V-2$ alors $\widetilde{\mathrm{C}}_{\mathrm{top}}^{\bullet}\hspace{1pt}X$ est acyclique.
\end{lem}

\medskip
Ce lemme est un cas particulier d'un lemme très général concernant la suite spectrale cohomologique associée à une filtration décroissante d'un espace (voir \cite[Lemma 5.6]{AFP1}). Cependant, pour le confort du lecteur, nous en explicitons une démonstration ci-après.

\medskip
\textit{Démonstration.} Le cas $\dim V=0$ est trivial~; le cas $\dim V=1$ est implicitement traité au début de la démonstration de l'implication (i)$\Rightarrow$(ii) de \ref{Ctop}. On suppose donc $\dim V\geq 2$~; on pose à nouveau $n:=\dim V$.

\medskip
On considère la suite spectrale ``cohomologique", $(\mathrm{E}_{r})_{r\geq1}$, convergeant vers $\mathrm{H}_{V}^{*}X$, associée à la filtration $V$-équivariante de $X$ par les $\mathrm{F}_{p}X$, introduite au début de cette section.

\medskip
On se convainc tout d'abord que cette suite spectrale dégénère au terme $\mathrm{E}_{2}$.

\medskip
On a par définition $\mathrm{E}_{2}^{p,*}=\mathrm{H}^{p}\hspace{1pt}\mathrm{C}_{\mathrm{top}}^{\bullet}\hspace{1pt}X$~; on a donc par hypothèse $\mathrm{E}_{2}^{p,*}=0$ pour $p\not=0,n-1,n$. Il en résulte que les seules différentielles qui peuvent être non nulles, pour $r\geq 2$, sont les $\mathrm{d}_{r}:\mathrm{E}_{r}^{0,q}\to\mathrm{E}_{r}^{r,q+1-r}$. On montre que celles-ci sont nulles en introduisant la suite spectrale $({}^{\mathrm{tr}}\mathrm{E}_{r})_{r\geq1}$ associée à la filtration triviale de $X$ également évoquée au début de cette section, et le morphisme canonique de suites spectrales, disons $\eta$, de $({}^{\mathrm{tr}}\mathrm{E}_{r})_{r\geq1}$ dans $(\mathrm{E}_{r})_{r\geq1}$. On observe que la suite spectrale  $({}^{\mathrm{tr}}\mathrm{E}_{r})_{r\geq1}$ est vraiment triviale~: ${}^{\mathrm{tr}}\mathrm{E}_{1}^{p,*}=0$ pour $p\not=0$ et ${}^{\mathrm{tr}}\mathrm{E}_{1}^{0,*}={}^{\mathrm{tr}}\mathrm{E}_{r}^{0,*}=\mathrm{H}_{V}^{*}X$ pour tout $r\geq 1$.

\medskip
On suppose $r\geq 2$ et on contemple le diagramme commutatif
$$
\begin{CD}
{}^{\mathrm{tr}}\mathrm{E}_{r}^{0,q}
@>\mathrm{d}_{r}>>
{}^{\mathrm{tr}}\mathrm{E}_{r}^{r,q+1-r}=0
\\ @V\eta VV @V\eta VV \\
\mathrm{E}_{r}^{0,q}
@>\mathrm{d}_{r}>>
\mathrm{E}_{r}^{r,q+1-r}
\end{CD}
$$
fourni par le morphisme de suites spectrales $\eta$. Pour $r=2$ on sait que  l'homomorphisme $\eta:{}^{\mathrm{tr}}\mathrm{E}_{2}^{0,q}\to\mathrm{E}_{2}^{0,q}$ est un isomorphisme : c'est une conséquence de l'égalité $\mathrm{H}^{0}\hspace{1pt}\widetilde{\mathrm{C}}_{\mathrm{top}}^{\bullet}\hspace{1pt}X=0$~; ceci montre que la différentielle $\mathrm{d}_{2}:\mathrm{E}_{2}^{0,q}\to\mathrm{E}_{2}^{2,q-1}$ est nulle et que l'homomorphisme $\eta:{}^{\mathrm{tr}}\mathrm{E}_{3}^{0,q}\to\mathrm{E}_{3}^{0,q}$ est un isomorphisme. En itérant l'argument, on montre que toutes les différentielles $\mathrm{d}_{r}:\mathrm{E}_{r}^{0,q}\to\mathrm{E}_{r}^{r,q+1-r}$ sont nulles pour $r\geq 2$.

\medskip
On sait maintenant que la suite spectrale $(\mathrm{E}_{r})_{r\geq1}$ dégénère au terme $\mathrm{E}_{2}$. En considérant à nouveau le morphisme $\eta$, on se convainc que l'homomorphisme canonique $\mathrm{E}_{\infty}^{0,*}\to\mathrm{H}_{V}^{*}X$ est un isomorphisme. Il en résulte $\mathrm{E}_{2}^{p,*}=0$ pour $p=n-1$ et $p=n$, c'est-à-dire $\mathrm{H}^{p}\hspace{1pt}\widetilde{\mathrm{C}}_{\mathrm{top}}^{\bullet}\hspace{1pt}X=0$ pour $p=n-1$ et $p=n$.
\hfill$\square$

\bigskip
\begin{cor}\label{instable}
Soient $V$ un $2$-groupe abélien élémentaire et $X$ un $V$-CW-complexe fini.  Si $\mathrm{H}^{*}_{V}X$ est libre comme $\mathrm{H}^{*}V$-module alors le $\mathrm{A}$-module
$$
\mathop{\Sigma^{-\dim V}}\mathrm{H}^{*}_{V}(X,\mathrm{Sing}_{V}X)
$$
est instable.
\end{cor}

\bigskip
\textit{Démonstration.} On procède toujours par récurrence sur l'entier $n=\dim V$. Là encore le cas $n=0$ est trivial~; on suppose donc $n\geq 1$ et on montre que si l'énoncé est vrai pour $\dim V=n-1$ alors il l'est aussi pour $\dim V=n$.

\medskip
Compte tenu du point (a) de la proposition \ref{Ctopefini}, le théorème \ref{Ctop} dit en particulier que l'on a un épimorphisme de $\mathrm{A}$-modules
$$
\bigoplus_{\mathop{\mathrm{dim}}W=1}
\hspace{-4pt}
\mathrm{H}^{*}V\otimes_{\mathrm{H}^{*}V/W}
\hspace{2pt}
\mathop{\Sigma^{-(n-1)}}\mathrm{H}_{V/W}^{*}(X^{W},\mathrm{Sing}_{V/W}X^{W})\to
\mathop{\Sigma^{-n}}\mathrm{H}^{*}_{V}(X,\mathrm{Sing}_{V}X).
$$
Pour tous les $W$ ci-dessus, $X^{W}$ est un $V/W$-CW-complexe fini avec $\mathrm{H}_{V/W}^{*}X^{W}$ libre comme $\mathrm{H}^{*}V/W$-module (Proposition \ref{libre})~; la source de l'épimorphisme en question est donc un $\mathrm{A}$-module instable d'après l'hypothèse de récurrence. Du coup il en est de même pour le but.
\hfill$\square$

\medskip
\begin{cor}\label{instablebis}
Soient $V$ un $2$-groupe abélien élémentaire et $X$ un $V$-CW-complexe fini.  Si $\mathrm{H}^{*}_{V}X$ est libre comme $\mathrm{H}^{*}V$-module alors le $\mathrm{A}$-module $\mathrm{C}_{\mathrm{top}}^{p}\hspace{1pt}X$ est instable pour tout $p$.
\end{cor}

\bigskip
L'énoncé ci-dessous  précise l'énoncé \ref{instable}:

\begin{theo}\label{Cntop}
Soient $V$ un $2$-groupe abélien élémentaire et $X$ un $V$-CW-complexe fini.  Si $\mathrm{H}^{*}_{V}X$ est libre comme $\mathrm{H}^{*}V$-module alors on a un isomorphisme canonique de $\mathrm{H}^{*}V_{\mathrm{tf}}$-$\mathrm{A}$-modules instables
$$
\mathop{\Sigma^{-n}}\mathrm{H}^{*}_{V}(X,\mathrm{Sing}_{V}X)
\hspace{4pt}\cong\hspace{4pt}
\mathrm{R}^{n}\mathrm{Pf}\hspace{2pt}\mathrm{H}^{*}_{V}X
$$
($n:=\dim_{\mathbb{Z}/2}V$).
\end{theo}

\bigskip
\textit{Démonstration.} Le cas $n=0$ est trivial. Le cas $n=1$ sera traité séparément. On suppose $n\geq 2$ et on reprend la suite exacte de complexes
$$
\begin{CD}
0@>>>\mathrm{N}^{\bullet}
@>>>\mathop{\Sigma^{n}}\hspace{-1pt}\widetilde{\mathrm{C}}_{\mathrm{top}}^{\bullet}\hspace{1pt}X
@>\rho_{\hspace{2pt}\mathop{\Sigma^{n}}\hspace{-1pt}\widetilde{\mathrm{C}}_{\mathrm{top}}^{\bullet}\hspace{1pt}X}>>
\lim_{\mathcal{W}_{0}}\Psi(\mathop{\Sigma^{n}}\hspace{-1pt}\widetilde{\mathrm{C}}_{\mathrm{top}}^{\bullet}\hspace{1pt}X)
@>>>0
\end{CD}
$$
implicitement introduite à la fin de la première étape de la démonstration de l'implication (i)$\Rightarrow$(ii) de \ref{Ctop}. Comme le complexe $\mathop{\Sigma^{n}}\hspace{-1pt}\widetilde{\mathrm{C}}_{\mathrm{top}}^{\bullet}\hspace{1pt}X$ est acyclique si $\mathrm{H}^{*}_{V}X$ est libre comme $\mathrm{H}^{*}V$-module, le connectant
$$
\mathrm{H}^{n-1}\hspace{1pt}{\lim}_{\mathcal{W}_{0}}\Psi(\mathop{\Sigma^{n}}\widetilde{\mathrm{C}}_{\mathrm{top}}^{\bullet}\hspace{1pt}X)
\longrightarrow
\mathrm{H}^{n}\hspace{1pt}\mathrm{N}^{\bullet}=\mathrm{H}^{*}_{V}(X,\mathrm{Sing}_{V}X)
$$
est un $(V_{\mathrm{tf}}\text{-}\mathcal{U})$-isomorphisme. On a déjà vu que l'on a un $(V_{\mathrm{tf}}\text{-}\mathcal{U})$-isomorphisme
$\mathrm{H}^{n-1}\lim_{\mathcal{W}_{0}}\Psi(\mathop{\Sigma^{n}}\widetilde{\mathrm{C}}_{\mathrm{top}}^{\bullet}\hspace{1pt}X)
\cong\mathop{\Sigma^{n}}\lim_{\mathcal{W}_{0}}^{n-1}\Psi(\mathrm{H}_{V}^{*}X)$ et \ref{derivePf} dit que $\lim_{\mathcal{W}_{0}}^{n-1}\Psi(\mathrm{H}_{V}^{*}X)$ est $(V_{\mathrm{tf}}\text{-}\mathcal{U})$-isomorphe à $\mathrm{R}^{n}\mathrm{Pf}\hspace{2pt}\mathrm{H}^{*}_{V}X$.
\hfill$\square$

\medskip
Le cas $n=1$. On a une suite exacte
$$
0\to\mathrm{H}^{*}_{V}X\to\mathrm{H}^{*}V\otimes\mathrm{H}^{*}X^{V}\to{\mathop{\Sigma}}^{-1}\mathrm{H}^{*}_{V}(X,X^{V})\to 0
$$
dans la catégorie $V_{\mathrm{tf}}\text{-}\mathcal{U}$ et $\mathrm{H}^{*}V\otimes\mathrm{H}^{*}X^{V}$ est $\mathrm{Pf}$-acyclique d'après \ref{DS3}.
\hfill$\square$

\bigskip
La démonstration que nous venons de donner du théorème \ref{Cntop} dans le cas $n=1$ sera généralisée en \ref{theoefini}. Ce théorème sera revisité en \ref{CalgCtop}.

\vspace{0.75cm}
\textit{Bande-annonce}

\medskip
Comme expliqué dans l'introduction le théorème \ref{Ctop} admet une généralisation, à savoir le théorème \ref{genafp} dû à Allday, Franz et Puppe. Le théorème \ref{introgentop} est une version de ce théorème que nous démontrerons en section 6 (l'énoncé \ref{gentop} reproduit \l'énoncé \ref{introgentop}).

\pagebreak

\sect{Le complexe algébrique}
 
Dans cette section nous associons à un $\mathrm{H}^{*}V$-$\mathrm{A}$-module instable $M$ un complexe de cochaînes
$$
\hspace{24pt}
\mathrm{C}^{0}\hspace{1pt}M\to\mathrm{C}^{1}\hspace{1pt}M\to\mathrm{C}^{2}\hspace{1pt}M\to\ldots\mathrm{C}^{p}\hspace{1pt}M\to\ldots\to\mathrm{C}^{n}\hspace{1pt}M
\hspace{24pt},
$$
coaugmenté par $M$~; nous le notons $\widetilde{\mathrm{C}}^{\bullet}\hspace{1pt}M$ (on a donc $\widetilde{\mathrm{C}}^{-1}\hspace{1pt}M=M$). Nous étudions ce complexe quand $M$ est de type fini comme $\mathrm{H}^{*}V$-module. Dans la prochaine section nous comparerons le complexe $\widetilde{\mathrm{C}}_{\mathrm{top}}^{\bullet}\hspace{1pt}X$, associé à un $V$-CW-complexe fini $X$, au complexe $\widetilde{\mathrm{C}}_{\mathrm{alg}}^{\bullet}\hspace{1pt}X:=\widetilde{\mathrm{C}}^{\bullet}\hspace{1pt}\mathrm{H}_{V}X$ que nous appelons le complexe algébrique associé à $X$. Nous montrerons en particulier que si $\mathrm{H}_{V}^{*}X$ est libre comme $\mathrm{H}^{*}V$-module alors les deux complexes sont isomorphes.

\vspace{0.625cm}
\textsc{Filtration d'un $\mathrm{H}^{*}V$-$\mathrm{A}$-module instable}

\medskip
Soient $V$ un $2$-groupe abélien élémentaire de dimension $n$, $W\subset V$ un sous-groupe et $M$ un $\mathrm{H}^{*}V$-$\mathrm{A}$-module instable.

\medskip
On rappelle que l'on note
$$
{\eta_{(V,W)}}_{M}:M\to\mathrm{H}^{*}V\otimes_{\mathrm{H}^{*}V/W}\mathrm{Fix}_{(V,W)}M
$$
l'unité de l'adjonction de la paire de foncteurs adjoints $(\mathrm{e}_{(V,W)},\mathrm{Fix}_{(V,W)}$).

\medskip
On définit une filtration décroissante de $M$ par des sous-$\mathrm{H}^{*}V$-$\mathrm{A}$-modules instables
$$
M=\mathrm{F}^{0}M\supset\mathrm{F}^{1}M\supset\ldots\supset\mathrm{F}^{n}M\supset\mathrm{F}^{n+1}M=0
$$
en posant
$$
\mathrm{F}^{p}M
\hspace{4pt}:=\hspace{4pt}
\bigcap_{\mathop{\mathrm{codim}}W<p}
\ker\hspace{1pt}({\eta_{(V,W)}}_{M}:M\to\mathrm{H}^{*}V\otimes_{\mathrm{H}^{*}V/W}\mathrm{Fix}_{(V,W)}M)
$$
(l'égalité $M=\mathrm{F}^{0}M$ tient au fait que l'ensemble $\{W; \mathop{\mathrm{codim}}W<0\}$ est vide~!).  Compte tenu de la proposition \ref{commutativity}, on peut alternativement définir cette filtration en posant
$$
\hspace{24pt}
\mathrm{F}^{p}M
\hspace{4pt}:=\hspace{4pt}
\bigcap_{\mathop{\mathrm{codim}}W<p}
\ker\hspace{1pt}({\rho_{(V,W)}}_{M}:M\to\mathrm{EFix}_{(V,W)}M)
\hspace{24pt}.
$$
Soient $W_{0}$ et $W_{1}$ deux sous-groupes de $V$ avec $W_{0}\subset W_{1}$, la proposition \ref{Wfonctoriel} dit que l'on a $\rho_{(V,W_{1})}=\rho(W_{0},W_{1})\circ\rho_{(V,W_{0})}$~; on peut donc également dans  la définition de $\mathrm{F}^{p}$ remplacer $\bigcap_{\hspace{2pt}\mathop{\mathrm{codim}}W<p}$ par  $\bigcap_{\hspace{2pt}\mathop{\mathrm{codim}}W=p-1}$.

\medskip
Trois  commentaires concernant cette définition~:

\medskip
1) L'exactitude des foncteurs $\mathrm{EFix}$ implique~:

\begin{pro}\label{incfilt} Soient $M$ un $\mathrm{H}^{*}V$-$\mathrm{A}$-module instable, $M'\subset M$ un sous-objet, et $p$ un entier naturel. On a l'égalité $\mathrm{F}^{p}M'=M'\cap\mathrm{F}^{p}M$.
\end{pro}

\medskip
2) Le scholie \ref{etaFix} (ou la proposition \ref{cohEFix}) conduit à l'énoncé ci-dessous et fournit du même coup une motivation pour la définition ci-dessus (voir la partie de l'introduction qui suit l'intertitre ``Le complexe algébrique")~:

\begin{pro}\label{Filttop2} Soit $X$ un $V$-CW-complexe fini. Alors on a
$$
\hspace{24pt}
\mathrm{F}^{p}\hspace{1pt}\mathrm{H}_{V}^{*}X
\hspace{4pt}=\hspace{4pt}
\bigcap_{\mathop{\mathrm{codim}}W<p}
\ker\hspace{1pt}(\mathrm{H}_{V}^{*}X\to\mathrm{H}_{V}^{*}X^{W})
\hspace{24pt}.
$$
\end{pro}

3) La théorie de Smith algébrique \cite{DWsmith2}\cite[Théorème 0.5]{LZsmith} et \cite[\S2, Corollaire]{Ser} impliquent que $\mathrm{F}^{1}M$ est la $\mathrm{H}^{*}V$-torsion de $M$ (sans aucune hypothèse sur le $\mathrm{H}^{*}V$-$\mathrm{A}$-module instable $M$,  voir \ref{Serre5}). Concernant $\mathrm{F}^{n}M$, on a~:

\begin{pro}\label{Fn=Pf} Soit $M$ un $\mathrm{H}^{*}V$-$\mathrm{A}$-module instable. Si l'on suppose que $M$ est de type fini comme $\mathrm{H}^{*}V$-module alors on a $\mathrm{F}^{n}M$ est la partie finie $\mathrm{Pf}M$ de $M$.
\end{pro}

\smallskip
\textit{Démonstration.} Soit $N\subset M$ un sous-$\mathrm{H}^{*}V$-$\mathrm{A}$-module instable fini~; l'implication $(iv)\Rightarrow(ii\text{-}bis)$ de  \ref{caracterisationfini} montre $N\subset\mathrm{F}^{n}M$. Puisque les foncteurs $\mathrm{EFix}_{(V,W)}$ sont exacts, on a  $\mathrm{EFix}_{(V,W)}(\mathrm{F}^{n}M)=0$ pour $\dim W=1$.  Si l'on suppose $M$ de type fini comme $\mathrm{H}^{*}V$-module alors il en est de même pour $\mathrm{F}^{n}M$ et l'implication $(ii\text{-}bis)\Rightarrow(iv)$ de  \ref{caracterisationfini} montre que $\mathrm{F}^{n}M$ est fini.
\hfill$\square$

\medskip
On voit donc que $\mathrm{F}^{1}M$ et $\mathrm{F}^{n}M$ peuvent être définis en termes de la structure de $\mathrm{H}^{*}V$-module de $M$, au moins sous l'hypothèse que celle-ci est de type fini. On verra dans la section \ref{Appendice} que c'est en fait le  cas de tous les $\mathrm{F}^{p}M$ (même sans l'hypothèse que $M$ est de type fini comme $\mathrm{H}^{*}V$-module).

\bigskip
On sera aussi amené à considérer le gradué de la filtration de $M$ par les $\mathrm{F}^{p}M$~: 
$$
\hspace{36pt}
\mathrm{Gr}^{p}M
\hspace{4pt}:=\hspace{4pt}
\mathrm{F}^{p}M\hspace{1pt}/\hspace{1pt}\mathrm{F}^{p+1}M
\hspace{12pt},\hspace{12pt}0\leq p\leq n\hspace{12pt}.
$$

\medskip
Il est clair que les applications $M\mapsto\mathrm{F}^{p}M$ et $M\mapsto\mathrm{Gr}^{p}M$ peuvent être considérées comme des endofoncteurs de $V\text{-}\mathcal{U}$ et comme des endofoncteurs de $V_{\mathrm{tf}}\text{-}\mathcal{U}$ si l'on suppose que $M$ est de type fini comme $\mathrm{H}^{*}V$-module. Ces endofoncteurs ``commutent'' aux sommes directes dans $V\text{-}\mathcal{U}$ et aux sommes directes finies dans $V_{\mathrm{tf}}\text{-}\mathcal{U}$ (qui coïncident avec les produit finis).

\bigskip
On dégage maintenant quelques propriétés de la filtration introduite ci-dessus (et de son gradué) lorsque $M$ est de type fini comme $\mathrm{H}^{*}V$-module.

\medskip
On explicite d'abord cette filtration dans le cas où $M$ est un objet injectif indécomposable de $V_{\mathrm{tf}}\text{-}\mathcal{U}$, c'est-à-dire, d'après le théorème \ref{injectifs3}, dans le cas $M=\mathrm{H}^{*}V\otimes_{\mathrm{H}^{*}V/U}\mathrm{J}_{V/U}(k)$ avec $(U,k)\in\mathcal{W}\times\mathbb{N}$. Le $\mathrm{H}^{*}V/U$-$\mathrm{A}$-module instable $\mathrm{J}_{V/U}(k)$ est fini et il ne coûte en fait pas plus cher de supposer plus généralement $M=\mathrm{H}^{*}V\otimes_{\mathrm{H}^{*}V/U}N$ avec $N$ un $\mathrm{H}^{*}V/U$-$\mathrm{A}$-module instable fini.

\begin{pro}\label{filtparticulier} Soient $U\subset V$ un sous-groupe et $N$ un $\mathrm{H}^{*}V/U$-$\mathrm{A}$-module instable fini. Alors on a
$$
\mathrm{F}^{p}\hspace{2pt}(\mathrm{H}^{*}V\otimes_{\mathrm{H}^{*}V/U}N)
\hspace{4pt}=\hspace{4pt}
\begin{cases}
\mathrm{H}^{*}V\otimes_{\mathrm{H}^{*}V/U}N & \text{pour $p\leq\mathop{\mathrm{codim}}U$}, \\
0 & \text{pour $p>\mathop{\mathrm{codim}}U$}.
\end{cases}
$$
\end{pro}

\medskip
\textit{Démonstration.} Conséquence de la proposition \ref{proclefbis}.
\hfill$\square$

\bigskip
L'énoncé ci-dessus ouvre la voie au suivant~:

\medskip
\begin{pro-def}\label{grad} Soient $M$ un $\mathrm{H}^{*}V_{\mathrm{tf}}$-$\mathrm{A}$-module instable et $p$ un entier avec $0\leq p\leq n$. Soit
$$
\Pi^{p}_{M}:M\longrightarrow\bigoplus_{\mathop{\mathrm{codim}}W=p}
\mathrm{H}^{*}V\otimes_{\mathrm{H}^{*}V/W}\mathrm{Fix}_{(V,W)}M
$$
le produit des homomorphismes $\eta_{(V,W),M}$ pour $\mathop{\mathrm{codim}}W=p$.

\medskip
{\em (a)} La restriction de $\Pi^{p}_{M}$ à $\mathrm{F}^{p+1}M$ est nulle.

\medskip
{\em (b)} L'image de la restriction de $\Pi^{p}_{M}$ à $\mathrm{F}^{p}M$ est contenue dans le sous-$\mathrm{H}^{*}V$-$\mathrm{A}$-module instable
$$
\hspace{24pt}
\bigoplus_{\mathop{\mathrm{codim}}W=p}
\mathrm{H}^{*}V\otimes_{\mathrm{H}^{*}V/W}\mathrm{Pf}_{V/W}\hspace{1pt}(\mathrm{Fix}_{(V,W)}M)
\hspace{24pt}.
$$
{\em (Décodons la notation~: $\mathrm{Fix}_{(V,W)}M$ est d'après \ref{typefiniFix} un $\mathrm{H}^{*}V/W$-$\mathrm{A}$-module instable de type fini comme $\mathrm{H}^{*}V/W$-module et $\mathrm{Pf}_{V/W}$ désigne l'endofoncteur ``partie finie'' de $V/W_{\mathrm{tf}}\text{-}\mathcal{U}$. La somme directe ci-dessus s'identifie bien à un sous-module de la précédente car $\mathrm{H}^{*}V$ est plat sur $\mathrm{H}^{*}V/W$.)}

\medskip
On note
$$
\overline{\Pi}^{\hspace{1pt}p}_{M}
\hspace{1pt}:\hspace{1pt}
\mathrm{Gr}^{p}M
\longrightarrow
\bigoplus_{\mathop{\mathrm{codim}}W=p}
\mathrm{H}^{*}V\otimes_{\mathrm{H}^{*}V/W}\mathrm{Pf}_{V/W}\hspace{1pt}(\mathrm{Fix}_{(V,W)}M)
$$
l'homomorphisme de $\mathrm{H}^{*}V$-$\mathrm{A}$-modules instables induit par $\Pi^{p}_{M}$.

\medskip
{\em (c)} L'homomorphisme $\overline{\Pi}^{\hspace{1pt}p}_{M}$ est un isomorphisme si $M$ est un $V_{\mathrm{tf}}\text{-}\mathcal{U}$-injectif.

\end{pro-def}

\bigskip
\textit{Démonstration du (a).} Conséquence de la définition même de $\mathrm{F}^{p+1}M$.
\hfill$\square$

\bigskip
\textit{Démonstration du (b).}

\medskip
1) On commence par le cas $M=\mathrm{H}^{*}V\otimes_{\mathrm{H}^{*}V/U}N$ avec $U$ un sous-groupe de $V$ et $N$ un $\mathrm{H}^{*}V/U$-$\mathrm{A}$-module instable fini~:

\smallskip
-- Si l'on a $\mathop{\mathrm{codim}}U<p$ alors on a  $\mathrm{F}^{p}M=0$ d'après \ref{filtparticulier} et  (b) est trivialement vérifié.

\smallskip
-- Si l'on a $\mathop{\mathrm{codim}}U>p$ alors on a $W\not\subset U$ pour tout $W$ avec $\mathop{\mathrm{codim}}W=p$. Il en résulte $\Pi^{p}_{M}=0$ d'après \ref{proclefbis} et  (b) est encore trivialement vérifié.

\smallskip
-- Si l'on a $\mathop{\mathrm{codim}}U=p$ alors il existe un seul $W$ avec $\mathop{\mathrm{codim}}W=p$ et $W\subset U$ à savoir $U$ lui-même si bien que, compte tenu de \ref{proclef} et \ref{EFix}, le seul facteur non nul de $\bigoplus_{\mathop{\mathrm{codim}}W=p}
\mathrm{H}^{*}V\otimes_{\mathrm{H}^{*}V/W}\mathrm{Fix}_{(V,W)}M$ est $\mathrm{H}^{*}V\otimes_{\mathrm{H}^{*}V/U}\mathrm{Fix}_{(V,U)}M$. On a $\mathrm{Fix}_{(V,U)}M\cong N$ d'après \ref{corproclef} et donc $\mathrm{Pf}_{V/U}\hspace{1pt}(\mathrm{Fix}_{(V,U)}M)=\mathrm{Fix}_{(V,U)}M$~; le point (b) est bien vérifié.

\medskip
2) D'après \ref{injectifs3} un $V_{\mathrm{tf}}\text{-}\mathcal{U}$-injectif est une somme directe finie de $\mathrm{H}^{*}V$-$\mathrm{A}$-modules instables du type considéré ci-dessus~; le point (b) est donc vérifié si $M$ est $V_{\mathrm{tf}}\text{-}\mathcal{U}$-injectif.

\medskip
3) On en vient enfin au cas général. D'après \ref{injectifs2} il existe un homomorphisme injectif de $\mathrm{H}^{*}V$-$\mathrm{A}$-modules instables $i: M\to I$ tel que $I$ est $V_{\mathrm{tf}}\text{-}\mathcal{U}$-injectif. Le diagramme suivant~:
$$
\hspace{-42pt}
\begin{CD}
\mathrm{F}^{p}M@>\text{restriction de}\hspace{2.5pt}\Pi^{p}_{M}>>
\bigoplus_{\mathop{\mathrm{codim}}W=p}\hspace{4pt}
\mathrm{H}^{*}V\otimes_{\mathrm{H}^{*}V/W}\mathrm{Fix}_{(V,W)}M \\
@V\mathrm{F}^{p}\hspace{1pt}iVV
@VV{\bigoplus_{\mathop{\mathrm{codim}}W=p}\hspace{2pt}\mathrm{H}^{*}V\otimes_{\mathrm{H}^{*}V/W}\mathrm{Fix}_{(V,W)}\hspace{1pt}i}V \\
\mathrm{F}^{p}I@>\text{restriction de}\hspace{2.5pt}\Pi^{p}_{I}>>
\bigoplus_{\mathop{\mathrm{codim}}W=p}\hspace{4pt}
\mathrm{H}^{*}V\otimes_{\mathrm{H}^{*}V/W}\mathrm{Fix}_{(V,W)}I
\end{CD}
$$
est commutatif (fonctorialité). La contemplation de ce diagramme montre que le point (b) est vérifié pour $M$. Détaillons un peu. On pose $P=\mathrm{Fix}_{(V,W)}M$, $J=\mathrm{Fix}_{(V,W)}I$ et $j=\mathrm{Fix}_{(V,W)}\hspace{1pt}i$. L'homomorphisme $j$ est injectif puisque le foncteur $\mathrm{Fix}_{(V,W)}$ est exact~; il en résulte $j^{-1}(\mathrm{Pf}_{V/W}(J))=\mathrm{Pf}_{V/W}(P)$. Soient $x$ un élément de $\mathrm{F}^{p}M$ et $y$ son image dans $\mathrm{H}^{*}V\otimes_{\mathrm{H}^{*}V/W}P$~; le fait que le point~(b) soit vérifié pour $I$ entraîne que $y$ appartient à l'image inverse
$$
\hspace{24pt}
(\mathrm{H}^{*}V\otimes_{\mathrm{H}^{*}V/W}j)^{-1}(\mathrm{H}^{*}V\otimes_{\mathrm{H}^{*}V/W}\mathrm{Pf}_{V/W}(J))
\hspace{24pt}.
$$
Or cette image inverse s'identifie à
$$
\hspace{24pt}
\mathrm{H}^{*}V\otimes_{\mathrm{H}^{*}V/W}j^{-1}(\mathrm{Pf}_{V/W}(J))
\hspace{4pt}=\hspace{4pt}
\mathrm{H}^{*}V\otimes_{\mathrm{H}^{*}V/W}\mathrm{Pf}_{V/W}(P)
\hspace{24pt};
$$
en effet, puisque $\mathrm{H}^{*}V$ est un $\mathrm{H}^{*}V/W$-module plat, le produit tensoriel par $\mathrm{H}^{*}V$ ``commute aux produits fibrés''.
\hfill$\square$

\bigskip
\textit{Démonstration du (c).} A nouveau il suffit de vérifier que $\overline{\Pi}^{\hspace{1pt}p}_{M}$ est un isomorphisme pour $M=\mathrm{H}^{*}V\otimes_{\mathrm{H}^{*}V/U}N$ avec $U$ un sous-groupe de $V$ et $N$ un $\mathrm{H}^{*}V/U$-$\mathrm{A}$-module instable fini ; pour cela on utilise encore \ref{proclefbis}.
\hfill$\square$

\vspace{0.75cm}
\textsc{La notion de $\mathrm{H}^{*}V$-$\mathrm{A}$-module instable e-fini}

\medskip
La  proposition \ref{filtparticulier} et la démonstration des points (b) et (c) de la proposition \ref{grad} suggère d'introduire la définition \textit{ad hoc} suivante~:

\bigskip
\begin{defi}\label{defefini} Nous dirons qu'un $\mathrm{H}^{*}V$-$\mathrm{A}$-module instable $M$ est {\em e-fini} s'il est isomorphe à une somme directe de la forme
$$
\bigoplus_{U\in\mathcal{W}}
\mathrm{H}^{*}V\otimes_{\mathrm{H}^{*}V/U}N_{U}
$$
avec $N_{U}$ un $\mathrm{H}^{*}V/U$-$\mathrm{A}$-module instable fini.
\end{defi}

\bigskip
\begin{exple}\label{injefini} Un $V_{\mathrm{tf}}\text{-}\mathcal{U}$-injectif est e-fini.
\end{exple}

\bigskip
Ce que l'on a appris en démontrant \ref{filtparticulier} et \ref{grad} conduit au scholie suivant~:

\begin{scho}\label{schoefini} Soient $N_{U}$, $U\in\mathcal{W}$, des $\mathrm{H}^{*}V/U$-$\mathrm{A}$-modules instable finis~; on pose
$$
\hspace{24pt}
M
\hspace{4pt}=\hspace{4pt}
\bigoplus_{U\in\mathcal{W}}
\mathrm{H}^{*}V\otimes_{\mathrm{H}^{*}V/U}N_{U}
\hspace{24pt}.
$$

\medskip {\em (a)} Les $\mathrm{H}^{*}V/U$-$\mathrm{A}$-modules instables finis $N_{U}$ sont uniquement déterminés, à isomorphisme près, en fonction du $\mathrm{H}^{*}V$-$\mathrm{A}$-module instable $M$. Plus précisé\-ment on a~:
$$
\hspace{24pt}
N_{U}
\hspace{4pt}\cong\hspace{4pt}
\mathrm{Pf}_{V/U}\hspace{1pt}(\mathrm{Fix}_{(V,U)}M)
\hspace{24pt}.
$$

\medskip {\em (b)} Soit $p$ un entier. On a~:
$$
\hspace{8pt}
\mathrm{F}^{p}\hspace{1pt}M
\hspace{4pt}\cong\hspace{4pt}
\bigoplus_{\mathop{\mathrm{codim}}U\geq p}
\mathrm{H}^{*}V\otimes_{\mathrm{H}^{*}V/U}N_{U}
\hspace{8pt}\text{et}\hspace{8pt}
\mathrm{Gr}^{p}\hspace{1pt}M
\hspace{4pt}\cong\hspace{4pt}
\bigoplus_{\mathop{\mathrm{codim}}U=p}
\mathrm{H}^{*}V\otimes_{\mathrm{H}^{*}V/U}N_{U}
\hspace{8pt}.
$$
\end{scho}

\bigskip
Les propositions \ref{annderivePf2} et \ref{annderivePf3} fournissent l'énoncé suivant~:

\medskip
\begin{pro}\label{proefini}Soit $M$ un $\mathrm{H}^{*}V$-$\mathrm{A}$-module instable e-fini. Alors on a $\mathrm{R}^{p}\mathrm{Pf}M=0$ pour $p>0$~; en d'autres termes $M$ est {\em $\mathrm{Pf}$-acyclique}.
\end{pro}

\medskip
\begin{cor-def}\label{cordefefini}Soient $M$ un $\mathrm{H}^{*}V_{\mathrm{tf}}$-$\mathrm{A}$-module instable et
$$
0\to M\to C^{0}\to C^{1}\to\ldots\to C^{p}\to\ldots
$$
une résolution de $M$ dans la catégorie $V_{\mathrm{tf}}$-$\mathcal{U}$, telle que  $C^{p}$ est e-fini pour tout~$p\geq 0$. Alors on a $\mathrm{R}^{p}\mathrm{Pf}M\cong\mathrm{H}^{p}\mathrm{Pf}C^{\bullet}$ pour tout $p\geq 0$.

\smallskip
On dira qu'une résolution comme ci-dessus est une {\em résolution e-finie} de $M$.
\end{cor-def}

\begin{cor}\label{corefini} Soient $M$ un $\mathrm{H}^{*}V_{\mathrm{tf}}$-$\mathrm{A}$-module instable et $S$ un $\mathrm{A}$-module instable fini~; on munit le produit tensoriel
$S\otimes M$ de sa structure évidente de $\mathrm{H}^{*}V$-$\mathrm{A}$-module instable. Alors on a un isomorphisme canonique de $\mathrm{H}^{*}V$-$\mathrm{A}$-modules instables (finis)
$$
\mathrm{R}^{p}\mathrm{Pf}\hspace{2pt}(S\otimes M)
\hspace{4pt}\cong\hspace{4pt}
S\otimes
\mathrm{R}^{p}\mathrm{Pf}\hspace{1pt}M
$$
pour tout $p\geq 0$ (en particulier les foncteurs dérivés $\mathrm{R}^{p}\mathrm{Pf}$ ``commutent à la suspension'').
\end{cor}

\bigskip
\textit{Démonstration.} Soit $M\to I^{\bullet}$ une résolution injective de $M$ dans la catégorie $V_{\mathrm{tf}}$-$\mathcal{U}$. La résolution $S\otimes M\to S\otimes I^{\bullet}$ n'est pas en général injective mais elle est toujours e-finie.
\hfill$\square$

\bigskip
Pour clore notre discussion sur la notion de $\mathrm{H}^{*}V$-$\mathrm{A}$-module instable e-fini, nous montrons que les énoncés \ref{cordefefini} et \ref{corefini} conduisent à une démonstration alternative du théorème \ref{Cntop}~:

\bigskip
\begin{exple}\label{theoefini} Soit $X$ un $V$-CW-complexe fini avec $\mathrm{H}^{*}_{V}X$ libre comme $\mathrm{H}^{*}V$-module.  Le théorème \ref{Ctop} dit que
$$
0\to\Sigma^{n}\hspace{1pt}\hspace{1pt}\mathrm{H}^{*}_{V}X
\to\Sigma^{n}\hspace{1pt}\mathrm{C}_{\mathrm{top}}^{0}\hspace{1pt}X
\to\ldots
\to\Sigma^{n}\hspace{1pt}\mathrm{C}_{\mathrm{top}}^{p}\hspace{1pt}X
\to\ldots
\to\Sigma^{n}\hspace{1pt}\mathrm{C}_{\mathrm{top}}^{n}\hspace{1pt}X\to 0
$$
est une résolution de $\Sigma^{n}\hspace{1pt}\hspace{1pt}\mathrm{H}^{*}_{V}X$ dans la catégorie $V_{\mathrm{tf}}$-$\mathcal{U}$. On fait les constatations suivantes~:

\smallskip
-- $\Sigma^{n}\hspace{1pt}\mathrm{C}_{\mathrm{top}}^{p}\hspace{1pt}X$ est e-fini pour tout $p$~;

\smallskip
-- $\mathrm{Pf}\hspace{1pt}(\Sigma^{n}\hspace{1pt}\mathrm{C}_{\mathrm{top}}^{p}\hspace{1pt}X$)=0 pour tout $p<n$~;

\smallskip
-- $\Sigma^{n}\hspace{1pt}\mathrm{C}_{\mathrm{top}}^{}\hspace{1pt}X\cong\mathrm{H}^{*}_{V}(X,\mathrm{Sing}_{V}X)$ est fini.

\smallskip
On a donc $\mathrm{R}^{n}\mathrm{Pf}\hspace{1pt}(\Sigma^{n}\hspace{1pt}\mathrm{H}^{*}_{V}X)\cong\mathrm{H}^{*}_{V}(X,\mathrm{Sing}_{V}X)$ d'après \ref{cordefefini}. D'autre part on~a $\mathrm{R}^{n}\mathrm{Pf}\hspace{1pt}(\Sigma^{n}\hspace{1pt}\mathrm{H}^{*}_{V}X)\cong\Sigma^{n}\hspace{1pt}\mathrm{R}^{n}\mathrm{Pf}\hspace{1pt}\mathrm{H}^{*}_{V}X$ d'après \ref{corefini}.
\end{exple}

\bigskip
\footnotesize
\begin{rem}\label{vanishing}
La démonstration ci-dessus dit aussi que l'on a $\mathrm{R}^{p}\mathrm{Pf}\hspace{1pt}\mathrm{H}^{*}_{V}X=0$ pour $p<n$ ce qui est en accord avec \ref{DS3}.
\end{rem}
\normalsize

\bigskip
On reviendra sur le théorème \ref{Cntop} dans la prochaine section et on montrera en particulier sa relation avec \cite{Henntop}.

\vspace{0.75cm}
\textsc{Complexe associé à un $\mathrm{H}^{*}V$-$\mathrm{A}$-module instable}

\bigskip
Soit $M$ un $\mathrm{H}^{*}V$-$\mathrm{A}$-module instable~; on définit un complexe de cochaînes dans la catégorie $V\text{-}\mathcal{U}$
$$
\mathrm{C}^{\bullet}M
\hspace{4pt}=\hspace{4pt}  (\hspace{2pt}\mathrm{C}^{0}M\to\ldots\to\mathrm{C}^{p}M\to\mathrm{C}^{p+1}M\to\ldots\to\mathrm{C}^{n}M\to 0\to 0\to\ldots\hspace{2pt})
$$
de la façon suivante~:

\medskip
On choisit une résolution injective $M\to I^{\bullet}$ de $M$ dans la catégorie $V\text{-}\mathcal{U}$. On pose
$$
\mathrm{C}^{p}M
\hspace{4pt}=\hspace{4pt}
\mathrm{H}^{p}\mathrm{Gr}^{p}I^{\bullet}
$$
et on prend pour cobord $\mathrm{C}^{p}M\to\mathrm{C}^{p+1}M$ l'homomorphisme connectant $\mathrm{H}^{p}\mathrm{Gr}^{p}I^{\bullet}\to\mathrm{H}^{p+1}\mathrm{Gr}^{p+1}I^{\bullet}$ associé à la suite exacte de complexes de cochaînes
$$
\hspace{24pt}
0\to\mathrm{Gr}^{p+1}I^{\bullet}\to
\mathrm{F}^{p}I^{\bullet}/\mathrm{F}^{p+2}I^{\bullet}
\to\mathrm{Gr}^{p}I^{\bullet}\to 0
\hspace{24pt}.
$$
Les mantras habituels de la théorie des résolutions injectives montrent que $\mathrm{C}^{\bullet}M$  est indépendant du choix de la résolution injective de $M$  et que la correspondance $M\mapsto\mathrm{C}^{\bullet}M$ est fonctorielle en $M$.

\medskip
\footnotesize
\begin{rem}
Si l'on a, comme Bourbaki, des scrupules de théorie des ensembles, on peut choisir pour résolution injective de $M$ une résolution ``standard". Précisons un peu. Soient $\mathcal{E}^{*}$ la catégorie des $\mathbb{F}_{2}$-espaces vectoriels $\mathbb{N}$-gradués et $\mathcal{O}:V\text{-}\mathcal{U}\to\mathcal{E}^{*}$ le foncteur oubli (notation locale~!). On observe que $\mathcal{O}$ admet un adjoint à droite  $\widetilde{\mathcal{O}}$, que $\widetilde{\mathcal{O}}\mathcal{O}M$ est un $V\text{-}\mathcal{U}$-injectif, et que l'unité de l'adjonction $M\to\widetilde{\mathcal{O}}\mathcal{O}M$ est injective~; la suite est routine.
\end{rem}
\normalsize

\medskip
Le complexe $\mathrm{C}^{\bullet}M$ est naturellement coaugmenté par $M$. On peut s'en con\-vaincre  en procèdant comme dans la section 4. On considère la filtration décroissante ``triviale'' de $M$ définie par
$$
\mathrm{F}^{p}_{\mathrm{tr}}M
\hspace{4pt}=\hspace{4pt}
\begin{cases}
M & \text{pour $p=0$}, \\
0 & \text{pour $1\leq p\leq n+1$}.
\end{cases}
$$
Le complexe obtenu en remplaçant dans la définition de $\mathrm{C}^{\bullet}M$ les $\mathrm{F}^{p}$ par les $\mathrm{F}^{p}_{\mathrm{tr}}$ est $\mathrm{C}^{\bullet}_{\mathrm{tr}}M:=(M\to 0\to 0\to\ldots)$. L'homomorphisme de complexes $\mathrm{C}^{\bullet}_{\mathrm{tr}}M\to\mathrm{C}^{\bullet}M$ fournit la coaugmentation évoquée plus haut~; le complexe coaugmenté $M\to\mathrm{C}^{\bullet}M$ est noté $\widetilde{\mathrm{C}}^{\bullet}M$.

\begin{pro}\label{complexealg} Soit $M$ un $\mathrm{H}^{*}V_{\mathrm{tf}}$-$\mathrm{A}$-module instable.

\medskip
{\em (a)} Le complexe $\mathrm{C}^{\bullet}M$ est un complexe dans la catégorie $V_{\mathrm{tf}}\text{-}\mathcal{U}$.

\medskip
{\em (b)} Pour $0\leq p\leq n$ on a un isomorphisme canonique de $\mathrm{H}^{*}V$-$\mathrm{A}$-modules instables
$$
\mathrm{C}^{p}M
\hspace{4pt}\cong\hspace{4pt}
\bigoplus_{\mathop{\mathrm{codim}}W=p}
\mathrm{H}^{*}V\otimes_{\mathrm{H}^{*}V/W}\mathrm{R}^{p}\mathrm{Pf}_{V/W}\hspace{1pt}(\mathrm{Fix}_{(V,W)}M)
$$
(la notation $\mathrm{R}^{p}\mathrm{Pf}_{V/W}$ ci-dessus désigne le $p$-ième dérivé de l'endofoncteur $\mathrm{Pf}_{V/W}$ de $V/W_{\mathrm{tf}}\text{-}\mathcal{U}$).
\end{pro}

\medskip
\textit{Démonstration.} Soit $M\to I^{\bullet}$ une résolution injective de $M$ dans la catégorie $V_{\mathrm{tf}}\text{-}\mathcal{U}$. C'est aussi une résolution injective dans la catégorie $V\text{-}\mathcal{U}$ (voir section~2)~; cette observation démontre le point (a). Passons au point (b). Le scholie \ref{schoefini} dit que l'on~a
$$
\hspace{24pt}
\mathrm{Gr}^{p}I^{\bullet}
\hspace{4pt}\cong\hspace{4pt}
\bigoplus_{\mathop{\mathrm{codim}}W=p}
\mathrm{H}^{*}V\otimes_{\mathrm{H}^{*}V/W}\mathrm{Pf}_{V/W}\hspace{1pt}(\mathrm{Fix}_{(V,W)}I^{\bullet})
\hspace{24pt}.
$$
D'après le point (b) de la proposition \ref{Fixinjectifs},  $\mathrm{Fix}_{(V,W)}M\to\mathrm{Fix}_{(V,W)}I^{\bullet}$ est une résolution injective dans la catégorie $V/W_{\mathrm{tf}}\text{-}\mathcal{U}$~; on a donc
$$
\hspace{24pt}
\mathrm{H}^{p}\hspace{1pt}\mathrm{Pf}_{V/W}\hspace{1pt}(\mathrm{Fix}_{(V,W)}I^{\bullet})
\hspace{4pt}\cong\hspace{4pt}
\mathrm{R}^{p}\mathrm{Pf}_{V/W}\hspace{1pt}(\mathrm{Fix}_{(V,W)}M)
\hspace{24pt}.
$$
Comme $\mathrm{H}^{*}V$ est plat sur $\mathrm{H}^{*}V/W$ on obtient bien au bout du compte un isomorphisme canonique
$$
\mathrm{H}^{p}\mathrm{Gr}^{p}I^{\bullet}
\hspace{4pt}\cong\hspace{4pt}
\bigoplus_{\mathop{\mathrm{codim}}W=p}
\mathrm{H}^{*}V\otimes_{\mathrm{H}^{*}V/W}
\mathrm{R}^{p}\mathrm{Pf}_{V/W}\hspace{1pt}(\mathrm{Fix}_{(V,W)}M)
$$
c'est-à-dire le point (b) (qui implique le point (a)~!).
\hfill$\square$

\pagebreak

\vspace{0.75cm}
\textsc{Pendant algébrique du théorème \ref{Ctop}}

\bigskip
Il s'agit de l'énoncé suivant~:

\begin{theo}\label{pendantalg} Soit $M$ un $\mathrm{H}^{*}V_{\mathrm{tf}}$-$\mathrm{A}$-module instable. Les deux propriétés suivantes sont équivalentes~:
\begin{itemize}
\item[(i)] Le $\mathrm{H}^{*}V$-module sous-jacent à $M$ est libre.
\item[(ii)] Le complexe coaugmenté $\widetilde{\mathrm{C}}^{\bullet}M$ est acyclique.
\end{itemize}
\end{theo}

\medskip
\textit{Démonstration de l'implication $(ii)\Rightarrow(i)$.} Elle est identique à celle de l'implication $(ii)\Rightarrow(i)$ du théorème \ref{Ctop}.
\hfill$\square$

\medskip
\textit{Démonstration de l'implication $(i)\Rightarrow(ii)$.} L'ingrédient essentiel de cette démonstration est le corollaire \ref{DS3}.

\medskip
Soient $M$ un $\mathrm{H}^{*}V_{\mathrm{tf}}$-$\mathrm{A}$-module instable et $M\to I^{\bullet}$ une résolution injective de~$M$ dans la catégorie $V_{\mathrm{tf}}\text{-}\mathcal{U}$. On considère la suite spectrale ``cohomologique'' $(\mathrm{E}_{r})_{r\geq 1}$  définie par la filtration du complexe $I^{\bullet}$ par les $\mathrm{F}^{p}I^{\bullet}$. Cette suite spectrale, à valeurs dans la catégorie $V_{\mathrm{tf}}\text{-}\mathcal{U}$, est indépendante du choix de $I^{\bullet}$ et fonctorielle en $M$.

\medskip
--\hspace{8pt}On a $\mathrm{E}_{1}^{p,q}=\mathrm{H}^{p+q}\mathrm{Gr}^{p}I^{\bullet}$ et la différentielle $\mathrm{d}_{1}:\mathrm{E}_{1}^{p,q}\to\mathrm{E}_{1}^{p+1,q}$ est l'homomorphisme connectant associé à la suite exacte de complexes $0\to\mathrm{Gr}^{p+1}I^{\bullet}\to
\mathrm{F}^{p}I^{\bullet}/\mathrm{F}^{p+2}I^{\bullet}
\to\mathrm{Gr}^{p}I^{\bullet}\to 0$ si bien que $\mathrm{C}^{\bullet}M$ n'est rien d'autre que la ligne~$q=0$ de $\mathrm{E}_{1}$ vue comme un complexe de cochaînes.

\medskip
--\hspace{8pt}Comme $\mathrm{Gr}^{p}I^{\bullet}$ est nul pour $p<0$ ou $p>n$,  $\mathrm{E}_{1}^{p,q}\not=0$ implique $p\geq 0$ et $p\leq n$. La suite spectrale dégénère donc \textit{a priori} au terme $\mathrm{E}_{n+1}$. On observera que $\mathrm{E}_{1}^{p,q}\not=0$ implique aussi $p+q\geq 0$.

\medskip
--\hspace{8pt}L'isomorphisme $\mathrm{Gr}^{p}I^{\bullet}
\cong\bigoplus_{\mathop{\mathrm{codim}}W=p}
\mathrm{H}^{*}V\otimes_{\mathrm{H}^{*}V/W}\mathrm{Pf}_{V/W}\hspace{1pt}(\mathrm{Fix}_{(V,W)}I^{\bullet})$ déjà utilisé dans la démonstration du point (b) de \ref{complexealg} donne pareillement l'isomorphisme
$$
\hspace{24pt}
\mathrm{E}_{1}^{p,q}
\hspace{4pt}\cong\hspace{4pt}
\bigoplus_{\mathop{\mathrm{codim}}W=p}
\mathrm{H}^{*}V\otimes_{\mathrm{H}^{*}V/W}
\mathrm{R}^{p+q}\mathrm{Pf}_{V/W}\hspace{1pt}(\mathrm{Fix}_{(V,W)}M)
\hspace{24pt};
$$
compte tenu de la proposition \ref{annderivePf1}, cet isomorphisme montre finalement que $\mathrm{E}_{1}^{p,q}\not=0$ implique $p\geq 0$, $p\leq n$, $p+q\geq 0$ et $q\leq 0$.

\medskip
-- Puisque $M\to I^{\bullet}$ est une résolution, l'aboutissement de la suite spectrale est connu~:
$$
\mathrm{E}_{\infty}^{p,q}
\hspace{4pt}=\hspace{4pt}
\begin{cases}
\mathrm{Gr}^{p}M & \text{pour $p+q=0$}, \\
0 & \text{pour $p+q>0$}. \\
\end{cases}
$$
(Pour se convaincre de l'égalité $\mathrm{E}_{\infty}^{p,-p}=\mathrm{Gr}^{p}M$, on vérifie que l'image de l'homomorphisme $\mathrm{H}^{0}\mathrm{F}^{p}I^{\bullet}\to\mathrm{H}^{0}I^{\bullet}$ s'identifie à $\mathrm{F}^{p}M$.)

\pagebreak
\bigskip
On suppose maintenant que  la condition (i) est vérifiée.  La proposition \ref{libre} et le corollaire \ref{DS3} montrent  que l'on a $\mathrm{E}_{1}^{p,q}=0$ pour $q\not=0$.  Il en résulte que la suite spectrale dégénère au terme $\mathrm{E}_{2}$. On en déduit 
$\mathrm{H}^{p}\mathrm{C}^{\bullet}M=\mathrm{E}_{2}^{p,0}=\mathrm{E}_{\infty}^{p,0}=0$ pour $p>0$. Pour $p=0$ on obtient $\mathrm{H}^{0}\mathrm{C}^{\bullet}M=\mathrm{E}_{2}^{0,0}=\mathrm{E}_{\infty}^{0,0}=\mathrm{Gr}^{0}M$. Comme le $\mathrm{H}^{*}V$-module sous-jacent à $M$ est libre et \textit{a fortiori} sans torsion on a  $\mathrm{F}^{1}M=0$ et donc $\mathrm{Gr}^{0}M\cong M$ (on peut contourner l'argument  $\mathrm{F}^{1}M=0$ en observant que la considération de la suite spectrale donne $\mathrm{Gr}^{p}M=0$ pour $p>0$). Pour se convaincre de ce que l'isomorphisme $\mathrm{H}^{0}\mathrm{C}^{\bullet}M\cong M$ obtenu ci-dessus est bien induit par la coaugmentation, considérer la suite spectrale, disons $({}_{\mathrm{tr}}\mathrm{E}_{r})_{r\geq 1}$, définie par la filtration de $I^{\bullet}$ par les $\mathrm{F}^{p}_{\mathrm{tr}}I^{\bullet}$ (suite spectrale particulièrement triviale~:  ${}_{\mathrm{tr}}\mathrm{E}_{1}^{0,0}\cong M$ et ${}_{\mathrm{tr}}\mathrm{E}_{1}^{p,q}=0$ pour $(p,q)\not=(0,0))$ et invoquer la commutativité du diagramme
$$
\begin{CD}
{}_{\mathrm{tr}}\mathrm{E}_{1}^{0,0}@>>>{}_{\mathrm{tr}}\mathrm{E}_{\infty}^{0,0} \\
@VVV @VVV \\
\mathrm{E}_{1}^{0,0}@>>>\mathrm{E}_{\infty}^{0,0} 
\end{CD}
$$
dans lequel les flèches verticales sont induites par les inclusions $\mathrm{F}^{p}_{\mathrm{tr}}\subset\mathrm{F}^{p}$.
\hfill$\square$

\bigskip
Maintenant que la démonstration de l'implication (i)$\Rightarrow$(ii) du théorème \ref{pendantalg} est achevée, nous dégageons, en vue d'une future utilisation, un énoncé qui résulte facilement de ce que nous avons appris, lors de cette démonstration, sur la suite spectrale $(\mathrm{E}_{r})_{r\geq 1}$ obtenue en filtrant une résolution injective $I^{\bullet}$ d'un $\mathrm{H}^{*}V_{\mathrm{tf}}$-$\mathrm{A}$-module instable $M$ \underline{arbitraire} par les $\mathrm{F}^{p}I^{\bullet}$~:

\begin{scho}\label{pre-n-2} Pour tout $\mathrm{H}^{*}V_{\mathrm{tf}}$-$\mathrm{A}$-module instable $M$, on a $\mathrm{H}^{p}\hspace{1pt}\mathrm{C}^{\bullet}M=0$ pour~$p\geq\sup\hspace{1pt}(1,\dim V-1)$.
\end{scho}

\medskip
\textit{Démonstration.} On pose comme d'habitude $n=\dim V$. On a vu notamment que l'on a~:

\smallskip
(1)\hspace{8pt}$\mathrm{E}_{1}^{p,q}=0$ pour $q>0$ et $p>n$~;

\smallskip
(2)\hspace{8pt}$\mathrm{E}_{\infty}^{p,q}=0$ pour $p+q>0$.

\smallskip
Le point (1) implique $\mathrm{E}_{2}^{p,0}=\mathrm{E}_{\infty}^{p,0}$ pour $p\geq n-1$~; le point (2) implique $\mathrm{E}_{\infty}^{p,0}=0$ pour~$p\geq 1$. Or on a par définition même $\mathrm{H}^{p}\hspace{1pt}\mathrm{C}^{\bullet}M=\mathrm{E}_{2}^{p,0}$.
\hfill$\square$

\pagebreak
\vspace{0.75cm}
\textsc{Généralisation du théorème \ref{pendantalg}}

\medskip
Le théorème \ref{genalg} ci-après (reproduisant \ref{introgenalg}) est le pendant algébrique d'une version du théorème 10.2 de \cite{AFP2} que nous avons évoquée dans l'intro\-duction (Théorème \ref{introgentop})~; nous verrons en section 6 que le théorème \ref{genalg} implique le théorème \ref{gentop} (reproduisant \ref{introgentop}) qui est la version en question. 

\medskip
L'énoncé \ref{genalg} fait intervenir les deux définitions suivantes (reproduisant \ref{intro-syzygie} et \ref{intro-syzygie-bis})~:

\bigskip
\begin{defi}\label{syzygie}
Soient $M$ un $\mathrm{H}^{*}V$-module $\mathbb{N}$-gradué et $j\geq 0$ un entier. On dit que $M$ est une  {\em $\mathrm{H}^{*}V$-$j$-syzygie} s'il existe une suite exacte de $\mathrm{H}^{*}V$-modules $\mathbb{N}$-gradués
$$
0\to M\to L^{0}\to L^{1}\to\ldots\to L^{j-1}
$$
avec $L^{0},L^{1},\ldots,L^{j-1}$ libres (on convient que tout $M$ est une $\mathrm{H}^{*}V$-$0$-syzygie).
\end{defi}

\bigskip
\begin{defi}\label{syzygie-bis}
Soient $M$ un $\mathrm{H}^{*}V_{\mathrm{tf}}$-$\mathrm{A}$-module instable et $j\geq 0$ un entier. Nous dirons que $M$ est une  {\em $(V_{\mathrm{tf}}\text{-}\mathcal{U})$-$j$-syzygie} s'il existe une suite exacte dans la catégorie $V_{\mathrm{tf}}\text{-}\mathcal{U}$
$$
0\to M\to L^{0}\to L^{1}\to\ldots\to L^{j-1}
$$
avec $L^{0},L^{1},\ldots,L^{j-1}$ libres comme $\mathrm{H}^{*}V$-modules (nous convenons que tout~$M$ est une $(V_{\mathrm{tf}}\text{-}\mathcal{U})$-$0$-syzygie).
\end{defi}

\bigskip
\begin{theo}\label{genalg} Soient $M$ un $\mathrm{H}^{*}V_{\mathrm{tf}}$-$\mathrm{A}$-module instable et $j\geq 0$ un entier. Les trois propriétés suivantes sont équivalentes~:
\begin{itemize}
\item[(i)] $M$ est une $(V_{\mathrm{tf}}\text{-}\mathcal{U})$-$j$-syzygie~;
\item[(ii)] le $\mathrm{H}^{*}V$-module sous-jacent à $M$ est une $\mathrm{H}^{*}V$-$j$-syzygie~;
\item[(iii)] on a $\mathrm{H}^{p}\hspace{1pt}\widetilde{\mathrm{C}}^{\bullet}M=0$ pour $p\leq j-2$.
\end{itemize}
\end{theo}

\medskip
L'implication (i)$\Rightarrow$(ii) est évidente. Nous vérifions (ii)$\Rightarrow$(iii) et (iii)$\Rightarrow$(i).

\bigskip
\textit{Démonstration de l'implication (ii)$\Rightarrow$(iii) de \ref{genalg}}

\medskip
On commence par deux énoncés bon marché qui concernent les $\mathrm{H}^{*}V$-modules $\mathbb{N}$-gradués.

\medskip
\begin{lem}\label{decalage} Soit $M$ un $\mathrm{H}^{*}V$-module $\mathbb{N}$-gradué. Si $M$ est une $\mathrm{H}^{*}V$-$j$-syzygie alors on a $\mathrm{Tor}_{k}^{\mathrm{H}^{*}V}(\mathbb{F}_{2},M)=0$ pour $k>\sup\hspace{1pt}(0,\dim V-j)$.
\end{lem}

\medskip
\textit{Démonstration.} On note $N$ le conoyau de la dernière flèche sur la droite dans la définition \ref{syzygie}~; on a donc une suite exacte de  $\mathrm{H}^{*}V$-modules $\mathbb{N}$-gradués de la forme
$$
0\to M\to L^{0}\to L^{1}\to\ldots\to L^{j-1}\to N\to 0
$$
avec $L_{0},L_{1},\ldots,L_{j-1}$ libres (on suppose $j\geq 1$). On constate que l'on a $\mathrm{Tor}_{k}^{\mathrm{H}^{*}V}(\mathbb{F}_{2},M)=\mathrm{Tor}_{k+j}^{\mathrm{H}^{*}V}(\mathbb{F}_{2},N)$ pour $k>0$~; or on a $\mathrm{Tor}_{k+j}^{\mathrm{H}^{*}V}(\mathbb{F}_{2},N)=0$ pour $k+j>\dim V$.
\hfill$\square$

\medskip
\begin{scho}\label{dpsyzygie} Soit $M$ un $\mathrm{H}^{*}V$-module $\mathbb{N}$-gradué. Si $M$ est une $\mathrm{H}^{*}V$-$j$-syzygie alors on a $\mathrm{dp}_{V}M\leq\sup\hspace{1pt}(0,\dim V-j)$.
\end{scho}

\bigskip
On revient maintenant aux $\mathrm{H}^{*}V$-$\mathrm{A}$-modules instables~:

\medskip
\begin{pro}\label{majdp} Soient $M$ un $\mathrm{H}^{*}V$-$\mathrm{A}$-module instable et $W$ un sous-groupe de $V$. Si le $\mathrm{H}^{*}V$-module sous-jacent à $M$  est une $\mathrm{H}^{*}V$-$j$-syzygie alors la dimension projective du $\mathrm{H}^{*}V/W$-module sous-jacent au $\mathrm{H}^{*}V/W$-$\mathrm{A}$-module instable $\mathrm{Fix}_{(V,W)}M$ vérifie  l'inégalité suivante~:
$$
\hspace{24pt}
\mathrm{dp}_{V/W}\hspace{2pt}\mathrm{Fix}_{(V,W)}M
\hspace{4pt}\leq\hspace{4pt}
\sup\hspace{1pt}(0,\mathop{\mathrm{codim}}W-j)
\hspace{24pt}.
$$
\end{pro}

\textit{Démonstration.} Il s'agit d'une généralisation de celle de la proposition \ref{libre}. On considère à nouveau l'isomorphisme de $\mathrm{H}^{*}V/W$-$\mathrm{A}$-modules instables
$$
\mathrm{Fix}_{(V,W)}\hspace{1pt}\mathrm{Tor}_{k}^{\mathrm{H}^{*}V}(\mathrm{H}^{*}W,M)
\hspace{4pt}\cong\hspace{4pt}
\mathrm{Tor}_{k}^{\mathrm{H}^{*}V/W}(\mathbb{F}_{2},\mathrm{Fix}_{(V,W)}M)
$$
du corollaire \ref{sptorFix}. On observe que l'on a un isomorphisme $\mathrm{Tor}_{k}^{\mathrm{H}^{*}V}(\mathrm{H}^{*}W,M)\cong\mathrm{Tor}_{k}^{\mathrm{H}^{*}V/W}(\mathbb{F}_{2},M)$, disons de $\mathbb{F}_{2}$-espaces vectoriels gradués (invoquer par exem\-ple l'isomorphisme canonique  $\mathrm{H}^{*}W\cong\mathbb{F}_{2}\otimes_{\mathrm{H}^{*}V/W}\mathrm{H}^{*}V$). Le $\mathrm{H}^{*}V$-module $M$ étant une $\mathrm{H}^{*}V$-$j$-syzygie est \textit{a fortiori} une $\mathrm{H}^{*}V/W$-$j$-syzygie si bien que le lemme \ref{decalage} implique $\mathrm{Tor}_{k}^{\mathrm{H}^{*}V}(\mathrm{H}^{*}W,M)=0$ pour $k>\sup\hspace{1pt}(0,\mathop{\mathrm{codim}}W-j)$. Il en résulte $\mathrm{Tor}_{k}^{\mathrm{H}^{*}V/W}(\mathbb{F}_{2},\mathrm{Fix}_{(V,W)}M)=0$ pour $k>\sup\hspace{1pt}(0,\mathop{\mathrm{codim}}W-j)$.
\hfill$\square$

\medskip
\begin{cor}\label{annderivePf5} Soient $M$ un $\mathrm{H}^{*}V_{\mathrm{tf}}$-$\mathrm{A}$-module instable et $W$ un sous-groupe de $V$~; soit $p$ la codimension de $W$ dans $V$. Si le $\mathrm{H}^{*}V$-module sous-jacent à $M$  est une $\mathrm{H}^{*}V$-$j$-syzygie alors on a $\mathrm{R}^{p+q}\mathrm{Pf}_{V/W}\hspace{1pt}(\mathrm{Fix}_{(V,W)}M)=0$ pour $\sup(0,p-j)+q<0$.
\end{cor}

\medskip
\textit{Démonstration.} Conséquence des propositions \ref{annderivePf4} et \ref{majdp}.
\hfill$\square$

\pagebreak

\bigskip
\textit{Fin de la démonstration de l'implication (ii)$\Rightarrow$(iii) de \ref{genalg}}

\smallskip
Le résultat d'annulation ci-dessus est exactement ce qu'il nous faut pour adapter la démonstration que nous avons donnée plus haut de l'implication (i)$\Rightarrow$(ii) du théorème \ref{pendantalg}.

\smallskip
On suppose  $j\geq 1$ (il n'y a rien à démontrer pour $j=0$~!). Dans ce cas le $\mathrm{H}^{*}V$-module sous-jacent à $M$ est sans torsion si bien que l'on a $\mathrm{Gr}^{0}M=M$ et $\mathrm{Gr}^{p}M=0$ pour $p>0$. En reprenant la suite spectrale considérée dans la démonstration du théorème \ref{pendantalg} on obtient $\mathrm{H}^{p}\hspace{1pt}\mathrm{C}^{\bullet}M=0$ pour $0<p\leq j-2$ et $\mathrm{H}^{0}\hspace{1pt}\mathrm{C}^{\bullet}M\cong M$ (cet isomorphisme étant induit par la coaugmentation).
\hfill$\square$

\bigskip
\textit{Démonstration de l'implication (iii)$\Rightarrow$(i) de \ref{genalg}}

\medskip
Il est clair que l'on peut supposer $j\leq\dim V$. Dans le cas $j=0$, il n'y a rien à démontrer. Le cas $j=1$ est trivial~: $\mathrm{H}^{-1}\hspace{1pt}\widetilde{\mathrm{C}}^{\bullet}M=0$ signifie que la coaugmentation $M\to\mathrm{H}^{*}V\otimes\mathrm{Fix}_{V}M\hspace{4pt}(=\mathrm{EFix}_{(V,V)}M)$ est injective.

\medskip
\begin{pro}\label{n-2} Soit $V$ un $2$-groupe abélien élémentaire avec $\dim V\geq 2$. Soit $M$ un $\mathrm{H}^{*}V_{\mathrm{tf}}$-$\mathrm{A}$-module instable. Si $\mathrm{H}^{p}\hspace{1pt}\widetilde{\mathrm{C}}^{\bullet}\hspace{1pt}M$ est nul pour $p\leq\dim V-2$, alors~:

\medskip
{\em (a)} le complexe $\widetilde{\mathrm{C}}^{\bullet}\hspace{1pt}M=0$ est acyclique~;

\medskip
{\em (b)} le $\mathrm{H}^{*}V$-module sous-jacent à $M$ est libre.
\end{pro}

\medskip
\textit{Démonstration.} Le point (a) résulte trivialement du scholie \ref{pre-n-2}. Le point (b) résulte ensuite de l'implication (ii)$\Rightarrow$(i) de \ref{pendantalg}.
\hfill$\square$

\bigskip
La proposition ci-dessous est la version algébrique de la proposition \ref{FixCtop}~:

\begin{pro}\label{FixCalg} Soient $M$ un $\mathrm{H}^{*}V_{\mathrm{tf}}$-$\mathrm{A}$-module instable et $W$ un sous-groupe de $V$. On a  des isomorphismes canoniques de complexes de $\mathrm{H}^{*}V/W_{\mathrm{tf}}$-$\mathrm{A}$-modules instables
$$
\hspace{8pt}
\mathrm{Fix}_{(V,W)}\hspace{1pt}\mathrm{C}^{\bullet}_{V}\hspace{1pt}M
\hspace{1pt}\cong\hspace{1pt}
\mathrm{C}^{\bullet}_{V/W}\hspace{1pt}\mathrm{Fix}_{(V,W)}M
\hspace{8pt},\hspace{8pt}
\mathrm{Fix}_{(V,W)\hspace{1pt}}\widetilde{\mathrm{C}}^{\bullet}_{V}\hspace{1pt}M
\hspace{1pt}\cong\hspace{1pt}
\widetilde{\mathrm{C}}^{\bullet}_{V/W}\hspace{1pt}\mathrm{Fix}_{(V,W)}M
\hspace{8pt}.
$$
\end{pro}

\textit{Démonstration.} Nous avons précisé  ci-dessus la notation $\mathrm{C}^{\bullet}$, $\widetilde{\mathrm{C}}^{\bullet}$ en $\mathrm{C}^{\bullet}_{V}$, $\widetilde{\mathrm{C}}^{\bullet}_{V}$ (resp. $\mathrm{C}^{\bullet}_{V/W}$, $\widetilde{\mathrm{C}}^{\bullet}_{V/W}$) pour souligner que dans le membre de gauche (resp. de droite) le foncteur $\mathrm{C}^{\bullet}$, $\widetilde{\mathrm{C}}^{\bullet}$ est défini sur la catégorie $V_{\mathrm{tf}}\text{-}\mathcal{U}$ (resp. $V/W_{\mathrm{tf}}\text{-}\mathcal{U}$). Pareillement, nous précisons ci-après la notation $\mathrm{F}^{p}M$, introduite au début de cette section,  en $\mathrm{F}^{p}_{V}M$.

\medskip
La proposition \ref{FixCalg} est essentiellement conséquence du lemme suivant~:

\begin{lem}\label{FixFp} Soient $M$ un $\mathrm{H}^{*}V_{\mathrm{tf}}$-$\mathrm{A}$-module instable, $W$ un sous-groupe de $V$ et $p$ un entier naturel. Les deux $\mathrm{H}^{*}V/W_{\mathrm{tf}}$-$\mathrm{A}$-modules instables $\mathrm{F}^{p}_{V/W}\mathrm{Fix}_{(V,W)}M$ et  $\mathrm{Fix}_{(V,W)}\mathrm{F}^{p}_{V}M$, vus comme des sous-objets de $\mathrm{Fix}_{(V,W)}M$ {\em (le second peut être vu ainsi car le foncteur $\mathrm{Fix}_{(V,W)}$ est exact)}, coïncident.
\end{lem}

\medskip
\textit{Démonstration.} Soit $I$ un $(V_{\mathrm{tf}}\text{-}\mathcal{U})$-injectif contenant $M$~; comme $\mathrm{Fix}_{(V,W)}$ est exact,  on peut également identifier $\mathrm{Fix}_{(V,W)}M$ à un sous-objet de $\mathrm{Fix}_{(V,W)}I$. Les égalités

\smallskip
--\hspace{8pt}$\mathrm{F}^{p}_{V}M=M\cap\mathrm{F}^{p}_{V}I$ (Proposition \ref{incfilt}),

\smallskip
--\hspace{8pt}$\mathrm{Fix}_{(V,W)}(M\cap\mathrm{F}^{p}_{V}I)=(\mathrm{Fix}_{(V,W)}M)\cap(\mathrm{Fix}_{(V,W)}\mathrm{F}^{p}_{V}I)$ (exactitude de $\mathrm{Fix}_{(V,W)}$),

\smallskip
--\hspace{8pt}$\mathrm{F}^{p}_{V/W}\mathrm{Fix}_{(V,W)}M=\mathrm{Fix}_{(V,W)}M\cap\mathrm{F}^{p}_{V/W}\mathrm{Fix}_{(V,W)}I$ (Proposition \ref{incfilt}),

\smallskip
montrent qu'il suffit de vérifier le lemme pour $M$ injectif (dans la catégorie $V_{\mathrm{tf}}\text{-}\mathcal{U}$).  Compte tenu de la classification des $(V_{\mathrm{tf}}\text{-}\mathcal{U})$-injectifs (voir \ref{injectifs3}) on est ramené à vérifier le lemme pour $M=\mathrm{H}^{*}V\otimes_{\mathrm{H}^{*}V/U}N$ avec $U\subset V$ et $N$ un $\mathrm{H}^{*}V/U$-$\mathrm{A}$-module instable fini.

\bigskip
On suppose donc $M$ de cette forme. La proposition \ref{filtparticulier} donne
$$
\mathrm{F}^{p}_{V}M
\hspace{4pt}=\hspace{4pt}
\begin{cases}
M & \text{pour $p\leq\mathop{\mathrm{codim}}U$} \\
0 & \text{pour $p>\mathop{\mathrm{codim}}U$}
\end{cases}
$$
et par conséquent
$$
\mathrm{Fix}_{(V,W)}\mathrm{F}^{p}_{V}M
\hspace{4pt}=\hspace{4pt}
\begin{cases}
\mathrm{Fix}_{(V,W)}M & \text{pour $p\leq\mathop{\mathrm{codim}}U$} \\
0 & \text{pour $p>\mathop{\mathrm{codim}}U$}
\end{cases}
$$
Le corollaire \ref{corproclef} donne quant à lui
$$
\mathrm{Fix}_{(V,W)}M
\hspace{4pt}\cong\hspace{4pt}
\begin{cases}
\mathrm{H}^{*}V/W\otimes_{\mathrm{H}^{*}V/U}N
& \text{pour $W\subset U$,} \\
0
& \hspace{1pt}\text{pour $W\not\subset U$.}
\end{cases}
$$
On obtient
$$
\mathrm{F}^{p}_{V/W}\mathrm{Fix}_{(V,W)}M
\hspace{4pt}=\hspace{4pt}
\begin{cases}
\mathrm{Fix}_{(V,W)}M & \text{pour $p\leq\mathop{\mathrm{codim}}U$} \\
0 & \text{pour $p>\mathop{\mathrm{codim}}U$}
\end{cases}
$$
en invoquant à nouveau \ref{filtparticulier} ($V$ étant remplacé par $V/W$).
\hfill$\square$

\pagebreak

\medskip
Soit maintenant $M\to I^{\bullet}$ une résolution injective de $M$ dans la catégorie $V_{\mathrm{tf}}\text{-}\mathcal{U}$~; comme le foncteur $\mathrm{Fix}_{(V,W)}$ est exact et transforme les $(V_{\mathrm{tf}}\text{-}\mathcal{U})$-injectifs en $(V/W_{\mathrm{tf}}\text{-}\mathcal{U})$-injectifs (voir le point (b) de \ref{Fixinjectifs}), $\mathrm{Fix}_{(V,W)}M\to\mathrm{Fix}_{(V,W)}I^{\bullet}$ est une résolution injective de $\mathrm{Fix}_{(V,W)}M$ dans la catégorie $V/W_{\mathrm{tf}}\text{-}\mathcal{U}$, si bien que le lemme \ref{FixFp} conduit formellement à la proposition \ref{FixCalg}.
\hfill$\square\square$

\bigskip
Soit $M$ un $\mathrm{H}^{*}V_{\mathrm{tf}}$-$\mathrm{A}$-module instable. On note $\mathrm{E}_{1}M$ le terme  $\mathrm{E}_{1}$ de  la suite spectrale de la démonstration de l'implication (i)$\Rightarrow$(ii) de \ref{pendantalg}~; $\mathrm{E}_{1}M$  est donc un $\mathrm{H}^{*}V_{\mathrm{tf}}$-$\mathrm{A}$-module instable $(\mathbb{Z}\times\mathbb{Z})$-gradué, dépendant fonctoriellement de~$M$, dont la structure est rappelée ci-dessous~:

\medskip
-- Soit $\mathrm{Tg}_{n}$ le sous-ensemble de $\mathbb{Z}\times\mathbb{Z}$ constitué des couples $(p,q)$ vérifiant $p\geq 0$, $p\leq n$, $p+q\geq 0$ et $q\leq 0$ ($\mathrm{Tg}$ est pour triangle~!) ; on a $\mathrm{E}_{1}^{p,q}M=0$ pour $(p,q)\not\in\mathrm{Tg}_{n}$.

\medskip
-- Il existe des homomorpismes $\mathrm{d}_{1}:\mathrm{E}_{1}^{p,q}M\to\mathrm{E}_{1}^{p+1,q}M$ (de $\mathrm{H}^{*}V_{\mathrm{tf}}$-$\mathrm{A}$-modules instables) avec  $\mathrm{d}_{1}\circ\mathrm{d}_{1}=0$~; en d'autres termes $\mathrm{E}_{1}^{\bullet,q}M$ est un complexe de cochaînes dans la catégorie $V_{\mathrm{tf}}\text{-}\mathcal{U}$.

\medskip
--  Le $(V_{\mathrm{tf}}\text{-}\mathcal{U})$-complexe de cochaînes $\mathrm{C}^{\bullet}M$ est $\mathrm{E}_{1}^{\bullet,0}M$.

\medskip
Le point (b) de la proposition \ref{n-2} et la proposition \ref{FixCalg} conduisent à l'énoncé suivant~:

\medskip
\begin{pro}\label{pre-reciproqueAFP} Soient $M$ un $\mathrm{H}^{*}V_{\mathrm{tf}}$-$\mathrm{A}$-module instable et $j\geq 0$ un entier. Si l'on a $\mathrm{H}^{p}\hspace{1pt}\widetilde{\mathrm{C}}^{\bullet}\hspace{1pt}M=0$ pour $p\leq j-2$, alors~:

\medskip
{\em (a)} $\mathrm{Fix}_{(V,W)}M$ est un $\mathrm{H}^{*}V/W$-module libre pour tout sous-groupe $W$ de $V$ avec $\mathop{\mathrm{codim}}W\leq j$~;

\medskip
{\em (b)} on a $\mathrm{E}_{1}^{p,q}M=0$ pour $p\leq j$ et $q\not=0$.
\end{pro}

\medskip
\textit{Démonstration du (a).} Soit $W$ un sous-groupe de $V$ avec $\mathop{\mathrm{codim}}W\leq j$. D'après la proposition \ref{FixCalg}, on a $\mathrm{Fix}_{(V,W)}\hspace{1pt}\widetilde{\mathrm{C}}^{\bullet}\hspace{1pt}M\cong
\widetilde{\mathrm{C}}^{\bullet}_{V/W}\hspace{1pt}\mathrm{Fix}_{(V,W)}M$~; puisque le foncteur $\mathrm{Fix}_{(V,W)}$ est exact, la condition  $\mathrm{H}^{p}\hspace{1pt}\widetilde{\mathrm{C}}^{\bullet}\hspace{1pt}M=0$ pour $p\leq j-2$ implique $\mathrm{H}^{p}\hspace{1pt}\widetilde{\mathrm{C}}_{V/W}^{\bullet}\hspace{1pt}M=0$ pour $p\leq j-2$. Le point (b) de la proposition \ref{n-2} dit alors que $\mathrm{Fix}_{(V,W)}M$ est un $\mathrm{H}^{*}V/W$-module libre.
\hfill$\square$

\medskip
\textit{Démonstration du (b).} On a  vu que l'on a
$$
\hspace{24pt}
\mathrm{E}_{1}^{p,q}
\hspace{4pt}\cong\hspace{4pt}
\bigoplus_{\mathop{\mathrm{codim}}W=p}
\mathrm{H}^{*}V\otimes_{\mathrm{H}^{*}V/W}
\mathrm{R}^{p+q}\mathrm{Pf}_{V/W}\hspace{1pt}(\mathrm{Fix}_{(V,W)}M)
\hspace{24pt};
$$
la proposition \ref{DS3} dit que l'on a $\mathrm{R}^{p+q}\mathrm{Pf}_{V/W}\hspace{1pt}(\mathrm{Fix}_{(V,W)}M)=0$ pour $q\not=0$, compte tenu du point (a).
\hfill$\square$

\pagebreak
\bigskip
La proposition suivante paraît au premier abord assez \textit{ad hoc}. Elle affirme l'existence d'un mystérieux bicomplexe $\mathrm{B}^{\bullet,\bullet}M$ dont certaines propriétés nous permettront d'achever la démonstration en cours. La construction de ce bicomplexe sera effectuée dans la sous-section \ref{Bpq}~; les propriétés évoquées ci-dessus y seront vérifiées. En fait l'introduction de $\mathrm{B}^{\bullet,\bullet}M$ n'est pas si artificielle : nous verrons en \ref{Bpq} qu'elle apporte un éclairage un peu plus concret sur le complexe $\mathrm{C}^{\bullet}M$ dont la définition dans la présente section est plutôt formelle.

\begin{pro}\label{bicomplexe} Soit $M$ un $\mathrm{H}^{*}V_{\mathrm{tf}}$-$\mathrm{A}$-module instable.

\medskip
Il existe un bicomplexe $\mathrm{B}^{\bullet,\bullet}M$ de cochaînes dans la catégorie $V_{\mathrm{tf}}\text{-}\mathcal{U}$, dépendant fonctoriellement de $M$, qui possède les propriétés $\mathcal{B}_{i}$, $i=0,1,2$, décrites ci-après.

\medskip
$(\mathcal{B}_{0})$ On a $\mathrm{B}^{p,q}M=0$ pour $(p,q)\not\in\mathrm{Tg}_{n}$.

\medskip
Soit ${}_\mathrm{B}\mathrm{E}_{1}M$ le terme $\mathrm{E}_{1}$ de la suite spectrale associée à $\mathrm{B}^{\bullet,\bullet}M$ obtenue en dérivant ``verticalement''.

\medskip
$(\mathcal{B}_{1})$ On a pour tout $(p,q)$ un isomorphisme canonique de $\mathrm{H}^{*}V_{\mathrm{tf}}$-$\mathrm{A}$-modules instables $\beta:{}_\mathrm{B}\mathrm{E}_{1}^{p,q}M\cong\mathrm{E}_{1}^{p,q}M$ tel que le diagramme de $\mathrm{H}^{*}V_{\mathrm{tf}}$-$\mathrm{A}$-modules instables
$$
\begin{CD}
{}_{\mathrm{B}}\mathrm{E}_{1}^{p,q}M@>\mathrm{d}_{1}>>{}_{\mathrm{B}}\mathrm{E}_{1}^{p+1,q}M  \\
@V\beta V\cong V @V\beta V\cong V \\
\mathrm{E}_{1}^{p,q}M@>\mathrm{d}_{1}>>\mathrm{E}_{1}^{p+1,q}M
\end{CD}
$$
est commutatif.

\medskip
$(\mathcal{B}_{2})$ Si $j\geq 2$ est un entier tel que $\mathrm{Fix}_{(V,W)}M$ est un $\mathrm{H}^{*}V/W$-module libre pour tout sous-groupe $W$ de $V$ avec $\mathop{\mathrm{codim}}W\leq j$, alors $\mathrm{B}^{p,q}M$ est un $\mathrm{H}^{*}V$-module libre pour $p\leq j$.
\end{pro}

\bigskip
\textit{Fin de la démonstration de l'implication (iii)$\Rightarrow$(i) du théorème \ref{genalg} à l'aide de \ref{pre-reciproqueAFP} et \ref{bicomplexe}}

\medskip
On suppose toujours  $j\geq 2$. On définit un bicomplexe $\mathrm{B}^{\bullet,\bullet}_{j}M$ de cochaînes dans la catégorie $V_{\mathrm{tf}}\text{-}\mathcal{U}$, dépendant fonctoriellement de $M$, de la manière suivante~:

\smallskip
On pose
$$
\mathrm{B}^{p,q}_{j}M
\hspace{4pt}=\hspace{4pt}
\begin{cases}
\mathrm{B}^{p,q}M & \text{pour $p\leq j-1$},
\\
0 & \text{pour $p\geq j$}.
\end{cases}
$$
Les cobords, horizontaux et verticaux, de $\mathrm{B}^{\bullet,\bullet}_{j}M$ sont induits par ceux du bicomplexe $\mathrm{B}^{\bullet,\bullet}M$ ($\mathrm{B}^{\bullet,\bullet}_{j}M$ est un quotient de $\mathrm{B}^{\bullet,\bullet}M$).

\smallskip
On note $({{}_{\mathrm{B}}}_{j}\mathrm{E}_{r}M)_{r\geq 1}$ la suite spectrale associée à $\mathrm{B}^{\bullet,\bullet}_{j}M$ obtenue en dérivant ``verticalement''. Par construction on a ${{}_{\mathrm{B}}}_{j}\mathrm{E}_{1}^{p,q}M={}_{\mathrm{B}}\mathrm{E}_{1}^{p,q}M$ pour $p\leq j-1$ et donc ${{}_{\mathrm{B}}}_{j}\mathrm{E}_{1}^{p,q}M\cong\mathrm{E}_{1}^{p,q}M$ pour $p\leq j-1$ (propriété $(\mathcal{B}_{1})$ de \ref{bicomplexe}).

\smallskip
On suppose maintenant $\mathrm{H}^{p}\hspace{1pt}\widetilde{\mathrm{C}}^{\bullet}M=0$ pour $p\leq j-2$. On note $\mathrm{C}^{\bullet}_{j}\hspace{1pt}M$ le complexe
$$
\mathrm{C}^{0}M\to\mathrm{C}^{1}M\to\ldots\to\mathrm{C}^{j-1}M\to 0\to 0\to\ldots
$$
($\mathrm{C}^{\bullet}_{j}\hspace{1pt}M$ est un quotient de $\mathrm{C}^{\bullet}M$). La définition même de $\mathrm{C}^{\bullet}M$ et le point (b) de \ref{pre-reciproqueAFP} montrent ${{}_{\mathrm{B}}}_{j}\mathrm{E}_{1}^{\bullet,0}M\cong\mathrm{C}^{\bullet}_{j}\hspace{1pt}M$ et ${{}_{\mathrm{B}}}_{j}\mathrm{E}_{1}^{\bullet,q}M=0$ pour $q\not=0$. On en déduit que la suite spectrale $({{}_{\mathrm{B}}}_{j}\mathrm{E}_{r}M)_{r\geq 1}$ dégénère au terme $\mathrm{E}_{2}$ et que l'homologie du complexe $\mathrm{Tot}\hspace{1pt}\mathrm{B}^{\bullet,\bullet}_{j}M$ (le totalisé du bicomplexe $\mathrm{B}^{\bullet,\bullet}_{j}M$) est canoniquement isomorphe à celle du complexe $\mathrm{C}^{\bullet}_{j}\hspace{1pt}M$. Il en résulte que l'on a un isomorphisme canonique $\mathrm{H}^{0}\hspace{1pt}\mathrm{Tot}\hspace{1pt}\mathrm{B}^{\bullet,\bullet}_{j}M\cong M$ et l'égalité $\mathrm{H}^{k}\hspace{1pt}\mathrm{Tot}\hspace{1pt}\mathrm{B}^{\bullet,\bullet}_{j}M=0$ pour $1\leq k\leq j-2$. D'autre part le point (a) de \ref{pre-reciproqueAFP} et la propriété $(\mathcal{B}_{2})$ de~\ref{bicomplexe} montrent que tous les termes du complexe $\mathrm{Tot}\hspace{1pt}\mathrm{B}^{\bullet,\bullet}_{j}M$ sont libres comme $\mathrm{H}^{*}V$-modules.

\smallskip
On constate au bout du compte que l'on dispose de l'énoncé suivant~:

\begin{pro}\label{AFPfonctoriel} Soient $M$ un $\mathrm{H}^{*}V_{\mathrm{tf}}$-$\mathrm{A}$-module instable et $j\geq 2$ un entier. Si l'on a $\mathrm{H}^{p}\hspace{1pt}\widetilde{\mathrm{C}}^{\bullet}M=0$ pour $p\leq j-2$ alors~:

\smallskip
{\em (a)} Les $\mathrm{H}^{*}V_{\mathrm{tf}}$-$\mathrm{A}$-modules instables $\mathrm{Tot}^{k}\hspace{1pt}\mathrm{B}^{\bullet,\bullet}_{j}M$, $0\leq k\leq j-1$, sont libres comme $\mathrm{H}^{*}V$-modules. 

\smallskip
{\em (b) } Il existe un homomorphisme canonique $\eta_{j}: M\to\mathrm{Tot}^{0}\hspace{1pt}\mathrm{B}^{\bullet,\bullet}_{j}M$ tel que
$$
0\longrightarrow M
\overset{\eta_{j}}{\longrightarrow}\mathrm{Tot}^{0}\hspace{1pt}\mathrm{B}^{\bullet,\bullet}_{j}M
\overset{\mathrm{d}}{\longrightarrow}\mathrm{Tot}^{1}\hspace{1pt}\mathrm{B}^{\bullet,\bullet}_{j}M
\overset{\mathrm{d}}{\longrightarrow}\ldots
\overset{\mathrm{d}}{\longrightarrow}\mathrm{Tot}^{j-1}\hspace{1pt}\mathrm{B}^{\bullet,\bullet}_{j}M
$$
est une suite exacte de $\mathrm{H}^{*}V_{\mathrm{tf}}$-$\mathrm{A}$-modules instables.

\smallskip
{\em (On observera que les $\mathrm{Tot}^{k}\hspace{1pt}\mathrm{B}^{\bullet,\bullet}_{j}$ sont des endofoncteurs de $V_{\mathrm{tf}}\text{-}\mathcal{U}$ que $\eta_{j}$ et les cobords  de $\mathrm{Tot}^{\bullet}\hspace{1pt}\mathrm{B}^{\bullet,\bullet}_{j}$ sont des transformations naturelles.)}
\end{pro}

\medskip
L'énoncé ci-dessus fournit une solution ``fonctorielle'' à la question de l'implication (iii)$\Rightarrow$(i) du théorème \ref{genalg}.
\hfill$\square$

\medskip
\begin{rem}\label{j=1} La condition $j\geq 2$ qui apparaît dans l'énoncé \ref{AFPfonctoriel} peut être remplacée par  $j\geq 1$ (qui n'est pas vraiment une restriction~!). On explique  pourquoi ci-après. On verra en \ref{Bpq} que l'on a $\mathrm{B}^{0,0}M=\mathrm{H}^{*}V\otimes\mathrm{Fix}_{V}M$, soit encore $\mathrm{B}^{0,0}M=\mathrm{C}^{0}M$ (en fait cette égalité est forcée par les propriétés $(\mathcal{B}_{0})$ et $(\mathcal{B}_{1})$ de \ref{bicomplexe}). On a donc $\mathrm{B}^{0,0}_{1}M=\mathrm{C}^{0}M$,  $\mathrm{B}^{p,q}_{1}M=0$ pour $(p,q)\not=(0,0)$ et $\mathrm{Tot}\hspace{1pt}\mathrm{B}^{\bullet,\bullet}_{1}M=\mathrm{C}^{\bullet}_{1}\hspace{1pt}M$. Pour avoir le point (b) on prend pour $\eta_{1}$ l'unité d'adjonction $M\to\mathrm{H}^{*}V\otimes\mathrm{Fix}_{V}M$.
\end{rem}

\pagebreak

\bigskip
\begin{exple}\label{babycase}\textit{La proposition \ref{AFPfonctoriel} pour $j=2$ (``baby case'')}

\medskip
On extrait de \ref{Bpq} les informations suivantes~:

\smallskip
--\hspace{8pt}$\mathrm{B}^{0,0}M=\mathrm{EFix}_{(V,V)}M\hspace{4pt}(=\mathrm{H}^{*}V\otimes\mathrm{Fix}_{(V,V)}M=\mathrm{H}^{*}V\otimes\mathrm{Fix}_{V}M)$~;

\smallskip
--\hspace{8pt}$\mathrm{B}^{1,0}M=\bigoplus_{\mathop{\mathrm{codim}}W=1}\mathrm{EFix}_{(V,V)}M$~;

\smallskip
--\hspace{8pt}$\mathrm{B}^{1,-1}M=\bigoplus_{\mathop{\mathrm{codim}}W=1}\mathrm{EFix}_{(V,W)}M$~;

\smallskip
--\hspace{8pt}$\mathrm{d}^{0,0}:\mathrm{B}^{0,0}M\to\mathrm{B}^{1,0}M$ l'homomorphisme diagonal évident~;

\smallskip
--\hspace{8pt}$\mathrm{d}^{1,-1}:\mathrm{B}^{0,0}M\to\mathrm{B}^{1,0}M$ la somme directe indexée par $W$ des homomorphismes $\rho(W,V)_{M}:\mathrm{EFix}_{(V,W)}M\to\mathrm{EFix}_{(V,V)}M$ (introduits en  \ref{Wfonctoriel}).

\medskip
Le diagramme
$$
\begin{CD}
\mathrm{B}^{0,0}(M)@>\mathrm{d}^{0,0}>>\mathrm{B}^{1,0}M \\
& & @A\mathrm{d}^{1,-1}AA \\
& & \mathrm{B}^{1,-1}M
\end{CD}
$$
est le bicomplexe $\mathrm{B}^{\bullet,\bullet}_{2}M$ dont il est question dans \ref{AFPfonctoriel}.

\medskip
On oublie maintenant la machinerie des quatre dernières sous-sections de la section 8, on prend ce qui précède comme définition de $\mathrm{B}^{\bullet,\bullet}_{2}M$ et on esquisse (en petits caractères) une démonstration \textit{ab initio}  de \ref{AFPfonctoriel} pour $j=2$.

\footnotesize
\medskip
On commence par faire les observations suivantes~:

\medskip
1) On a $\mathrm{B}^{0,0}M=\mathrm{C}^{0}M$.

\medskip
2) On a $\mathop{\mathrm{coker}}\mathrm{d}^{1,-1}\cong\mathrm{C}^{1}M$. Précisons un peu. Le cas $\dim V=1$ est facile~; il peut être vu comme une conséquence de la première partie de \ref{derivePf}. Le cas général en résulte grâce à la remarque \ref{Wfonctoriel-4}.

\medskip
3) Les $\mathrm{H}^{*}V_{\mathrm{tf}}$-$\mathrm{A}$-instables $\mathrm{B}^{0,0}M$ et $\mathrm{B}^{1,0}M$ sont manifestement des $\mathrm{H}^{*}V$-modules libres. Il~en est de même pour  $\mathrm{B}^{1,-1}M$ si l'on a $\mathrm{H}^{-1}\widetilde{\mathrm{C}}^{\bullet}M=0$ . En effet cette condition équivaut au fait que le $\mathrm{H}^{*}V$-module $M$ est sans torsion (théorie de Smith algébrique).  Soit $W$ un sous-groupe de $V$, compte tenu de la proposition \ref{FixCalg} et de l'exactitude des foncteurs $\mathrm{Fix}$,  on~a $\mathrm{H}^{-1}\widetilde{\mathrm{C}}_{V/W}^{\bullet}\mathrm{Fix}_{(V,W)}M=0$ et donc $\mathrm{Fix}_{(V,W)}M$ est aussi un $\mathrm{H}^{*}V/W$-module sans torsion. Si $W$ est de codimension $1$, alors $\mathrm{Fix}_{(V,W)}M$ est un $\mathrm{H}^{*}V/W$-module libre puisque $\mathrm{H}^{*}V/W$ est principal et que $\mathrm{Fix}_{(V,W)}M$ est un $\mathrm{H}^{*}V/W$-module de type fini (Proposition \ref{typefiniFix}). L'isomorphisme $\mathrm{EFix}_{(V,W)}M\cong\mathrm{H}^{*}V\otimes_{\mathrm{H}^{*}V/W}\mathrm{Fix}_{(V,W)}M$ montre enfin que $\mathrm{EFix}_{(V,W)}M$ est un $\mathrm{H}^{*}V$-module libre.

\medskip
4)  Soit $\pi:\mathrm{B}^{1,0}M\to\mathrm{C}^{1}M$ l'épimorphisme induit par l'isomorphisme de la deuxième observation~; on constate que $\pi\circ\mathrm{d}^{0,0}$ s'identifie à $\mathrm{d}^{0}:\mathrm{C}^{0}M\to\mathrm{C}^{1}M$.

\medskip
5) On constate que les deux complexes de cochaînes
$$
\begin{CD}
\mathrm{C}^{0}M @>\mathrm{d}^{0}>> \mathrm{C}^{1}M
\end{CD}
$$
et
$$
\begin{CD}
\mathrm{B}^{0,0}M\oplus\mathrm{B}^{1,0}M
@>{\begin{bmatrix}
\mathrm{d}^{0,0} & \mathrm{d}^{1,-1}
\end{bmatrix}}
>> \mathrm{B}^{1,0}M
\end{CD}
$$
(le totalisé du bicomplexe introduit plus haut) ont même cohomologie.

\pagebreak

\bigskip
On en déduit que l'on dispose d'une suite exacte de $\mathrm{H}^{*}V_{\mathrm{tf}}$-$\mathrm{A}$-instables de la forme
$$
\begin{CD}
0@>>> M @>>>
\mathrm{B}^{0,0}M\oplus\mathrm{B}^{1,0}M
@>{\begin{bmatrix}
\mathrm{d}^{0,0} & \mathrm{d}^{1,-1}
\end{bmatrix}}
>> \mathrm{B}^{1,0}M
\end{CD}
$$
et donc que $M$ est une $(V_{\mathrm{tf}}\text{-}\mathcal{U})$-$2$-syzygie (fonctoriellement en $M$).
\normalsize
\end{exple}

\pagebreak

\sect{Comparaison entre les complexes algébrique et topologique}

Le premier objet de cette section est de démontrer la proposition \ref{introCalgCtop} qui compare les complexes $\widetilde{\mathrm{C}}^{\bullet}_{\mathrm{alg}}X$ et  $\widetilde{\mathrm{C}}^{\bullet}_{\mathrm{top}}X$ asociés à un $V$-CW-complexe fini $X$. Le second est d'utiliser cette comparaison pour obtenir la version topologique du théorème \ref{genalg} à savoir le théorème  \ref{gentop} (reproduisant \ref{introgentop}).

\medskip
Notre démonstration va comporter une courte incursion dans la théorie des catégories dérivées, aussi nous commençons par rappeler quelques énoncés classiques d'algèbre homologique qui apparaissent dans cette théorie.

\bigskip
Soit $\mathcal{A}$ une catégorie abélienne~; soit $\mathrm{Ch}^{\geq 0}(\mathcal{A})$ la catégorie (abélienne) des complexes de cochaînes dont les termes sont des objets de $\mathcal{A}$. Soient $Z^{\bullet}$ et $K^{\bullet}$ deux objets de $\mathrm{Ch}^{\geq 0}(\mathcal{A})$, on note $[Z^{\bullet},K^{\bullet}]$ le groupe abélien des classes d'homotopie de morphismes de  $Z^{\bullet}$ dans $K^{\bullet}$.

\medskip
L'énoncé clef (voir par exemple \cite[I, Theorem 6.2]{Ive}) est le suivant~:

\begin{theo}\label{catder1}
Soit $K^{\bullet}$ un complexe de cochaînes dont les termes sont des objets injectifs de $\mathcal{A}$, alors le foncteur contravariant $Z^{\bullet}\mapsto [Z^{\bullet},K^{\bullet}]$, défini sur la catégorie $\mathrm{Ch}^{\geq 0}(\mathcal{A})$ et à valeurs dans la catégorie des groupes abéliens, envoie un quasi-isomorphisme sur un isomorphisme.
\end{theo}

(On rappelle qu'un  quasi-isomorphisme est un morphisme de complexes qui induit un isomorphisme en homologie.)

\medskip
On suppose  maintenant que $\mathcal{A}$ a assez d'injectifs. Soit $C^{\bullet}$ un complexe de cochaînes dont les termes sont des objets de $\mathcal{A}$~; un ``remplacement injectif" de $C^{\bullet}$ est la donnée d'un  complexe de cochaînes $I^{\bullet}$ dont les termes sont des objets injectifs de $\mathcal{A}$ et d'un quasi-isomorphisme $\phi:C^{\bullet}\to I^{\bullet}$. L'hypothèse faite sur $\mathcal{A}$ garantit l'existence de remplacements injectifs pour tout $C^{\bullet}$.

\begin{cor}\label{catder2}
On considère un diagramme dans la catégorie  $\mathrm{Ch}^{\geq 0}(\mathcal{A})$
$$
\begin{CD}
C^{\bullet}@>\phi>>I^{\bullet} \\
@VfVV & & \\
D^{\bullet}@>\psi>>J^{\bullet} \\
\end{CD}
$$
avec $\phi$ et $\psi$ des remplacements injectifs.

\medskip
Il existe un homomorphisme de complexes $F:I^{\bullet}\to J^{\bullet}$, unique à homotopie près, tel que le diagramme suivant
$$
\begin{CD}
C^{\bullet}@>\phi>>I^{\bullet} \\
@VfVV @VFVV \\
D^{\bullet}@>\psi>>J^{\bullet} 
\end{CD}
$$
est commutatif à homotopie près.
\end{cor}

\bigskip
\textit{Démonstration.} On prend $K^{\bullet}=J^{\bullet}$ dans \ref{catder1}. La classe de $F$ est l'image inverse de $\psi\circ f$ par l'isomorphisme $[I^{\bullet},J^{\bullet}]\to[C^{\bullet},J^{\bullet}]$ induit par $\phi$.
\hfill$\square$

\bigskip
En prenant ci-dessus pour $f$ l'identité de $C^{\bullet}$ on obtient un ``énoncé d'unicité" pour les remplacements injectifs~:

\begin{cor}\label{catder3}
Soient $\phi:C^{\bullet}\to I^{\bullet}$ et $\phi':C^{\bullet}\to I'^{\bullet}$ deux remplacements injectifs dans la catégorie  $\mathrm{Ch}^{\geq 0}(\mathcal{A})$. Alors il existe une équivalence d'homotopie $H:I^{\bullet}\to I'^{\bullet}$, unique à homotopie près, telle que $\phi'$ est homotope à $H\circ\phi$. 
\end{cor}

\begin{scho}\label{catder4}
On reprend les hypothèses de \ref{catder2} et on suppose en outre que $f$ est un quasi-isomorphisme. Alors $F$ (qui apparaît dans la conclusion de \ref{catder2}) est une équivalence d'homotopie (toujours unique à homotopie près).
\end{scho}

\bigskip
Ces rappels étant faits, on prend pour $\mathcal{A}$ la catégorie abélienne $V_{\mathrm{tf}}\text{-}\mathcal{U}$.

\medskip
On définit un foncteur $\mathrm{G}:\mathrm{Ch}^{\geq 0}(V_{\mathrm{tf}}\text{-}\mathcal{U})\to\mathrm{Ch}^{\geq 0}(V_{\mathrm{tf}}\text{-}\mathcal{U})$ par le procédé déjà utilisé dans la section 5. On pose $(\mathrm{G}\hspace{1pt}C^{\bullet})^{\hspace{1pt}p}:=\mathrm{H}^{p}\mathrm{Gr}^{p}\hspace{1pt}C^{\bullet}$ et on prend pour cobord l'homomorphisme connectant $\mathrm{H}^{p}\mathrm{Gr}^{p}\hspace{1pt}C^{\bullet}\to\mathrm{H}^{p+1}\mathrm{Gr}^{p+1}\hspace{1pt}C^{\bullet}$ associé à la suite exacte $0\to\mathrm{Gr}^{p+1}\hspace{1pt}C^{\bullet}\to\mathrm{F}^{p}\hspace{1pt}C^{\bullet}/\mathrm{F}^{p+2}\hspace{1pt}C^{\bullet}
\to\mathrm{Gr}^{p}\hspace{1pt}C^{\bullet}\to 0$. Le foncteur $\mathrm{G}$ donne naissance à un foncteur $\mathbf{G}:\mathrm{Ch}^{\geq 0}(V_{\mathrm{tf}}\text{-}\mathcal{U})\to\mathrm{Ch}^{\geq 0}(V_{\mathrm{tf}}\text{-}\mathcal{U})$ de la façon suivante~:  on choisit un remplacement injectif $\phi:C^{\bullet}\to I^{\bullet}$ et on pose $\mathbf{G}\hspace{1pt}C^{\bullet}:=\mathrm{G}\hspace{1pt}I^{\bullet}$. On observera que l'on dispose d'une transformation naturelle $\mathrm{G}\hspace{1pt}C^{\bullet}\to\mathbf{G}\hspace{1pt}C^{\bullet}$. On observera également que si $f:C^{\bullet}\to D^{\bullet}$ est un quasi-isomorphisme alors $\mathbf{G}\hspace{1pt}f$ est un isomorphisme~; en d'autres termes le foncteur $\mathbf{G}$  se factorise à travers la catégorie dérivée.

\medskip
\footnotesize
\begin{rem}
Les remplacements injectifs de $C^{\bullet}$ tels que les termes de $I^{\bullet}$ sont des sommes directes du type de celles qui apparaissent dans le theorème \ref{injectifs3} forment un ensemble. Cette observation permet de se débarasser d'éventuels scrupules ``ensemblistes'' concernant la définition du foncteur $\mathbf{G}$.
\end{rem}
\normalsize

\bigskip
\begin{exple}\label{CpointMbis} Soit $M$ un objet de $V_{\mathrm{tf}}\text{-}\mathcal{U}$. Soit $\mathrm{c}^{\bullet}M$ le complexe de cochaînes
$$
\hspace{24pt}
M\to 0\to 0\to\ldots
\hspace{24pt};
$$
alors $\mathbf{G}\hspace{1pt}\mathrm{c}^{\bullet}M$ est le complexe $\mathrm{C}^{\bullet} M$ de la section 5. On constate que l'on a $\mathrm{G}\hspace{1pt}\mathrm{c}^{\bullet}M=\mathrm{c}^{\bullet}(M/\mathrm{F}^{1}M)$ et que la coaugmentation de $\mathrm{C}^{\bullet}M$ est fournie par le composé de l'homomorphisme canonique $\mathrm{c}^{\bullet}M\to\mathrm{c}^{\bullet}(M/\mathrm{F}^{1}M)$ et de l'homomorphisme naturel $\mathrm{G}\hspace{1pt}\mathrm{c}^{\bullet}M\to\mathbf{G}\hspace{1pt}\mathrm{c}^{\bullet}M$.
\end{exple}

\bigskip
Soit $C^{\bullet}$ l'un des deux complexes $\mathop{\Sigma^{\dim V}}\hspace{-1pt}\mathrm{C}^{\bullet}_{\mathrm{top}}X$ ($X$ un $V$-CW-com\-plexe fini) ou  $\mathrm{C}^{\bullet}M$ ($M$ un $\mathrm{H}^{*}V_{\mathrm{tf}}$-$\mathrm{A}$-module instable), introduits respectivement dans les sections 4 et 5~;  $C^{\bullet}$ est un complexe de cochaînes dans la catégorie $V_{\mathrm{tf}}\text{-}\mathcal{U}$ qui possède la propriété suivante~: pour tout entier $p\geq 0$, $C^{p}$ est isomorphe à une somme directe de la forme
$$
\bigoplus_{\mathop{\mathrm{codim}}W=p}
\mathrm{H}^{*}V\otimes_{\mathrm{H}^{*}V/W} N_{W}
$$
($W$ sous-groupe de $V$) avec $N_{W}$ un $\mathrm{H}^{*}V/W$-$\mathrm{A}$-module instable fini (ce qui implique $C^{p}=0$ pour $p>\dim V$). Nous dirons qu'un complexe de ce type est {\em e-fini spécial}.

\medskip
\begin{pro}\label{efinispecial1}
Soit $C^{\bullet}$ un complexe de cochaînes dans la catégorie $V_{\mathrm{tf}}\text{-}\mathcal{U}$ que l'on suppose e-fini spécial.

\medskip
{\em (a)} Les complexes $C^{\bullet}$ et $\mathrm{G}\hspace{1pt}C^{\bullet}$ sont canoniquement isomorphes.

\medskip
{\em (b)} L'homomorphisme naturel $\mathrm{G}\hspace{1pt}C^{\bullet}\to\mathbf{G}\hspace{1pt}C^{\bullet}$ est un isomorphisme.

\end{pro}

\textit{Démonstration du (a).} On commence par rappeler la définition des troncations brutales d'un complexe. Soient
$$
C^{\bullet}
\hspace{4pt}=\hspace{4pt}(C^{0}\to C^{1}\to\ldots\to C^{p}\to C^{p+1}\to\ldots)
$$
un complexe de cochaînes dans une catégorie abélienne et $p\geq 0$ un entier~; on note $\sigma^{\geq p}\hspace{1pt}C^{\bullet}$ le sous-complexe
$$
0\to 0\to\ldots\to 0\to C^{p}\to C^{p+1}\to C^{p+2}\to\ldots
$$
de $C^{\bullet}$ ($p$-ième troncation ``brutale'' de $C^{\bullet}$) . Ces sous-complexes définissent une filtration décroissante de $C^{\bullet}$~:
$$
\hspace{24pt}
C^{\bullet}=\sigma^{\geq 0}\hspace{2pt}C^{\bullet}\supset\sigma^{\geq 1}\hspace{2pt}C^{\bullet}\supset\ldots\supset\sigma^{\geq p}\hspace{2pt}C^{\bullet}\supset\ldots
\hspace{24pt}.
$$
Soit maintenant $C^{\bullet}$ un complexe de cochaînes dans la catégorie $V_{\mathrm{tf}}\text{-}\mathcal{U}$~; si $C^{\bullet}$ est e-fini spécial alors, d'après \ref{filtparticulier}, la filtration de $C^{\bullet}$ par les $\mathrm{F}^{p}\hspace{1pt}C^{\bullet}$ coïncide avec celle par les $\sigma^{\geq p}\hspace{1pt}C^{\bullet}$. Le point (a) de la proposition résulte du point (a) de l'énoncé ci-dessous dont la vérification est laissée au lecteur~:

\medskip
\begin{lem}\label{E1troncbrutal} Soit $C^{\bullet}=(C^{0}\to C^{1}\to  C^{2}\to\ldots)$ un complexe de cochaînes dans une catégorie abélienne. Le terme $\mathrm{E}_{1}^{\bullet,\bullet}$ de la suite spectrale associée à la filtration par les troncations brutales vérifie les deux propriétés suivantes~:

\medskip
{\em (a)} $\mathrm{E}_{1}^{\bullet,0}$ (vu comme un complexe grâce à $\mathrm{d}_{1}$) s'identifie à $C^{\bullet}$~;

\medskip
{\em (b)}  $\mathrm{E}_{1}^{p,q}$ est nul pour $q\not=0$.

\end{lem}

(Le point (b) est là en vue d'une future référence.)
\hfill$\square$

\bigskip
\textit{Démonstration du (b) de \ref{efinispecial1}.} On va utiliser une résolution injective, au sens de Cartan-Eilenberg, de $C^{\bullet}$. Rappelons un peu de quoi il s'agit.

\medskip
On note $C^{\bullet,\bullet}$ le bicomplexe de cochaînes du premier quadrant défini par~:
$$
C^{\bullet,q}
\hspace{4pt}=\hspace{4pt}
\begin{cases}
C^{\bullet} & \text{pour}\hspace{6pt}q=0,
\\
0 & \text{pour}\hspace{6pt}q>0.
\end{cases}
$$
Une résolution injective de Cartan-Eilenberg de $C^{\bullet}$ est la donnée d'un bicomplexe de cochaînes du premier quadrant $I^{\bullet,\bullet}$, dont tous les termes sont injectifs, et d'un homomorphisme de bicomplexes $\eta:C^{\bullet,\bullet}\to I^{\bullet,\bullet}$, telle que six propriétés sont vérifiées (voir \cite[Chap. XVII]{CE}). La première, la seule que nous utiliserons ci-après, est que $\eta:C^{p}=C^{p,0}\to I^{p,\bullet}$ est une résolution injective. Cette propriété implique, par un argument trivial de suites spectrales (considérer les ``cohomologies verticales''), que
$$
\begin{CD}
C^{\bullet}=\mathrm{Tot}\hspace{1pt}C^{\bullet,\bullet}
@>\mathrm{Tot}\hspace{1pt}\eta>>
\mathrm{Tot}\hspace{1pt}I^{\bullet,\bullet}
\end{CD}
$$
est un remplacement injectif de $C^{\bullet}$.

\medskip
On fixe maintenant un entier $k\geq 0$ et on considère les deux suites spectrales du premier quadrant, disons ${}^{\mathrm{II}}\mathrm{E}_{r}^{p,q}(0,k)$ et ${}^{\mathrm{II}}\mathrm{E}_{r}^{p,q}(1,k)$ (la décoration~${}^{\mathrm{II}}$ indique que l'on commence par mettre en {\oe}uvre le cobord  vertical), associées respectivement aux bicomplexes $\mathrm{Gr}^{k}\hspace{1pt}C^{\bullet,\bullet}$ et $\mathrm{Gr}^{k}\hspace{1pt}I^{\bullet,\bullet}$, convergeant vers les cohomologies de $\mathrm{Tot}\hspace{1pt}\mathrm{Gr}^{k}\hspace{1pt}C^{\bullet,\bullet}$ et $\mathrm{Tot}\hspace{1pt}\mathrm{Gr}^{k}\hspace{1pt}I^{\bullet,\bullet}$.

\medskip
Le lemme clef est le suivant~:

\begin{lem}\label{lemmevert} 
L'homomorphisme $\eta:C^{\bullet,\bullet}\to I^{\bullet,\bullet}$ induit  pour tout couple $(p,q)$ un isomorphisme ${}^{\mathrm{II}}\mathrm{E}_{1}^{p,q}(0,k)\cong{}^{\mathrm{II}}\mathrm{E}_{1}^{p,q}(1,k)$.
\end{lem}

\medskip
\textit{Démonstration.}

\medskip
0) On constate, compte tenu de  \ref{filtparticulier} et de la définition même de ${}^{\mathrm{II}}\mathrm{E}_{1}^{p,q}(0,k)$, que l'on a
$$
{}^{\mathrm{II}}\mathrm{E}_{1}^{p,q}(0,k)
\hspace{4pt}=\hspace{4pt}
\begin{cases}
C^{k} & \text{pour}\hspace{6pt}(p,q)=(k,0),
\\
0 & \text{pour}\hspace{6pt}(p,q)\not=(k,0).
\end{cases}
$$

\medskip
1) Comme on l'a rappelé plus haut $\eta:C^{p}=C^{p,0}\to I^{p,\bullet}$ est une résolution injective et par définition encore ${}^{\mathrm{II}}\mathrm{E}_{1}^{p,q}(1,k)$ est $\mathrm{H}^{q}\hspace{1pt}\mathrm{Gr}^{k}\hspace{1pt}I^{p,\bullet}$. On va déterminer cette cohomologie en remplaçant la résolution injective $\eta:C^{p}\to I^{p,\bullet}$ par une résolution injective $\eta':C^{p}\to I'^{\bullet}$ adaptée à la forme particulière de $C^{p}$.

\medskip
On a par hypothèse
$$
C^{p}
\hspace{4pt}\cong\hspace{4pt}
\bigoplus_{\mathop{\mathrm{codim}}W=p}
\mathrm{H}^{*}V\otimes_{\mathrm{H}^{*}V/W}N_{W}
$$
avec $N_{W}$ un  un $\mathrm{H}^{*}V/W$-$\mathrm{A}$-module instable fini~; soit $N_{W}\to K_{W}^{\bullet}$ une résolution injective dans la catégorie $V/W_{\mathrm{tf}}\text{-}\mathcal{U}$ dont tous les termes sont finis (une telle résolution existe d'après \ref{resolutionfini}), alors
$$
C^{p}\to
\bigoplus_{\mathop{\mathrm{codim}}W=p}
\mathrm{H}^{*}V\otimes_{\mathrm{H}^{*}V/W} K_{W}^{\bullet}
$$
est une résolution injective dans la catégorie $V_{\mathrm{tf}}\text{-}\mathcal{U}$ (invoquer la platitude du $\mathrm{H}^{*}V/W$-module $\mathrm{H}^{*}V$ et la proposition \ref{exactFixbis})~: c'est la résolution injective $\eta':C^{p}\to I'^{\bullet}$ évoquée plus haut.

\medskip
D'après \ref{filtparticulier}, on a cette fois
$$
\mathrm{Gr}^{k}\hspace{1pt}I'^{\bullet}
\hspace{4pt}=\hspace{4pt}
\begin{cases}
I'^{\bullet} & \text{pour}\hspace{6pt}p=k,
\\
0 & \text{pour}\hspace{6pt}p\not=k.
\end{cases}
$$
On en déduit
$$
\mathrm{H}^{q}\hspace{1pt}\mathrm{Gr}^{k}\hspace{1pt}I'^{\bullet}
\hspace{4pt}\cong\hspace{4pt}
\begin{cases}
C^{k} & \text{pour}\hspace{6pt}(p,q)=(k,0),
\\
0 & \text{pour}\hspace{6pt}(p,q)\not=(k,0),
\end{cases}
$$
le premier isomorphisme étant induit par $\eta'$, et donc
$$
\mathrm{H}^{q}\hspace{1pt}\mathrm{Gr}^{k}\hspace{1pt}I^{p,\bullet}
\hspace{4pt}\cong\hspace{4pt}
\begin{cases}
C^{k} & \text{pour}\hspace{6pt}(p,q)=(k,0),
\\
0 & \text{pour}\hspace{6pt}(p,q)\not=(k,0).
\end{cases}
$$
le premier isomorphisme étant induit par $\eta$.

\medskip
La démonstration du lemme est achevée. Pour compléter celle du point (b) de la proposition on observe que l'on a $\mathrm{Tot}\hspace{1pt}\mathrm{Gr}^{k}\hspace{1pt}C^{\bullet,\bullet}=\mathrm{Gr}^{k}\hspace{1pt}\mathrm{Tot}\hspace{1pt}C^{\bullet,\bullet}$ et $\mathrm{Tot}\hspace{1pt}\mathrm{Gr}^{k}\hspace{1pt}I^{\bullet,\bullet}=\mathrm{Gr}^{k}\hspace{1pt}\mathrm{Tot}\hspace{1pt}I^{\bullet,\bullet}$ et on considère l'homomorphisme de complexes de chaînes $\mathrm{G}\hspace{1pt}\mathrm{Tot}\hspace{1pt}C^{\bullet,\bullet}\to\mathrm{G}\hspace{1pt}\mathrm{Tot}\hspace{1pt}I^{\bullet,\bullet}$ induit par $\eta$. C'est un isomorphisme car le lemme implique que l'homomorphisme $(\mathrm{G}\hspace{1pt}\mathrm{Tot}\hspace{1pt}C^{\bullet,\bullet})^{\hspace{1pt}p}\to(\mathrm{G}\hspace{1pt}\mathrm{Tot}\hspace{1pt}I^{\bullet,\bullet})^{\hspace{1pt}p}$, induit par $\eta$, est un isomorphisme pour tout $p$.
\hfill$\square$

\bigskip
\begin{rem}\label{degenerescence} La démonstration que nous avons donnée du lemme \ref{lemmevert} montre en fait que les deux suites spectrales ${}^{\mathrm{II}}\mathrm{E}_{r}^{p,q}(0,k)$ et ${}^{\mathrm{II}}\mathrm{E}_{r}^{p,q}(1,k)$ (qui sont isomorphes) dégénèrent au terme $\mathrm{E}_{1}$.
\end{rem}

\bigskip
Soit $X$ un $V$-CW-complexe fini~; le complexe $\mathrm{C}^{\bullet}_{\mathrm{top}}X$ n'est pas en général un objet de $\mathrm{Ch}^{\geq 0}(V_{\mathrm{tf}}\text{-}\mathcal{U})$ à cause de l'apparition du foncteur $\Sigma^{-p}$ dans la définition du terme $\mathrm{C}^{p}_{\mathrm{top}}X$ (considérer le cas où l'action de $V$ sur $X$ est libre). Aussi, pour pouvoir appliquer la proposition \ref{efinispecial1}, nous allons commencer par ``instabiliser'' $\mathrm{C}^{\bullet}_{\mathrm{top}}X$ à l'aide de l'endofoncteur $\widetilde{\Sigma}$ dont nous rappelons ci-dessous la définition et les propriétés dont nous aurons besoin.

\bigskip
On note $\widetilde{\Sigma}:\mathcal{U}\to\mathcal{U}$ l'adjoint à droite du foncteur $\Sigma:\mathcal{U}\to\mathcal{U}$ (voir \cite[\S\hspace{2pt}2.3]{LZens}). En fait $\widetilde{\Sigma}$ peut être explicité très concrètement~: soit $M$ un $\mathrm{A}$-module instable, $\Sigma\widetilde{\Sigma}M$ s'identifie (\textit{via} la co-unité de l'adjonction) au sous-$\mathrm{A}$-module de $M$ constitué des éléments $x$ vérifiant $\mathrm{Sq}_{0}\hspace{1pt}x=0$ ($\mathrm{Sq}_{0}\hspace{1pt}x:=\mathrm{Sq}^{\vert x\vert}x$, $\vert x\vert$ désignant le degré de $x$). Ceci montre en particulier, \cite[Proposition 8.1.1]{LZens}, que l'homomorphisme $L\otimes\widetilde{\Sigma}M\to\widetilde{\Sigma}(L\otimes M)$, naturels en les $\mathrm{A}$-modules instables $L$ et $M$, est un isomorphisme si l'on a $\widetilde{\Sigma}L=0$ (on dit alors que $L$ est {\em réduit})~; il en résulte que l'endofoncteur $\widetilde{\Sigma}:\mathcal{U}\to\mathcal{U}$ induit des endofoncteurs, toujours noté $\widetilde{\Sigma}$, des catégories abéliennes $V\text{-}\mathcal{U}$ et $V_{\mathrm{tf}}\text{-}\mathcal{U}$ (dans ce dernier cas invoquer le fait que $\mathrm{H}^{*}V$ est noethérien). On observera que les endofoncteurs composés $\widetilde{\Sigma}\Sigma$ sont les foncteurs identité (l'unité de l'adjonction est un isomorphisme).

\begin{pro-def}\label{sigmatilde1} Soit $p\geq 0$ un entier~; on note $\widetilde{\Sigma}^{\hspace{1pt}p}$ le $p$-ième itéré de $\widetilde{\Sigma}$.

\medskip
{\em (a)} Soit $M$ un $\mathrm{A}$-module instable~;  $\widetilde{\Sigma}^{\hspace{1pt}p}M$ s'identifie au plus grand sous-$\mathrm{A}$-module instable du  $\mathrm{A}$-module $\Sigma^{-p}M$.

\medskip
{\em (b)} Soit $M$ un $\mathrm{H}^{*}V$-$\mathrm{A}$-module instable~;  le plus grand sous-$\mathrm{A}$-module instable du  $\mathrm{H}^{*}V$-$\mathrm{A}$-module $\Sigma^{-p}M$ est stable sous l'action de $\mathrm{H}^{*}V$ et s'identifie au $\mathrm{H}^{*}V$-$\mathrm{A}$-module instable $\widetilde{\Sigma}^{\hspace{1pt}p}M$.

\medskip
{\em (c)} Soient $W\subset V$ un sous-groupe et $N$ un $\mathrm{H}^{*}V/W$-$\mathrm{A}$-module instable~;  le $\mathrm{H}^{*}V$-$\mathrm{A}$-module instable $\widetilde{\Sigma}^{\hspace{1pt}p}(\mathrm{H}^{*}V\otimes_{\mathrm{H}^{*}V/W}N)$ est naturellement isomorphe au $\mathrm{H}^{*}V$-$\mathrm{A}$-module instable $\mathrm{H}^{*}V\otimes_{\mathrm{H}^{*}V/W}\widetilde{\Sigma}^{\hspace{1pt}p}N$.
\end{pro-def}

\bigskip
\textit{Démonstration.} Le cas $p=0$ est trivial. Le cas $p\geq 1$ se démontre par une récurrence immédiate à partir du cas $p=1$. Dans ce dernier cas, on se convainc grâce aux rappels précédents. Donnons un peu plus de détails pour le point (c). Soit $\gamma:\mathrm{H}^{*}V\otimes_{\mathrm{H}^{*}V/W}\widetilde{\Sigma}N\to\widetilde{\Sigma}(\mathrm{H}^{*}V\otimes_{\mathrm{H}^{*}V/W}N)$ l'adjoint de l'homomorphisme
$$
\hspace{24pt}
\Sigma\hspace{1pt}(\mathrm{H}^{*}V\otimes_{\mathrm{H}^{*}V/W}\widetilde{\Sigma}N)
=
\mathrm{H}^{*}V\otimes_{\mathrm{H}^{*}V/W}\Sigma\widetilde{\Sigma}N
\to
\mathrm{H}^{*}V\otimes_{\mathrm{H}^{*}V/W}N
\hspace{24pt}.
$$
On a vu lors de la démonstration du lemme \ref{pre-annPf} qu'une rétraction de groupes $r:V\to W$ fournit en particulier un isomorphisme de $\mathrm{A}$-modules instables $\mathrm{H}^{*}V\otimes_{\mathrm{H}^{*}V/W}N\cong\mathrm{H}^{*}W\otimes N$, naturel en le $\mathrm{H}^{*}V/W$-$\mathrm{A}$-module instable $N$. On en déduit que l'on dispose d'un diagramme commutatif dans la catégorie $\mathcal{U}$
$$
\begin{CD}
\mathrm{H}^{*}V\otimes_{\mathrm{H}^{*}V/W}\widetilde{\Sigma}N
@>\gamma>>
\widetilde{\Sigma}(\mathrm{H}^{*}V\otimes_{\mathrm{H}^{*}V/W}N)
\\
@V\cong VV @V\cong VV
\\
\mathrm{H}^{*}W\otimes\widetilde{\Sigma}N
@>\gamma_{\mathcal{U}}>>
\widetilde{\Sigma}(\mathrm{H}^{*}W\otimes N)
\end{CD}
$$
dans lequel $\gamma_{\mathcal{U}}$ est l'isomorphisme évoqué plus haut (celui de \cite[Proposition 8.1.1]{LZens}). Le $\mathcal{U}$-homomorphisme sous-jacent à $\gamma$ est donc un isomorphisme si bien qu'il en est de même pour $\gamma$.
\hfill$\square$

\bigskip
Nous voici armés pour définir et décrire ``l'instabilisé'' de $\mathrm{C}^{\bullet}_{\mathrm{top}}X$ que nous notons $\mathrm{C}^{\bullet}_{\mathrm{utop}}X$. Nous posons
$$
\mathrm{C}^{\bullet}_{\mathrm{utop}}X
\hspace{4pt}:=\hspace{4pt}
\widetilde{\Sigma}^{\hspace{1pt}n}(\Sigma^{\hspace{1pt}n}\mathrm{C}^{\bullet}_{\mathrm{top}}X)
$$
(ci-dessus $n$ désigne, comme convenu dans ce mémoire, la dimension de $V$). Par construction $\mathrm{C}^{\bullet}_{\mathrm{utop}}X$ est un complexe (de cochaînes) dont les termes sont des $\mathrm{H}^{*}V$-$\mathrm{A}$-modules instables et qui est un sous-complexe de $\mathrm{C}^{\bullet}_{\mathrm{top}}X$, dans la catégorie des complexes dont les termes sont des $\mathrm{H}^{*}V$-$\mathrm{A}$-modules. Précisons. Les isomorphismes
$$
\mathrm{C}^{p}_{\mathrm{utop}}X
\hspace{4pt}\cong\hspace{4pt}
\widetilde{\Sigma}^{\hspace{1pt}p}(\widetilde{\Sigma}^{\hspace{1pt}n-p}\Sigma^{\hspace{1pt}n-p})\Sigma^{\hspace{1pt}p}\mathrm{C}^{p}_{\mathrm{top}}X
\hspace{4pt}\cong\hspace{4pt}
\widetilde{\Sigma}^{\hspace{1pt}p}\Sigma^{\hspace{1pt}p}\mathrm{C}^{p}_{\mathrm{top}}X
$$
et la proposition \ref{sigmatilde1} montrent~:

\smallskip
-- que $\mathrm{C}^{p}_{\mathrm{utop}}X$ est le plus grand sous-$\mathrm{A}$-module instable de $\mathrm{C}^{p}_{\mathrm{top}}X$

\smallskip
-- que $\mathrm{C}^{p}_{\mathrm{utop}}X$ un sous-$\mathrm{H}^{*}V$-$\mathrm{A}$-module de $\mathrm{C}^{p}_{\mathrm{top}}X$

\smallskip
-- et que l'on a un isomorphisme de $\mathrm{H}^{*}V$-$\mathrm{A}$-modules instables
$$
\hspace{24pt}
\mathrm{C}^{p}_{\mathrm{utop}}X
\hspace{4pt}\cong\hspace{4pt}
\bigoplus_{\mathop{\mathrm{codim}}W=p}
\mathrm{H}^{*}V
\otimes_{\mathrm{H}^{*}V/W}
\widetilde{\Sigma}^{\hspace{1pt}p}\hspace{1pt}\mathrm{H}_{V/W}^{*}(X^{W},\mathrm{Sing}_{V/W}X^{W})
\hspace{24pt}.
$$

\begin{pro}\label{efinispecial2}
Le complexe $\mathrm{C}^{\bullet}_{\mathrm{utop}}X$ est e-fini spécial.
\end{pro}

\bigskip
\textit{Démonstration.} Il suffit de ce convaincre de ce que le $\mathrm{H}^{*}V/W$-$\mathrm{A}$-module instable $\widetilde{\Sigma}^{\hspace{1pt}p}\hspace{1pt}\mathrm{H}_{V/W}^{*}(X^{W},\mathrm{Sing}_{V/W}X^{W})$ est fini. Or celui-ci est en particulier un sous-module du module $\Sigma^{-p}\hspace{1pt}\mathrm{H}_{V/W}^{*}(X^{W},\mathrm{Sing}_{V/W}X^{W})$ qui est fini.
\hfill$\square$

\bigskip
On rappelle que l'on a posé $\mathrm{C}^{\bullet}_{\mathrm{alg}}X:=\mathrm{C}^{\bullet}\hspace{1pt}\mathrm{H}^{*}_{V}X$ (resp. $\widetilde{\mathrm{C}}^{\bullet}_{\mathrm{alg}}X:=\widetilde{\mathrm{C}}^{\bullet}\hspace{1pt}\mathrm{H}^{*}_{V}X$), $\mathrm{C}^{\bullet}$~désignant le foncteur de $V_{\mathrm{tf}}\text{-}\mathcal{U}$ dans $\mathrm{Ch}^{\geq 0}(V_{\mathrm{tf}}\text{-}\mathcal{U})$ introduit dans la section~5. Nous sommes maintenant en mesure d'exhiber un homomorphisme de complexes de cochaînes (resp. complexes de cochaînes coaugmentés), de $\mathrm{H}^{*}V$-$\mathrm{A}$-modules, $\varkappa:\mathrm{C}^{\bullet}_{\mathrm{alg}}X\to\mathrm{C}^{\bullet}_{\mathrm{top}}X$ (resp. $\varkappa:\widetilde{\mathrm{C}}^{\bullet}_{\mathrm{alg}}X\to\widetilde{\mathrm{C}}^{\bullet}_{\mathrm{top}}X$) qui sera un isomorphisme si l'on suppose $\mathrm{H}^{*}_{V}X$ libre comme $\mathrm{H}^{*}V$ module.

\bigskip
On note $\mathrm{c}^{\bullet}X$ le complexe $\mathrm{c}^{\bullet}\hspace{1pt}\mathrm{H}^{*}_{V}X=(\mathrm{H}^{*}_{V}X\to 0\to 0\to\ldots)$~; $\mathrm{c}^{\bullet}X$ est un objet de $\mathrm{Ch}^{\geq 0}(V_{\mathrm{tf}}\text{-}\mathcal{U})$. Nous avons vu dans la section 4 que l'on dispose d'un homomorphisme de complexes de cochaînes de $\mathrm{H}^{*}V$-$\mathrm{A}$-modules $\eta:\mathrm{c}^{\bullet}X\to\mathrm{C}^{\bullet}_{\mathrm{top}}X$ (qui fournit la coaugmentation de $\widetilde{\mathrm{C}}^{\bullet}_{\mathrm{top}}X$). On note $\eta_{\mathrm{u}}:\mathrm{c}^{\bullet}X\to\mathrm{C}^{\bullet}_{\mathrm{utop}}X$ l'homomorphisme $\widetilde{\Sigma}^{\hspace{1pt}n}(\Sigma^{\hspace{1pt}n}\eta)$~; par construction $\eta_{\mathrm{u}}$ est un morphisme de $\mathrm{Ch}^{\geq 0}(V_{\mathrm{tf}}\text{-}\mathcal{U})$ et l'on a $\eta=\iota\circ\eta_{\mathrm{u}}$, $\iota$ désignant l'inclusion $\mathrm{C}^{\bullet}_{\mathrm{utop}}X\hookrightarrow\mathrm{C}^{\bullet}_{\mathrm{top}}X$. On note $\widetilde{\mathrm{C}}^{\bullet}_{\mathrm{utop}}X$ le complexe coaugmenté grâce à $\eta_{\mathrm{u}}$~; il est clair que $\iota$ se prolonge en un homomorphisme de complexes de cochaînes coaugmentés, prolongement que nous notons encore $\iota$. On considère le $\mathrm{Ch}^{\geq 0}(V_{\mathrm{tf}}\text{-}\mathcal{U})$-morphisme
$$
\hspace{24pt}
\mathbf{G}\hspace{1pt}\eta_{u}:
\mathbf{G}\hspace{1pt}\mathrm{c}^{\bullet}X
\to
\mathbf{G}\hspace{1pt}\mathrm{C}^{\bullet}_{\mathrm{utop}}X
\hspace{24pt}.
$$
Par définition $\mathbf{G}\hspace{1pt}\mathrm{c}^{\bullet}X$ est $\mathrm{C}^{\bullet}_{\mathrm{alg}}X$~; les propositions \ref{efinispecial1} et \ref{efinispecial2} fournissent un isomorphisme canonique $\nu:\mathrm{C}^{\bullet}_{\mathrm{utop}}X\to\mathbf{G}\hspace{1pt}\mathrm{C}^{\bullet}_{\mathrm{utop}}X$. On prend pour
$$
\varkappa:\mathrm{C}^{\bullet}_{\mathrm{alg}}X
\longrightarrow
\mathrm{C}^{\bullet}_{\mathrm{top}}X
$$
l'homomorphisme $\iota\circ\nu^{-1}\circ\mathbf{G}\hspace{1pt}\eta_{u}$. En contemplant le diagramme commutatif
$$
\begin{CD}
\mathrm{G}\hspace{1pt}\mathrm{c}^{\bullet}X
@>\mathrm{G}\hspace{1pt}\eta_{\mathrm{u}}>>
\mathrm{G}\hspace{1pt}\mathrm{C}^{\bullet}_{\mathrm{utop}}X \\
@VVV @ VVV \\
\mathbf{G}\hspace{1pt}\mathrm{c}^{\bullet}X
@>\mathbf{G}\hspace{1pt}\eta_{\mathrm{u}}>>
\mathbf{G}\hspace{1pt}\mathrm{C}^{\bullet}_{\mathrm{utop}}X
\end{CD}
$$
on se convainc que $\nu^{-1}\circ\mathbf{G}\hspace{1pt}\eta_{\mathrm{u}}$ et $\varkappa$ se prolongent respectivement en des homomorphismes $\widetilde{\mathrm{C}}^{\bullet}_{\mathrm{alg}}X\to\widetilde{\mathrm{C}}^{\bullet}_{\mathrm{utop}}X$ et $\widetilde{\mathrm{C}}^{\bullet}_{\mathrm{alg}}X\to\widetilde{\mathrm{C}}^{\bullet}_{\mathrm{top}}X$. Ce dernier est encore noté~$\varkappa$~; par construction il est l'identité pour $\bullet=-1$. Ceci achève la vérification du point (a) de \ref{introCalgCtop}.

\medskip
La proposition \ref{kappanaturel} ci-après dit que l'homorphisme de complexes de $\mathrm{H}^{*}V$-$\mathrm{A}$-modules $\varkappa$ que nous venons de définir est une transformation naturelle de foncteurs (contravariants en $X$) en un sens facile à deviner. Précisons tout de même. Soient  $X$ et $Y$ deux $V$-CW-complexes finis et $f:X\to Y$ une application $V$-équivariante. On a $f(X^{W})\subset f(Y^{W})$ pour tout sous-groupe $W$ de $V$ et \textit{a fortiori} $f(\mathrm{F}_{p}X)\subset f(\mathrm{F}_{p}Y)$ (la filtration $\mathrm{F}_{p}$ est introduite en section 4) si bien que $f$ induit un homomorphisme de complexes de $\mathrm{H}^{*}V$-$\mathrm{A}$-modules $\widetilde{\mathrm{C}}_{\mathrm{top}}^{\bullet}\hspace{1pt}Y\to\widetilde{\mathrm{C}}_{\mathrm{top}}^{\bullet}\hspace{1pt}X$, disons $\widetilde{\mathrm{C}}_{\mathrm{top}}^{\bullet}(f)$. Notons d'autre part $\widetilde{\mathrm{C}}_{\mathrm{alg}}^{\bullet}(f):\widetilde{\mathrm{C}}_{\mathrm{alg}}^{\bullet}\hspace{1pt}Y\to\widetilde{\mathrm{C}}_{\mathrm{alg}}^{\bullet}\hspace{1pt}X$ l'homomorphisme de complexes de $\mathrm{H}^{*}V$-$\mathrm{A}$-modules instables induit par l'homomorphisme de $\mathrm{H}^{*}V$-$\mathrm{A}$-modules instables $f^{*}:\mathrm{H}^{*}_{V}Y\to\mathrm{H}^{*}_{V}X$ (le foncteur $V_{\mathrm{tf}}\text{-}\mathcal{U}\to\mathrm{Ch}^{\geq-1}(V_{\mathrm{tf}}\text{-}\mathcal{U}),M\mapsto\widetilde{\mathrm{C}}^{\bullet}M$, introduit en section 5, est covariant en $M$). En reprenant chaque étape de la définition de $\varkappa$ on constate que l'énoncé suivant est vérifié~:

\begin{pro}\label{kappanaturel} Soient  $X$ et $Y$ deux $V$-CW-complexes finis et $f:X\to Y$ une application $V$-équivariante. Le diagramme de complexes de $\mathrm{H}^{*}V$-$\mathrm{A}$-modules
$$
\begin{CD}
\widetilde{\mathrm{C}}_{\mathrm{alg}}^{\bullet}\hspace{1pt}Y
@>\varkappa_{Y}>>
\widetilde{\mathrm{C}}_{\mathrm{top}}^{\bullet}\hspace{1pt}Y \\
@V\widetilde{\mathrm{C}}_{\mathrm{alg}}^{\bullet}(f)VV @V\widetilde{\mathrm{C}}_{\mathrm{top}}^{\bullet}(f)VV \\
\widetilde{\mathrm{C}}_{\mathrm{alg}}^{\bullet}\hspace{1pt}X
@>\varkappa_{X}>>
\widetilde{\mathrm{C}}_{\mathrm{top}}^{\bullet}\hspace{1pt}X
\end{CD}
$$
est commutatif.

\smallskip
\footnotesize
{\em On a précisé ci-dessus la notation $\varkappa:\widetilde{\mathrm{C}}_{\mathrm{alg}}^{\bullet}\hspace{1pt}X\to
\widetilde{\mathrm{C}}_{\mathrm{top}}^{\bullet}\hspace{1pt}X$ en $\varkappa_{X}$ pour souligner la naturalité en $X$ de $\varkappa$ \ldots qui est le propos de l'énoncé.
\normalsize}
\end{pro}

\bigskip
On en vient maintenant au point (b) de \ref{introCalgCtop}~:

\begin{pro}\label{CalgCtop} Soient $V$ un $2$-groupe abélien élémentaire et $X$ un $V$-CW-complexe fini. Si $\mathrm{H}^{*}_{V}X$ est libre comme $\mathrm{H}^{*}V$-module alors l'homomorphisme de complexes de cochaînes de $\mathrm{H}^{*}V$-$\mathrm{A}$-modules (coaugmentés par $\mathrm{H}^{*}_{V}X$)
$$
\varkappa:
\hspace{4pt}\widetilde{\mathrm{C}}_{\mathrm{alg}}^{\bullet}\hspace{1pt}X
\longrightarrow
\widetilde{\mathrm{C}}_{\mathrm{top}}^{\bullet}\hspace{1pt}X
$$
est un isomorphisme.
\end{pro}

\bigskip
\textit{Démonstration.} D'apès \ref{instablebis} tous les termes du complexe $\mathrm{C}_{\mathrm{top}}^{\bullet}\hspace{1pt}X$ sont des $\mathrm{A}$-modules instables, en d'autres termes on a $\mathrm{C}_{\mathrm{utop}}^{\bullet}\hspace{1pt}X=\mathrm{C}_{\mathrm{top}}^{\bullet}\hspace{1pt}X$ et $\eta_{\mathrm{u}}=\eta$~; l'homomorphisme $\varkappa:\mathrm{C}^{\bullet}_{\mathrm{alg}}X\to\mathrm{C}^{\bullet}_{\mathrm{top}}X$ est le composé $\nu^{-1}\circ\mathbf{G}\hspace{1pt}\eta$. Le théorème \ref{Ctop} dit que $\eta:\mathrm{c}^{\bullet}X\to\mathrm{C}^{\bullet}_{\mathrm{top}}X$ est un quasi-isomorphisme, ce qui entraîne que $\mathbf{G}\hspace{1pt}\eta$ est un isomorphisme.
\hfill$\square$

\bigskip
La proposition ci-dessus dit en particulier que si $\mathrm{H}^{*}_{V}X$ est libre comme $\mathrm{H}^{*}V$-module alors l'homomorphisme de $\mathrm{H}^{*}V$-$\mathrm{A}$-modules
$$
\mathrm{R}^{n}\mathrm{Pf}\hspace{1pt}\mathrm{H}^{*}_{V}X
\cong\mathrm{C}_{\mathrm{alg}}^{n}\hspace{1pt}X
\overset{\varkappa_{X}^{n}}{\longrightarrow}
\mathrm{C}_{\mathrm{top}}^{n}\hspace{1pt}X
:=
\Sigma^{-n}\hspace{1pt}\mathrm{H}^{*}_{V}(X,\mathrm{Sing}_{V}X)
$$
est un isomorphisme~; on retrouve le résultat de \ref{Cntop}.

\bigskip
La proposition ci-dessous montre que l'homomorphisme $\varkappa:\mathrm{C}_{\mathrm{alg}}^{\bullet}\hspace{1pt}X\to\mathrm{C}_{\mathrm{top}}^{\bullet}\hspace{1pt}X$ peut s'exprimer en fonction des $\varkappa^{p}:\mathrm{C}_{\mathrm{alg},V/W}^{p}\hspace{1pt}X^{W}\to\mathrm{C}_{\mathrm{top,V/W}}^{p}\hspace{1pt}X^{W}$ (la présence en indice de $V/W$ est là pour signaler que l'on considère $X^{W}$ comme un $V/W$-espace) pour $p=\dim V/W$.

\medskip
\begin{pro}\label{decompkappa} Soient $V$ un $2$-groupe abélien élémentaire, $X$ un $V$-CW-complexe fini et $p$ un entier avec $0\leq p\leq n$.

\smallskip
Soit
$$
\mathrm{d\acute{e}c}_{\mathrm{alg}}^{p}:
\mathrm{C}_{\mathrm{alg}}^{p}\hspace{1pt}X
\longrightarrow
\bigoplus_{\mathop{\mathrm{codim}}W=p}
\mathrm{H}^{*}V\otimes_{\mathrm{H}^{*}V/W}\mathrm{C}_{\mathrm{alg,V/W}}^{p}\hspace{1pt}X^{W}
$$
l'isomorphisme de $\mathrm{H}^{*}V$-$\mathrm{A}$-modules instables induit par l'isomorphisme du point (b) de \ref{complexealg} et celui de \ref{cohFix}.

\smallskip
Soit
$$
\mathrm{d\acute{e}c}_{\mathrm{top}}^{p}:
\mathrm{C}_{\mathrm{top}}^{p}\hspace{1pt}X
\longrightarrow
\bigoplus_{\mathop{\mathrm{codim}}W=p}
\mathrm{H}^{*}V\otimes_{\mathrm{H}^{*}V/W}\mathrm{C}_{\mathrm{top,V/W}}^{p}\hspace{1pt}X^{W}
$$
l'isomorphisme de $\mathrm{H}^{*}V$-$\mathrm{A}$-modules donné par le point (b) de \ref{Ctopefini}.

\medskip
Alors le diagramme de de $\mathrm{H}^{*}V$-$\mathrm{A}$-modules
$$
\begin{CD}
\mathrm{C}_{\mathrm{alg}}^{p}\hspace{1pt}X
@>\mathrm{d\acute{e}c}_{\mathrm{alg}}^{p}>>
\bigoplus_{\mathop{\mathrm{codim}}W=p}\hspace{4pt}
\mathrm{H}^{*}V\otimes_{\mathrm{H}^{*}V/W}\mathrm{C}_{\mathrm{alg,V/W}}^{p}\hspace{1pt}X^{W} \\
@VV\varkappa_{X}^{p}V
@VV\bigoplus_{\mathop{\mathrm{codim}}W=p}\hspace{1pt}
\mathrm{H}^{*}{V}\otimes_{\mathrm{H}^{*}V/W}\hspace{3pt}\varkappa_{X^{W},V/W}^{p} V \\
\mathrm{C}_{\mathrm{top}}^{p}\hspace{1pt}X
@>\mathrm{d\acute{e}c}_{\mathrm{top}}^{p}>>
\bigoplus_{\mathop{\mathrm{codim}}W=p}\hspace{4pt}
\mathrm{H}^{*}V\otimes_{\mathrm{H}^{*}V/W}\mathrm{C}_{\mathrm{top,V/W}}^{p}\hspace{1pt}X^{W}
\end{CD}
$$
est commutatif.
\end{pro}

\medskip
\textit{Démonstration.} On propose une démonstration en kit. Les pièces de ce kit sont les énoncés \ref{kit-1}, \ref{kit-2} et \ref{kit-3} ci-après~; les mots-clef de la notice de montage sont ``fonctorialité'' et ``naturalité''.

\medskip
\begin{pro-def}\label{kit-1} Soit $W$ un sous-groupe de $V$. On a  un isomorphisme canonique de complexes de $\mathrm{H}^{*}V$-$\mathrm{A}$-modules
$$
\iota_{\mathrm{top},W}\hspace{2pt}:\hspace{2pt}
\mathrm{C}_{\mathrm{top},V}^{\bullet}\hspace{1pt}X^{W}
\hspace{4pt}\cong\hspace{4pt}
\mathrm{H}^{*}V\otimes_{\mathrm{H}^{*}V/W}\mathrm{C}_{\mathrm{top},V/W}^{\bullet}\hspace{1pt}X^{W}
$$
(à gauche $X^{W}$ est vu comme un $V$-espace, à droite comme un $V/W$-espace).
\end{pro-def}

\pagebreak

\textit{Démonstration.} On constate que l'on a $\mathrm{F}_{p,V}X^{W}=\mathrm{F}_{p,V/W}X^{W}$.
\hfill$\square$

\medskip
\begin{pro-def}\label{kit-2} Soient $W$ un sous-groupe de $V$ et $N$ un $\mathrm{H}^{*}V/W$-$\mathrm{A}$-module instable. On a  un isomorphisme canonique de complexes de $\mathrm{H}^{*}V$-$\mathrm{A}$-modules instables
$$
\hspace{24pt}
\mathrm{C}^{\bullet}_{V}\hspace{1pt}
(\mathrm{H}^{*}V\otimes_{\mathrm{H}^{*}V/W}N)
\hspace{4pt}\cong\hspace{4pt}
\mathrm{H}^{*}V\otimes_{\mathrm{H}^{*}V/W}\mathrm{C}^{\bullet}_{V/W}\hspace{1pt}N
\hspace{24pt}.
$$
Dans le cas $N=\mathrm{H}_{V/W}^{*}X^{W}$ (et donc $\mathrm{H}^{*}V\otimes_{\mathrm{H}^{*}V/W}N=\mathrm{H}_{V}^{*}X^{W}$) cet isomorphisme est noté
$$
\hspace{24pt}
\iota_{\mathrm{alg},W}\hspace{2pt}:\hspace{2pt}
\mathrm{C}_{\mathrm{alg},V}^{\bullet}\hspace{1pt}X^{W}
\hspace{4pt}\cong\hspace{4pt}
\mathrm{H}^{*}V\otimes_{\mathrm{H}^{*}V/W}\mathrm{C}_{\mathrm{alg},V/W}^{\bullet}\hspace{1pt}X^{W}
\hspace{24pt}.
$$
\end{pro-def}

\textit{Démonstration.} On observe que si $N\to I^{\bullet}$ est une résolution injective dans la catégorie $V/W\text{-}\mathcal{U}$ alors $\mathrm{H}^{*}V\otimes_{\mathrm{H}^{*}V/W}N\to\mathrm{H}^{*}V\otimes_{\mathrm{H}^{*}V/W}I^{\bullet}$ est une résolution injective dans la catégorie $V\text{-}\mathcal{U}$ (variante de l'argument utilisé dans la démonstration de \ref{annderivePf3}).
\hfill$\square$

\medskip
\begin{pro}\label{kit-3} Soit $W$ un sous-groupe de $V$. Le diagramme
$$
\begin{CD}
\mathrm{C}_{\mathrm{alg},V}^{\bullet}\hspace{1pt}X^{W}
@>\iota_{\mathrm{alg},W}>\cong>
\mathrm{H}^{*}V\otimes_{\mathrm{H}^{*}V/W}\mathrm{C}_{\mathrm{alg},V/W}^{\bullet}\hspace{1pt}X^{W} \\
@VV\varkappa_{V,X^{ W}}V
@VV{\mathrm{H}^{*}V\otimes_{\mathrm{H}^{*}V/W}\hspace{3pt}\varkappa_{V/W,X^{ W}}}V \\
\mathrm{C}_{\mathrm{top},V}^{\bullet}\hspace{1pt}X^{W}
@>\iota_{\mathrm{top},W}>\cong>
\mathrm{H}^{*}V\otimes_{\mathrm{H}^{*}V/W}\mathrm{C}_{\mathrm{top},V/W}^{\bullet}\hspace{1pt}X^{W}
\end{CD}
$$
est commutatif.
\end{pro}

\medskip
\textit{Démonstration.} On reprend la définition de la transformation naturelle $\varkappa$ et on utilise les points suivants~:

\smallskip
-- Soit $C^{\bullet}$ un complexe de cochaînes dans la catégorie $V/W_{\mathrm{tf}}\text{-}\mathcal{U}$~; on a un isomorphisme $\mathbf{G}_{V}(\mathrm{H}^{*}V\otimes_{\mathrm{H}^{*}V/W}C^{\bullet})\cong\mathrm{H}^{*}V\otimes_{\mathrm{H}^{*}V/W}\mathbf{G}_{V/W}\hspace{1pt}C^{\bullet}$ (l'endofoncteur~$\mathbf{G}$ de la catégorie $\mathrm{Ch}^{\geq 0}(V_{\mathrm{tf}}\text{-}\mathcal{U})$ est introduit au début de cette section, la signification de $V$ et $V/W$ en indice est transparente).

\smallskip
-- Soit $C^{\bullet}$ un complexe de cochaînes dans la catégorie $V/W_{\mathrm{tf}}\text{-}\mathcal{U}$~; si $C^{\bullet}$ est e-fini spécial (cette terminologie est introduite juste avant l'énoncé \ref{efinispecial1}) alors il en est de même pour $\mathrm{H}^{*}V\otimes_{\mathrm{H}^{*}V/W}C^{\bullet}$ (qui est un complexe de cochaînes dans la catégorie $V_{\mathrm{tf}}\text{-}\mathcal{U}$).

\smallskip
-- L'isomorphisme $\iota_{\mathrm{top},W}$ de \ref{kit-1} induit un isomorphisme de complexes de cochaînes $\mathrm{C}_{\mathrm{utop},V}^{\bullet}\hspace{1pt}X^{W}
\cong\mathrm{H}^{*}V\otimes_{\mathrm{H}^{*}V/W}\mathrm{C}_{\mathrm{utop},V/W}^{\bullet}\hspace{1pt}X^{W}$ dans la catégorie $V_{\mathrm{tf}}\text{-}\mathcal{U}$ (invoquer le point (c) de \ref{sigmatilde1}).

\medskip
\textit{Fin de la démonstration de la proposition \ref{decompkappa}}

\medskip
On constate (``par fonctorialité'') que les homomorphismes $\mathrm{d\acute{e}c}_{\mathrm{alg}}^{p}$, $\mathrm{d\acute{e}c}_{\mathrm{top}}^{p}$ coïncident respectivement avec les produits, indexés par les $W$ de codimension $p$, des composés
$$
\begin{CD}
\hspace{24pt}
\mathrm{C}_{\mathrm{alg},V}^{p}\hspace{1pt}X
@>>>
\mathrm{C}_{\mathrm{alg},V}^{p}\hspace{1pt}X^{W}
@>\iota_{\mathrm{alg},W}^{p}>\cong>
\mathrm{H}^{*}V\otimes_{\mathrm{H}^{*}V/W}\mathrm{C}_{\mathrm{alg},V/W}^{\bullet}\hspace{1pt}X^{W}
\hspace{24pt},
 \\
\hspace{24pt}
\mathrm{C}_{\mathrm{top},V}^{p}\hspace{1pt}X
@>>>
\mathrm{C}_{\mathrm{top},V}^{p}\hspace{1pt}X^{W}
@>\iota_{\mathrm{top},W}^{p}>\cong>
\mathrm{H}^{*}V\otimes_{\mathrm{H}^{*}V/W}\mathrm{C}_{\mathrm{top},V/W}^{\bullet}\hspace{1pt}X^{W}
\hspace{24pt},
\end{CD}
$$
et on observe que la naturalité de $\varkappa$ (voir \ref{kappanaturel}) dit que le diagramme
$$
\begin{CD}
\mathrm{C}_{\mathrm{alg},V}^{p}\hspace{1pt}X
@>>>
\mathrm{C}_{\mathrm{alg},V}^{p}\hspace{1pt}X^{W} \\
@VV\varkappa_{X}^{p}V @VV\varkappa_{X^{W}}^{p}V \\
\mathrm{C}_{\mathrm{top},V}^{p}\hspace{1pt}X
@>>>
\mathrm{C}_{\mathrm{top},V}^{p}\hspace{1pt}X^{W}
\end{CD}
$$
est commutatif
\hfill$\square\square$

\bigskip
La proposition \ref{decompkappa} conduit à l'énoncé ci-dessous qui est une généralisation (prendre $j=\dim V$) de la proposition \ref{CalgCtop}~:

\medskip
\begin{pro}\label{genCalgCtop} Soient $V$ un $2$-groupe abélien élémentaire, $X$ un $V$-CW-complexe fini et $j\geq 0$ un entier. Si le $\mathrm{H}^{*}V$-module sous-jacent à $\mathrm{H}_{V}^{*}X$ est une $\mathrm{H}^{*}V$-$j$-syzygie, alors l'homomorphisme de $\mathrm{H}^{*}V$-$\mathrm{A}$-modules 
$$
\varkappa_{X}^{p}\hspace{2pt}:
\hspace{2pt}\widetilde{\mathrm{C}}_{\mathrm{alg}}^{p}\hspace{1pt}X
\longrightarrow
\widetilde{\mathrm{C}}_{\mathrm{top}}^{p}\hspace{1pt}X
$$
est un isomorphisme pour $p\leq j$.
\end{pro}

\bigskip
\textit{Démonstration.} L'implication (ii)$\Rightarrow$(iii) du théorème \ref{genalg} montre que l'on a $\mathrm{H}^{p}\hspace{1pt}\widetilde{\mathrm{C}}^{\bullet}\mathrm{H}_{V}^{*}X=0$ pour $p\leq j-2$. Le point (a) de \ref{pre-reciproqueAFP} dit alors que $\mathrm{Fix}_{(V,W)}\mathrm{H}_{V}^{*}X$ est libre comme $\mathrm{H}^{*}V/W$-module pour $\mathop{\mathrm{codim}}W\leq j$. Comme l'on a $\mathrm{H}_{V/W}^{*}X\cong\mathrm{Fix}_{(V,W)}\mathrm{H}_{V}^{*}X$ (Proposition \ref{cohFix}) la proposition \ref{CalgCtop} dit que $\varkappa_{V/W,X^{W}}:\widetilde{\mathrm{C}}_{\mathrm{alg},V/W}^{\bullet}\hspace{1pt}X
\to\widetilde{\mathrm{C}}_{\mathrm{top},V/W}^{\bullet}\hspace{1pt}X$ est un isomorphisme pour $\mathop{\mathrm{codim}}W\leq j$. La proposition \ref{decompkappa} montre que $\varkappa_{X}^{p}:\widetilde{\mathrm{C}}_{\mathrm{alg}}^{p}\hspace{1pt}X\to\widetilde{\mathrm{C}}_{\mathrm{top}}^{p}\hspace{1pt}X$ est un isomorphisme pour $0\leq p\leq j$~; le cas particulier $p=-1$ est trivial.
\hfill$\square$

\bigskip
Nous sommes à présent en mesure d'obtenir la version topologique du théo\-rème \ref{genalg}~:

\pagebreak

\medskip
\begin{theo}\label{gentop} Soient $V$ un $2$-groupe abélien élémentaire, $X$ un $V$-CW-complexe fini et $j\geq 0$ un entier. Les trois propriétés suivantes sont équi\-valentes~:
\begin{itemize}
\item[(i)] $\mathrm{H}_{V}^{*}X$ est une $(V_{\mathrm{tf}}\text{-}\mathcal{U})$-$j$-syzygie~;
\item[(ii)] le $\mathrm{H}^{*}V$-module sous-jacent à $\mathrm{H}_{V}^{*}X$ est une $\mathrm{H}^{*}V$-$j$-syzygie~;
\item[(iii)] on a $\mathrm{H}^{p}\hspace{1pt}\widetilde{\mathrm{C}}_{\mathrm{top}}^{\bullet}X=0$ pour $p\leq j-2$.
\end{itemize}
\end{theo}

\medskip
L'implication  (i)$\Rightarrow$(ii) est toujours triviale~!

\medskip
\textit{Démonstration de (ii)$\Rightarrow$(iii).} Le théorème \ref{genalg} dit que la condition (ii) implique  $\mathrm{H}^{p}\hspace{1pt}\widetilde{\mathrm{C}}_{\mathrm{alg}}^{\bullet}X=0$ pour $p\leq j-2$. Compte tenu de \ref{genCalgCtop} cette dernière condition implique  $\mathrm{H}^{p}\hspace{1pt}\widetilde{\mathrm{C}}_{\mathrm{top}}^{\bullet}X=0$ pour $p\leq j-2$.
\hfill$\square$

\medskip
\textit{Démonstration de (iii)$\Rightarrow$(i).} Soit $W$ un sous-groupe de $V$~; l'exactitude du foncteur $\mathrm{Fix}_{(V,W)}$ et la proposition \ref{FixCtop} montrent que la condition (iii) implique $\mathrm{H}^{p}\hspace{1pt}\widetilde{\mathrm{C}}_{\mathrm{top},V/W}^{\bullet}X^{W}=0$ pour $p\leq j-2$. On suppose maintenant $\mathop{\mathrm{codim}}W\leq j$~; on a \textit{a fortiori} $\mathrm{H}^{p}\hspace{1pt}\widetilde{\mathrm{C}}_{\mathrm{top},V/W}^{\bullet}X^{W}=0$ pour $p\leq\mathop{\mathrm{codim}}W-2$ si bien que le lemme \ref{top-n-2} implique que $\widetilde{\mathrm{C}}_{\mathrm{top},V/W}^{\bullet}X^{W}$ est acyclique. L'implication (ii)$\Rightarrow$(i) du théorème \ref{Ctop} montre alors que $\mathrm{H}_{V/W}^{*}X^{W}$ est libre comme $\mathrm{H}^{*}V/W$-module. La proposition \ref{CalgCtop} dit ensuite que $\varkappa_{X^{W},V/W}:\widetilde{\mathrm{C}}_{\mathrm{alg},V/W}^{\bullet}\hspace{1pt}X\to\widetilde{\mathrm{C}}_{\mathrm{top},V/W}^{\bullet}\hspace{1pt}X^{W}$ est un isomorphisme. Du coup, la proposition \ref{decompkappa} montre que l'homomorphisme $\varkappa_{X}^{p}:\widetilde{\mathrm{C}}_{\mathrm{alg}}^{p}\hspace{1pt}X\to\widetilde{\mathrm{C}}_{\mathrm{top}}^{p}\hspace{1pt}X$ est un isomorphisme pour $p\leq j$ et donc que l'on a $\mathrm{H}^{p}\hspace{1pt}\widetilde{\mathrm{C}}_{\mathrm{alg}}^{\bullet}X=0$ pour $p\leq j-2$. On conclut en invoquant l'implication (iii)$\Rightarrow$(i) du théorème \ref{genalg}.
\hfill$\square$

\vspace{0,75cm}
Nous achevons cette section en expliquant comme promis le lien entre le théorème \ref{Cntop} et la suite spectrale du théorème 0.3 de \cite{Henntop}. Dans cet article Henn considère un $G$-CW-complexe $X$, disons avec $G$ un groupe de Lie compact, et étudie entre autres une suite spectrale convergeant vers la cohomologie équivariante modulo $\ell$ ($\ell$ un nombre premier fixé) de sa partie $\ell$-singulière $X_{\mathrm{s}}$ ($X_{\mathrm{s}}$ est constitué des points fixés par un élément d'ordre $\ell$). Le point de départ de \cite{Henntop} est \cite{JMcC}, cependant Henn mentionne dans son introduction que le cas où $G$ est abélien ne nécessite pas les résultats de Jackowski-McClure et est très élémentaire. Nous prenons $G=V$. Soit $X$ un $V$-CW-complexe~; Henn utilise implicitement le fait que l'application canonique
$$
X_{\mathrm{s}}=\mathrm{Sing}_{V}X={\mathop{\mathrm{colim}}}_{\mathcal{W}_{0}^{\mathrm{op}}}X^{W}
\longrightarrow
{\mathop{\mathrm{hocolim}}}_{\mathcal{W}_{0}^{\mathrm{op}}}X^{W}
$$
est une équivalence d'homotopie et donc qu'il en est de même pour l'application (tout aussi canonique)
$$
\hspace{24pt}
\mathrm{E}V\times_{V}\mathrm{Sing}_{V}X
\longrightarrow
{\mathop{\mathrm{hocolim}}}_{\mathcal{W}_{0}^{\mathrm{op}}}\hspace{4pt}\mathrm{E}V\times_{V}X^{W}
\hspace{24pt}.
$$
Les cas $n=0$ et $n=1$ ($n:=\dim V$) sont triviaux~; nous supposons $n\geq 2$. On dispose d'une suite spectrale cohomologique du premier quadrant
$$
\hspace{24pt}
\mathrm{E}_{2}^{p,q}
\hspace{4pt}=\hspace{4pt}
{\lim_{\mathcal{W}_{0}}}^{\hspace{1pt}p}\hspace{2pt}
\mathrm{H}^{q}_{V}X^{W}
\hspace{4pt}\Rightarrow\hspace{4pt}
\mathrm{H}^{p+q}_{V}\hspace{1pt}\mathrm{Sing}_{V}X
\hspace{24pt}.
$$
Compte tenu de la proposition \ref{derivePf}, on a sous l'hypothèse ``$\mathrm{H}^{*}_{V}X$ libre comme $\mathrm{H}^{*}V$-module"~:
$$
\mathrm{E}_{2}^{p,*}
\hspace{4pt}=\hspace{4pt}
\begin{cases}
\mathrm{H}^{*}_{V}X & \text{pour}\hspace{6pt}p=0,
\\
0 & \text{pour}\hspace{6pt} 1\leq p\leq n-2,
\\
\mathrm{R}^{n}\mathrm{Pf}\hspace{1pt}\mathrm{H}^{*}_{V}X & \text{pour}\hspace{6pt}p=n-1.
\end{cases}
$$
Soit $\mathrm{b}:\mathrm{H}^{*}_{V}\mathrm{Sing}_{V}X\to\mathrm{E}_{2}^{0,*}=\mathrm{H}^{*}_{V}X$ l'homomorphisme de bord~; le composé de l'homomorphisme de restriction $\mathrm{H}^{*}_{V}X\to\mathrm{H}^{*}_{V}\mathrm{Sing}_{V}X$ et de $\mathrm{b}$ est l'identité. Pour s'en convaincre considérer l'application $\mathop{\mathrm{hocolim}}_{\mathcal{W}_{0}^{\mathrm{op}}}\hspace{2pt}\mathrm{E}V\times_{V}X^{W}\to\mathop{\mathrm{hocolim}}_{\mathcal{W}_{0}^{\mathrm{op}}}\hspace{2pt}\mathrm{E}V\times_{V}X$ (la colimite homotopique de droite étant celle du foncteur ``constant'', $W\mapsto\mathrm{E}V\times_{V}X$) et utiliser un argument de fonctorialité. Nous laissons au lecteur le soin d'achever l'analyse du cas $n=2$ et nous supposons $n\geq 3$. Dans ce cas il existe \textit{a priori} une seule différentielle non triviale à savoir $\mathrm{d}_{n-1}$ et  l'on a une suite exacte canonique (de $\mathrm{H}^{*}V$-$\mathrm{A}$-modules instables)
\begin{multline*}
\hspace{24pt}
\Sigma\hspace{1pt}\mathrm{H}^{*}_{V}X
\overset{\Sigma\hspace{1pt}\mathrm{d}_{n-1}}{\longrightarrow}
\Sigma^{n-1}\hspace{1pt}\mathrm{R}^{n}\mathrm{Pf}\hspace{1pt}\mathrm{H}^{*}_{V}X
\longrightarrow\mathrm{H}^{*}_{V}\mathrm{Sing}_{V}X
\overset{\mathrm{b}}{\longrightarrow}
\mathrm{H}^{*}_{V}X \\
\overset{\mathrm{d}_{n-1}}{\longrightarrow}
\Sigma^{n-2}\hspace{1pt}\mathrm{R}^{n}\mathrm{Pf}\hspace{1pt}\mathrm{H}^{*}_{V}X
\hspace{24pt};
\end{multline*}
comme $\mathrm{b}$ est surjectif la différentielle $\mathrm{d}_{n-1}$ est triviale.

\bigskip
La discussion précédente conduit à l'énoncé suivant~:

\bigskip
\begin{pro}\label{Henn} Soit $X$ un $V$-CW-complexe fini. Si $\mathrm{H}^{*}_{V}X$ est libre comme $\mathrm{H}^{*}V$-module alors on a des isomorphismes canoniques de $\mathrm{H}^{*}V$-$\mathrm{A}$-modules instables~:

\medskip
\hspace{12pt}
$\mathrm{H}^{*}_{V}\mathrm{Sing}_{V}X
\hspace{4pt}\cong\hspace{4pt}
\mathrm{H}^{*}_{V}X\oplus\Sigma^{n-1}\mathrm{R}^{n}\mathrm{Pf}\hspace{1pt}\mathrm{H}^{*}_{V}X
\hspace{12pt};$

\medskip
\hspace{12pt}
$\mathrm{H}^{*}_{V}\mathrm{Sing}_{V}X
\hspace{4pt}\cong\hspace{4pt}
\mathrm{H}^{*}_{V}X\oplus\Sigma^{-1}\mathrm{H}^{*}_{V}(X,\mathrm{Sing}_{V}X)
\hspace{12pt};$

\medskip
\hspace{12pt}
$\mathrm{H}^{*}_{V}(X,\mathrm{Sing}_{V}X)
\hspace{4pt}\cong\hspace{4pt}
\Sigma^{n}\hspace{1pt}\mathrm{R}^{n}\mathrm{Pf}\hspace{1pt}\mathrm{H}^{*}_{V}X
\hspace{12pt}.$
\end{pro}

\pagebreak

\sect{Illustrations}

Dans cette section nous illustrons par quelques exemples certains des énoncés des sections précédentes. Toutes ces illustrations dérivent peu ou prou de la considération de représentations linéaires réelles d'un $2$-groupe abélien élémentaire $V$.

\bigskip
\begin{pro}\label{representation} Soit $E$ un $\mathbb{R}$-espace vectoriel de dimension finie muni d'une action linéaire de $V$, en d'autres termes une représentation linéaire réelle de dimension finie du groupe $V$.

\smallskip
On pose $n=\dim_{\mathbb{Z}/2}V$, $m=\dim_{\mathbb{R}}E$ et $f=\dim_{\mathbb{R}}E^{V}$ ($E^{V}$ désigne le sous-espace de $E$ constitué des vecteurs invariants par $V$).

\smallskip
On note $\mathrm{w}_{k}(E)$ la $k$-ième classe de Stiefel-Whitney, appartenant à $\mathrm{H}^{k}V$, de la représentation linéaire $E$.

\smallskip
On note enfin $E_{\mathrm{r\acute{e}g}}$ le plus grand ouvert de $E$ sur lequel l'action de $V$ est libre et $V\backslash E_{\mathrm{r\acute{e}g}}$ le quotient de cette action ($V\backslash E_{\mathrm{r\acute{e}g}}$ est une variété de classe $\mathrm{C}^{\infty}$).

\medskip
{\em (a)} La classe $\mathrm{w}_{m-f}(E)$ est produit de classes de  $\mathrm{H}^{1}V-\{0\}$.

\medskip
{\em (b)} On a un isomorphisme canonique de $\mathrm{H}^{*}V$-$\mathrm{A}$-modules instables
$$
\hspace{24pt}
\mathrm{H}^{*}_{\mathrm{c}}(V\backslash E_{\mathrm{r\acute{e}g}})
\hspace{4pt}\cong\hspace{4pt}
\Sigma^{f+n}\hspace{2pt}\mathrm{R}^{n}\mathrm{Pf}_{V}\hspace{1pt}(\hspace{1pt}\mathrm{w}_{m-f}(E)\hspace{1pt}\mathrm{H}^{*}V\hspace{1pt})
\hspace{24pt}.
$$

\smallskip
{\em (La notation $\mathrm{H}^{*}_{\mathrm{c}}$ désigne la cohomologie modulo $2$ à support compact, pour ce qui est de la structure de $\mathrm{H}^{*}V$-$\mathrm{A}$-module instable de $\mathrm{H}^{*}_{\mathrm{c}}(V\backslash E_{\mathrm{r\acute{e}g}})$ voir la parenthèse qui suit l'énoncé \ref{finitude3}.)}
\end{pro}

\bigskip
\textit{Démonstration.} Comme $V$ est un groupe fini on peut supposer, ce que nous faisons ci-après, que la représentation linéaire $E$ est orthogonale, en clair que le $\mathbb{R}$-espace vectoriel $E$ est euclidien et que $V$ agit par isométries.

\medskip
 On considère le fibré vectoriel euclidien $\mathrm{E}\times_{V}E\to\mathrm{B}V$ que l'on note $\xi(E)$~; les classes de Stiefel-Whitney de la représentation $E$ sont par définition celles de~$\xi(E)$. Soit $u$ un élément du dual $V^{*}$ de $V$~; $\mathbb{R}$ muni de l'action de $V$ définie par $v.x=(-1)^{u(v)}x$ est une représentation orthogonale notée~$\mathbb{R}_{u}$. Si l'on identifie $\mathrm{H}^{1}V$ avec $V^{*}$ alors on a $\mathrm{w}_{1}(\mathbb{R}_{u})=u$. La représentation orthogonale $E$ est isomorphe à une somme directe $\bigoplus_{i\in I}\mathbb{R}_{u_{i}}$ ($I$ ensemble fini à $m$ éléments). On note $J$ (resp. $K$) le sous-ensemble de $I$ constitué des éléments $i$ de $I$ avec $u_{i}=0$ (resp. $u_{i}\not=0$)~; par définition $J$ a $f$ éléments. Puisque l'on a $\xi(E)\cong\bigoplus_{i\in I}\xi(\mathbb{R}_{u_{i}})$, la classe de Stiefel-Whitney ``totale'' est le produit $\prod_{i\in I}(1+u_{i})=\prod_{i\in K}(1+u_{i})$~; on a donc en particulier $\mathrm{w}_{m-f}(E)=\prod_{i\in K}u_{i}$ ce qui donne le point (a).

\medskip
Passons au point (b). L'isomorphisme de Thom conduit à l'énoncé suivant~:

\begin{pro}\label{Thom} Soit $\mathrm{S}(E\oplus\mathbb{R}_{0})$ la sphère unité de la représentation orthogonale $E\oplus\mathbb{R}_{0}$ de $V$. On a un isomorphisme de $\mathrm{H}^{*}V$-$\mathrm{A}$-modules instables
$$
\hspace{24pt}
\mathrm{H}^{*}_{V}\hspace{1pt}\mathrm{S}(E\oplus\mathbb{R}_{0})
\hspace{4pt}\cong\hspace{4pt}
\Sigma^{f}\hspace{1pt}\mathrm{w}_{m-f}(E)\hspace{1pt}\mathrm{H}^{*}V
\hspace{2pt}\oplus\hspace{2pt}
\mathrm{H}^{*}V
\hspace{24pt}.
$$
\end{pro}

\textit{Démonstration.} Pour alléger la notation on pose $\mathrm{S}:=\mathrm{S}(E\oplus\mathbb{R}_{0})$. On note respectivement $\mathrm{N}_{+}$ et $\mathrm{N}_{-}$ les points $(0,1)$ et $(0,-1)$ de $\mathrm{S}$~; comme ces points sont fixes sous l'action de $V$, les ouverts $\mathrm{S}-\mathrm{N}_{+}$ et $\mathrm{S}-\mathrm{N}_{-}$ sont stables sous l'action de $V$. On note respectivement $\mathrm{st}_{+}$ et  $\mathrm{st}_{-}$ les projections stéréographiques de $\mathrm{S}-\mathrm{N}_{+}$ et $\mathrm{S}-\mathrm{N}_{-}$  sur $E$~; on vérifie que ces deux homéomorphismes sont  $V$-équivariants.

\smallskip
On considère maintenant le diagramme commutatif
$$
\begin{CD}
\mathrm{H}^{*}_{V}\hspace{1pt}\mathrm{S}@>>>\mathrm{H}^{*}_{V}(\mathrm{S}-\mathrm{N}_{-}) \\
@AAA @AAA \\
\mathrm{H}^{*}V@>=>>\mathrm{H}^{*}V
\end{CD}
$$
dans lequel les flèches verticales sont les homomorphismes de $\mathrm{A}$-algèbres instables $\mathrm{H}^{*}_{V}(\mathrm{point})\to\mathrm{H}^{*}_{V}(-)$~; celle de droite s'identifie \textit{via} $\mathrm{H}^{*}_{V}(\mathrm{st}_{-})$ à l'homomorphisme $\mathrm{H}^{*}V\to\mathrm{H}^{*}_{V}\hspace{1pt}E$ et est donc un isomorphisme. On en déduit que la longue suite exacte de la paire $(\mathrm{S},\mathrm{S}-\mathrm{N}_{-})$ en cohomologie équivariante donne en fait une suite exacte courte canoniquement scindée
$$
0\to\mathrm{H}^{*}_{V}(\mathrm{S},\mathrm{S}-\mathrm{N}_{-})\to\mathrm{H}^{*}_{V}\hspace{1pt}\mathrm{S}\to\mathrm{H}^{*}_{V}(\mathrm{S}-\mathrm{N}_{-})\to 0
$$
soit encore un isomorphisme de  $\mathrm{H}^{*}V$-$\mathrm{A}$-modules instables
$$
\hspace{24pt}
\mathrm{H}^{*}_{V}\hspace{1pt}\mathrm{S}
\hspace{4pt}\cong\hspace{4pt}
\mathrm{H}^{*}_{V}(\mathrm{S},\mathrm{S}-\mathrm{N}_{-})
\hspace{2pt}\oplus\hspace{2pt}
\mathrm{H}^{*}V
\hspace{24pt}.
$$
Or on a par excision $\mathrm{H}^{*}_{V}(\mathrm{S},\mathrm{S}-\mathrm{N}_{-})\cong\mathrm{H}^{*}_{V}(\mathrm{S}-\mathrm{N}_{+},(\mathrm{S}-\mathrm{N}_{+})-\mathrm{N}_{-})$ et le membre de droite s'identifie \textit{via} $\mathrm{H}^{*}_{V}(\mathrm{st}_{+})$ à $\mathrm{H}^{*}_{V}(E,E-\{0\})$ qui d'après l'isomorphisme de Thom pour le fibré $\xi(E)$ est un $\mathrm{H}^{*}V$-module libre de dimension~$1$ (engendré par la classe de Thom). On dispose dans notre contexte d'un énoncé plus précis concernant $\mathrm{H}^{*}_{V}(E,E-\{0\})$ qui permet d'achever la démonstration de \ref{Thom}~:

\medskip
\begin{lem}\label{Thomsp} On a un isomorphisme de $\mathrm{H}^{*}V$-$\mathrm{A}$-modules instables
$$
\hspace{24pt}
\mathrm{H}^{*}_{V}(E,E-\{0\})
\hspace{4pt}\cong\hspace{4pt}
\Sigma^{f}\hspace{1pt}\mathrm{w}_{m-f}(E)\hspace{1pt}\mathrm{H}^{*}V
\hspace{24pt}.
$$
\end{lem}

\textit{Démonstration.} Soit $\widetilde{E}$ l'orthogonal de $E^{V}$ dans $E$~; $E^{V}$ et $\widetilde{E}$ sont stables par~$V$ et l'on a un isomorphisme de représentations orthogonales $E\cong E^{V}\oplus\widetilde{E}$. La paire d'espaces $(\mathrm{E}V\times_{V}E,\mathrm{E}V\times_{V}(E-{0}))$ s'écrit
$$
\hspace{24pt}
(E^{V},E^{V}-\{0\})\times(\mathrm{E}V\times_{V}\widetilde{E},\mathrm{E}V\times_{V}(\widetilde{E}-{0}))
\hspace{24pt};
$$
on en déduit que l'on a un isomorphisme de $\mathrm{H}^{*}V$-$\mathrm{A}$-modules instables
$$
\hspace{24pt}
\mathrm{H}^{*}_{V}(E,E-\{0\})
\hspace{4pt}\cong\hspace{4pt}
\Sigma^{f}\hspace{1pt}\mathrm{H}^{*}_{V}(\widetilde{E},\widetilde{E}-\{0\})
\hspace{24pt}.
$$
Par définition même de la classe d'Euler de $\xi(\widetilde{E})$ (ici modulo 2), l'homomorphisme composé
$$
\begin{CD}
\mathrm{H}^{*}V@>\text{isomorphisme de Thom}>>
\mathrm{H}^{*}_{V}(\widetilde{E},\widetilde{E}-\{0\})@>>>
\mathrm{H}^{*}_{V}\hspace{1pt}\widetilde{E}
\cong
\mathrm{H}^{*}V
\end{CD}
$$
est la multiplication par cette classe. Comme la classe d'Euler de $\xi(\widetilde{E})$ coïncide avec $\mathrm{w}_{\dim\xi(\widetilde{E})}(\xi(\widetilde{E}))=\mathrm{w}_{m-f}(\xi(\widetilde{E}))=\mathrm{w}_{m-f}(\xi(E)):=\mathrm{w}_{m-f}(E)$, que $\mathrm{w}_{m-f}(E)$ est non nulle (point (a)) et que $\mathrm{H}^{*}V$ est intègre, on constate que l'homomorphisme de $\mathrm{H}^{*}V$-$\mathrm{A}$-modules instables $\mathrm{H}^{*}_{V}(\widetilde{E},\widetilde{E}-\{0\})\to\mathrm{H}^{*}_{V}\hspace{1pt}\widetilde{E}\cong\mathrm{H}^{*}V$ est injectif et que son image s'identifie à $\mathrm{w}_{m-f}(E)\hspace{1pt}\mathrm{H}^{*}V$.
\hfill$\square\square$

\bigskip
On achève maintenant la démonstration du point (b) à l'aide du théorème~\ref{Cntop}. La stratégie est la suivante~:

\medskip
1) On vérifie tout d'abord que le $V$-espace $\mathrm{S}:=\mathrm{S}(E\oplus\mathbb{R}_{0})$ admet une structure de $V$-CW-complexe fini. On pourrait invoquer \cite{Il}, mais nous effectuerons cette vérification de manière explicite ci-après.

\medskip
2) On constate que la proposition \ref{Thom} dit en particulier que $\mathrm{H}^{*}_{V}\hspace{1pt}\mathrm{S}$ est un $\mathrm{H}^{*}V$-module libre (de dimension $2$).

\medskip
3) Compte tenu de \ref{Thom}, le théorème \ref{Cntop} donne un isomorphisme de $\mathrm{H}^{*}V$-$\mathrm{A}$-modules instables
$$
\hspace{24pt}
\mathrm{H}^{*}_{V}(\mathrm{S},\mathrm{Sing}_{V}\mathrm{S})
\hspace{4pt}\cong\hspace{4pt}
\Sigma^{n}\hspace{1pt}\mathrm{R}^{n}\mathrm{Pf}_{V}\hspace{1pt}(\hspace{1pt}\Sigma^{f}\hspace{2pt}\mathrm{w}_{m-f}(E)\hspace{1pt}\mathrm{H}^{*}V\hspace{2pt}\oplus\hspace{2pt}\mathrm{H}^{*}V\hspace{1pt})
\hspace{24pt}.
$$

\medskip
4) On observe que l'on a $\mathrm{R}^{n}\mathrm{Pf}_{V}(\Sigma^{f}\mathrm{w}_{m-f}(E)\mathrm{H}^{*}V)\hspace{-1pt}\cong\hspace{-1pt}\Sigma^{f}\mathrm{R}^{n}\mathrm{Pf}_{V}(\mathrm{w}_{m-f}(E)\mathrm{H}^{*}V)$ d'après \ref{corefini} et $\mathrm{R}^{n}\mathrm{Pf}_{V}\hspace{1pt}\mathrm{H}^{*}V=0$ car $\mathrm{H}^{*}V$ est un $V_{\mathrm{tf}}\text{-}\mathcal{U}$-injectif (le cas $n=0$ est trivial).

\medskip
5) On a un isomorphisme de $\mathrm{H}^{*}V$-$\mathrm{A}$-modules instables
$$
\mathrm{H}^{*}_{V}(\mathrm{S},\mathrm{Sing}_{V}\mathrm{S})
\hspace{4pt}\cong\hspace{4pt}
\mathrm{H}^{*}_{\mathrm{c}}(V\backslash(\mathrm{S}-\mathrm{Sing}_{V}\mathrm{S}))
$$
d'après le scholie \ref{finitude3}.

\medskip
6) On observe que la projection stéréographique $\mathrm{st}_{+}$ induit un homéomor\-phisme $V$-équivariant de $\mathrm{S}-\mathrm{Sing}_{V}\mathrm{S}$ sur $E_{\mathrm{r\acute{e}g}}$.

\medskip
On s'acquitte enfin de la promesse faite en 1)~:

\bigskip
\textit{Triangulation hyperoctaèdrique de la sphère unité d'une représentation orthogonale de dimension finie d'un $2$-groupe abélien élémentaire}

\bigskip
Soient $R$ une représentation orthogonale de dimension $d$ de $V$ et $u_{i}:V\to\nolinebreak\mathbb{Z}/2$, $1\leq i\leq d$, des formes linéaires telles que l'on a $R\cong\bigoplus_{1\leq i\leq d}\hspace{2pt}\mathbb{R}_{u_{i}}$. On note~:

\smallskip
--\hspace{8pt}$\varrho:V\to(\mathbb{Z}/2)^{d}$ le produit des homomorphismes $u_{i}$,

\smallskip
--\hspace{8pt}$\{\varepsilon_{1},\varepsilon_{2},\ldots,\varepsilon_{d}\}$ la base canonique de $(\mathbb{Z}/2)^{d}$,

\smallskip
--\hspace{8pt}$\{\mathrm{e}_{1},\mathrm{e}_{2},\ldots,\mathrm{e}_{d}\}$ la base canonique de $\mathbb{R}^{d}$,

\smallskip
--\hspace{8pt}$\mathrm{U}^{d}$ la représentation orthogonale $\bigoplus_{1\leq i\leq d}\hspace{2pt}\mathbb{R}_{\varepsilon_{i}}$ de $(\mathbb{Z}/2)^{d}$, en clair $\mathrm{U}^{d}$ est l'espace euclidien  $\mathbb{R}^{d}$ muni de l'action de $(\mathbb{Z}/2)^{d}$ définie par
$$
\varepsilon_{i}.\mathrm{e}_{j}
\hspace{4pt}=\hspace{4pt}
\begin{cases}
-\mathrm{e}_{j} & \text{pour $i=j$}, \\
\mathrm{e}_{j} & \text{pour $i\not=j$}.
\end{cases}
$$
Soit $\Delta^{d-1}$ le $(d-1)$-simplexe standard de $\mathbb{R}^{d}$, c'est-à-dire l'enveloppe convexe des points $\mathrm{e}_{1},\mathrm{e}_{2},\ldots,\mathrm{e}_{d}$. On pose
$$
\hspace{24pt}
\mathrm{O}^{d-1}
\hspace{4pt}:=\hspace{4pt}
\bigcup_{\varepsilon\in(\mathbb{Z}/2)^{d}}
\varepsilon\hspace{1pt}.\hspace{1pt}\Delta^{d-1}
\hspace{24pt};
$$
$\mathrm{O}^{d-1}$ est {\em l'hyperoctaèdre standard} de $\mathbb{R}^{d}$ ($\mathrm{O}^{2}$ est l'octaèdre standard de $\mathbb{R}^{3}$). La projection radiale
$$
\mathrm{O}^{d-1}\to\mathrm{S}^{d-1}=\mathrm{S}(\mathrm{U}^{d})=\mathrm{S}(R)
\hspace{24pt},\hspace{24pt}
x\mapsto\frac{x}{\Vert x\Vert}
$$
est un homéomorphisme, disons $\tau$, qui est la triangulation hyperoctaèdrique mentionnée dans l'intertitre. Par construction le groupe $(\mathbb{Z}/2)^{d}$ agit sur $\mathrm{O}^{d-1}$ et $\tau:\mathrm{O}^{d-1}\to\mathrm{S}(\mathrm{U}^{d})$ est $(\mathbb{Z}/2)^{d}$-équivariante~; l'action de $(\mathbb{Z}/2)^{d}$ sur $\mathrm{O}^{d-1}$ induit \textit{via} $\varrho$ une action de $V$ et $\tau:\mathrm{O}^{d-1}\to\mathrm{S}(R)$ est $V$-équivariante. 

\medskip
Soit $\Sigma_{m}$ l'ensemble des $m$-simplexes de $\mathrm{O}^{d-1}$~; les $m$-cellules de la structure de $V$-CW-complexe de $\mathrm{S}(R)$ sont indexées par $V\backslash\Sigma_{m}$.
\hfill$\square$

\bigskip
\begin{rem}\label{questionsuspension} Le point (b) de la proposition \ref{representation} dit en particulier que le\linebreak $\mathrm{A}$-module instable $\mathrm{H}^{*}_{\mathrm{c}}(V\backslash E_{\mathrm{r\acute{e}g}})$ est une suspension $(f+n)$-ième d'un $\mathrm{A}$-module instable. Existe-t-il une explication ``topologique'' de ce phénomène~? (On peut supposer $f=0$ puisque l'espace $V\backslash E_{\mathrm{r\acute{e}g}}$ est homéomorphe au produit $E^{V}\times V\backslash\widetilde{E}_{\mathrm{r\acute{e}g}}\simeq\mathbb{R}^{f}\times\widetilde{E}_{\mathrm{r\acute{e}g}}$\hspace{2pt}.)
\end{rem}

\bigskip
La proposition ci-dessous dit notamment que tous les $\mathrm{H}^{*}V$-$\mathrm{A}$-modules instables qui sont libres de dimension $1$ comme $\mathrm{H}^{*}V$-modules, sont ceux qui apparaissent dans le lemme \ref{Thomsp}.

\bigskip
\begin{pro}\label{libredim1} Soit $M$ un $\mathrm{H}^{*}V$-$\mathrm{A}$-module instable qui est libre de dimension $1$ comme $\mathrm{H}^{*}V$-module.

\medskip
{\em (a)} Il existe un entier naturel $f$ et un élément $e$ de $\mathrm{H}^{*}V$, produit d'éléments de $\mathrm{H}^{1}V-\{0\}$, tels que l'on a un isomorphisme de $\mathrm{H}^{*}V$-$\mathrm{A}$-modules instables
$$
\hspace{24pt}
M
\hspace{4pt}\cong\hspace{4pt}
\Sigma^{f}\hspace{1pt}e\hspace{1pt}\mathrm{H}^{*}V
\hspace{24pt}.
$$

\medskip
{\em (b)} Le couple $(f,e)$ est uniquement déterminé par la classe d'isomorphisme de $M$.
\end{pro}

\bigskip
\textit{Démonstration.} Le théorème \ref{pendantalg} dit entre autres que  l'unité d'adjonction $\eta_{M}:M\to\mathrm{H}^{*}V\otimes\mathrm{Fix}_{V}M$ ($\mathrm{Fix}_{V}:=\mathrm{Fix}_{(V,V)}$) est injective et que son conoyau est de torsion~; en fait la théorie de Smith algébrique \cite{DWsmith2}\cite{LZsmith} dit que $\eta_{M}[\mathrm{c}_{V}^{-1}]$ est un isomorphisme pour tout $\mathrm{H}^{*}V_{\mathrm{tf}}$-$\mathrm{A}$-module instable $M$, $\mathrm{c}_{V}$
désignant le produit de tous les éléments de $\mathrm{H}^{1}V-\{0\}$. Il en résulte $\dim_{\mathrm{H}^{*}V}M=\dim_{\mathrm{H}^{*}V}(\mathrm{H}^{*}V\otimes\mathrm{Fix}_{V}M)=\dim_{\mathbb{F}_{2}}\mathrm{Fix}_{V}M$ (pour la première égalité oublier la graduation et tensoriser par le corps des fractions de $\mathrm{H}^{*}V$)~; on a donc ici $\dim_{\mathbb{F}_{2}}\mathrm{Fix}_{V}M=1$. Cette égalité montre qu'il existe un entier naturel $f$, uniquement déterminé, tel que le $A$-module instable $\mathrm{Fix}_{V}M$ est isomorphe à $\Sigma^{f}\hspace{1pt}\mathbb{F}_{2}$ si bien que l'on peut identifier $\mathrm{H}^{*}V\otimes\mathrm{Fix}_{V}M$ à $\Sigma^{f}\hspace{1pt}\mathrm{H}^{*}V$. Soit $g$ le générateur de $M$ (l'article est défini parce que le groupe des unités de $\mathrm{H}^{*}V$ est trivial~!)~; on pose $\eta_{M}(g)=\Sigma^{f}\hspace{1pt}e$. On constate que l'idéal $e\hspace{1pt}\mathrm{H}^{*}V$ de $\mathrm{H}^{*}V$ est stable sous l'action de $\mathrm{A}$ et que $\eta_{M}$ induit un isomorphisme $M\cong\Sigma^{f}\hspace{1pt}e\hspace{1pt}\mathrm{H}^{*}V$. Le fait que $e$ est un produit d'éléments de $\mathrm{H}^{1}V-\{0\}$ est conséquence d'un résultat de Serre \cite[\S 2, Corollaire]{Ser} (voir \ref{Serre4}). En effet ce résultat dit que l'idéal $e\hspace{1pt}\mathrm{H}^{*}V$ contient un produit d'éléments de $\mathrm{H}^{1}V-\{0\}$, disons $\pi$, comme $e$ divise $\pi$, $e$ est également un produit d'éléments de $\mathrm{H}^{1}V-\{0\}$.
\hfill$\square$

\pagebreak

\bigskip
\begin{scho}\label{libredim1bis} La correspondance
$$
E
\hspace{4pt}\mapsto\hspace{4pt}
\mathrm{H}^{*}_{V}(E,E-\{0\})
$$
induit une bijection entre classes d'isomorphisme de représentations linéaires réelles de $V$ et classes d'isomorphisme de $\mathrm{H}^{*}V$-$\mathrm{A}$-modules instables qui sont libres de dimension $1$ comme $\mathrm{H}^{*}V$-module.
\end{scho}

\bigskip
\begin{rem}\label{Fixdirect} Soit $(e,f)$ comme dans la proposition \ref{libredim1}. Il est facile de se convaincre directement de l'isomorphisme de $A$-modules instables
$$
\hspace{24pt}
\mathrm{Fix}_{V}\hspace{1pt}(\Sigma^{f}\hspace{1pt}e\hspace{1pt}\mathrm{H}^{*}V)
\hspace{4pt}\cong\hspace{4pt}
\Sigma^{f}\hspace{1pt}\mathbb{F}_{2}
\hspace{24pt}.
$$
Soit $\iota:\Sigma^{f}\hspace{1pt}e\hspace{1pt}\mathrm{H}^{*}V\to\Sigma^{f}\hspace{1pt}\mathrm{H}^{*}V$ l'inclusion évidente, on considère la suite exacte de $\mathrm{H}^{*}V$-$\mathrm{A}$-modules instables $0\to\Sigma^{f}\hspace{1pt}e\hspace{1pt}\mathrm{H}^{*}V\overset{\iota}{\to}\Sigma^{f}\hspace{1pt}\mathrm{H}^{*}V\to\mathop{\mathrm{coker}}\iota\to 0$. Comme $\mathop{\mathrm{coker}}\iota$ est un $\mathrm{H}^{*}V$-module de torsion, $\mathrm{Fix}_{V}\hspace{1pt}(\mathop{\mathrm{coker}}\iota)$ est nul, comme $\mathrm{Fix}_{V}$ est exact, $\mathrm{Fix}_{V}(\iota)$ est un isomorphisme. D'après le point (b) de \ref{suspensionFix} on~a $\mathrm{Fix}_{V}(\Sigma^{f}\hspace{1pt}\mathrm{H}^{*}V)\cong\Sigma^{f}\hspace{1pt}\mathrm{Fix}_{V}(\mathrm{H}^{*}V)$ et on constate que l'on a $\mathrm{Fix}_{V}(\mathrm{H}^{*}V)\cong\mathbb{F}_{2}$.
\end{rem}

\bigskip
\begin{rem}\label{librespecial} Posons
$$
M
\hspace{4pt}:=\hspace{4pt}
\bigoplus_{i=1}^{d}\hspace{2pt}\Sigma^{f_{i}}\hspace{1pt}e_{i}\hspace{1pt}\mathrm{H}^{*}V
$$
avec $(f_{i},e_{i})$ comme dans la proposition \ref{libredim1}. Le $\mathrm{A}$-module instable $\mathrm{Fix}_{V}M$ est isomorphe à la somme directe $\bigoplus_{i=1}^{d}\hspace{2pt}\Sigma^{f_{i}}\hspace{1pt}\mathbb{F}_{2}$ et est donc un $\mathrm{A}$-module instable trivial. Cette observation montre qu'un $\mathrm{H}^{*}V$-$\mathrm{A}$-module instable de la forme $\mathrm{H}^{*}V\otimes N$ avec $N$ un $\mathrm{A}$-module instable fini ($\mathrm{H}^{*}V$-$\mathrm{A}$-module instable qui est libre de dimension finie comme $\mathrm{H}^{*}V$-module) ne peut être isomorphe à une somme directe du type ci-dessus si $N$ n'est pas trivial. En effet on a $\mathrm{Fix}_{V}(\mathrm{H}^{*}V\otimes N)=N$, à nouveau d'après le point (b) de \ref{suspensionFix}.
\end{rem}

\bigskip
\begin{rem}\label{illustrationcohFix}Reprenons les notations de la proposition \ref{Thom}. Cette proposition et la remarque précédente montrent que l'on a un isomorphisme de $\mathrm{A}$-modules instables
$$
\hspace{24pt}
\mathrm{Fix}_{V}\hspace{1pt}\mathrm{H}^{*}_{V}\hspace{1pt}\mathrm{S}(E\oplus\mathbb{R}_{0})
\hspace{4pt}\cong\hspace{4pt}
\Sigma^{f}\hspace{1pt}\mathbb{F}_{2}
\hspace{2pt}\oplus\hspace{2pt}\mathbb{F}_{2}
\hspace{24pt}.
$$
Cet isomorphisme illustre la proposition \ref{cohFix} (avec $W=V$). En effet l'espace des points fixes $(\mathrm{S}(E\oplus\mathbb{R}_{0}))^{V}$ est la sphère $\mathrm{S}(E^{V}\oplus\mathbb{R}_{0})$ qui est de dimension~$f$.
\end{rem}

\pagebreak

\vspace{0.75cm}
\textsc{Etude du complexe $\widetilde{\mathrm{C}}^{\bullet}\hspace{1pt}(\mathrm{c}_{V}\hspace{1pt}\mathrm{H}^{*}V)$}

\bigskip
On prend maintenant pour $E$ la représentation régulière réelle que l'on note $\mathbb{R}[V]$~;  il sera commode de considérer les vecteurs de $\mathbb{R}[V]$ comme des fonctions $f:V\to\mathbb{R}$, l'action d'un élément $v_{0}$ de $V$ étant définie par la formule $(v_{0}.f)(v)=f(v+v_{0})$. On munit $\mathbb{R}[V]$ du produit scalaire $V$-équivariant
$$
(f,g)
\hspace{4pt}\mapsto\hspace{4pt}
\frac{1}{\vert V\vert}\sum_{v\in V}f(v)g(v)
$$
($\vert V\vert$ désignant le cardinal de $V$). Le sous-espace invariant $(\mathbb{R}[V])^{V}$ est consti\-tué des fonctions constantes, il est de dimension $1$ engendré par la fonction constante $1$.
On pose $\widetilde{\mathbb{R}}[V]:={((\mathbb{R}[V])^{V})}^{\perp}$. En clair $\widetilde{\mathbb{R}}[V]$ est le sous-espace de~$\mathbb{R}[V]$, stable sous l'action de $V$, constitué des fonctions $f$ vérifiant $\sum_{v\in V}f(v)=0$~; $\widetilde{\mathbb{R}}[V]$ est appelée la {\em représentation réelle régulière réduite} de~$V$. On dispose de deux isomorphismes de représentations orthogonales de~$V$~:
$$
\hspace{24pt}
\mathbb{R}[V]
\hspace{4pt}\cong\hspace{4pt}
\bigoplus_{u\in V^{*}}
\mathbb{R}_{u}
\hspace{24pt},\hspace{24pt}
\widetilde{\mathbb{R}}[V]
\hspace{4pt}\cong\hspace{4pt}
\bigoplus_{u\in V^{*}-\{0\}}
\mathbb{R}_{u}
\hspace{24pt};
$$
le second montre que la classe de Stiefel-Whitney $\mathrm{w}_{2^{n}-1}(\widetilde{\mathbb{R}}[V])$ (qui est aussi la classe d'Euler modulo $2$ de $\widetilde{\mathbb{R}}[V]$) est le produit  $\mathrm{c}_{V}:=\prod_{u\in\mathrm{H}^{1}V-\{0\}}u$. Le lemme \ref{Thomsp} se spécialise en $\mathrm{H}^{*}_{V}\hspace{1pt}(\widetilde{\mathbb{R}}[V],\widetilde{\mathbb{R}}[V]-\{0\})\cong\mathrm{c}_{V}\hspace{1pt}\mathrm{H}^{*}V$. Nous allons étudier ci-après le complexe de $\mathrm{H}^{*}V$-$\mathrm{A}$-modules instables $\widetilde{\mathrm{C}}^{\bullet}\hspace{1pt}(\mathrm{c}_{V}\hspace{1pt}\mathrm{H}^{*}V)$ (voir section 5). D'après le théorème  \ref{pendantalg} ce complexe est acyclique, si bien que la coaugmentation $\mathrm{c}_{V}\hspace{1pt}\mathrm{H}^{*}V\to\mathrm{C}^{\bullet}\hspace{1pt}(\mathrm{c}_{V}\hspace{1pt}\mathrm{H}^{*}V)$ peut être vue comme une résolution du $\mathrm{H}^{*}V$-$\mathrm{A}$-module instable $\mathrm{c}_{V}\hspace{1pt}\mathrm{H}^{*}V$~; c'est d'ailleurs la recherche d'une résolution de ce type qui a été le point de départ de ce mémoire. Notre stratégie va être d'utiliser la proposition \ref{CalgCtop}.

\bigskip
On commence par quelques observations~:

\medskip
($\mathrm{O}_{1}$)\hspace{8pt}On a $\mathrm{S}(\mathbb{R}[V])=\mathrm{S}(\widetilde{\mathbb{R}}[V]\oplus\mathbb{R}_{0})$ (la représentation $\mathbb{R}_{0}$ au second membre est engendrée par la fonction constante~$1$ qui est un vecteur unitaire de $\mathbb{R}[V]$ fixe sous l'action de $V$).

\medskip
($\mathrm{O}_{2}$)\hspace{8pt}Le sous-espace $(\mathbb{R}[V])^{W}$, $W$ sous-groupe de $V$, s'identifie à $\mathbb{R}[V/W]$~; plus précisément $(\mathbb{R}[V])^{W}$ vu comme une représentation orthogonale de $V/W$ s'identifie à $\mathbb{R}[V/W]$ et $(\mathrm{S}(\mathbb{R}[V]))^{W}$ s'identifie à $\mathrm{S}(\mathbb{R}[V/W])$

\medskip
($\mathrm{O}_{3}$)\hspace{8pt} Compte tenu de ($\mathrm{O}_{1}$) et de  la proposition \ref{Thom}, on a un isomorphisme canonique de $\mathrm{H}^{*}V$-$\mathrm{A}$-modules instables
$$
\hspace{24pt}
\mathrm{H}^{*}_{V}\hspace{1pt}\mathrm{S}(\mathbb{R}[V])
\hspace{4pt}\cong\hspace{4pt}
\mathrm{c}_{V}\hspace{1pt}\mathrm{H}^{*}V\oplus\mathrm{H}^{*}V
\hspace{24pt};
$$
l'inclusion $\mathrm{H}^{*}V\to\mathrm{H}^{*}_{V}\hspace{1pt}\mathrm{S}(\mathbb{R}[V])$ et la projection $\mathrm{H}^{*}_{V}\hspace{1pt}\mathrm{S}(\mathbb{R}[V])\to\mathrm{H}^{*}V$ s'identifient respectivement aux homomorphismes de $\mathrm{H}^{*}V$-$\mathrm{A}$-modules instables $\mathrm{H}^{*}V=\mathrm{H}^{*}_{V}(\mathrm{point})\to\mathrm{H}^{*}_{V}\hspace{1pt}\mathrm{S}(\mathbb{R}[V])$ et $\mathrm{H}^{*}_{V}\hspace{1pt}\mathrm{S}(\mathbb{R}[V])\to\mathrm{H}^{*}_{V}\hspace{1pt}\{1\}=\mathrm{H}^{*}V$.

\medskip
($\mathrm{O}_{4}$)\hspace{8pt} Compte tenu de ($\mathrm{O}_{3}$), on a un isomorphisme canonique de complexes de $\mathrm{H}^{*}V$-$\mathrm{A}$-modules instables
$$
\widetilde{\mathrm{C}}^{\bullet}\hspace{1pt}(\mathrm{H}^{*}_{V}\hspace{1pt}\mathrm{S}(\mathbb{R}[V]))
\hspace{4pt}\cong\hspace{4pt}
\widetilde{\mathrm{C}}^{\bullet}\hspace{1pt}(\mathrm{c}_{V}\hspace{1pt}\mathrm{H}^{*}V)\oplus\widetilde{\mathrm{C}}^{\bullet}\hspace{1pt}(\mathrm{H}^{*}V)
$$
et les homomorphismes canoniques $\widetilde{\mathrm{C}}^{\bullet}\hspace{1pt}(\mathrm{H}^{*}V)\leftrightarrow\widetilde{\mathrm{C}}^{\bullet}\hspace{1pt}(\mathrm{H}^{*}_{V}\hspace{1pt}\mathrm{S}(\mathbb{R}[V]))$ s'identifient aux homomorphisme évidents.

\medskip
($\mathrm{O}_{5}$)\hspace{8pt} On a $\mathrm{C}^{\bullet}\hspace{1pt}(\mathrm{H}^{*}V)=(\mathrm{H}^{*}V\to 0\to 0\to\ldots)$ et la coaugmentation $\mathrm{H}^{*}V\to\mathrm{C}^{0}\hspace{1pt}(\mathrm{H}^{*}V)$ est l'identité (observer que $\mathrm{H}^{*}V$ est un $V_{\mathrm{tf}}\text{-}\mathcal{U}$-injectif avec\linebreak $\mathrm{F}^{1}\mathrm{H}^{*}V=0$). Même chose pour $\mathrm{C}^{\bullet}_{\mathrm{top}}\hspace{1pt}(\{1\})$~;  l'isomorphisme $\varkappa:\widetilde{\mathrm{C}}^{\bullet}\hspace{1pt}(\mathrm{H}^{*}V)\to\widetilde{\mathrm{C}}^{\bullet}_{\mathrm{top}}\hspace{1pt}(\{1\})$ fourni par la proposition \ref{CalgCtop} ``est'' l'identité.

\bigskip
Nous sommes maintenant en mesure de décrire le complexe $\widetilde{\mathrm{C}}^{\bullet}\hspace{1pt}(\mathrm{c}_{V}\hspace{1pt}\mathrm{H}^{*}V)$. La proposition \ref{CalgCtop} et les observations précédentes montrent que l'on a un isomorphisme canonique de complexes de $\mathrm{H}^{*}V$-$\mathrm{A}$-modules instables
$$
\hspace{24pt}
\widetilde{\mathrm{C}}^{\bullet}\hspace{1pt}(\mathrm{c}_{V}\hspace{1pt}\mathrm{H}^{*}V)
\hspace{4pt}\cong\hspace{4pt}
\ker
\hspace{2pt}(\hspace{2pt}
\widetilde{\mathrm{C}}^{\bullet}_{\mathrm{top}}\hspace{1pt}\mathrm{S}(\mathbb{R}[V])
\to\widetilde{\mathrm{C}}^{\bullet}_{\mathrm{top}}\hspace{1pt}\{1\}
\hspace{2pt})
\hspace{18pt}.
$$
On pose
$$
\hspace{24pt}
\mathrm{M}(V)
\hspace{4pt}:=\hspace{4pt}
\ker
\hspace{2pt}(
\hspace{2pt}\mathrm{H}^{*}_{V}(\hspace{1pt}\mathrm{S}(\mathbb{R}[V]),\mathrm{Sing}_{V}\mathrm{S}(\mathbb{R}[V])\hspace{1pt})
\to
\mathrm{H}^{*}_{V}(\hspace{1pt}\{1\},\mathrm{Sing}_{V}\{1\}\hspace{1pt})
\hspace{2pt})
\hspace{18pt};
$$
on observera que l'on a $\mathrm{H}^{*}_{V}(\hspace{1pt}\{1\},\mathrm{Sing}_{V}\{1\}\hspace{1pt})=0$ pour $V\not=0$~!

\medskip
Le $\mathrm{H}^{*}V$-$\mathrm{A}$-module instable $\mathrm{M}(V)$ a deux avatars~: un avatar ``topologique'' $\mathrm{M}(V)\cong\Sigma^{-n}\hspace{1pt}\mathrm{H}^{*}_{\mathrm{c}}(V\backslash\widetilde{\mathbb{R}}[V]_{\mathrm{r\acute{e}g}})$ et  l'autre ``algébrique'' $\mathrm{M}(V)\cong\mathrm{R}^{n}\mathrm{Pf}_{V}\hspace{1pt}(\hspace{1pt}\mathrm{c}_{V}\hspace{1pt}\mathrm{H}^{*}V\hspace{1pt})$ ($n:=\dim V$). On constate que l'on a $\mathrm{M}(V)=\mathbb{F}_{2}$ pour $n=0$ et $n=1$, et $\mathrm{M}(V)\cong\mathrm{J}_{V}(1)$ pour $n=2$ ($\mathrm{J}_{V}(1)$ est l'un des $V\text{-}\mathcal{U}$-injectifs ``tautologiques'' introduits au début de la section 2).

\medskip
Voici la forme du complexe (acyclique) $\widetilde{\mathrm{C}}^{\bullet}\hspace{1pt}(\mathrm{c}_{V}\hspace{1pt}\mathrm{H}^{*}V)$~:
\begin{multline*}
\hspace{8pt}0\to\mathrm{c}_{V}\hspace{1pt}\mathrm{H}^{*}V\to\mathrm{H}^{*}V\to
\bigoplus_{\mathop{\mathrm{codim}}W=1}\mathrm{H}^{*}W\to\bigoplus_{\mathop{\mathrm{codim}}W=2}\mathrm{H}^{*}V\otimes_{\mathrm{H}^{*}V/W}\mathrm{M}(V/W) \\ 
\to\ldots\to 
\bigoplus_{\mathop{\mathrm{codim}}W=n-1}\mathrm{H}^{*}V\otimes_{\mathrm{H}^{*}V/W}\mathrm{M}(V/W)
\to\mathrm{M}(V)\to 0\hspace{8pt}.
\end{multline*}
(On constate par inspection que $\mathrm{c}_{V}\hspace{1pt}\mathrm{H}^{*}V\to\mathrm{C}^{\bullet}\hspace{1pt}(\mathrm{c}_{V}\hspace{1pt}\mathrm{H}^{*}V)$ est une résolution injective dans la catégorie $V_{\mathrm{tf}}\text{-}\mathcal{U}$ pour $n\leq 2$.)

\bigskip
\begin{rem} On aurait pu définir \textit{mutatis mutandis} $\widetilde{\mathrm{C}}^{\bullet}_{\mathrm{top}}\hspace{1pt}(X,Y)$ pour une paire $(X,Y)$ de $V$-CW-complexes finis et montrer que $\widetilde{\mathrm{C}}^{\bullet}\hspace{1pt}(\mathrm{H}^{*}_{V}(X,Y))$ et $\widetilde{\mathrm{C}}^{\bullet}_{\mathrm{top}}\hspace{1pt}(X,Y)$ sont naturellement isomorphes si $\mathrm{H}^{*}_{V}(X,Y)$ est libre comme $\mathrm{H}^{*}V$-module. L'isomorphisme $\widetilde{\mathrm{C}}^{\bullet}\hspace{1pt}(\mathrm{c}_{V}\hspace{1pt}\mathrm{H}^{*}V)
\cong\ker\hspace{2pt}(\hspace{2pt}
\widetilde{\mathrm{C}}^{\bullet}_{\mathrm{top}}\hspace{1pt}\mathrm{S}(\mathbb{R}[V])
\to\widetilde{\mathrm{C}}^{\bullet}_{\mathrm{top}}\hspace{1pt}\{1\}
\hspace{2pt})$ aurait pu être remplacé par $\widetilde{\mathrm{C}}^{\bullet}\hspace{1pt}(\mathrm{c}_{V}\hspace{1pt}\mathrm{H}^{*}V)
\cong\widetilde{\mathrm{C}}^{\bullet}_{\mathrm{top}}\hspace{1pt}(\mathrm{S}(\mathbb{R}[V]),\{1\})$ ce qui nous aurait évité quelques contorsions disgracieuses.
\end{rem}

\vspace{0.5cm}
\textit{Action de $\mathrm{GL}(V)$ sur $\widetilde{\mathrm{C}}^{\bullet}\hspace{1pt}(\mathrm{c}_{V}\hspace{1pt}\mathrm{H}^{*}V)$}

\medskip
La spécificité du complexe $\widetilde{\mathrm{C}}^{\bullet}\hspace{1pt}(\mathrm{c}_{V}\hspace{1pt}\mathrm{H}^{*}V)$ (et plus généralement des complexes $\widetilde{\mathrm{C}}^{\bullet}\hspace{1pt}(\mathrm{c}^{h}_{V}\hspace{1pt}\mathrm{H}^{*}V)$, $h\in\mathbb{N}$) parmi les complexes $\widetilde{\mathrm{C}}^{\bullet}\hspace{1pt}(e\hspace{1pt}\mathrm{H}^{*}V)$, avec $e$ un produit arbitraire d'éléments de $\mathrm{H}^{1}V-\{0\}$, est la suivante~:

\medskip
L'action à gauche tautologique de $\mathrm{GL}(V)$ sur $V$ induit une action à droite de $\mathrm{GL}(V)$ sur  $\mathrm{H}^{*}V$ qui préserve la structure de $\mathrm{A}$-module instable. En ce qui concerne la structure de $\mathrm{H}^{*}V$-module, cette action est ``tordue''~:

\medskip
\begin{defi}\label{GL-tordu} Soit $M$ un $\mathrm{H}^{*}V$-$\mathrm{A}$-module instable muni d'une action à droite de $\mathrm{GL}(V)$ qui préserve la structure de $\mathrm{A}$-module instable de $M$. Nous dirons que cette action est {\em tordue} si l'on a
$$
(\hspace{1pt}c\hspace{1pt}x\hspace{1pt})\hspace{1pt}.\hspace{1pt}\alpha
\hspace{4pt}=\hspace{4pt}
(\alpha^{*}c)\hspace{1pt}(\hspace{1pt}x\hspace{1pt}.\hspace{1pt}\alpha)
$$
pour tout $\alpha$ dans $\mathrm{GL}(V)$, tout $x$ dans $M$ et tout $c$ dans $\mathrm{H}^{*}V$.
\end{defi}

\bigskip
Comme $\mathrm{c}_{V}$ est invariant sous l'action à droite de $\mathrm{GL}(V)$ sur $\mathrm{H}^{*}V$, le $\mathrm{H}^{*}V$-$\mathrm{A}$-module instable $\mathrm{c}_{V}\hspace{1pt}\mathrm{H}^{*}V$ est également muni d'une action à droite tordue de~$\mathrm{GL}(V)$. On va voir ci-après que celle-ci se prolonge en une action à droite tordue de $\mathrm{GL}(V)$ sur $\widetilde{\mathrm{C}}^{\bullet}\hspace{1pt}(\mathrm{c}_{V}\hspace{1pt}\mathrm{H}^{*}V)$. Précisons lourdement. Soit $\mathrm{GL}(V)_{\mathrm{td}}\text{-}V\text{-}\mathcal{U}$ la catégorie dont les objets sont les $\mathrm{H}^{*}V$-$\mathrm{A}$-modules instables munis d'une action à droite tordue de $\mathrm{GL}(V)$ et dont les morphismes sont les $V\text{-}\mathcal{U}$-morphismes $\mathrm{GL}(V)$-équivariants~; $\mathrm{C}^{\bullet}\hspace{1pt}(\mathrm{c}_{V}\hspace{1pt}\mathrm{H}^{*}V)$ est un complexe dans cette catégorie et la coaugmentation $\mathrm{c}_{V}\hspace{1pt}\mathrm{H}^{*}V\to\mathrm{C}^{0}\hspace{1pt}(\mathrm{c}_{V}\hspace{1pt}\mathrm{H}^{*}V)$ est un $\mathrm{GL}(V)_{\mathrm{td}}\text{-}V\text{-}\mathcal{U}$-morphisme.

\vspace{0.5cm}
\textit{Définition ``topologique'' de l'action de $\mathrm{GL}(V)$ sur $\widetilde{\mathrm{C}}^{\bullet}\hspace{1pt}(\mathrm{c}_{V}\hspace{1pt}\mathrm{H}^{*}V)$}

\medskip
Soit $X$ un $V$-espace muni d'une action à gauche (continue) de $\mathrm{GL}(V)$. Nous pourrions imiter la définition \ref{GL-tordu} et dire que cette action est tordue si l'on a
$$
\alpha.(x+v)
\hspace{4pt}=\hspace{4pt}
\alpha.x+\alpha(v)
$$
pour tout $\alpha$ dans $\mathrm{GL}(V)$ et tout $v$ dans $V$
(les actions de  $\mathrm{GL}(V)$ et $V$ sur $X$ sont respectivement notées $(\alpha,x)\mapsto\alpha.x$ et $(v,x)\mapsto v+x$). Cependant nous n'introduirons pas cette terminologie car il est équivalent  de dire que $X$ est muni d'une action à gauche (continue) du produit semi-direct $V\rtimes\mathrm{GL}(V)$ (le {\em groupe affine} de $V$) qui prolonge celle de $V$.

\medskip
\begin{pro}\label{GL-action1}
Soit $X$ un $V$-espace muni d'une action à gauche (continue) de $V\rtimes\mathrm{GL}(V)$ qui prolonge celle de $V$.

\medskip
{\em (a)} Le $\mathrm{H}^{*}V$-$\mathrm{A}$-module instable $\mathrm{H}^{*}_{V}X$ est naturellement muni d'une action à droite tordue de $\mathrm{GL}(V)$.

\medskip
{\em (b)} Soit $Y$ un sous-espace de $X$ stable sous l'action de $V\rtimes\mathrm{GL}(V)$ (et donc de~$V$) alors le $\mathrm{H}^{*}V$-$\mathrm{A}$-module instable $\mathrm{H}^{*}_{V}(X,Y)$ est naturellement muni d'une action à droite tordue de $\mathrm{GL}(V)$ et les $V\text{-}\mathcal{U}$-morphismes
$$
\begin{CD}
\mathrm{H}^{*}_{V}(X,Y)@>>>\mathrm{H}^{*}_{V}X@>>>\mathrm{H}^{*}_{V}Y@>\partial>>\Sigma\hspace{1pt}\mathrm{H}^{*}_{V}(X,Y)
\end{CD}
$$
sont $\mathrm{GL}(V)$-équivariants.
\end{pro}

\bigskip
\textit{Démonstration.} On vérifie le point (a)~; la vérification du point (b) est laissée au lecteur. L'action tautologique de $V\rtimes\mathrm{GL}(V)$ sur $V$ induit une action à gauche de $V\rtimes\mathrm{GL}(V)$ sur $\mathrm{E}V$ et on constate que l'action diagonale de $V\rtimes\mathrm{GL}(V)$ sur $\mathrm{E}V\times X$ induit une action à gauche de $\mathrm{GL}(V)$ sur le quotient topologique $V\backslash(EV\times X)$, naturelle en $X$. En considérant l'application $(V\rtimes\mathrm{GL}(V))$-équivariante $X\to\mathrm{point}$, on obtient un diagramme commutatif
$$
\begin{CD}
\mathrm{GL}(V)\times(V\backslash(\mathrm{E}V\times X))
@>>> V\backslash(\mathrm{E}V\times X) \\
@VVV @VVV \\
\mathrm{GL}(V)\times\mathrm{B}V @>>>\mathrm{B}V
\end{CD}
$$
dans lequel la flèche horizontale du haut est l'action introduite ci-dessus  et celle du bas l'action de $\mathrm{GL}(V)$ sur $\mathrm{B}V$ induite par l'action tautologique de $\mathrm{GL}(V)$ sur $V$. Le point (a) de la proposition en résulte.
\hfill$\square$

\medskip
\begin{cor}
\label{GL-action2} Soit $X$ un $V$-CW-complexe fini. Si l'action de $V$ sur $X$ se prolonge en une action à gauche (continue) de $V\rtimes\mathrm{GL}(V)$ alors la structure de complexe de cochaînes dans la catégorie $V\text{-}\mathcal{U}$-complexe de $\Sigma^{n}\hspace{1pt}\widetilde{\mathrm{C}}^{\bullet}_{\mathrm{top}}\hspace{1pt}X$ se raffine naturellement en une structure de complexe de cochaînes  dans la catégorie $\mathrm{GL}(V)_{\mathrm{td}}\text{-}V\text{-}\mathcal{U}$.

\smallskip
\em{(La suspension $n$-ième, $n=\dim V$, n'est là que pour garantir que les termes du complexe en question sont ``instables''.)}
\end{cor}

\medskip
\textit{Démonstration.} On constate que les sous-espaces $\mathrm{F}_{p}\hspace{1pt}X$ sont stables sous l'action de $V\rtimes\mathrm{GL}(V)$ (observer que si $W$ est  un sous-groupe de $V$ et $\alpha$ un élément de $\mathrm{GL}(V)$ alors on a $\alpha.X^{W}=X^{\alpha(W)}$). Il en résulte que le connectant
$$
\mathrm{H}^{*}_{V}(\mathrm{F}_{p}\hspace{1pt}X,\mathrm{F}_{p-1}\hspace{1pt}X)
\longrightarrow
\Sigma\hspace{2pt}\mathrm{H}^{*}_{V}(\mathrm{F}_{p+1}\hspace{1pt}X,\mathrm{F}_{p}\hspace{1pt}X)$$
est $\mathrm{GL}(V)$-équivariant et donc que $\Sigma^{n}\hspace{1pt}\widetilde{\mathrm{C}}^{\bullet}_{\mathrm{top}}\hspace{1pt}X$ est un $\mathrm{GL}(V)_{\mathrm{td}}\text{-}V\text{-}\mathcal{U}$-complexe de cochaînes.
\hfill$\square$

\bigskip
Armés de l'énoncé \ref{GL-action2} nous pouvons maintenant vérifier ce que nous avons affirmé plus haut concernant le complexe $\widetilde{\mathrm{C}}^{\bullet}\hspace{1pt}(\mathrm{c}_{V}\hspace{1pt}\mathrm{H}^{*}V)$~:

\medskip
-- L'action tautologique de $V\rtimes\mathrm{GL}(V)$ sur $V$ induit une action à gauche de $V\rtimes\mathrm{GL}(V)$ sur $\mathrm{S}(\mathbb{R}[V])$ qui prolonge celle de $V$ et qui fixe le point $1$~;  $\widetilde{\mathrm{C}}^{\bullet}_{\mathrm{top}}\hspace{1pt}\mathrm{S}(\mathbb{R}[V])$ et $\widetilde{\mathrm{C}}^{\bullet}_{\mathrm{top}}\hspace{1pt}\{1\}
$ sont donc des $\mathrm{GL}(V)_{\mathrm{td}}\text{-}V\text{-}\mathcal{U}$-complexes de cochaînes (on peut oublier ici la suspension $n$-ième compte tenu de \ref{instablebis}).

\medskip
--  L'homomorphisme $\widetilde{\mathrm{C}}^{\bullet}_{\mathrm{top}}\hspace{1pt}\mathrm{S}(\mathbb{R}[V])
\to\widetilde{\mathrm{C}}^{\bullet}_{\mathrm{top}}\hspace{1pt}\{1\}
$ est un morphisme dans la catégorie de ces complexes~; son noyau  $\widetilde{\mathrm{C}}^{\bullet}\hspace{1pt}(\mathrm{c}_{V}\hspace{1pt}\mathrm{H}^{*}V)$ acquiert donc la structure de $\mathrm{GL}(V)_{\mathrm{td}}\text{-}V\text{-}\mathcal{U}$-complexe de cochaînes.

\vspace{0.5cm}
\textit{Définition ``algébrique'' de l'action de $\mathrm{GL}(V)$ sur $\widetilde{\mathrm{C}}^{\bullet}\hspace{1pt}(\mathrm{c}_{V}\hspace{1pt}\mathrm{H}^{*}V)$}

\medskip
On peut munir $\widetilde{\mathrm{C}}^{\bullet}\hspace{1pt}(\mathrm{c}_{V}\hspace{1pt}\mathrm{H}^{*}V)$ d'une structure de $\mathrm{GL}(V)_{\mathrm{td}}\text{-}V\text{-}\mathcal{U}$-complexe de cochaînes de façon purement algébrique~:

\medskip
\begin{pro}\label{GL-action3} Soit $M$ un $\mathrm{H}^{*}V$-$\mathrm{A}$-module instable. Si $M$ est muni d'une action à droite tordue de $\mathrm{GL}(V)$ alors celle-ci peut être naturellement prolongée à $\widetilde{\mathrm{C}}^{\bullet}M$.
\end{pro}

\medskip
Avant d'attaquer la démonstration de cette proposition on introduit une définition et on dégage un énoncé (dont la vérification est immédiate) qui seront utiles à notre exposition.

\medskip
\begin{pro-def}\label {thetaalphaGL-action} Soit $M$ un $\mathrm{H}^{*}V\text{-}\mathrm{A}$-module instable.

\smallskip
Soit $\alpha$ un élément de $\mathrm{GL}(V)$~; on note $\theta_{\alpha}M$ le $\mathrm{A}$-module instable $M$ muni de la structure de $\mathrm{H}^{*}V$-module définie \textit{via} l'isomorphisme $\alpha^{*}:\mathrm{H}^{*}V\to\mathrm{H}^{*}V$ et $\theta_{\alpha}:V\text{-}\mathcal{U}\to V\text{-}\mathcal{U}$ le foncteur $M\mapsto\theta_{\alpha}M$.

\smallskip
Si $M$ est muni d'une action à droite tordue de $\mathrm{GL}(V)$ alors les applications $M\to M, x\mapsto x.\alpha$, $\alpha\in\mathrm{GL}(V)$, sont  des $V\text{-}\mathcal{U}$-morphismes, disons $a_{\alpha}:M\to\theta_{\alpha}M$ vérifiant $a_{\alpha\beta}=(\theta_{\alpha}a_{\beta})\circ a_{\alpha}$, pour tous $\alpha$ et $\beta$ dans $\mathrm{GL}(V)$ et $a_{1_{\mathrm{GL}(V)}}=\mathrm{id}_{M}$.

\smallskip
Réciproquement,  une famille $(a_{\alpha}:M\to\theta_{\alpha}M)_{\alpha\in\mathrm{GL}(V)}$ de $V\text{-}\mathcal{U}$-morphismes, vérifiant les propriétés ci-dessus donne une action à droite tordue de $\mathrm{GL}(V)$ sur $M$ en posant $x.\alpha=a_{\alpha}(x)$.

\end{pro-def}

\medskip
\textit{Démonstration de \ref{GL-action3}.} Elle va faire intervenir les trois propositions \ref{thetaalpharesinj},  \ref{GL-actionresinj}, et  \ref{thetaalphafilt} ci-après.

\medskip
\begin{pro}\label{thetaalpharesinj} Soient $M$ un $\mathrm{H}^{*}V\text{-}\mathrm{A}$-module instable, $M\to I^{\bullet}$ une résolution injective de $M$ dans la catégorie $V\text{-}\mathcal{U}$ et $\alpha$ un élément de $\mathrm{GL}(V)$. Alors $\theta_{\alpha}M\to\theta_{\alpha}I^{\bullet}$ est une résolution injective de $\theta_{\alpha}M$ dans la catégorie~$V\text{-}\mathcal{U}$.
\end{pro}

\textit{Démonstration.} Les foncteurs $\theta_{\alpha}$ sont exacts et préservent les injectifs  (observer que l'on a $\mathrm{Hom}_{V\text{-}\mathcal{U}}(-,\theta_{\alpha}-)=\mathrm{Hom}_{V\text{-}\mathcal{U}}(\theta_{\alpha^{-1}}-,-)$).
\hfill$\square$

\medskip
\begin{pro}\label{GL-actionresinj} Soit $M$ un $\mathrm{H}^{*}V\text{-}\mathrm{A}$-module instable muni d'une action à droite tordue de $\mathrm{GL}(V)$, $(a_{\alpha}:M\to\theta_{\alpha}M)_{\alpha\in\mathrm{GL}(V)}$ la famille de $V\text{-}\mathcal{U}$-morphismes associée à cette action (voir \ref{thetaalphaGL-action}) et $M\to I^{\bullet}$ une résolution injective de $M$ dans la catégorie $V\text{-}\mathcal{U}$.

\smallskip
Alors on dispose d'une famille $(a_{\alpha}^{\bullet}:I^{\bullet}\to\theta_{\alpha}I^{\bullet})_{\alpha\in\mathrm{GL}(V)}$ d'homomorphismes de complexes dans la catégorie $V\text{-}\mathcal{U}$, unique à homotopie près, homomorphismes qui font commuter les diagrammes
$$
\begin{CD}
M@>>>I^{\bullet} \\
@Va_{\alpha}VV @Va_{\alpha}^{\bullet}VV \\
\theta_{\alpha}M@>>>\theta_{\alpha}I^{\bullet}
\end{CD}
$$
et tels que $a_{\alpha\beta}^{\bullet}$ est homotope à $(\theta_{\alpha}a_{\beta}^{\bullet})\circ a_{\alpha}^{\bullet}$ pour tous $\alpha$, $\beta$ dans $\mathrm{GL}(V)$ (et que $a_{1_{\mathrm{GL}(V)}}^{\bullet}$ est l'identité).
\end{pro}

\medskip
\textit{Démonstration.} Mantras de la théorie des résolutions injectives.\hfill$\square$

\begin{pro}\label{thetaalphafilt} Soient $M$ un $\mathrm{H}^{*}V\text{-}\mathrm{A}$-module instable et $\alpha$ un élément de $\mathrm{GL}(V)$. Alors les deux filtrations $(\mathrm{F}^{p}\theta_{\alpha}M)_{p\in\mathbb{N}}$ et $(\theta_{\alpha}\mathrm{F}^{p}M)_{p\in\mathbb{N}}$ de $\theta_{\alpha}M$ coïncident (la filtration $(\mathrm{F}^{p\hspace{1pt}}-)_{p\in\mathbb{N}}$ est introduite au début de la section 5).
\end{pro}

\medskip
La démonstration de \ref{thetaalphafilt} est reportée après la fin de celle de \ref{GL-action3}.

\medskip
\textit{Fin de la démonstration de \ref{GL-action3}.} Soit $M\to I^{\bullet}$ une $V\text{-}\mathcal{U}$-résolution injective dans la~catégorie $V\text{-}\mathcal{U}$. On a défini en section 5 le complexe $\mathrm{C}^{\bullet}\hspace{1pt}M$ en posant $\mathrm{C}^{p}\hspace{1pt}M:=\mathrm{H}^{p}\mathrm{Gr}^{p}I^{\bullet}$ (la notation $\mathrm{Gr}^{p}-$ désigne le quotient $\mathrm{F}^{p}-/\mathrm{F}^{p+1}-$) et en prenant pour cobord le connectant associé à la suite exacte de complexes $0\to\mathrm{Gr}^{p+1}I^{\bullet}\to\mathrm{F}^{p}I^{\bullet}/\mathrm{F}^{p+2}I^{\bullet}\to\mathrm{Gr}^{p}I^{\bullet}\to 0$. Les propositions \ref{thetaalpharesinj} et \ref{thetaalphafilt} conduisent à l'énoncé suivant (que nous archivons en vue de futures références)~:

\begin{scho}\label{thetaalpha0} Soient $M$ un $\mathrm{H}^{*}V\text{-}\mathrm{A}$-module instable et $\alpha$ un élément de $\mathrm{GL}(V)$. On a un isomorphisme naturel de complexes de $\mathrm{H}^{*}V\text{-}\mathrm{A}$-modules instables $\widetilde{\mathrm{C}}^{\bullet}\hspace{1pt}\theta_{\alpha}M\cong\theta_{\alpha}\hspace{1pt}\widetilde{\mathrm{C}}^{\bullet}M$.
\end{scho}

La proposition \ref{GL-actionresinj} implique quant à elle que la famille de $V\text{-}\mathcal{U}$-morphismes
$$
(\hspace{2pt}\mathrm{C}^{p}(a_{\alpha}):\mathrm{C}^{p}\hspace{1pt}M\to
\mathrm{C}^{p}\hspace{1pt}\theta_{\alpha}M=\theta_{\alpha}\mathrm{C}^{p}\hspace{1pt}M\hspace{2pt})_{\alpha\in\mathrm{GL}(V)}
$$
définit une action à droite tordue de $\mathrm{GL}(V)$ sur $\mathrm{C}^{p}\hspace{1pt}M$, uniquement déterminée par l'action à droite tordue de $\mathrm{GL}(V)$ sur $M$, que les cobords $\mathrm{C}^{p}\hspace{1pt}M\to\mathrm{C}^{p+1}\hspace{1pt}M$ sont $\mathrm{GL}(V)$-équivariants et qu'il en est de même pour la coaugmentation $M\to \mathrm{C}^{0}\hspace{1pt}M$. Expliquons par exemple la première implication~: les égalités à homotopie près vérifiées par la famille d'homomorphismes de $V\text{-}\mathcal{U}$-complexes
$$
(\hspace{2pt}\mathrm{Gr}^{p}(a_{\alpha}^{\bullet}):
\mathrm{Gr}^{p}I^{\bullet}\to\mathrm{Gr}^{p}\theta_{\alpha}I^{\bullet}=\theta_{\alpha}\mathrm{Gr}^{p}I^{\bullet}\hspace{2pt})\hspace{1pt}_{\alpha\in\mathrm{GL}(V)}
$$
deviennent des égalités après application du foncteur $\mathrm{H}^{p}$.
\hfill$\square$

\bigskip
\textit{Démonstration de \ref{thetaalphafilt}.} Comme l'on a par définition
$$
\hspace{24pt}
\mathrm{F}^{p}M
\hspace{4pt}:=\hspace{4pt}
\bigcap_{\mathop{\mathrm{codim}}W<p}
\ker\hspace{1pt}({\rho_{(V,W)}}_{M}:M\to\mathrm{EFix}_{(V,W)}M)
\hspace{24pt},
$$
la proposition \ref{thetaalphafilt} est conséquence de la suivante~:

\medskip
\begin{pro}\label{Psitildethetaalpha} Soit $M$ un $\mathrm{H}^{*}V$-$\mathrm{A}$-module instable. Soit $\alpha$ un élément de $\mathrm{GL}(V)$~; on note encore $\alpha$  l'automorphisme de la catégorie $\mathcal{W}$ qu'il induit. On a un isomorphisme dans la catégorie $(V\text{-}\mathcal{U})^{\mathcal{W}}$
$$
\hspace{24pt}
\widehat{\Psi}(\theta_{\alpha}M)
\hspace{4pt}\cong\hspace{4pt}
\theta_{\alpha}\circ\widehat{\Psi}(M)\circ\alpha^{-1}
\hspace{24pt},
$$
naturel en $M$.

\smallskip
{\em (La notation $\widehat{\Psi}(M)$, ou $\widehat{\Psi}_{M}$, désigne le foncteur $W\mapsto\mathrm{EFix}_{(V,W)}M$ de $\mathcal{W}$ dans $V\text{-}\mathcal{U}$.)}
\end{pro}

\medskip
\textit{Démonstration.} On commence  par le lemme suivant~:

\begin{lem}\label{Efixtheta} Soient $\alpha$ un élément de $\mathrm{GL}(V)$ et $W$ un sous-groupe de $V$. On a un isomorphisme fonctoriel entre endofoncteurs de $V\text{-}\mathcal{U}$~:
$$
\hspace{24pt}
\mathrm{EFix}_{(V,W)}\circ\theta_{\alpha}
\hspace{4pt}\cong\hspace{4pt}
\theta_{\alpha}\circ\mathrm{EFix}_{(V,\alpha^{-1}(W))}
\hspace{24pt}.
$$
\end{lem}

\textit{Démonstration.} Les endofoncteurs de $V\text{-}\mathcal{U}$, $\theta_{\alpha}$ et $\mathrm{EFix}_{(V,W)}$, sont respectivement adjoints  à gauche des  endofoncteurs de $V\text{-}\mathcal{U}$, $\theta_{\alpha^{-1}}$ et $\mathrm{E}_{(V,W)}$ (on rappelle que $\mathrm{E}_{(V,W)}$ est l'endofoncteur $M\mapsto\mathrm{H}^{*}V\otimes_{\mathrm{H}^{*}V/W}M$). Il est donc équivalent de montrer que l'on a un isomorphisme fonctoriel entre endofoncteurs de $V\text{-}\mathcal{U}$~:
$$
\hspace{24pt}
\theta_{\alpha^{-1}}\circ\mathrm{E}_{(V,W)}
\hspace{4pt}\cong\hspace{4pt}
\mathrm{E}_{(V,\alpha^{-1}(W))}\circ\theta_{\alpha^{-1}}
\hspace{24pt},
$$
ou encore~:
$$
\hspace{24pt}
\mathrm{E}_{(V,W)}\circ\theta_{\alpha}
\hspace{4pt}\cong\hspace{4pt}
\theta_{\alpha}\circ\mathrm{E}_{(V,\alpha^{-1}(W))}
\hspace{24pt}
$$
En clair il faut vérifier que l'on a pour tout $\mathrm{H}^{*}V$-$\mathrm{A}$-module instable $M$ un isomorphisme de $\mathrm{H}^{*}V$-$\mathrm{A}$-modules instables
$$
\mathrm{H}^{*}V\otimes_{\mathrm{H}^{*}V/W}\hspace{1pt}\theta_{\alpha}M
\hspace{4pt}\cong\hspace{4pt}
\theta_{\alpha\hspace{1pt}}(\mathrm{H}^{*}V\otimes_{\mathrm{H}^{*}V/\alpha^{-1}(W)}M)
$$
naturel en $M$. Il s'agit d'un phénomène très général~: soient $A$ un anneau, $B$ un sous-anneau de $A$, $\phi:A\to A$ un automorphisme d'anneau, $M$ un $A$-module à gauche et  $\theta_{\phi}M$ le $A$-module à gauche ``\hspace{2pt}obtenu en tordant l'action de $A$ sur $M$ par $\phi$\hspace{2pt}'', alors l'application
$$
A\times M\to A\otimes_{\phi(B)}M
\hspace{24pt},\hspace{24pt}
(a,x)\mapsto\phi(a)\otimes_{\phi(B)}x
$$
induit un isomorphisme
$$
A\otimes_{B}\theta_{\phi}M
\hspace{4pt}\cong\hspace{4pt}
\theta_{\phi}(A\otimes_{\phi(B)}M)
$$
de $A$-modules à gauche.
\hfill$\square$

\bigskip
On achève la démonstration de \ref{Psitildethetaalpha} à l'aide de la remarque \ref{Wfonctoriel-2}. Précisons un peu. Soit $\mathrm{i}(V,W): \mathrm{EFix}_{(V,W)}\circ\theta_{\alpha}\to\theta_{\alpha}\circ\mathrm{EFix}_{(V,\alpha^{-1}(W))}$ l'isomorphisme fonctoriel du lemme \ref{Efixtheta}. Soient $W_{0}$ et $W_{1}$ deux sous-groupes de $V$ avec $W_{0}\subset W_{1}$~; on doit vérifier que le diagramme
$$
\begin{CD}
\mathrm{EFix}_{(V,W_{0})}\hspace{1pt}\theta_{\alpha}M
@>\mathrm{i}(V,W_{0})_{M}>>
\theta_{\alpha}\mathrm{EFix}_{(V,\alpha^{-1}(W_{0}))}\hspace{1pt}M \\
@V\rho(W_{0},W_{1})_{\theta_{\alpha}M}VV 
@VV\theta_{\alpha}\rho\hspace{1pt}(\alpha^{-1}(W_{0}),\alpha^{-1}(W_{1}))_{M}V \\
\mathrm{EFix}_{(V,W_{1})}\hspace{1pt}\theta_{\alpha}M
@>\mathrm{i}(V,W_{1})_{M}>>
\theta_{\alpha}\mathrm{EFix}_{(V,\alpha^{-1}(W_{1}))}\hspace{1pt}M
\end{CD}
$$
est commutatif (la notation $\rho(W_{0},W_{1})$ est introduite dans la proposition-définition \ref{Wfonctoriel}). En invoquant \ref{Wfonctoriel-2}, on se ramène, par la même stratégie que précédemment, à vérifier qu'un diagramme ``naturel'' de la forme
$$
\begin{CD}
\theta_{\alpha^{-1}}\mathrm{E}_{(V,W_{0})}N
@<\cong<<
\mathrm{E}_{(V,\alpha^{-1}(W_{0}))}\hspace{1pt}\theta_{\alpha^{-1}}N \\
@AAA @AAA \\
\theta_{\alpha^{-1}}\mathrm{E}_{(V,W_{1})}N
@<\cong<<
\mathrm{E}_{(V,\alpha^{-1}(W_{1}))}\hspace{1pt}\theta_{\alpha^{-1}}N 
\end{CD}
$$
est commutatif, $N$ étant un $\mathrm{H}^{*}V$-$\mathrm{A}$-module instable arbitraire et les isomorphismes horizontaux étant les isomorphismes naturels considérés dans la démonstration du lemme \ref{Efixtheta}. Cette vérification est routine.
\hfill$\square$

\medskip
\begin{cor}\label{Psithetaalpha} Soit $M$ un $\mathrm{H}^{*}V$-$\mathrm{A}$-module instable. Soit $\alpha$ un élément de $\mathrm{GL}(V)$~; on note encore $\alpha$  l'automorphisme de la catégorie $\mathcal{W}_{0}$ qu'il induit. On a un isomorphisme dans la catégorie $(V\text{-}\mathcal{U})^{\mathcal{W}_{0}}$
$$
\hspace{24pt}
\Psi(\theta_{\alpha}M)
\hspace{4pt}\cong\hspace{4pt}
\theta_{\alpha}\circ\Psi(M)\circ\alpha^{-1}
\hspace{24pt},
$$
naturel en $M$.
\end{cor}

\medskip
\textit{Démonstration.} Par définition le foncteur $\Psi(M):\mathcal{W}_{0}\to V\text{-}\mathcal{U}$ est la res\-triction de $\widehat{\Psi}(M)$ à $\mathcal{W}_{0}$.
\hfill$\square$

\bigskip
La proposition \ref{GL-actionresinj} conduit également à l'énoncé suivant~:

\begin{pro}\label{GL-action3.5} Soit $M$ un $\mathrm{H}^{*}V_{\mathrm{tf}}\text{-}\mathrm{A}$-module instable. Si $M$ est muni d'une action à droite tordue de $\mathrm{GL}(V)$ alors sa partie finie $\mathrm{Pf}\hspace{1pt}M$ est stable sous cette action et  les $\mathrm{H}^{*}V_{\mathrm{tf}}\text{-}\mathrm{A}$-modules instables $\mathrm{R}^{k}\mathrm{Pf}\hspace{1pt}M$, $k\geq 1$, sont aussi naturellement munis d'une action à droite tordue de $\mathrm{GL}(V)$. En d'autres termes, les foncteurs $\mathrm{R}^{k}\mathrm{Pf}:V_{\mathrm{tf}}\text{-}\mathcal{U}\to V_{\mathrm{tf}}\text{-}\mathcal{U}$, $k\geq 0$, induisent canoniquement des foncteurs $\mathrm{R}^{k}\mathrm{Pf}:\mathrm{GL}(V)_{\mathrm{td}}\text{-}V_{\mathrm{tf}}\text{-}\mathcal{U}\to\mathrm{GL}(V)_{\mathrm{td}}\text{-}V_{\mathrm{tf}}\text{-}\mathcal{U}$.

\smallskip
{\em (Le lecteur décodera sans peine la notation $\mathrm{GL}(V)_{\mathrm{td}}\text{-}V_{\mathrm{tf}}\text{-}\mathcal{U}$.)}
\end{pro}

\medskip
\textit{Démonstration.}  La partie de l'énoncé concernant $\mathrm{Pf}\hspace{1pt}M$ est triviale~; elle permet d'ailleurs d'identifier $\mathrm{Pf}\hspace{1pt}(\theta_{\alpha}-)$ et $\theta_{\alpha}\mathrm{Pf}\hspace{1pt}(-)$. Passons à la partie concernant les $\mathrm{R}^{k}\mathrm{Pf}\hspace{1pt}M$. Soient $M\to I^{\bullet}$ et $(a_{\alpha}^{\bullet}:I^{\bullet}\to\theta_{\alpha}I^{\bullet})_{\alpha\in\mathrm{GL}(V)}$ la famille d'homomorphismes de complexes dans la catégorie $V_{\mathrm{tf}}\text{-}\mathcal{U}$ fournie par \ref{GL-actionresinj}~; la famille ($\mathrm{H}^{k}\mathrm{Pf}\hspace{1pt}(a_{\alpha}^{\bullet}):\mathrm{R}^{k}\mathrm{Pf}\hspace{1pt}M\to\theta_{\alpha}\mathrm{R}^{k}\mathrm{Pf}\hspace{1pt}M)_{\alpha\in\mathrm{GL}(V)}$ définit une action à droite tordue de $\mathrm{GL}(V)$ sur $\mathrm{R}^{k}\mathrm{Pf}\hspace{1pt}M$, uniquement déterminée par l'action à droite tordue de $\mathrm{GL}(V)$ sur $M$.
\hfill$\square$

\begin{scho}\label{thetaalpha0.5} Soient $M$ un $\mathrm{H}^{*}V_{\mathrm{tf}}\text{-}\mathrm{A}$-module instable et $\alpha$ un élément de $\mathrm{GL}(V)$. Pour tout entier $k\geq 0$, on a un isomorphisme naturel de $\mathrm{H}^{*}V_{\mathrm{tf}}\text{-}\mathrm{A}$-modules instables $\mathrm{R}^{k}\mathrm{Pf}\hspace{1pt}\theta_{\alpha}M\cong\theta_{\alpha}\mathrm{R}^{k}\mathrm{Pf}\hspace{1pt}M$.
\end{scho}

\vspace{0.5cm}
\textit{Coïncidence des actions de $\mathrm{GL}(V)$ sur $\widetilde{\mathrm{C}}^{\bullet}\hspace{1pt}(\mathrm{c}_{V}\hspace{1pt}\mathrm{H}^{*}V)$ définies topologiquement et algébriquement}

\medskip
Sans surprise les deux structures que l'on vient de définir sur $\widetilde{\mathrm{C}}^{\bullet}\hspace{1pt}(\mathrm{c}_{V}\hspace{1pt}\mathrm{H}^{*}V)$, respectivement de façon topologique et algébrique, coïncident. Ceci résulte de la proposition \ref{CalgCtop} et la proposition \ref{GL-action4} ci-dessous. Avant d'énoncer cette dernière, deux observations. Soit $X$ un $V$-CW-complexe fini muni d'une action à gauche  (continue) de $V\rtimes\mathrm{GL}(V)$ qui prolonge celle de $V$~:

\smallskip
-- Le $\mathrm{H}^{*}V$-$\mathrm{A}$-module instable $\mathrm{H}^{*}_{V}X$ est naturellement muni d'une action à droite tordue de $\mathrm{GL}(V)$ d'après le point (a) de \ref{GL-action1} si bien que le complexe $\widetilde{\mathrm{C}}^{\bullet}\hspace{1pt}\mathrm{H}^{*}_{V}X$ est naturellement muni d'une action à droite de $\mathrm{GL}(V)$ d'après \ref{GL-action3}.

\smallskip
-- Le complexe $\widetilde{\mathrm{C}}^{\bullet}_{\mathrm{top}}\hspace{1pt}X$ est naturellement muni d'une action à droite de $\mathrm{GL}(V)$ d'après \ref{GL-action2}.

\medskip
\begin{pro}\label{GL-action4} Soit $X$ un $V$-CW-complexe fini. Si l'action de $V$ sur $X$ se prolonge en une action à gauche (continue) de $V\rtimes\mathrm{GL}(V)$ alors l'homomorphisme de complexes
$$
\varkappa:
\hspace{4pt}\widetilde{\mathrm{C}}^{\bullet}\hspace{1pt}\mathrm{H}^{*}_{V}X
\longrightarrow
\widetilde{\mathrm{C}}_{\mathrm{top}}^{\bullet}\hspace{1pt}X
$$
introduit en section 6 est $\mathrm{GL}(V)$-équivariant.

\medskip
\em(On rappelle que l'on a posé $\widetilde{\mathrm{C}}^{\bullet}\hspace{1pt}\mathrm{H}^{*}_{V}X=:\widetilde{\mathrm{C}}_{\mathrm{alg}}^{\bullet}\hspace{1pt}X$ en section 6 pour faire pendant à la notation $\widetilde{\mathrm{C}}_{\mathrm{top}}^{\bullet}\hspace{1pt}X$.)
\end{pro}

\medskip
\textit{Démonstration.} Elle repose essentiellement sur la proposition \ref{kappanaturel}. On trouvera les détails ci-après (détails techniques, d'où les petits caractères).

\medskip
 \footnotesize
 On commence par dégager cinq énoncés, \ref{thetaalpha1}, \ref{thetaalpha2}, \ref{thetaalpha3},  \ref{thetaalpha4} et \ref{thetaalpha5}, dont la vérification ne présente pas de difficultés.

\begin{lem-def}\label{thetaalpha1} Soient $X$ un $V$-espace et $\alpha$ un élément de $\mathrm{GL}(V)$. L'espace $X$ muni de l'action de $V$ définie \textit{via} l'isomorphisme $\alpha:V\to V$ est un $V$-espace noté $\theta_{\alpha}X$. L'homéomor\-phisme $\mathrm{E}\alpha\times\mathrm{id}:\mathrm{E}V\times X\to\mathrm{E}V\times X$ induit un homéomorphisme
$$
\hspace{24pt}
V\backslash(\mathrm{E}V\times\theta_{\alpha}X))
\longrightarrow
V\backslash(\mathrm{E}V\times X)
\hspace{24pt},
$$
noté $\mathrm{h}_{\alpha}$, qui fait commuter le diagramme
$$
\begin{CD}
V\backslash(\mathrm{E}V\times\theta_{\alpha}X)@>\mathrm{h}_{\alpha}>>
V\backslash(\mathrm{E}V\times X) \\
@VVV @VVV \\
\mathrm{B}V@>\mathrm{B}\alpha>>
\hspace{24pt}\mathrm{B}V\hspace{23pt}.
\end{CD}
$$
\end{lem-def}

\begin{scho}\label{thetaalpha2} Soient $X$ un $V$-espace et $\alpha$ un élément de $\mathrm{GL}(V)$. L'homomorphisme
$$
\mathrm{h}_{\alpha}^{*}:\mathrm{H}^{*}_{V}X\longrightarrow
\theta_{\alpha}\hspace{1pt}\mathrm{H}^{*}_{V}\hspace{1pt}\theta_{\alpha}X
$$
est un isomorphisme de $\mathrm{H}^{*}V$-$\mathrm{A}$-modules instables.
\end{scho}

\medskip
Soient $X$ un $V$-CW-complexe fini et $\alpha$ un élément de $\mathrm{GL}(V)$~; il est clair que $\theta_{\alpha}X$ est encore un $V$-CW-complexe fini. En observant que les énoncés \ref{thetaalpha1} et \ref{thetaalpha2} s'étendent \textit{mutatis mutandis} aux paires de $V$-espaces, on obtient~:

\pagebreak

\begin{scho-def}\label{thetaalpha3} Soient $X$ un $V$-CW-complexe fini et $\alpha$ un élément de $\mathrm{GL}(V)$. L'isomorphisme $\mathrm{h}_{\alpha}^{*}:\mathrm{H}^{*}_{V}X\to
\theta_{\alpha}\hspace{1pt}\mathrm{H}^{*}_{V}\hspace{1pt}\theta_{\alpha}X$ se prolonge en un isomorphisme, toujours noté $\mathrm{h}_{\alpha}^{*}:\widetilde{\mathrm{C}}^{\bullet}_{\mathrm{top}}\hspace{1pt}X\to\theta_{\alpha}\hspace{1pt}\widetilde{\mathrm{C}}^{\bullet}_{\mathrm{top}}\hspace{1pt}\theta_{\alpha}X$, de complexes de $\mathrm{H}^{*}V$-$\mathrm{A}$-modules.
\end{scho-def}

\medskip
On note encore $\mathrm{h}_{\alpha}^{*}:\widetilde{\mathrm{C}}^{\bullet}\hspace{1pt}\mathrm{H}^{*}_{V}X\to\theta_{\alpha}\hspace{1pt}\widetilde{\mathrm{C}}^{\bullet}\hspace{1pt}\mathrm{H}^{*}_{V}\hspace{1pt}\theta_{\alpha}X$ l'isomorphisme naturel de complexes de $\mathrm{H}^{*}V$-$\mathrm{A}$-modules donné par \ref{thetaalpha0} et \ref{thetaalpha3}.

\begin{pro}\label{thetaalpha4} Soient $X$ un $V$-CW-complexe fini et $\alpha$ un élément de $\mathrm{GL}(V)$. Le diagramme de complexes de $\mathrm{H}^{*}V$-$\mathrm{A}$-modules
$$
\begin{CD}
\widetilde{\mathrm{C}}^{\bullet}\hspace{1pt}\mathrm{H}^{*}_{V}X
@>\varkappa_{X}>>
\widetilde{\mathrm{C}}_{\mathrm{top}}^{\bullet}\hspace{1pt}X \\
@V\mathrm{h}_{\alpha}^{*}V\cong V @V\mathrm{h}_{\alpha}^{*}V\cong V \\
\theta_{\alpha}\hspace{1pt}\widetilde{\mathrm{C}}^{\bullet}\hspace{1pt}\mathrm{H}^{*}_{V}\hspace{1pt}\theta_{\alpha}X
@>\theta_{\alpha}\hspace{1pt}\varkappa_{\hspace{1pt}\theta_{\alpha}X}>>
\theta_{\alpha}\hspace{1pt}\widetilde{\mathrm{C}}_{\mathrm{top}}^{\bullet}\hspace{1pt}\theta_{\alpha}X
\end{CD}
$$
est commutatif (par souci de clarté on  précise la notation $\varkappa:\widetilde{\mathrm{C}}^{\bullet}\hspace{1pt}\mathrm{H}^{*}_{V}X\to
\widetilde{\mathrm{C}}_{\mathrm{top}}^{\bullet}\hspace{1pt}X$ en $\varkappa_{X}:\widetilde{\mathrm{C}}^{\bullet}\hspace{1pt}\mathrm{H}^{*}_{V}X
\to
\widetilde{\mathrm{C}}_{\mathrm{top}}^{\bullet}\hspace{1pt}X$).
\end{pro}

\medskip
\begin{pro-def}\label{thetaalpha5} Soit $X$ un $V$-espace muni d'une action à gauche (continue) de $V\rtimes\mathrm{GL}(V)$ qui prolonge celle de $V$. Soit $\alpha$ un élément de $\mathrm{GL}(V)$~; on note $\mathrm{t}_{\alpha}:X\to\theta_{\alpha}\hspace{1pt}X$ l'homéomorphisme $V$-équivariant $x\mapsto\alpha.x$. La famille d'homomorphismes de $\mathrm{H}^{*}V$-$\mathrm{A}$-modules instables $(\hspace{1pt}\theta_{\alpha}\mathrm{t}_{\alpha}^{*}\circ\mathrm{h}_{\alpha}^{*}:\mathrm{H}^{*}_{V}X\to\theta_{\alpha}\mathrm{H}^{*}_{V}X\hspace{1pt})_{\alpha\in\mathrm{GL}(V)}$ coïncide avec celle qui donne (voir \ref{thetaalphaGL-action}) l'action à droite tordue de $\mathrm{GL}(V)$ sur $\mathrm{H}^{*}_{V}X$ (introduite dans le point~(a) de \ref{GL-action1}).
\end{pro-def}

\medskip
On en vient maintenant à la démonstration proprement dite de la proposition.

\medskip
Soit $V$ un $V$-CW-complexe fini muni d'une action à gauche (continue) de $V\rtimes\mathrm{GL}(V)$ qui prolonge celle de $V$. D'après la proposition \ref{kappanaturel} le diagramme de complexes de $\mathrm{H}^{*}V$-$\mathrm{A}$-modules
$$
\begin{CD}
\widetilde{\mathrm{C}}^{\bullet}\hspace{1pt}\mathrm{H}^{*}_{V}\hspace{1pt}\theta_{\alpha}X
@>\varkappa_{\hspace{1pt}\theta_{\alpha}X}>>
\widetilde{\mathrm{C}}_{\mathrm{top}}^{\bullet}\hspace{1pt}\theta_{\alpha}X \\
@V\widetilde{\mathrm{C}}^{\bullet}(\mathrm{t}_{\alpha}^{*})V\cong V
@V\widetilde{\mathrm{C}}_{\mathrm{top}}^{\bullet}(\mathrm{t}_{\alpha})V\cong V \\
\widetilde{\mathrm{C}}^{\bullet}\hspace{1pt}\mathrm{H}^{*}_{V}X
@>\varkappa_{X}>>
\widetilde{\mathrm{C}}_{\mathrm{top}}^{\bullet}\hspace{1pt}X
\end{CD}
\leqno{(\mathcal{N}_{\alpha})}
$$
est commutatif. Soit $(\theta_{\alpha}\hspace{1pt}\mathcal{N}_{\alpha})$ le diagramme obtenu en appliquant le foncteur $\theta_{\alpha}$ à $(\mathcal{N}_{\alpha})$~; en ``concaténant'' le diagramme commutatif $\theta_{\alpha}\hspace{1pt}(\mathcal{N}_{\alpha})$ et celui de \ref{thetaalpha4}, on obtient au bout du compte un diagramme commutatif de la forme suivante~:
$$
\begin{CD}
\widetilde{\mathrm{C}}^{\bullet}\hspace{1pt}\mathrm{H}^{*}_{V}X
@>\varkappa_{X}>>
\widetilde{\mathrm{C}}_{\mathrm{top}}^{\bullet}\hspace{1pt}X \\
@V\mathrm{a}_{\alpha}^{\mathrm{g}}V\cong V @V\mathrm{a}_{\alpha}^{\mathrm{d}}V\cong V \\
\theta_{\alpha}\hspace{1pt}\widetilde{\mathrm{C}}^{\bullet}\hspace{1pt}\mathrm{H}^{*}_{V}X
@>\theta_{\alpha}\hspace{1pt}\varkappa_{X}>>
\theta_{\alpha}\hspace{1pt}\widetilde{\mathrm{C}}_{\mathrm{top}}^{\bullet}\hspace{1pt}X
\end{CD}
$$
(les exposants $\mathrm{g}$ et $\mathrm{d}$ pour gauche et droite). On vérifie enfin, à l'aide de l'énoncé \ref{thetaalpha5} (et de son extension aux paires de $(V\rtimes\mathrm{GL}(V))$-espaces) et de l'énoncé \ref{thetaalpha0}, que les familles d'homomorphismes $(\hspace{1pt}\mathrm{a}_{\alpha}^{\mathrm{g}}\hspace{1pt})_{\alpha\in\mathrm{GL}(V)}$ et $(\hspace{1pt}\mathrm{a}_{\alpha}^{\mathrm{d}}\hspace{1pt})_{\alpha\in\mathrm{GL}(V)}$ correspondent (voir \ref{thetaalphaGL-action}) aux actions à droite tordues respectivement introduites en \ref{GL-action3} et \ref{GL-action2}.
\normalsize
\hfill$\square$

\vspace{0.5cm}
\textit{Quelques informations sur  le  $\mathrm{GL}(V)_{\mathrm{td}}\text{-}\mathrm{H}^{*}V\text{-}\mathrm{A}$-module instable $\mathrm{M}(V)$}
 
\bigskip
La théorie que nous venons de développer concernant $\widetilde{\mathrm{C}}^{\bullet}\hspace{1pt}(\mathrm{c}_{V}\hspace{1pt}\mathrm{H}^{*}V)$, où plus directement la proposition \ref{GL-action3.5}, montre que le $\mathrm{H}^{*}V\text{-}\mathrm{A}$-module instable $\mathrm{M}(V)$ est canoniquement muni d'une action à droite tordue de $\mathrm{GL}(V)$~; le sous-espace $\mathrm{M}^{0}(V)$ des éléments de degré zéro est donc en particulier muni d'une action à droite de $\mathrm{GL}(V)$.

\begin{pro}\label{MV} Le $\mathrm{GL}(V)_{\mathrm{td}}\text{-}\mathrm{H}^{*}V\text{-}\mathrm{A}$-module instable $\mathrm{M}(V)$ possède les propriétés suivantes~:

\medskip
{\em (a)}  $\mathrm{M}(V)$ est engendré comme $\mathrm{H}^{*}V$-module par $\mathrm{M}^{0}(V)$~;

\medskip
{\em (b)}  comme $\mathbb{F}_{2}[\mathrm{GL}(V)]$-module à droite $\mathrm{M}^{0}(V)$ est isomorphe au dual $\mathrm{St}_{V}^{*}$ de la représentation de Steinberg (modulo $2$) $\mathrm{St}_{V}$ de $\mathrm{GL}(V)$ (on convient de poser $\mathrm{St}_{V}=\mathbb{F}_{2}$ pour $\dim V=0$ et $\dim V=1$).
\end{pro}

\bigskip
\textit{Démonstration du (a).} On procède par récurrence sur l'entier $n=\dim V$. Le cas $n=0$ est trivial~; on franchit le pas de récurrence en observant que le cobord
$$
\begin{CD}
\bigoplus_{\dim W=1}\mathrm{H}^{*}V\otimes_{\mathrm{H}^{*}V/W}\mathrm{M}(V/W)
@>\mathrm{d}^{n-1}>>
\mathrm{M}(V)
\end{CD}
$$
de  $\widetilde{\mathrm{C}}^{\bullet}\hspace{1pt}(\mathrm{c}_{V}\hspace{1pt}\mathrm{H}^{*}V)$ est surjectif.
\hfill$\square$

\bigskip
\textit{Démonstration du (b).} On suppose $n\geq 2$. On va utiliser l'égalité $\mathrm{M}(V):=\mathrm{R}^{n}\mathrm{Pf}(\mathrm{c}_{V}\hspace{1pt}\mathrm{H}^{*}V)$ et l'isomorphisme $\mathrm{R}^{n}\mathrm{Pf}(\mathrm{c}_{V}\hspace{1pt}\mathrm{H}^{*}V)\cong\lim_{\mathcal{W}_{0}}^{n-1}\Psi(\mathrm{c}_{V}\hspace{1pt}\mathrm{H}^{*}V)$ fourni par la proposition \ref{derivePf}.

\medskip
On détermine ci-après le foncteur $\Psi(\mathrm{c}_{V}\hspace{1pt}\mathrm{H}^{*}V)$ et plus généralement le foncteur $\Psi(e\hspace{1pt}\mathrm{H}^{*}V)$, $e$ désignant un produit d'éléments de $\mathrm{H}^{1}V-\{0\}$. 

\medskip
\begin{pro-def}\label{spePsi} Soit $e$ un produit d'éléments de $\mathrm{H}^{1}V-\{0\}$~; soit $W$ un sous-groupe de $V$.

\medskip
{\em (a)} La classe $e$ s'écrit de façon unique comme un produit $e'e''$ avec $e'$ (resp. $e''$) un produit d'éléments de $\mathrm{H}^{1}V-\{0\}$ dont la restriction à $W$ est nulle (resp. non nulle). La classe $e'$ est notée $e_{W}$ (on a en particulier $e_{0}=e$ et $e_{V}=1$).

\medskip
{\em (b)} Si $W'$ est un sous-groupe de $V$ avec $W'\subset W$ alors $e_{W'}$ est un multiple de $e_{W}$ et l'on a une inclusion canonique $e_{W'}\hspace{1pt}\mathrm{H}^{*}V\subset e_{W}\hspace{1pt}\mathrm{H}^{*}V$.

\medskip
{\em (c)} Le foncteur $\Psi(e\hspace{1pt}\mathrm{H}^{*}V)$ est le foncteur
$$
\hspace{24pt}
\mathcal{W}_{0}\to V_{\mathrm{tf}}\text{-}\mathcal{U}
\hspace{12pt},\hspace{12pt}
W\mapsto
e_{W}\hspace{1pt}\mathrm{H}^{*}V
\hspace{24pt}.
$$
\end{pro-def}

\medskip
\textit{Démonstration.} Les points (a) et (b) sont triviaux~; on démontre le point~(c). On écrit $e=u_{1}u_{2}\ldots u_{r}$, avec $(u_{1}u_{2},\ldots,u_{r})$  une suite finie d'éléments de $\mathrm{H}^{1}V-\{0\}$~; il est clair que l'on a un isomorphisme canonique de $\mathrm{H}^{*}V\text{-}\mathrm{A}$-modules instables~:
$$
\hspace{24pt}
e\hspace{1pt}\mathrm{H}^{*}V
\hspace{4pt}\cong\hspace{4pt}
u_{1}\hspace{1pt}\mathrm{H}^{*}V\otimes_{\mathrm{H}^{*}V}
u_{2}\hspace{1pt}\mathrm{H}^{*}V\otimes_{\mathrm{H}^{*}V}\ldots \otimes_{\mathrm{H}^{*}V}u_{r}\hspace{1pt}\mathrm{H}^{*}V
\hspace{24pt}.
$$
La proposition \ref{tensEFix} montre que l'on a~:
\begin{multline*}
\hspace{24pt}
\mathrm{EFix}_{(V,W)}(e\hspace{1pt}\mathrm{H}^{*}V)
\hspace{4pt}\cong\hspace{4pt}
\mathrm{EFix}_{(V,W)}(u_{1}\hspace{1pt}\mathrm{H}^{*}V)\otimes_{\mathrm{H}^{*}V}
\mathrm{EFix}_{(V,W)}(u_{2}\hspace{1pt}\mathrm{H}^{*}V)\ldots \\\ldots\otimes_{\mathrm{H}^{*}V}\mathrm{EFix}_{(V,W)}(u_{r}\hspace{1pt}\mathrm{H}^{*}V)
\hspace{24pt}.
\end{multline*}
On est donc amené à déterminer $\mathrm{EFix}_{(V,W)}(u\hspace{1pt}\mathrm{H}^{*}V)$ pour  $u$ un élément de $\mathrm{H}^{1}V-\{0\}$. Pour cela on considère la suite exacte de $\mathrm{H}^{*}V\text{-}\mathrm{A}$-modules instables
$$
\begin{CD}
0@>>>u\hspace{1pt}\mathrm{H}^{*}V@>\subset>>\mathrm{H}^{*}V
@>>>\mathrm{H}^{*}\ker u@>>>0
\end{CD}
$$
(on rappelle que l'on identifie $\mathrm{H}^{1}V$ et $V^{*}$). Puisque le foncteur $\mathrm{EFix}_{(V,W)}$  est exact on a aussi une suite exacte
de $\mathrm{H}^{*}V\text{-}\mathrm{A}$-modules instables~:
$$
\hspace{12pt}
0\to\mathrm{EFix}_{(V,W)}(u\hspace{1pt}\mathrm{H}^{*}V)\to
\mathrm{EFix}_{(V,W)}(\mathrm{H}^{*}V)\to
\mathrm{EFix}_{(V,W)}(\mathrm{H}^{*}\ker u)\to 0
\hspace{12pt}.
$$
En observant que l'on a $\mathrm{H}^{*}V\cong\mathrm{H}^{*}V\otimes_{\mathrm{H}^{*}V/V}\mathbb{F}_{2}$ et $\mathrm{H}^{*}\ker u\cong\mathrm{H}^{*}V\otimes_{\mathrm{H}^{*}V/\ker u}\mathbb{F}_{2}$ et en utilisant la proposition \ref{proclef} on obtient des isomorphismes canoniques de $\mathrm{H}^{*}V\text{-}\mathrm{A}$-modules instables~:
$$
\mathrm{EFix}_{(V,W)}(\mathrm{H}^{*}V)\cong\mathrm{H}^{*}V
\hspace{6pt}\text{et}\hspace{6pt}
\mathrm{EFix}_{(V,W)}(\mathrm{H}^{*}\ker u)\cong
\begin{cases} \mathrm{H}^{*}\ker u & \text{pour $u_{\vert W}=0$,} \\
0 & \text{pour $u_{\vert W}\not=0.$}
\end{cases}
$$
Dans le cas $u_{\vert W}=0$ la proposition \ref{proclef} montre également que l'homomorphisme $\mathrm{EFix}_{(V,W)}(\mathrm{H}^{*}V)\to\mathrm{EFix}_{(V,W)}(\mathrm{H}^{*}\ker u)$ s'identifie à l'épimorphisme canonique $\mathrm{H}^{*}V\to\mathrm{H}^{*}\ker u$.
 
\smallskip
Il en résulte que le monomorphisme $\mathrm{EFix}_{(V,W)}(u\hspace{1pt}\mathrm{H}^{*}V)\to\mathrm{EFix}_{(V,W)}(\mathrm{H}^{*}V)$ s'identifie à l'inclusion $u\hspace{1pt}\mathrm{H}^{*}V\hookrightarrow\mathrm{H}^{*}V$ dans le cas $u_{\vert W}=0$ et à l'identité de $\mathrm{H}^{*}V$ dans le cas $u_{\vert W}\not=0$. Le monomorphisme $\mathrm{EFix}_{(V,W)}(e\hspace{1pt}\mathrm{H}^{*}V)\to\mathrm{EFix}_{(V,W)}(\mathrm{H}^{*}V)$ s'identifie donc à l'inclusion $e_{W}\hspace{1pt}\mathrm{H}^{*}V\hookrightarrow\mathrm{H}^{*}V$.

\smallskip
On détermine la valeur de $\Psi(e\hspace{1pt}\mathrm{H}^{*}V)$ sur les morphismes de $\mathcal{W}_{0}$ en observant que $\Psi(\mathrm{H}^{*}V)$ est le ``foncteur constant'' $W\mapsto\mathrm{H}^{*}V$ et en considérant la transformation naturelle de foncteurs de $\mathcal{W}_{0}$ dans $V_{\mathrm{tf}}\text{-}\mathcal{U}$ induite par l'inclusion $e\hspace{1pt}\mathrm{H}^{*}V\hookrightarrow\mathrm{H}^{*}V$.
\hfill$\square$

\medskip
\begin{exple}\label{spespePsi} Soit $W$ un sous-groupe de $V$, on pose $\mathrm{c}_{(V,W)}:=q^{*}\mathrm{c}_{V/W}$, $q$~désignant la surjection canonique $V\to V/W$~; on a donc
$$
\mathrm{c}_{(V,W)}
\hspace{4pt}:=\hspace{4pt}
\prod_{u\in V^{*}\hspace{2pt},\hspace{2pt}u_{\vert W}=0\hspace{2pt}, \hspace{2pt}u\not=0} u
$$
(en particulier $\mathrm{c}_{(V,0)}=\mathrm{c}_{V}$ et $\mathrm{c}_{(V,V)}=1$). Soit $h\geq 0$ un entier, la proposition précédente dit que l'on a $\Psi(\mathrm{c}_{V}^{h}\hspace{1pt}\mathrm{H}^{*}V)(W)=\mathrm{c}_{(V,W)}^{h}\hspace{1pt}\mathrm{H}^{*}V$ et que l'image par $\Psi(\mathrm{c}_{V}^{h}\hspace{1pt}\mathrm{H}^{*}V)$ d'une inclusion $W'\subset W$ (vue comme morphisme de $\mathcal{W}_{0}$) est l'inclusion $\mathrm{c}_{(V,W')}^{h}\hspace{1pt}\hspace{1pt}\mathrm{H}^{*}V\subset\mathrm{c}_{(V,W)}^{h}\hspace{1pt}\mathrm{H}^{*}V$.
\end{exple}

\medskip
\textit{Suite de la démonstration du (b) de \ref{MV}.} Soient $\mathcal{E}$ la catégorie des $\mathbb{F}_{2}$-espaces vectoriels, et $M$ un  $\mathrm{H}^{*}V_{\mathrm{tf}}$-$\mathrm{A}$-module instable~; on note $\Psi^{0}(M):\mathcal{W}_{0}\to\mathcal{E}$ le~foncteur $W\mapsto(\Psi(M)(W))^{0}$ (le $\mathbb{F}_{2}$-espace vectoriel constitué des éléments de degré $0$ du  $\mathrm{H}^{*}V$-$\mathrm{A}$-module instable $\Psi(M)(W)$). Comme l'on a pour tout $k\geq 0$ l'égalité $\lim_{\mathcal{W}_{0}}^{k}\Psi^{0}(M)=(\lim_{\mathcal{W}_{0}}^{k}\Psi(M))^{0}$ (invoquer par exemple le rappel (R.2) de \ref{rappelslim}), la proposition \ref{derivePf} implique $(\mathrm{R}^{n}\mathrm{Pf}\hspace{1pt}M)^{0}\cong\lim_{\mathcal{W}_{0}}^{n-1}\Psi^{0}(M)$. On a donc en particulier $\mathrm{M}^{0}(V)\cong\lim_{\mathcal{W}_{0}}^{n-1}\Psi^{0}(\mathrm{c}_{V}\hspace{1pt}\mathrm{H}^{*}V)$. On constate que l'on~a :
$$
\Psi^{0}(\mathrm{c}_{V}\hspace{1pt}\mathrm{H}^{*}V)(W)
\hspace{4pt}=\hspace{4pt}
\begin{cases}
\mathbb{F}_{2} & \text{pour $W=V$}, \\
0 & \text{pour $W\not=V$}~;
\end{cases}
$$
en effet la classe $\mathrm{c}_{(V,W)}$ (notation de \ref{spespePsi}) est de degré strictement positif pour $W\not=V$ et est égale à $1$ pour $W=V$.

\smallskip
Récapitulons~: on a $\mathrm{M}^{0}(V)\cong\lim_{\mathcal{W}_{0}}^{n-1}\mathrm{F}$, $\mathrm{F}$ désignant le foncteur de $\mathcal{W}_{0}$ dans~$\mathcal{E}$ défini par
$$
\mathrm{F}(W)
\hspace{4pt}=\hspace{4pt}
\begin{cases}
\mathbb{F}_{2} & \text{pour $W=V$,} \\
0 & \text{pour $W\not=V$,}
\end{cases}
$$
(la valeur de $\mathrm{F}$ sur les objets force la valeur de $\mathrm{F}$ sur les morphismes).

\smallskip
Soit $\mathrm{E}:\mathcal{W}_{0}\to\mathcal{E}$ le foncteur ``constant à valeur $\mathbb{F}_{2}$'' ($\mathrm{E}:=\Delta_{\mathbb{F}_{2}}$ avec la notation du (R.1) de \ref{rappelslim}). On observe que l'on a $\mathrm{E}=\mathbb{F}_{2}^{\mathrm{Hom}_{\mathcal{W}_{0}}(-,V)}$, égalité qui entraîne $\mathrm{Hom}_{\mathcal{E}^{\mathcal{W}_{0}}}(G,\mathrm{E})=\mathrm{Hom}_{\mathcal{E}}(G(V),\mathbb{F}_{2})$ pour tout objet $G$ de $\mathcal{E}^{\mathcal{W}_{0}}$ et donc que $\mathrm{E}$ est un injectif de $\mathcal{E}^{\mathcal{W}_{0}}$. Soit $\mathrm{i}:\mathrm{F}\to\mathrm{E}$ l'homomorphisme correspondant à l'isomorphisme $\mathrm{F}(V)\cong\mathbb{F}_{2}$~; $\mathrm{i}$ est injectif $((\mathrm{E},\mathrm{i})$ est en fait une enveloppe injective de $\mathrm{F}$ dans la catégorie $\mathcal{E}^{\mathcal{W}_{0}}$). On note $\mathrm{Q}$ le conoyau de $\mathrm{i}$~; on a ~:

\smallskip
--\hspace{8pt}$\lim_{\mathcal{W}_{0}}^{n-1}\mathrm{F}\cong\lim_{\mathcal{W}_{0}}^{n-2}\mathrm{Q}$ pour $n>2$,

 \smallskip
--\hspace{8pt}$\lim_{\mathcal{W}_{0}}^{n-1}\mathrm{F}\cong\mathrm{coker}\hspace{1pt}(\hspace{1pt}\lim_{\mathcal{W}_{0}}^{n-2}\mathrm{E}\longrightarrow\lim_{\mathcal{W}_{0}}^{n-2}\mathrm{Q}\hspace{1pt})$ pour $n=2$.

\smallskip
On a
$$
\mathrm{Q}(W)\hspace{4pt}=\hspace{4pt}
\begin{cases}
\mathbb{F}_{2} & \text{pour $W\not=V$} \\
0 & \text{pour $W=V$}
\end{cases}
$$
et $\mathrm{Q}(W'\subset W)=\mathrm{id}_{\mathbb{F}_{2}}$ pour $W\not=V$. Soit $\mathcal{W}_{0,n}$ le sous-ensemble ordonné de~$\mathcal{W}$ constitué des $W$ avec $W\not=0$ et $W\not=V$~; soit $\Vert\mathcal{W}_{0,n}\Vert$ le polyèdre associé. On constate que le complexe $\mathrm{L}^{\bullet}(\mathrm{Q})$  (voir le rappel (R.2) de \ref{rappelslim}) s'identifie au complexe de cochaînes polyédrales modulo~$2$ de $\Vert\mathcal{W}_{0,n}\Vert$ si bien que l'on a
$$
\hspace{24pt}
{\lim}_{\mathcal{W}_{0}}^{n-2}\hspace{2pt}\mathrm{Q}
\hspace{4pt}\cong\hspace{4pt}
\mathrm{H}^{n-2}(\Vert\mathcal{W}_{0,n}\Vert;\mathbb{F}_{2})
\hspace{24pt}.
$$
Pareillement, le complexe $\mathrm{L}^{\bullet}(\mathrm{E})$ s'identifie au complexe de cochaînes poly\-édrales modulo $2$ du cône sur $\Vert\mathcal{W}_{0,n}\Vert$ et l'homomorphisme $\mathrm{L}^{\bullet}(\mathrm{E})\to\mathrm{L}^{\bullet}(\mathrm{Q})$ à~l'homomorphisme induit par l'inclusion de $\Vert\mathcal{W}_{0,n}\Vert$ dans ce cône si bien que l'on a pour tout $n\geq 2$~:
$$
{\lim}_{\mathcal{W}_{0}}^{n-1}\hspace{2pt}\mathrm{F}
\hspace{4pt}\cong\hspace{4pt}
\widetilde{\mathrm{H}}^{n-2}(\Vert\mathcal{W}_{0,n}\Vert;\mathbb{F}_{2})
$$
(voir \ref{Bob}). Comme l'on peut définir $\mathrm{St}_{V}$ par  $\mathrm{St}_{V}:=\widetilde{\mathrm{H}}_{n-2}(\Vert\mathcal{W}_{0,n}\Vert;\mathbb{F}_{2})$ on obtient bien au bout du compte un isomorphisme canonique
$$
\hspace{24pt}
\mathrm{M}^{0}(V)
\hspace{4pt}\cong\hspace{4pt}
{(\mathrm{St}_{V})}^{*}
\hspace{24pt}.
$$

\medskip
Il reste à se convaincre qu'il s'agit d'un isomorphisme de $\mathbb{F}_{2}[\mathrm{GL}(V)]$-modules à droite. Pour cela on réexamine la proposition \ref{derivePf} dans le contexte des $\mathrm{H}^{*}V_{\mathrm{tf}}$-$\mathrm{A}$-modules instables munis d'une action à droite tordue de $\mathrm{GL}(V)$.

\medskip
Soient $M$ un $\mathrm{H}^{*}V_{\mathrm{tf}}$-$\mathrm{A}$-module instable et $M\to I^{\bullet}$ une résolution injective dans la  catégorie $V_{\mathrm{tf}}\text{-}\mathcal{U}$. Les connectants associés à la suite exacte de complexes (voir le point (c) de \ref{proPsi})
$$
\begin{CD}
0@>>>\mathrm{Pf}(I^{\bullet})@>\subset>>
I^{\bullet}@>\rho_{I^{\bullet}}>>
\lim_{\mathcal{W}_{0}}\Psi(I^{\bullet})@>>>0
\end{CD}
$$
fournissent  des $V_{\mathrm{tf}}\text{-}\mathcal{U}$-morphismes naturels  $\partial:\lim_{\mathcal{W}_{0}}^{k}\Psi(M)\to\mathrm{R}^{k+1}\mathrm{Pf}\hspace{1pt}(M)$ pour tout $k\geq 0$ ($\partial$~est un isomorphisme pour $k>0$). Nous avons déjà montré que si $M$ est muni d'une action à droite tordue de $\mathrm{GL}(V)$ alors il en est de même pour $\mathrm{R}^{k+1}\mathrm{Pf}\hspace{1pt}(M)$ (Proposition \ref{GL-action3.5}). Nous allons d'abord montrer, de façon formelle, que $\lim_{\mathcal{W}_{0}}^{k}\Psi(M)$ est aussi canoniquement muni d'une action à droite tordue de $\mathrm{GL}(V)$ et que pour cette action l'homomorphisme~$\partial$ ci-dessus est $\mathrm{GL}(V)$-équivariant. Nous décrirons ensuite (Proposition \ref{GL-actionlimPsi}) cette action de $\mathrm{GL}(V)$ sur $\lim_{\mathcal{W}_{0}}^{k}\Psi(M)$ en termes de la définition ``combinatoire'' des dérivés  de $\lim$ (rappel (R.2) de \ref{rappelslim}). Passons à la réalisation de ce programme~:

\medskip
Supposons $M$ muni d'une action à droite tordue de $\mathrm{GL}(V)$. On dispose donc (Proposition \ref{thetaalphaGL-action}) d'une famille $(a_{\alpha}:M\to\theta_{\alpha}M)_{\alpha\in\mathrm{GL}(V)}$ d'homomorphismes dans la catégorie $V_{\mathrm{tf}}\text{-}\mathcal{U}$ avec $a_{\alpha\beta}=(\theta_{\alpha}a_{\beta})\circ a_{\alpha}$ pour tous $\alpha$ et $\beta$ dans $\mathrm{GL}(V)$ et $a_{1_{\mathrm{GL}(V)}}=\mathrm{id}_{M}$. Le fait que $\theta_{\alpha}M\to\theta_{\alpha}I^{\bullet}$ est encore une résolution injective dans la  catégorie $V_{\mathrm{tf}}\text{-}\mathcal{U}$ (Proposition \ref{thetaalpharesinj}) entraîne (Proposition \ref{GL-actionresinj}) que l'on dispose également d'une famille $(a_{\alpha}^{\bullet}:I^{\bullet}\to\theta_{\alpha}I^{\bullet})_{\alpha\in\mathrm{GL}(V)}$ d'homomorphismes de complexes dans la catégorie $V_{\mathrm{tf}}\text{-}\mathcal{U}$, unique à homotopie près, homomorphismes qui font commuter les diagrammes
$$
\begin{CD}
M@>>>I^{\bullet} \\
@Va_{\alpha}VV @Va_{\alpha}^{\bullet}VV \\
\theta_{\alpha}M@>>>\theta_{\alpha}I^{\bullet}
\end{CD}
$$
et tels que $a_{\alpha\beta}^{\bullet}$ est homotope à $(\theta_{\alpha}a_{\beta}^{\bullet})\circ a_{\alpha}^{\bullet}$ pour tous $\alpha$, $\beta$ dans $\mathrm{GL}(V)$ (et que $a_{1_{\mathrm{GL}(V)}}^{\bullet}$ est l'identité).

\medskip
Soit $\mathrm{R\acute{e}d}$ l'endofoncteur $M\mapsto M/\mathrm{Pf}(M)$ de la catégorie $V_{\mathrm{tf}}\text{-}\mathcal{U}$. La transformation naturelle, entre endofoncteurs de $V_{\mathrm{tf}}\text{-}\mathcal{U}$, $\rho:\mathrm{id}\to\lim_{\mathcal{W}_{0}}\Psi(-)$, qui apparaît dans la proposition \ref{derivePf}, induit d'après cette proposition une transformation naturelle $\bar{\rho}:\mathrm{R\acute{e}d}\to\lim_{\mathcal{W}_{0}}\Psi(-)$ telle que $\bar{\rho}_{I}$ est un isomorphisme si~$I$ est injectif (point (c) de \ref{proPsi}). On note $b_{\alpha}^{\bullet}$ l'homomorphisme de complexes composé
\begin{multline*}
\begin{CD} \lim_{\mathcal{W}_{0}}\Psi(I^{\bullet})
@>(\bar{\rho}_{I^{\bullet}})^{-1}>>
\mathrm{R\acute{e}d}(I^{\bullet})
@>\mathrm{R\acute{e}d}(a_{\alpha}^{\bullet})>>
\mathrm{R\acute{e}d}(\theta_{\alpha}I^{\bullet})=
\theta_{\alpha}\mathrm{R\acute{e}d}(I^{\bullet}) \end{CD} \\
\begin{CD} 
@>\theta_{\alpha}\bar{\rho}_{I^{\bullet}}>>
\theta_{\alpha}\lim_{\mathcal{W}_{0}}\Psi(I^{\bullet})\hspace{8pt}.
\end{CD}
\end{multline*}
Il est clair que $b_{\alpha\beta}^{\bullet}$ est homotope à $(\theta_{\alpha}b_{\beta}^{\bullet})\circ b_{\alpha}^{\bullet}$ pour tous $\alpha$, $\beta$ dans $\mathrm{GL}(V)$ (et que $b_{1_{\mathrm{GL}(V)}}^{\bullet}$ est l'identité)~: la famille $(b_{\alpha}^{\bullet}:\lim_{\mathcal{W}_{0}}\hspace{-1pt}\Psi(I^{\bullet})\to\theta_{\alpha}\lim_{\mathcal{W}_{0}}\hspace{-1pt}\Psi(I^{\bullet})_{\alpha\in\mathrm{GL}(V)}$ fournit une action à droite tordue de $\mathrm{GL}(V)$ sur les $\lim^{k}_{\mathcal{W}_{0}}\Psi(M)$, $k\geq 0$, uniquement déterminée par l'action à droite tordue de $\mathrm{GL}(V)$ sur $M$. La contemplation du diagramme commutatif de complexes (dont les deux lignes sont exactes) 
$$
\begin{CD}
0@>>>\mathrm{Pf}(I^{\bullet})
@>\subset>>I^{\bullet}
@>\rho_{I^{\bullet}}>>
\lim_{\mathcal{W}_{0}}\Psi(I^{\bullet})@>>>0 \\
&& @V\mathrm{Pf}(a_{\alpha}^{\bullet})VV @Va_{\alpha}^{\bullet}VV @Vb_{\alpha}^{\bullet}VV  \\
0@>>>\theta_{\alpha}\mathrm{Pf}(I^{\bullet})
@>\subset>>\theta_{\alpha}I^{\bullet}
@>\theta_{\alpha}\rho_{I^{\bullet}}>>
\theta_{\alpha}\lim_{\mathcal{W}_{0}}\Psi(I^{\bullet})@>>>0
\end{CD}
$$
montre que les homomorphismes naturels  $\partial:\lim_{\mathcal{W}_{0}}^{k}\Psi(M)\to\mathrm{R}^{k+1}\mathrm{Pf}\hspace{1pt}(M)$ sont $\mathrm{GL}(V)$-équivariants.

\bigskip
La façon dont nous avons obtenue ci-dessus une action à droite tordue de $\mathrm{GL}(V)$ sur les $\lim_{\mathcal{W}_{0}}^{k}\Psi(M)$ est très formelle~; nous allons ci-dessous décrire cette action en termes un peu plus concrets, dans la ligne du (R.2) de \ref{rappelslim}.

\pagebreak
\medskip
\begin{defi}\label{lambdaphi1} On reprend les notations du (R.2) de \ref{rappelslim}.

\smallskip
Soient $F$ un foncteur de $\mathcal{C}$ dans $\mathcal{A}$ et $\phi$ un automorphisme de $\mathcal{C}$ (en clair une permutation qui respecte l'ordre).

\smallskip
On note $\lambda_{\phi}^{\bullet}:\mathrm{L}^{\bullet}(F)\to\mathrm{L}^{\bullet}(F\circ\phi)$ l'homomorphisme de complexes tel que le $\mathcal{A}$-morphisme
$$
\lambda_{\phi}^{k}:\mathrm{L}^{k}(F):=\prod_{\sigma\in\mathrm{S}_{k}\mathcal{C}}F(\sup\sigma)
\longrightarrow
\prod_{\sigma\in\mathrm{S}_{k}\mathcal{C}}F(\phi\sup\sigma)
=:\mathrm{L}^{k}(F\circ\phi)
$$
est le produit  $\prod_{\sigma\in\mathrm{S}_{k}\mathcal{C}}\pi_{\phi\sigma}$,  $\pi_{\phi\sigma}$ désignant la projection du produit à la source sur le facteur d'indice $\phi\sigma$ à savoir $F(\sup\phi\sigma)=F(\phi\sup\sigma)$.

\smallskip
On note encore $\lambda_{\phi}^{k}:\lim_{\mathcal{C}}^{k}F\to\lim_{\mathcal{C}}^{k}\hspace{1pt}(F\circ\phi)$ la transformation naturelle $\mathrm{H}^{k}(\lambda_{\phi}^{\bullet})$~; la transformation naturelle $\lambda_{\phi}^{0}:\lim_{\mathcal{C}}F\to\lim_{\mathcal{C}}\hspace{1pt}(F\circ\phi)$ est aussi simplement notée $\lambda_{\phi}$.
\end{defi}

\bigskip
La vérification de l'énoncé suivant est immédiate~:

\begin{pro}\label{lambdaphi2} Soient $F$ un foncteur de $\mathcal{C}$ dans $\mathcal{A}$ (notations du (R.2) de \ref{rappelslim}).

\medskip
{\em (a)} Soient $\phi$ et $\psi$ deux automorphismes de $\mathcal{C}$. Alors le composé des homomorphismes $\lambda_{\phi}^{\bullet}:\mathrm{L}^{\bullet}(F)\to\mathrm{L}^{\bullet}(F\circ\phi)$ et $\lambda_{\psi}^{\bullet}:\mathrm{L}^{\bullet}(F\circ\phi)\to\mathrm{L}^{\bullet}(F\circ(\phi\circ\psi))$ est~$\lambda_{\phi\circ\psi}^{\bullet}:\mathrm{L}^{\bullet}(F)\to\mathrm{L}^{\bullet}(F\circ(\phi\circ\psi))$.

\medskip
{\em (b)} Soit $\Gamma$ un sous-groupe d'automorphismes de $\mathcal{C}$ avec $F\circ\phi=F$ pour tout $\phi$ dans $\Gamma$~; soit $k\geq 0$ un entier. Alors les $\lambda_{\phi}^{k}:\lim_{\mathcal{C}}^{k}F\to\lim_{\mathcal{C}}^{k}\hspace{1pt}(F\circ\phi)$, $\phi$~parcourant $\Gamma$, définissent une action à droite de $\Gamma$ sur $\lim_{\mathcal{C}}^{k}F$.
\end{pro}

\medskip
\begin{pro}\label{lambdaphi3} On reprend les notations du (R.2) de \ref{rappelslim}.

\smallskip
Soient $F$ un foncteur de $\mathcal{C}$ dans $\mathcal{A}$, $\phi$ un automorphisme de $\mathcal{C}$ et $F\to I^{\bullet}$ une résolution injective dans la catégorie $\mathcal{A}^{\mathcal{C}}$. Alors les deux propriétés suivantes sont vérifiées~:

\smallskip
--\hspace{8pt}$F\circ\phi\to I^{\bullet}\circ\phi$ une résolution injective dans la catégorie $\mathcal{A}^{\mathcal{C}}$~;

\smallskip
--\hspace{8pt}le $\mathcal{A}$-morphisme induit entre $k$-ième objets de cohomologie par l'homomorphisme de $\mathcal{A}$-complexes $\lambda_{\phi}:\lim_{\mathcal{C}}I^{\bullet}\to\lim_{\mathcal{C}}\hspace{1pt}(I^{\bullet}\circ\phi)$ s'identifie au $\mathcal{A}$-morphisme $\lambda_{\phi}^{k}:\lim_{\mathcal{C}}^{k}F\to\lim_{\mathcal{C}}^{k}\hspace{1pt}(F\circ\phi)$.
\end{pro}

\medskip
\textit{Démonstration.} La première propriété est évidente~: l'endofoncteur $F\mapsto F\circ\phi$ de $\mathcal{A}^{\mathcal{C}}$ est exact et préserve les injectifs. Passons à la seconde. La famille d'homomorphismes de complexes $(\lambda_{\phi}^{\bullet}:\mathrm{L}^{\bullet}(I^{q})\to\mathrm{L}^{\bullet}(I^{q}\circ\phi))_{q\in\mathbb{N}}$ définit un homomorphisme de bicomplexes $\mathrm{L}^{\bullet}(I^{\bullet})\to\mathrm{L}^{\bullet}(I^{\bullet}\circ\phi)$ que l'on note $\lambda_{\phi}^{\bullet,\bullet}$. La seconde propriété est obtenue en appliquant le foncteur $\mathrm{H}^{k}$ au diagramme commutatif de complexes
$$
\begin{CD}
\lim_{\mathcal{C}}I^{\bullet}@>\eta_{\mathrm{v}}>>
\mathrm{Tot}\hspace{2pt}\mathrm{L}^{\bullet}(I^{\bullet})
@<\eta_{\mathrm{h}}<< \mathrm{L}^{\bullet}(F) \\
@V\lambda_{\phi}VV @V\mathrm{Tot}\hspace{1pt}\lambda_{\phi}^{\bullet,\bullet}VV @V\lambda_{\phi}^{\bullet}VV \\
\lim_{\mathcal{C}}\hspace{1pt}(I^{\bullet}\circ\phi)@>\eta_{\mathrm{v}}>>
\mathrm{Tot}\hspace{2pt}\mathrm{L}^{\bullet}(I^{\bullet}\circ\phi)
@<\eta_{\mathrm{h}}<< \mathrm{L}^{\bullet}(F\circ\phi)
\end{CD}
$$
(pour les notations se reporter au dernier point du (R.2) de \ref{rappelslim}).
\hfill$\square$

\bigskip
On revient maintenant  au cas $\mathcal{A}=V_{\mathrm{tf}}\text{-}\mathcal{U}$, $\mathcal{C}=\mathcal{W}_{0}$ et $F=\Psi(M)$. On constate que l'énoncé (technique~!) suivant est vérifié~:

\begin{lem}\label{factthetaalpharho} Soient $M$ un $\mathrm{H}^{*}V_{\mathrm{tf}}\text{-}\mathrm{A}$-module instable et $\alpha$ un élément de $\mathrm{GL}(V)$ (que l'on peut aussi considérer comme un automorphisme de $\mathcal{W}_{0}$). L'homomorphisme $\theta_{\alpha}\hspace{1pt}\rho_{M}:\theta_{\alpha}M\to\theta_{\alpha}\lim_{\mathcal{W}_{0}}\Psi(M)$
admet la factorisation suivante~:
\begin{multline*}
\hspace{12pt}
\begin{CD}
\theta_{\alpha}M@>\rho\hspace{1pt}_{\theta_{\alpha}M}>>
\lim_{\mathcal{W}_{0}}\Psi(\theta_{\alpha}M)
@>\cong>>
\theta_{\alpha}\lim_{\mathcal{W}_{0}}(\Psi(M)\circ\alpha^{-1})
\end{CD} \\
\begin{CD}
@>\theta_{\alpha}\hspace{1pt}\lambda_{\alpha}>>
\theta_{\alpha}\lim_{\mathcal{W}_{0}}\Psi(M)
\end{CD}
\hspace{12pt},
\end{multline*}
la deuxième flèche étant induite par l'isomorphisme de \ref{Psithetaalpha}.
\end{lem}

\medskip
Ce lemme conduit à la proposition ci-dessous.

\medskip
\begin{pro}\label{GL-actionlimPsi} Soient $M$ un $\mathrm{H}^{*}V_{\mathrm{tf}}\text{-}\mathrm{A}$-module instable  muni d'une action à droite tordue de $\mathrm{GL}(V)$, $k\geq 0$ un entier  et $\alpha$ un élément de $\mathrm{GL}(V)$~; soit $a_{\alpha}:M\to\theta_{\alpha}M$, $\alpha\in\mathrm{GL}(V)$, l'homomorphisme $x\mapsto x.\alpha$.

\smallskip
Alors l'homomorphisme $\lim_{\mathcal{W}_{0}}^{k}\Psi(M)\to\theta_{\alpha}\lim_{\mathcal{W}_{0}}^{k}\Psi(M)$ donnée par l'action à droite tordue de $\mathrm{GL}(V)$ sur $\lim_{\mathcal{W}_{0}}^{k}\Psi(M)$ admet la factorisation suivante~:
\begin{multline*}
\hspace{12pt}
\begin{CD}
\lim_{\mathcal{W}_{0}}^{k}\Psi(M)
@>\lim_{\mathcal{W}_{0}}^{k}\Psi(a_{\alpha})>>
\lim_{\mathcal{W}_{0}}^{k}\Psi(\theta_{\alpha}M)
@>\cong>>
\theta_{\alpha}\lim_{\mathcal{W}_{0}}^{k}(\Psi(M)\circ\alpha^{-1})
\end{CD} \\
\begin{CD}
@>\theta_{\alpha}\hspace{1pt}\lambda_{\alpha}^{k}>>
\theta_{\alpha}\lim_{\mathcal{W}_{0}}^{k}\Psi(M)
\end{CD}
\hspace{12pt},
\end{multline*}
la deuxième flèche étant induite par l'isomorphisme de \ref{Psithetaalpha}.

\end{pro}

\medskip
\textit{Démonstration.} On reprend le diagramme commutatif
$$
\begin{CD}
I^{\bullet}
@>\rho_{I^{\bullet}}>>
\lim_{\mathcal{W}_{0}}\Psi(I^{\bullet}) \\
@Va_{\alpha}^{\bullet}VV @Vb_{\alpha}^{\bullet}VV  \\
\theta_{\alpha}I^{\bullet}
@>\theta_{\alpha}\rho_{I^{\bullet}}>>
\theta_{\alpha}\lim_{\mathcal{W}_{0}}\Psi(I^{\bullet})
\end{CD}
$$
considéré plus haut. Soit $c_{\alpha}^{\bullet}$ l'homomorphisme de complexes composé
\begin{multline*}
\hspace{12pt}
\begin{CD}
\lim_{\mathcal{W}_{0}}\Psi(I^{\bullet})
@>\lim_{\mathcal{W}_{0}}\Psi(a_{\alpha}^{\bullet})>>
\lim_{\mathcal{W}_{0}}\Psi(\theta_{\alpha}I^{\bullet})
@>\cong>>
\theta_{\alpha}\lim_{\mathcal{W}_{0}}(\Psi(I^{\bullet})\circ\alpha^{-1})
\end{CD} \\
\begin{CD}
@>\theta_{\alpha}\hspace{1pt}\lambda_{\alpha}>>
\theta_{\alpha}\lim_{\mathcal{W}_{0}}\Psi(I^{\bullet})
\end{CD}
\hspace{12pt}.
\end{multline*}
En utilisant la commutativité du diagramme
$$
\begin{CD}
I^{\bullet}
@>\rho_{I^{\bullet}}>>
\lim_{\mathcal{W}_{0}}\Psi(I^{\bullet}) \\
@Va_{\alpha}^{\bullet}VV @Vlim_{\mathcal{W}_{0}}\Psi(a_{\alpha}^{\bullet})VV  \\
\theta_{\alpha}I^{\bullet}
@>\rho\hspace{1pt}_{\theta_{\alpha}I^{\bullet}}>>
\lim_{\mathcal{W}_{0}}\Psi(\theta_{\alpha}I^{\bullet})
\end{CD}
$$
(naturalité de $\rho$) et le lemme \ref{factthetaalpharho} on constate que le diagramme
$$
\begin{CD}
I^{\bullet}
@>\rho_{I^{\bullet}}>>
\lim_{\mathcal{W}_{0}}\Psi(I^{\bullet}) \\
@Va_{\alpha}^{\bullet}VV @Vc_{\alpha}^{\bullet}VV  \\
\theta_{\alpha}I^{\bullet}
@>\theta_{\alpha}\rho_{I^{\bullet}}>>
\theta_{\alpha}\lim_{\mathcal{W}_{0}}\Psi(I^{\bullet})
\end{CD}
$$
est aussi commutatif. On en déduit  $b_{\alpha}^{\bullet}=c_{\alpha}^{\bullet}$ (car $\rho_{I^{\bullet}}$ est surjectif).
\hfill$\square$

\bigskip
On prend enfin $M=\mathrm{c}_{V}\mathrm{H}^{*}V$, $k=n-1$ ($n=\dim V$) et on se focalise sur le degré zéro. On rappelle que l'on a posé $\mathrm{F}:=\Psi^{0}(\mathrm{c}_{V}\mathrm{H}^{*}V)$.

\smallskip
On constate que l'on a $\mathrm{End}_{\mathcal{E}^{\mathcal{W}_{0}}}(\mathrm{F})=\{0,\mathrm{id}\}$ et donc que tout automorphisme de $\mathrm{F}$ est l'identité. On constate également que l'on a $\mathrm{F}\circ\alpha=\mathrm{F}$ pour tout~$\alpha$ dans $\mathrm{GL}(V)$ et donc que les $\lambda_{\alpha}^{n-1}:\lim_{\mathcal{W}_{0}}^{n-1}\mathrm{F}\to\lim_{\mathcal{W}_{0}}^{n-1}\mathrm{F}$ définissent une action à droite de $\mathrm{GL}(V)$ sur $\lim_{\mathcal{W}_{0}}^{n-1}\mathrm{F}$ (point (b) de \ref{lambdaphi2}). Ces constations faites, la proposition \ref{GL-actionlimPsi} permet de se convaincre que si l'on munit $\lim_{\mathcal{W}_{0}}^{n-1}\mathrm{F}$ de cette action alors l'homomorphisme $\lim_{\mathcal{W}_{0}}^{n-1}\mathrm{F}\to\mathrm{M}^{0}(V)$, considéré plus haut, est $\mathrm{GL}(V)$-équivariant.

\smallskip
On reprend ensuite la suite exacte $0\to\mathrm{F}\to\mathrm{E}\to\mathrm{Q}\to 0$. On observe que l'on a aussi $\mathrm{E}\circ\alpha=\mathrm{E}$ et $\mathrm{Q}\circ\alpha=\mathrm{Q}$ pour tout $\alpha$ dans $\mathrm{GL}(V)$. En contemplant les diagrammes commutatifs de complexes
$$
\begin{CD}
0@>>>\mathrm{L}^{\bullet}(\mathrm{F})
@>>>\mathrm{L}^{\bullet}(\mathrm{E})
@>>>\mathrm{L}^{\bullet}(\mathrm{Q})@>>>0 \\
&& @V\lambda_{\alpha}^{\bullet}VV
@V\lambda_{\alpha}^{\bullet}VV.
@V\lambda_{\alpha}^{\bullet}VV \\
0@>>>\mathrm{L}^{\bullet}(\mathrm{F})
@>>>\mathrm{L}^{\bullet}(\mathrm{E})
@>>>\mathrm{L}^{\bullet}(\mathrm{Q})@>>>0
\end{CD}
$$
dont les lignes sont exactes (et en invoquant, si l'on en éprouve le besoin, le rappel (R.3) de \ref{rappelslim}) on voit que la suite exacte
$$
\begin{CD}
\lim_{\mathcal{W}_{0}}^{n-2}\mathrm{E}
@>>>\lim_{\mathcal{W}_{0}}^{n-2}\mathrm{Q}
@>>>\lim_{\mathcal{W}_{0}}^{n-1}\mathrm{F}
@>>> 0
\end{CD}
$$
est une suite exacte de $\mathbb{F}_{2}[\mathrm{GL}(V)]$-modules à droite. Ceci achève (finalement~!) la démonstration du point (b) de la proposition \ref{MV}.
\hfill$\square\square$

\pagebreak

\vspace{0.75cm}
\textsc{Les complexes $\widetilde{\mathrm{C}}^{\bullet}\hspace{1pt}(\mathrm{c}_{V}^{h}\hspace{1pt}\mathrm{H}^{*}V)$, $h\in\mathbb{N}$}

\bigskip
Les résultats concernant le complexe $\widetilde{\mathrm{C}}^{\bullet}\hspace{1pt}(\mathrm{c}_{V}\hspace{1pt}\mathrm{H}^{*}V)$ que nous venons de présenter s'étendent \textit{mutatis mutandis} aux complexes $\widetilde{\mathrm{C}}^{\bullet}\hspace{1pt}(\mathrm{c}_{V}^{h}\hspace{1pt}\mathrm{H}^{*}V)$, $h\in\mathbb{N}$ (on aura intérêt à traiter séparément le cas $h=0$\ldots qui est vraiment trivial).

\medskip
Posons
$$
\mathrm{M}(V;h):=\mathrm{R}^{n}\mathrm{Pf}(\mathrm{c}_{V}^{h}\hspace{1pt}\mathrm{H}^{*}V)\cong\Sigma^{-n}\hspace{1pt}\mathrm{H}^{*}_{\mathrm{c}}(V\backslash(\widetilde{\mathbb{R}}[V]^{\oplus h})_{\mathrm{r\acute{e}g}})
$$
(on a donc $\mathrm{M}(V;1)=\mathrm{M}(V)$).

\smallskip
\hspace{8pt}-- Le complexe $\mathrm{C}^{\bullet}\hspace{1pt}(\mathrm{c}_{V}^{h}\hspace{1pt}\mathrm{H}^{*}V)$ est obtenu en remplaçant $\mathrm{M}(-)$ par $\mathrm{M}(-;h)$.

\smallskip
\hspace{8pt}-- L'énoncé obtenu en remplaçant $\mathrm{M}(-)$ par $\mathrm{M}(-;h)$ dans l'énoncé \ref{MV} reste valable pour $h\geq 1$.

\smallskip
Précisons la seconde affirmation ci-dessus et tirons-en les conséquences~:

\medskip
\begin{pro}\label{surjectiviteM} Soient $h$ et $h'$ deux entiers avec $1\leq h\leq h'$. L'inclusion de $\mathrm{H}^{*}V_{\mathrm{tf}}$-$\mathrm{A}$-modules instables, $\mathrm{c}_{V}^{h'}\hspace{1pt}\mathrm{H}^{*}V\hookrightarrow\mathrm{c}_{V}^{h}\hspace{1pt}\mathrm{H}^{*}V$, induit un isomorphisme $\mathrm{M}^{0}(V;h')\cong\mathrm{M}^{0}(V;h)$ et un homomorphisme surjectif $\mathrm{M}(V;h')\twoheadrightarrow\mathrm{M}(V;h)$ de $\mathrm{GL}(V)_{\mathrm{td}}$-$\mathrm{H}^{*}V_{\mathrm{tf}}$-$\mathrm{A}$-modules instables.

\smallskip
{\em (La terminologie {\em  $\mathrm{GL}(V)_{\mathrm{td}}$-$\mathrm{H}^{*}V_{\mathrm{tf}}$-$\mathrm{A}$-module instable} est transparente.)}

\end{pro}

\textit{Démonstration.} Compte tenu du point (a) de \ref{MV}, la seconde partie de la proposition résulte de la première. Pour démontrer celle-ci on considère le diagramme commutatif
$$
\begin{CD}
\lim_{\mathcal{W}_{0}}^{n-1}(\Psi(\mathrm{c}_{V}^{h'}\hspace{1pt}\mathrm{H}^{*}V))^{0}
@>\cong>> \mathrm{M}^{0}(V;h') \\
@VVV @VVV \\
\lim_{\mathcal{W}_{0}}^{n-1}(\Psi(\mathrm{c}_{V}^{h}\hspace{1pt}\mathrm{H}^{*}V))^{0}
@>\cong>> \mathrm{M}^{0}(V;h)
\end{CD}
$$
dans lequel les isomorphismes horizontaux sont fournis par la proposition \ref{proPsi} (on suppose $V\not=0$) et les flèches verticales sont induites par  l'inclusion $\mathrm{c}_{V}^{h'}\hspace{1pt}\mathrm{H}^{*}V\hookrightarrow\mathrm{c}_{V}^{h}\hspace{1pt}\mathrm{H}^{*}V$, disons $\iota$~;  on conclut en observant que la transformation naturelle $\Psi(\mathrm{c}_{V}^{h'}\hspace{1pt}\mathrm{H}^{*}V))^{0}\to\Psi(\mathrm{c}_{V}^{h}\hspace{1pt}\mathrm{H}^{*}V))^{0}$ induite par $\iota$ est un isomorphisme fonctoriel d'après \ref{spePsi}.
\hfill$\square$

\vspace{0.5cm} \textit{Relation entre $\mathrm{M}(V;h)$ et $\mathrm{M}(V;2h)$}

\bigskip
Pour décrire cette relation il nous faut introduire l'endofoncteur ``double'' de la catégorie $\mathcal{U}$, foncteur noté $\Phi$  (voir par exemple \cite[\S 2.3]{LZens}). Rappelons sa définition. Soit $M$ un $\mathrm{A}$-module instable, on note $\Phi M$ le $\mathrm{A}$-module instable défini par
$$
(\Phi M)^{k}=
\begin{cases}
M^{\frac{k}{2}} & \text{si $k$ est pair} \\
0 &  \text{sinon}
\end{cases}
\hspace{12pt}\text{et}\hspace{12pt}
\mathrm{Sq}^{i}(\Phi x)=
\begin{cases}
\Phi(\mathrm{Sq}^{\frac{i}{2}}x) & \text{si $i$ est pair} \\
0 &  \text{sinon},
\end{cases}
$$
$\Phi x$ désignant l'élément de $(\Phi M)^{2k}$ associé à un élément $x$ de $M^{k}$. L'application, $\Phi M\to M,\Phi x\mapsto\mathrm{Sq}^{\vert x\vert}x$ ($\vert x\vert$ désignant le degré de $x$) est $\mathrm{A}$-linéaire~; elle fournit une transformation naturelle, que l'on note $\lambda$, de $\Phi$ dans $1_{\mathcal{U}}$. Soit $K$ une $\mathrm{A}$-algèbre instable~; $\Phi K$ est muni d'une structure naturelle de $\mathrm{A}$-algèbre instable et $\lambda_{K}:\Phi K\to K$ est un homomorphisme de $\mathrm{A}$-algèbres instables qui par définition même d'une $\mathrm{A}$-algèbre instable coïncide avec l'élévation au carré. Si $M$ est muni d'une structure de $K$-$\mathrm{A}$-module instable alors $\Phi M$ est naturellement muni d'une structure de $\Phi K$-$\mathrm{A}$-module instable. On note $\mathrm{E}\Phi:V\text{-}\mathcal{U}\to V\text{-}\mathcal{U}$ le foncteur
$$
\hspace{24pt}
M
\hspace{4pt}\mapsto\hspace{4pt}
\mathrm{H}^{*}V\otimes_{\Phi\mathrm{H}^{*}V}\Phi M
\hspace{24pt},
$$
$\mathrm{H}^{*}V$ étant ci-dessus un $\Phi\mathrm{H}^{*}V$-module \textit{via} $\lambda_{\mathrm{H}^{*}V}$~;  la transformation naturelle $\lambda:\Phi\to 1_{\mathcal{U}}$ induit une transformation naturelle $\mathrm{E}\Phi\to 1_{V\text{-}\mathcal{U}}$ que l'on note~$\mathrm{E}\lambda$.

\begin{pro}\label{double1} Le foncteur $\mathrm{E}\Phi$ possède les cinq propriétés suivantes~:

\medskip
{\em (a)} Le foncteur $\mathrm{E}\Phi$ est exact.

\medskip
{\em (b)} Si $M$ est un $\mathrm{H}^{*}V$-$\mathrm{A}$-module instable de type fini comme $\mathrm{H}^{*}V$-module alors il en est de même pour $\mathrm{E}\Phi(M)$.

\medskip
{\em (c)} Si $M$ est un $\mathrm{H}^{*}V$-$\mathrm{A}$-module instable e-fini (terminologie introduite en \ref{defefini}) alors il en est de même pour $\mathrm{E}\Phi(M)$.

\medskip
{\em (d)} Soient $M$ est un $\mathrm{H}^{*}V$-$\mathrm{A}$-module instable et $\alpha$ un élément de $\mathrm{GL}(V)$, on a un isomorphisme de  $\mathrm{H}^{*}V$-$\mathrm{A}$-modules instables, naturel en $M$~:
$$
\hspace{24pt}
\mathrm{E}\Phi(\hspace{1pt}\theta_{\alpha}M)
\hspace{4pt}\cong\hspace{4pt}
\theta_{\alpha}\hspace{1pt}\mathrm{E}\Phi(M)
\hspace{24pt}.
$$

\medskip
{\em (e)} Si $M$ est un $\mathrm{H}^{*}V$-$\mathrm{A}$-module instable muni d'une action à droite tordue de $\mathrm{GL}(V)$ alors $\mathrm{E}\Phi(M)$ est aussi naturellement muni d'une telle action. En d'autres termes l'endofoncteur $\mathrm{E}\Phi$ de $V\text{-}\mathcal{U}$ induit un endofoncteur de $\mathrm{GL}(V)_{\mathrm{td}}\text{-}V\text{-}\mathcal{U}$.
\end{pro}

\medskip
\textit{Démonstration de (a), (b) et (c).} On note tout d'abord que $\mathrm{H}^{*}V$ est un $\Phi\mathrm{H}^{*}V$-module libre de dimension finie (avec $\mathbb{F}_{2}\otimes_{\Phi\mathrm{H}^{*}V}\mathrm{H}^{*}V$ isomorphe à l'algèbre extérieure $\Lambda^{*}V$).

\medskip
(a) Le foncteur $\Phi$ est exact et $\mathrm{H}^{*}V$ est un $\Phi\mathrm{H}^{*}V$-module plat.

\medskip
(b) Le $\Phi\mathrm{H}^{*}V$-module $\Phi M$ est de type fini.

\medskip
(c) Soient $W$ un sous-groupe de $V$ et $N$ un $\mathrm{H}^{*}V/W$-$\mathrm{A}$-module instable. On constate que l'on a un isomorphisme de $\mathrm{H}^{*}V$-$\mathrm{A}$-modules instables
$$
\mathrm{H}^{*}V\otimes_{\Phi\mathrm{H}^{*}V}\Phi(\mathrm{H}^{*}V\otimes_{\mathrm{H}^{*}V/W}N)
\hspace{4pt}\cong\hspace{4pt}
\mathrm{H}^{*}V\otimes_{\mathrm{H}^{*}V/W}(\mathrm{H}^{*}V/W\otimes_{\Phi\mathrm{H}^{*}V/W}\Phi N)
$$
(naturel en $N$) et que si le $\mathrm{H}^{*}V/W$-$\mathrm{A}$-module instable $N$ est fini alors il en est de même pour le $\mathrm{H}^{*}V/W$-$\mathrm{A}$-module instable $\mathrm{H}^{*}V/W\otimes_{\Phi\mathrm{H}^{*}V/W}\Phi N$ (parce que $\mathrm{H}^{*}V/W$ est un $\Phi\mathrm{H}^{*}V/W$-module libre de dimension finie).
\hfill$\square$

\medskip
\textit{Démonstration de (d) et (e).} Le point (e) résulte formellement du point (d) et de \ref{thetaalphaGL-action}. L'isomorphisme en question dans le point (d) est induit par l'application $\mathrm{H}^{*}V\times\Phi M\to\mathrm{E}\Phi(M)\hspace{2pt},\hspace{2pt}(a,\Phi x)\mapsto\alpha^{*}a\otimes_{\Phi\mathrm{H}^{*}V}\Phi x$.
\hfill$\square$

\bigskip
\begin{pro}\label{double2} Soit $M$ un  $\mathrm{H}^{*}V_{\mathrm{tf}}$-$\mathrm{A}$-module instable.

\medskip
{\em (a)} On a un isomorphisme de $\mathrm{H}^{*}V_{\mathrm{tf}}$-$\mathrm{A}$-modules instables
$$
\mathrm{R}^{p}\mathrm{Pf}\hspace{2pt}(\mathrm{E}\Phi(M))
\hspace{4pt}\cong\hspace{4pt}
\mathrm{E}\Phi\hspace{1pt}(\mathrm{R}^{p}\mathrm{Pf}\hspace{1pt}M)
$$
(naturel en $M$) pour tout entier naturel $p$.

\medskip
{\em (b)} Si $M$ est muni d'une action à droite tordue de $\mathrm{GL}(V)$ alors l'isomorphisme ci-dessus préserve les actions à droite tordues données par  \ref{GL-action3.5} et le point (e) de \ref{double1}.

\end{pro}

\medskip
\textit{Démonstration du (a).} Soit $r:M\to I^{\bullet}$ une résolution injective dans la caté\-gorie $V_{\mathrm{tf}}\text{-}\mathcal{U}$. Comme les injectifs de cette catégorie sont e-finis, les trois premiers points de la proposition \ref{double1} montre que $\mathrm{E}\Phi(r):\mathrm{E}\Phi(M)\to\mathrm{E}\Phi(I^{\bullet})$ est une résolution e-finie au sens de la partie définition de \ref{cordefefini}. La partie corollaire de \ref{cordefefini} permet de conclure.

\medskip
\textit{Démonstration du (b).} On note (localement) $\nu$ l'isomorphisme naturel du point (d) de \ref{double1} et on reprend les notations de \ref{GL-actionresinj}.

\smallskip
Soit $s:\mathrm{E}\Phi(I^{\bullet})\to J^{\bullet}$ un remplacement injectif dans la catégorie $V_{\mathrm{tf}}\text{-}\mathcal{U}$ (terminologie introduite au début de la section 6)~; on constate que $s\circ\mathrm{E}\Phi(r):\mathrm{E}\Phi(M)\to J^{\bullet}$ est une résolution injective. Comme $\theta_{\alpha}s:\theta_{\alpha}\mathrm{E}\Phi(I^{\bullet})\to\theta_{\alpha}J^{\bullet}$ est également un remplacement injectif, il existe un homomorphisme de complexes $b_{\alpha}^{\bullet}$, unique à homotopie près, tel que le diagramme
$$
\begin{CD}
\mathrm{E}\Phi(I^{\bullet})
@>s>>
J^{\bullet} \\
@V\nu_{I^{\bullet}}\circ\mathrm{E}\Phi(a_{\alpha}^{\bullet})VV 
@Vb_{\alpha}^{\bullet}VV \\
\theta_{\alpha}\mathrm{E}\Phi(I^{\bullet})
@>\theta_{\alpha}s>>
\theta_{\alpha}J^{\bullet} 
\end{CD}
$$
est commutatif à homotopie près (voir \ref{catder2}). Il en résulte que le diagramme
$$
\begin{CD}
\mathrm{E}\Phi(M)
@>s\circ\mathrm{E}\Phi(r)>>
J^{\bullet} \\
@V\nu_{M}\circ\mathrm{E}\Phi(a_{\alpha})VV 
@Vb_{\alpha}^{\bullet}VV \\
\theta_{\alpha}\mathrm{E}\Phi(M)
@>\theta_{\alpha}(s\circ\mathrm{E}\Phi(r))>>
\theta_{\alpha}J^{\bullet} 
\end{CD}
$$
est commutatif. Comme les actions de $\mathrm{GL}(V)$ sur $\mathrm{R}^{p}\mathrm{Pf}\hspace{1pt}M)$ et $\mathrm{R}^{p}\mathrm{Pf}\hspace{2pt}(\mathrm{E}\Phi(M)$) sont respectivement définies par les homomorphismes $\mathrm{H}^{p}(a_{\alpha}^{\bullet})$ et $\mathrm{H}^{p}(b_{\alpha}^{\bullet})$, les deux commutativités ci-dessus impliquent le point (b).
\hfill$\square$

\medskip
\begin{cor}\label{double3} On a un isomorphisme canonique de $\mathrm{GL}(V)_{\mathrm{td}}$-$\mathrm{H}^{*}V_{\mathrm{tf}}$-$\mathrm{A}$-modules instables~:
$$
\mathrm{M}(V,2h)
\hspace{4pt}\cong\hspace{4pt}
\mathrm{E}\Phi\hspace{1pt}(\mathrm{M}(V,h))
$$
pour tout entier naturel $h$.
\end{cor}

\medskip
\textit{Démonstration.} On prend $M=\mathrm{c}_{V}^{h}\hspace{1pt}\mathrm{H}^{*}V$, $p=n$, dans l'énoncé précédent et l'on observe que l'homomorphisme $\mathrm{E}\lambda :\mathrm{E}\Phi(\mathrm{c}_{V}^{h}\hspace{1pt}\mathrm{H}^{*}V)\to\mathrm{c}_{V}^{h}\hspace{1pt}\mathrm{H}^{*}V$ induit un isomorphisme de $\mathrm{GL}(V)_{\mathrm{td}}$-$\mathrm{H}^{*}V_{\mathrm{tf}}$-$\mathrm{A}$-modules instables $\mathrm{E}\Phi(\mathrm{c}_{V}^{h}\hspace{1pt}\mathrm{H}^{*}V)\cong\mathrm{c}_{V}^{2h}\hspace{1pt}\mathrm{H}^{*}V$.
\hfill$\square$

\pagebreak

\sect{Modules de Steinberg et algèbre homologique}

L'objectif de cette section est double~:

\smallskip
-- Dégager diverses connexions entre notre travail et l'article \cite{Ol} de Bob Oliver.

\smallskip
-- Expliciter la construction du bicomplexe $\mathrm{B}^{\bullet,\bullet}M$ évoqué dans la proposition~\ref{bicomplexe}, bicomplexe qui joue un rôle crucial dans la démonstration de l'implication (iii)$\Rightarrow$(i) du théorème \ref{genalg} (et par ricochet dans celle  de l'implication (iii)$\Rightarrow$(i) du théorème \ref{gentop}).

\medskip

\subsect{Spécialisation dans notre contexte des résultats de \cite{Ol}}\label{speOl}

\medskip
Tout d'abord un commentaire~:

\medskip
\begin{com}\label{Bob} Soient $\mathcal{A}b$ la catégorie des groupes abéliens et $\mathrm{F}_{\mathbb{Z}}$ le foncteur de $\mathcal{W}_{0}$ dans $\mathcal{A}b$ défini par
$$
\mathrm{F}_{\mathbb{Z}}(W)
\hspace{4pt}=\hspace{4pt}
\begin{cases}
\mathbb{Z} & \text{pour $W=V$,} \\
0 & \text{pour $W\not=V$.}
\end{cases}
$$
Bob Oliver signale dans l'introduction de \cite{Ol} qu'il est facile de déduire des isomorphismes
$$
\hspace{24pt}
{\lim}_{\mathcal{W}_{0}}^{*}\hspace{1pt}\mathrm{F}_{\mathbb{Z}}
\hspace{4pt}\cong\hspace{4pt}
\mathrm{Ext}_{\mathbb{Z}}^{*}(\Delta_{\mathbb{Z}},\mathrm{F}_{\mathbb{Z}})
\hspace{4pt}\cong\hspace{4pt}
\mathrm{H}^{*}\hspace{1pt}\mathrm{Hom}_{\mathcal{A}b^{\mathcal{W}_{0}}}(\mathrm{P}_{\bullet},\mathrm{F}_{\mathbb{Z}})
\hspace{24pt},
$$
$\mathrm{P}_{\bullet}$ désignant une résolution projective ``standard'' de $\Delta_{\mathbb{Z}}$ dans $\mathcal{A}b^{\mathcal{W}_{0}}$, que les groupes abéliens $\lim_{\mathcal{W}_{0}}^{k}\mathrm{F}_{\mathbb{Z}}$ sont isomorphes, de façon naturelle, aux groupes de cohomologie $\widetilde{\mathrm{H}}^{k-1}(\Vert\mathcal{W}_{0,n}\Vert;\mathbb{Z})$ et  en particulier que $\lim_{\mathcal{W}_{0}}^{n-1}\mathrm{F}_{\mathbb{Z}}$ est isomorphe au dual du $\mathbb{Z}$-module de Steinberg de $\mathrm{GL}(V)$.

\smallskip
La catégorie $V_{\mathrm{tf}}\text{-}\mathcal{U}$ n'ayant pas assez de projectifs (voir \ref{proj}) nous avons choisi en \ref{rappelslim} de présenter une variante (fort classique) de ce qui précède à savoir une définition des foncteurs dérivés de $\lim$ comme groupe de cohomologie d'un complexe ``combinatoire''~$\mathrm{L}^{\bullet}(-)$. Dans le cas  où l'on considère des foncteurs de $\mathcal{W}_{0}$ dans $\mathcal{A}b$ on a
$\mathrm{L}^{\bullet}(-)=\mathrm{Hom}_{\mathcal{A}b^{\mathcal{W}_{0}}}(\mathrm{P}_{\bullet},-)$.
\end{com}

\medskip
Ce commentaire étant fait, nous expliquons ci-après ce que donne le résultat principal de \cite{Ol} dans notre contexte.

\smallskip
On spécialise la proposition 5 de \cite{Ol} en prenant $p=2$, $\mathcal{A}=\mathcal{A}_{2}(V)(:=\mathcal{W}_{0})$. Soit $F$ un foncteur de $\mathcal{W}_{0}$ dans $\mathcal{E}$ (catégorie des $\mathbb{F}_{2}$-espaces vectoriels)~; cette proposition dit que l'on dispose d'un  $\mathcal{E}$-complexe de la forme
$$
\mathrm{C}^{0}_{\mathrm{St}}(F)\to\mathrm{C}^{1}_{\mathrm{St}}(F)\to
\ldots\to\mathrm{C}^{k}_{\mathrm{St}}(F)
\overset{\mathrm{d}^{k}}{\to}
\mathrm{C}^{k+1}_{\mathrm{St}}(F)\to
\ldots\to\mathrm{C}^{n-1}_{\mathrm{St}}(F)
$$
dont le $k$-ième groupe de cohomologie est isomorphe à $\lim_{\mathcal{W}_{0}}^{k}F$. Précisons~: On a
$$
\mathrm{C}^{k}_{\mathrm{St}}(F)
\hspace{4pt}=\hspace{4pt}
\prod_{\dim W=k+1}
\mathrm{St}_{W}^{*}\otimes F(W)
\hspace{4pt}=\hspace{4pt}
\bigoplus_{\dim W=k+1}
\mathrm{St}_{W}^{*}\otimes F(W)
$$
(la notation $\mathrm{St}_{W}^{*}$ désigne ci-dessus le $\mathbb{F}_{2}$-espace vectoriel dual de l'espace de la représentation de Steinberg modulo $2$ de $\mathrm{GL}(W)$, le produit et la somme directe sont indexées par l'ensemble des sous-groupes $W\subset V$ de dimension $k+1$)~; le cobord
$$
\mathrm{d}^{k}:
\bigoplus_{\dim W=k+1}
\mathrm{St}_{W}^{*}\otimes F(W)
\longrightarrow
\bigoplus_{\dim W'=k+2}
\mathrm{St}_{W'}^{*}\otimes F(W')
$$
a pour matrice $[\hspace{1pt}\mathrm{d}_{W',W}^{k}:\mathrm{St}_{W}^{*}\otimes F(W)\to\mathrm{St}_{W'}^{*}\otimes F(W')\hspace{1pt}]$~avec
$$
\mathrm{d}_{W',W}^{k}=
\begin{cases}
(\mathrm{r}_{W',W})^{*}\otimes F(W\subset W') & \text{pour $W\subset W'$}  \\
0 & \text{pour $W\not\subset W'$} ,
\end{cases}
$$
$\mathrm{r}_{W',W}$ désignant l'homomorphisme canonique de $\mathrm{St}_{W'}$ dans $\mathrm{St}_{W}$ (pour une définition de cet homomorphisme on pourra se reporter à \ref{rapSt}).

\bigskip
\begin{exple}\label{Lusztig} On suppose $V\not=0$ et on prend pour $F$ le foncteur $\Delta_{\mathbb{F}_{2}}$, à savoir le foncteur constant à valeur $\mathbb{F}_{2}$. On a $\lim_{\mathcal{W}_{0}}\Delta_{\mathbb{F}_{2}}=\mathbb{F}_{2}$ et $\lim_{\mathcal{W}_{0}}^{k}\Delta_{\mathbb{F}_{2}}=0$ pour $k>0$~; en effet on a $\Delta_{\mathbb{F}_{2}}=\mathbb{F}_{2}^{\mathrm{Hom}_{\mathcal{W}_{0}}(-,V)}$, égalité qui montre que $\Delta_{\mathbb{F}_{2}}$ et un objet injectif de $\mathcal{E}^{\mathcal{W}_{0}}$. On constate donc que l'on dispose d'un complexe acyclique, que nous notons $\mathrm{Lu}^{\bullet}=(\mathrm{Lu}^{0}\to\mathrm{Lu}^{1}\to\ldots\to\mathrm{Lu}^{n})$, avec $\mathrm{Lu}^{k}=\bigoplus_{\dim W=k}
\mathrm{St}_{W}^{*}$ (on convient que l'on a $\mathrm{St}_{0}=\mathbb{F}_{2}$), le cobord $\mathrm{d}^{k}:\bigoplus_{\dim W=k}
\mathrm{St}_{W}^{*}\to\bigoplus_{\dim W'=k+1}
\mathrm{St}_{W}^{*}$ étant donné par la matrice d'homomorphismes $[\hspace{1pt}\mathrm{d}_{W',W}^{k}\hspace{1pt}]$ avec
$$
\mathrm{d}_{W',W}^{k}=
\begin{cases}
(\mathrm{r}_{W',W})^{*} & \text{pour $W\subset W'$}  \\
0 & \text{pour $W\not\subset W'$} 
\end{cases}
$$
(on convient que l'on a $\mathrm{r}_{U,0}=\mathrm{id}_{\mathbb{F}_{2}}$ pour $\dim U=1$). Le complexe $\mathrm{Lu}^{\bullet}$ est le dual du complexe (b) de la page 22 de \cite{Lu} (avec $A=\mathbb{F}_{2}$) dont nous reparlerons en \ref{rapSt}~; ceci explique notre notation.
 \end{exple}

\bigskip
Soit  maintenant $F$ un foncteur de $\mathcal{W}_{0}$ dans $V_{\mathrm{tf}}\text{-}\mathcal{U}$. Par les mêmes formules que ci-dessus on définit un complexe $\mathrm{C}_{\mathrm{St}}^{\bullet}(F)$ dans la catégorie $V_{\mathrm{tf}}\text{-}\mathcal{U}$. On se convainc aisément que  $\lim_{\mathcal{W}_{0}}^{k}F$ est isomorphe à $\mathrm{H}^{k}\hspace{1pt}\mathrm{C}_{\mathrm{St}}^{\bullet}(F)$ comme $\mathrm{H}^{*}V_{\mathrm{tf}}$-$\mathrm{A}$-module instable.

\bigskip
Soit $M$ un $\mathrm{H}^{*}V_{\mathrm{tf}}$-$\mathrm{A}$-module instable. Le complexe $\mathrm{C}_{\mathrm{St}}^{\bullet}(\Psi_{M})$ est muni d'une coaugmentation ``tautologique'' $\lim_{\mathcal{W}_{0}}\Psi_{M}\to\mathrm{C}_{\mathrm{St}}^{\bullet}(\Psi_{M})$~; en composant à la source par l'homomorphisme $\rho_{M}$ de \ref{derivePf} on obtient une nouvelle coaugmentation que l'on note $\mathrm{d}^{-1}:M\to\mathrm{C}_{\mathrm{St}}^{\bullet}(\Psi_{M})$. Le complexe coaugmenté associé
$$
M=:\mathrm{C}_{\mathrm{St}}^{-1}(\Psi_{M})
\overset{\mathrm{d}^{-1}}{\to}
\mathrm{C}_{\mathrm{St}}^{0}(\Psi_{M})
\to\mathrm{C}_{\mathrm{St}}^{1}(\Psi_{M})
\to\ldots\to
\mathrm{C}_{\mathrm{St}}^{n-1}(\Psi_{M})
$$
est noté $\widetilde{\mathrm{C}}_{\mathrm{St}}^{\bullet}(\Psi_{M})$. Compte tenu de ce qui précède, la proposition suivante est une simple reformulation de la proposition \ref{derivePf}~:

\medskip
\begin{pro}\label{derivePf-bis} Soit $M$ un $\mathrm{H}^{*}V_{\mathrm{tf}}$-$\mathrm{A}$-module instable~; on a pour tout entier $k\geq 0$ un isomorphisme de $\mathrm{H}^{*}V_{\mathrm{tf}}$-$\mathrm{A}$-modules instables~:
$$
\hspace{24pt}
\mathrm{R}^{k}\mathrm{Pf}\hspace{1pt}M
\hspace{4pt}\cong\hspace{4pt}
\mathrm{H}^{k}\hspace{1pt}(\widetilde{\mathrm{C}}_{\mathrm{St}}^{\bullet}(\Psi_{M})[1])
\hspace{24pt}.
$$
{\em Décodons la notation. Le $k$-ième terme de $\widetilde{\mathrm{C}}_{\mathrm{St}}^{\bullet}(\Psi_{M})[1]$ est le $(k-1)$-ième terme de $\widetilde{\mathrm{C}}_{\mathrm{St}}^{\bullet}(\Psi_{M})$~: on a ``décalé'' $\widetilde{\mathrm{C}}_{\mathrm{St}}^{\bullet}(\Psi_{M})$ ``de $1$ vers la droite''.}

\end{pro}

\medskip
\begin{exple}\label{resproj} Cet exemple fait suite à l'exemple \ref{Lusztig}.

\medskip
Soit $N$ un $\mathrm{A}$-module instable fini. La proposition \ref{proclef} (avec $U=V$), ou le point (b) du corollaire \ref{suspensionEFix}, implique que $\Psi_{\mathrm{H}^{*}V\otimes N}$ est le foncteur constant $\Delta_{\mathrm{H}^{*}V\otimes N}$~; on a donc $\widetilde{\mathrm{C}}_{\mathrm{St}}^{\bullet}(\Psi_{\mathrm{H}^{*}V\otimes N})[1]=\mathrm{Lu}^{\bullet}\otimes(\mathrm{H}^{*}V\otimes N)$. Ce complexe est acyclique ce qui est bien en accord avec les propositions \ref{derivePf-bis} et \ref{annderivePf3}.

\medskip
Considérons plus généralement un $\mathrm{H}^{*}V$-$\mathrm{A}$-module instable $M$ qui est libre de dimension finie comme $\mathrm{H}^{*}V$-module. La théorie de Smith algébrique dit que le $\mathrm{A}$-module instable $\mathrm{Fix}_{V}M$ est fini et que l'unité d'adjonction $\eta:M\to\mathrm{H}^{*}V\otimes\mathrm{Fix}_{V}M$ est injective ($\eta$ s'identifie d'ailleurs au cobord $\mathrm{d}^{-1}$ du complexe $\widetilde{\mathrm{C}}^{\bullet}M$ qui est acyclique d'après \ref{pendantalg}). Puisque le foncteur $\Psi:V_{\mathrm{tf}}\text{-}\mathcal{U}\to(V_{\mathrm{tf}}\text{-}\mathcal{U})^{\mathcal{W}_{0}}$ est exact, $\eta$ induit un monomorphisme $\Psi_{M}\hookrightarrow\Psi_{\mathrm{H}^{*}V\otimes\mathrm{Fix}_{V}M}$ si bien que le complexe $\widetilde{\mathrm{C}}_{\mathrm{St}}^{\bullet}(\Psi_{M})[1]$ s'identifie à un sous-complexe de $\mathrm{Lu}^{\bullet}\otimes(\mathrm{H}^{*}V\otimes\mathrm{Fix}_{V}M)$.

\smallskip
On observe que tous les termes du complexe $\widetilde{\mathrm{C}}_{\mathrm{St}}^{\bullet}(\Psi_{M})[1]$ sont des $\mathrm{H}^{*}V$-$\mathrm{A}$-modules instables qui sont libres de dimension finie comme $\mathrm{H}^{*}V$-modules. En effet le $k$-ième terme de ce complexe, $0\leq k\leq n$, est par définition une somme directe finie de $\mathrm{H}^{*}V$-$\mathrm{A}$-modules instables de la forme $\mathrm{EFix}_{(V,W)}M$ avec $\dim W=k$ et $\mathrm{EFix}_{(V,W)}M\cong\mathrm{H}^{*}V\otimes_{\mathrm{H}^{*}V/W}\mathrm{Fix}_{(V,W)}M$ est libre comme $\mathrm{H}^{*}V$-module d'après~\ref{libre} (et de type fini comme  $\mathrm{H}^{*}V$-module d'après \ref{typefiniFix}).

\pagebreak

\smallskip
La proposition \ref{derivePf-bis} et le corollaire \ref{DS3} montrent que $\mathrm{H}^{k}\hspace{1pt}(\widetilde{\mathrm{C}}_{\mathrm{St}}^{\bullet}(\Psi_{M})[1])$ est nul pour $0\leq k\leq n-1$ et isomorphe à $\mathrm{R}^{n}\mathrm{Pf}\hspace{1pt}M$ pour $k=n$~; on note $\varepsilon$ l'homomorphisme de $\mathrm{H}^{*}V_{\mathrm{tf}}$-$\mathrm{A}$-modules instables de $\mathrm{St}_{V}^{*}\otimes(\mathrm{H}^{*}V\otimes\mathrm{Fix}_{V}M)$ (le~$n$-ième terme de $\widetilde{\mathrm{C}}_{\mathrm{St}}^{\bullet}(\Psi_{M})[1])$ dans $\mathrm{R}^{n}\mathrm{Pf}\hspace{1pt}M$ induit par cet isomorphisme. On dispose donc d'une suite exacte dans la catégorie $V_{\mathrm{tf}}\text{-}\mathcal{U}$ de la forme
\begin{multline*}
0\to M\to\bigoplus_{\dim W=1}
\mathrm{St}_{W}^{*}\otimes\mathrm{EFix}_{(V,W)}M\to\ldots\to
\bigoplus_{\dim W=k}
\mathrm{St}_{W}^{*}\otimes\mathrm{EFix}_{(V,W)}M \\
\to\ldots\to\bigoplus_{\dim W=n-1}
\mathrm{St}_{W}^{*}\otimes\mathrm{EFix}_{(V,W)}M\to
\mathrm{St}_{V}^{*}\otimes(\mathrm{H}^{*}V\otimes\mathrm{Fix}_{V}M)
\overset{\varepsilon}{\to}
\mathrm{R}^{n}\mathrm{Pf}\hspace{1pt}M
\end{multline*}
telle que la suite exacte sous-jacente dans la catégorie des $\mathrm{H}^{*}V$-modules de type fini est une résolution libre du $\mathrm{H}^{*}V$-module fini $\mathrm{R}^{n}\mathrm{Pf}\hspace{1pt}M$.

\smallskip
Voici une spécialisation concrète de ce qui précède~: on prend $M=\mathrm{c}_{V}\mathrm{H}^{*}V$. Le complexe $\widetilde{\mathrm{C}}_{\mathrm{St}}^{\bullet}(\Psi_{\mathrm{c}_{V}\mathrm{H}^{*}V})[1]$ est le sous-complexe de $\mathrm{Lu}^{\bullet}\otimes\mathrm{H}^{*}V$ dont le $k$-ième terme est $\bigoplus_{\dim W=k}
\mathrm{St}_{W}^{*}\otimes\mathrm{c}_{(V,W)}\mathrm{H}^{*}V$. On rappelle que nous avons posé $\mathrm{M}(V):=\mathrm{R}^{n}\mathrm{Pf}\hspace{1pt}(\mathrm{c}_{V}\mathrm{H}^{*}V)$. L'homomorphisme $\varepsilon:\mathrm{St}_{V}^{*}\otimes\mathrm{H}^{*}V\to\mathrm{M}(V)$ s'identifie à l'application $\mathrm{M}^{0}(V)\otimes\mathrm{H}^{*}V\to\mathrm{M}(V)$ induite par la structure de $\mathrm{H}^{*}V$-module de $\mathrm{M}(V)$. Comme $\varepsilon$ est $\mathrm{H}^{*}V$-linéaire il suffit de vérifier cette affirmation en degré zéro~: c'est essentiellement ce que nous avons fait au tout début de la démonstration du point (b) de la proposition \ref{MV}.
\end{exple}

\medskip
\subsect{Quelques rappels sur les modules de Steinberg}\label{rapSt}

\medskip
Soit $V$ un $\mathbb{Z}/2$-espace vectoriel de dimension finie.

\medskip
Soient $d$ et $c$ deux entiers naturels~; on note $\mathcal{B}_{d,c}(V)$ l'ensemble des sous-espaces $W$ de $V$ avec $\dim W\geq d$~et $\mathop{\mathrm{codim}}W\geq c$ ($d\leq\dim W\leq\dim V-c$). L'ensemble $\mathcal{B}_{d,c}(V)$ est ordonné par inclusion. On observera que l'on a en particulier $\mathcal{B}_{0,0}(V)=\mathcal{W}$ et $\mathcal{B}_{1,0}(V)=\mathcal{W}_{0}$.

\medskip
Soit $n$ la dimension de $V$~;  pour $n\geq 1$, on pose
$$
\hspace{24pt}
\mathrm{St}_{V}
\hspace{4pt}:=\hspace{4pt}
\widetilde{\mathrm{H}}_{n-2}\hspace{1pt}(\mathcal{B}_{1,1}(V);\mathbb{F}_{2})
\hspace{24pt}.
$$
La notation $\widetilde{\mathrm{H}}_{n-2}\hspace{1pt}(\mathcal{B}_{1,1}(V);\mathbb{F}_{2})$ désigne le $(n-2)$-ième groupe d'homologie réduite, à coefficients $\mathbb{F}_{2}$, du polyèdre associé à l'ensemble ordonné $\mathcal{B}_{1,1}(V)$. Tous les groupes d'homologie apparaissant ci-après seront à coefficients $\mathbb{F}_{2}$~; ce choix des coefficients est justifié par  les applications que nous avons en vue\ldots il nous permettra  par ailleurs de ne pas nous préoccuper de questions de signes~!

\medskip
Cette définition donne $\mathrm{St}_{V}=\mathbb{F}_{2}$ pour $n=1$. On convient que l'on a $\mathrm{St}_{V}=\mathbb{F}_{2}$ pour $n=0$ (cette convention sera justifiée en \ref{convSt0}).

\bigskip
On note respectivement $\mathrm{S}_{\bullet}I$ l'ensemble $\mathbb{N}$-gradué des simplexes  et  $\mathrm{C}_{\bullet}I$ le complexe des chaînes polyédrales, à coefficients $\mathbb{F}_{2}$, d'un ensemble ordonné~$I$ ($\mathrm{C}_{\bullet}I=\mathbb{F}_{2}[\mathrm{S}_{\bullet}I]$). On note $\widetilde{\mathrm{S}}_{\bullet}I$ et $\widetilde{\mathrm{C}}_{\bullet}I$ les versions augmentées ($\widetilde{\mathrm{S}}_{-1}I:=\{\mathrm{point}\}$). Soit $I'$ un sous-ensemble de $I$~; on pose  $\mathrm{C}_{\bullet}(I,I'):=\mathrm{C}_{\bullet}I/\mathrm{C}_{\bullet}I'$ (ou ce qui revient au même $\mathrm{C}_{\bullet}(I,I'):=\widetilde{\mathrm{C}}_{\bullet}I/\widetilde{\mathrm{C}}_{\bullet}I'=:\widetilde{\mathrm{C}}_{\bullet}(I,I')$). 

\bigskip
On rappelle maintenant la définition de l'homomorphisme $\mathrm{r}_{V,H}:\mathrm{St}_{V}\to\mathrm{St}_{H}$, $H$ hyperplan de $V$, évoqué en \ref{speOl}.\hspace{-1pt}1.

\medskip
On constate que l'on a
$$
\mathrm{S}_{\bullet}\mathcal{B}_{1,1}(V)-\mathrm{S}_{\bullet}\mathcal{B}_{1,2}(V)
\hspace{4pt}=\hspace{4pt}
\coprod_{\mathop{\mathrm{codim}}H=1}(\widetilde{\mathrm{S}}_{\bullet}\mathcal{B}_{1,1}(H))[-1]
$$
(le coproduit est indexé par l'ensemble des hyperplans $H$ de $V$, la notation $[-1]$ désigne le décalage de $1$ vers la droite). Il en résulte~:
$$
\hspace{24pt}
\mathrm{C}_{\bullet}(\mathcal{B}_{1,1}(V),\mathcal{B}_{1,2}(V))
\hspace{4pt}\cong\hspace{4pt}
\bigoplus_{\mathop{\mathrm{codim}}H=1}(\widetilde{\mathrm{C}}_{\bullet}\mathcal{B}_{1,1}(H))[-1]
\hspace{24pt}.
$$
La composition
$$
\widetilde{\mathrm{H}}_{n-2}\hspace{1pt}\mathcal{B}_{1,1}(V)\to\mathrm{H}_{n-2}\hspace{1pt}(\mathcal{B}_{1,1}(V),\mathcal{B}_{1,2}(V))
\cong
\bigoplus_{\mathop{\mathrm{codim}}H=1}\widetilde{\mathrm{H}}_{n-3}\hspace{1pt}\mathcal{B}_{1,1}(H) 
$$
donne un homomorphisme $\mathrm{St}_{V}\to\bigoplus_{\mathop{\mathrm{codim}}H=1}\mathrm{St}_{H}$. Le composé de cet homomorphisme et de la projection sur le facteur indexé par $H$ est noté $\mathrm{r}_{V,H}$.

\bigskip
\textit{Remarque.} Dans le cas $\dim V=2$, l'homomorphisme $\mathrm{St}_{V}\to\bigoplus_{\mathop{\mathrm{codim}}H=1}\mathrm{St}_{H}$ s'identifie à l'homomorphisme $\widetilde{\mathrm{H}}_{n-2}\hspace{1pt}\mathcal{B}_{1,1}(V)\to\mathrm{H}_{n-2}\hspace{1pt}\mathcal{B}_{1,1}(V).$

\bigskip
\textit{Convention.} Dans le cas $\dim V=1$,  on convient de prendre $\mathrm{r}_{V,H}=\mathrm{id}_{\mathbb{F}_{2}}$.

\bigskip
On définit pareillement un homomorphisme $\mathrm{St}_{V}\to\mathrm{St}_{V/D}$ pour $D$ une droite de $V$. On constate que l'on a
$$
\mathrm{S}_{\bullet}\mathcal{B}_{1,1}(V)-\mathrm{S}_{\bullet}\mathcal{B}_{2,1}(V)
\hspace{4pt}=\hspace{4pt}
\coprod_{\dim D=1}(\widetilde{\mathrm{S}}_{\bullet}\mathcal{B}_{1,1}(V/D))[-1]
$$
(le coproduit est indexé par l'ensemble des droites $D$ de $E$). Il en résulte~:
$$
\hspace{24pt}
\mathrm{C}_{\bullet}(\mathcal{B}_{1,1}(V),\mathcal{B}_{2,1}(V))
\hspace{4pt}\cong\hspace{4pt}
\bigoplus_{\dim D=1}(\widetilde{\mathrm{C}}_{\bullet}\mathcal{B}_{1,1}(V/D))[-1]
\hspace{24pt}.
$$
La composition
$$
\widetilde{\mathrm{H}}_{n-2}\hspace{1pt}\mathcal{B}_{1,1}(V)\to\mathrm{H}_{n-2}\hspace{1pt}(\mathcal{B}_{1,1}(V),\mathcal{B}_{2,1}(V))
\cong
\bigoplus_{\dim D=1}\widetilde{\mathrm{H}}_{n-3}\hspace{1pt}\mathcal{B}_{1,1}(V/D) 
$$
donne un homomorphisme $\mathrm{St}_{V}\to\bigoplus_{\dim D=1}\mathrm{St}_{V/D}$. Le composé de cet homomorphisme et de la projection sur le facteur indexé par $D$ est noté $\mathrm{s}_{V,V/D}$.

\bigskip
\textit{Remarque.} Dans le cas $\dim V=2$, l'homomorphisme $\mathrm{St}_{V}\to\bigoplus_{\dim D=1}\mathrm{St}_{V/D}$ s'identifie à l'homomorphisme $\widetilde{\mathrm{H}}_{n-2}\hspace{1pt}\mathcal{B}_{1,1}(V)\to\mathrm{H}_{n-2}\hspace{1pt}\mathcal{B}_{1,1}(V).$

\bigskip
\textit{Convention.} Dans le cas $\dim V=1$,  on convient de prendre $\mathrm{s}_{V,V/D}=\mathrm{id}_{\mathbb{F}_{2}}$.

\bigskip
\begin{pro}\label{Stcom}Soient respectivement $D$ et $H$ une droite et un hyperplan de $V$ avec $D\subset H$. Alors le diagramme
$$
\begin{CD}
\mathrm{St}_{V}@>\mathrm{r}_{V,H}>>\mathrm{St}_{H} \\
@V\mathrm{s}_{V,V/D}VV @V\mathrm{s}_{H,H/D}VV \\
\mathrm{St}_{V/D}@>\mathrm{r}_{V/D,V/D}>>\mathrm{St}_{V/D}
\end{CD}
$$
est commutatif.
\end{pro}

\bigskip
\textit{Démonstration.} Par hypothèse, on a $\dim V\geq 2$. Le cas $\dim V=2$ découle des remarques et conventions ci-dessus. On suppose $\dim V>2$. On a
\begin{multline*}
\hspace{4pt}\mathrm{S}_{\bullet}\mathcal{B}_{1,1}(V)-(\mathrm{S}_{\bullet}\mathcal{B}_{1,2}(V)\cup\mathrm{S}_{\bullet}\mathcal{B}_{2,1}(V))=\\
(\mathrm{S}_{\bullet}\mathcal{B}_{1,1}(V)-\mathrm{S}_{\bullet}\mathcal{B}_{1,2}(V))\cap(\mathrm{S}_{\bullet}\mathcal{B}_{1,1}(V)-\mathrm{S}_{\bullet}\mathcal{B}_{2,1}(V))=\\
(\coprod_{\mathop{\mathrm{codim}}H=1} \{\sigma;\sup\sigma=H\})\cap
(\coprod_{\dim D=1} \{\sigma;\inf\sigma=D\})=\\
\coprod_{D\subset H} \{\sigma;\inf\sigma=D, \sup\sigma=H\}\hspace{4pt}.
\end{multline*}
Il en résulte~:
$$
\hspace{24pt}
\mathrm{C}_{\bullet}(\mathcal{B}_{1,1}(V),\mathcal{B}_{1,2}(V))\cup\mathcal{B}_{2,1}(V))
\hspace{4pt}\cong\hspace{4pt}
\bigoplus_{D\subset H}(\widetilde{\mathrm{C}}_{\bullet}\mathcal{B}_{1,1}(H/D))[-2]
\hspace{24pt}.
$$
La commutativité du diagramme de complexes
$$
\begin{CD}
\widetilde{\mathrm{C}}_{\bullet}\mathcal{B}_{1,1}(V)@>>>\mathrm{C}_{\bullet}(\mathcal{B}_{1,1}(V),\mathcal{B}_{1,2}(V)) \\
@VVV @VVV \\
\mathrm{C}_{\bullet}(\mathcal{B}_{1,1}(V),\mathcal{B}_{2,1}(V))@>>>\mathrm{C}_{\bullet}(\mathcal{B}_{1,1}(V),\mathcal{B}_{1,2}(V))\cup\mathcal{B}_{2,1}(V))
\end{CD}
$$
permet de conclure.
\hfill$\square$

\pagebreak

\vspace{0.75cm}
\textit{Dualité}

\bigskip
L'application $W\mapsto W^{\perp}$ fournit un isomorphisme canonique d'ensembles ordonnés $\mathcal{B}_{d,c}(V)\to\mathcal{B}_{c,d}^{\hspace{1pt}\mathrm{op}}(V^{*})$. On en déduit un isomorphisme de $\mathbb{F}_{2}$-espaces vectoriels, disons $\delta:\mathrm{St}_{V}\overset{\cong}{\longrightarrow}\mathrm{St}_{V^{*}}$ ($\delta$ est $\mathrm{GL}(V)$-équivariant si l'on fait agir $\mathrm{GL}(V)$ sur $\mathrm{St}_{V^{*}}$ \textit{via} l'isomorphisme $\mathrm{GL}(V)\to\mathrm{GL}(V^{*}),\alpha\mapsto(\alpha^{*})^{-1}$).

\bigskip
La vérification de la proposition suivante est laissée au lecteur~:

\begin{pro}\label{dualrs} Soit $D$ une droite de $E$~; le diagramme suivant
$$
\begin{CD}
\mathrm{St}_{E}@>\mathrm{s}_{E,E/D}>>\mathrm{St}_{E/D} \\
@V\delta VV @V\delta VV \\
\mathrm{St}_{E^{*}}@>\mathrm{r}_{E,D^{\perp}}>>\mathrm{St}_{D^{\perp}}
\end{CD}
$$
est commutatif (on identifie l'hyperplan $D^{\perp}$ de $E^{*}$ à $(E/D)^{*}$).
\end{pro}

\vspace{0.75cm}
\textit{Retour sur le complexe de Lusztig}

\bigskip
On commence par observer que $\Vert\mathcal{B}_{1,0}(V)\Vert$ (le polyèdre associé à l'ensemble ordonné $\mathcal{B}_{1,0}(V)$) est contractile. En effet, puisque $\mathcal{B}_{1,0}(V)$ possède un plus grand élément, à savoir $V$, $\Vert\mathcal{B}_{1,0}(V)\Vert$ est un cône de sommet le $0$-simplexe $\{V\}$ et de base $\Vert\mathcal{B}_{1,1}(V)\Vert$ ($\mathcal{B}_{1,0}(V)-\{V\}=\mathcal{B}_{1,1}(V)$).

\bigskip
On considère ensuite la filtration de $\Vert\mathcal{B}_{1,0}(V)\Vert$ induite par la filtration croissante $(\mathrm{F}_{p}\hspace{1pt}\mathcal{B}_{1,0}(V))_{-1\leq p\leq n}$ de $\mathcal{B}_{1,0}(V)$ définie par $\mathrm{F}_{p}\hspace{1pt}\mathcal{B}_{1,0}(V)=\mathcal{B}_{1,n-p}(V)$ (on a donc $\mathrm{F}_{p}\hspace{1pt}\mathcal{B}_{1,0}(V)=\emptyset$ pour $p=-1,0$ et $\mathrm{F}_{n}\hspace{1pt}\mathcal{B}_{1,0}(V)=\mathcal{B}_{1,0}(V)$) et la suite spectrale associée $(\mathrm{E}_{p,q}^{r})_{r\geq 1}$ en homologie modulo $2$. On a
$$
\hspace{24pt}
\mathrm{E}_{p,q}^{1}
\hspace{4pt}=\hspace{4pt}
\mathrm{H}_{p+q}(\mathcal{B}_{1,n-p}(V),\mathcal{B}_{1,n+1-p}(V))
\hspace{24pt}.
$$
L'égalité (pour $p\geq 1$)
$$
\mathrm{S}_{\bullet}\mathcal{B}_{1,n-p}(V)-\mathrm{S}_{\bullet}\mathcal{B}_{1,n+1-p}(V)
\hspace{4pt}=\hspace{4pt}
\coprod_{\dim W=p}(\widetilde{\mathrm{S}}_{\bullet}\mathcal{B}_{1,1}(W))[-1]
$$
implique
$$
\mathrm{E}_{p,q}^{1}
\hspace{4pt}\cong\hspace{4pt}
\begin{cases}
\bigoplus_{\dim W=p} \mathrm{St}_{W} & \text{pour $p\geq 1$ et $q=-1$}, \\
0 & \text{sinon}.
\end{cases}
$$

\pagebreak

\medskip
\footnotesize
\begin{rem}\label{Folkman} On utilise ici que l'on a $\widetilde{\mathrm{H}}_{k}\hspace{1pt}\mathcal{B}_{1,1}(V)=0$ pour $k<n-2$ ($n:=\dim V\geq 3$). Pour le confort du lecteur, nous reproduisons la démonstration de cette propriété, que l'on trouve à la  page 12 de \cite{Lu}, par une méthode que Lusztig attribue à Jon Folkman~\cite{Fo}. On considère le ``recouvrement'', disons $\mathcal{F}$, de $\mathcal{B}_{1,1}(V)$ par les $\mathcal{B}_{1,0}(H)$, $H$ décrivant l'ensemble, disons~$\mathcal{H}$, des hyperplans de~$V$. Soit $\mathcal{P}$ une partie non vide de $\mathcal{H}$~; l'intersection $\bigcap_{H\in\mathcal{P}}\mathcal{B}_{1,0}(H)$ possède un plus grand élément (si $\bigcap_{H\in\mathcal{P}}H\not=0$) ou est vide (si $\bigcap_{H\in\mathcal{P}}H=0$) si bien que $\mathcal{B}_{1,1}(V)$ a le type d'homotopie du nerf de $\mathcal{F}$, disons $\mathrm{N}(\mathcal{F})$. Si le cardinal de $\mathcal{P}$ est strictement inférieur à $n$ alors $\bigcap_{H\in\mathcal{P}}\mathcal{B}_{1,0}(H)$ est non vide~;  le $(n-2)$-squelette de $\mathrm{N}(\mathcal{F})$ coïncide donc avec le $(n-2)$-squelette du simplexe ``standard'' sur l'ensemble $\mathcal{H}$. L'égalité $\widetilde{\mathrm{H}}_{k}\hspace{1pt}\mathrm{N}(\mathcal{F})=0$ pour $k<n-2$ en résulte.
\end{rem}
\normalsize

\medskip
L'isomorphisme précédant la remarque et le fait que l'aboutissement de la suite spectrale est l'homologie modulo $2$ d'un espace contractile montrent que la suite
$$
0\to\mathrm{E}_{n,-1}^{1}\overset{\mathrm{d}^{1}}{\to}\mathrm{E}_{n-1,-1}^{1}\overset{\mathrm{d}^{1}}{\to}\ldots\overset{\mathrm{d}^{1}}{\to}
\mathrm{E}_{p,-1}^{1}\overset{\mathrm{d}^{1}}{\to}\ldots\overset{\mathrm{d}^{1}}{\to}\mathrm{E}_{1,-1}^{1}
$$
est exacte et que le conoyau de la dernière flèche est isomorphe à $\mathbb{F}_{2}$. On dispose donc d'un complexe, acyclique pour $V\not=0$, {\em le complexe de Lustig} \cite[p. 22]{Lu},
$$
\mathrm{Lu}_{\bullet}
\hspace{4pt}:=\hspace{4pt}
(\hspace{2pt}\mathrm{Lu}_{n}\overset{\mathrm{d}_{n}}{\to}\mathrm{Lu}_{n-1}\ldots\to\mathrm{Lu}_{p}\overset{\mathrm{d}_{p}}{\to}\mathrm{Lu}_{p-1}\to\ldots\to\mathrm{Lu}_{1}\overset{\mathrm{d}_{1}}{\to}\mathrm{Lu}_{0}\hspace{2pt})
$$
avec $\mathrm{Lu}_{p}=\bigoplus_{\dim W=p}\mathrm{St}_{W}$ (on rappelle que l'on a posé $\mathrm{St}_{0}=\mathbb{F}_{2}$) dont les opérateurs de bord sont explicités ci-dessous.

\medskip
\begin{pro}\label{bordsLu} Le $p$-ième opérateur de bord du complexe de Lusztig
$$
\mathrm{d}_{p}:
\mathrm{Lu}_{p}:=\bigoplus_{\dim W=p}\mathrm{St}_{W}\longrightarrow
\bigoplus_{\dim W'=p-1}\mathrm{St}_{W'}=:\mathrm{Lu}_{p-1}
$$
est donné par
$$
\pi_{W'}\circ\mathrm{d}_{p}\circ\iota_{W}
\hspace{4pt}=\hspace{4pt}
\begin{cases}
\mathrm{r}_{W,W'} & \text{pour\hspace{6pt}$W'\subset W$}, \\
0 & \text{pour\hspace{6pt}$W'\not\subset W$,}
\end{cases}
$$
$\iota_{W}$ et $\pi_{W'}$ désignant respectivement l'inclusion de $\mathrm{St}_{W}$ dans $\mathrm{Lu}_{p}$ et la projection de $\mathrm{Lu}_{p-1}$ sur $\mathrm{St}_{W'}$.
\end{pro}

\medskip
\textit{Démonstration.} On traite d'abord le cas $p=n$. L'opérateur de bord $\mathrm{d}_{n}$ est le composé
$$
\hspace{12pt}
\mathrm{H}_{n-1}(\mathcal{B}_{1,0}(V),\mathcal{B}_{1,1}(V))\overset{\partial}{\longrightarrow}
\widetilde{\mathrm{H}}_{n-2}\hspace{1pt}\mathcal{B}_{1,1}(V)
\longrightarrow
\mathrm{H}_{n-2}(\mathcal{B}_{1,1}(V),\mathcal{B}_{1,2}(V))
\hspace{12pt}.
$$
Or la flèche de gauche est un isomorphisme car $\Vert\mathcal{B}_{1,0}(V)\Vert$ est contractile et celle de droite est par définition le produit des $\mathrm{r}_{V,H}$, $H$ décrivant l'ensemble des hyperplans de $V$. On passe ensuite au cas général. Soit $W$ un sous-groupe de $V$ avec $\dim W=p$~; on considère l'application de triples évidente
$$
\hspace{9pt}
(\mathcal{B}_{1,0}(W),\mathcal{B}_{1,1}(W),\mathcal{B}_{1,2}(W))\to
(\mathcal{B}_{1,n-p}(V),\mathcal{B}_{1,n+1-p}(V),\mathcal{B}_{1,n+2-p}(V))
\hspace{9pt}.
$$
En remplaçant $V$ par $W$ dans l'argument que l'on vient de donner pour $\mathrm{d}_{n}$, on constate que $\mathrm{d}_{p}\circ\iota_{W}$ est le produit des $\mathrm{r}_{W,W'}$, $W'$ décrivant l'ensemble des hyperplans de $W$. 
\hfill$\square$

\bigskip
Quand cela nous paraîtra utile nous préciserons la notation $\mathrm{Lu}_{\bullet}$ en $\mathrm{Lu}_{\bullet,V}$.

\medskip
Nous exploitons maintenant l'isomorphisme $\mathcal{B}_{d,c}(V)\cong\mathcal{B}_{c,d}^{\hspace{1pt}\mathrm{op}}(V^{*})$, dont nous avons traité plus haut sous l'intertitre ``Dualité'', pour décrire $\mathrm{Lu}_{\bullet,V^{*}}$ en termes des sous-groupes $W$ de $V$.

\medskip
On a
$$
\mathrm{Lu}_{p,V^{*}}
\hspace{4pt}=\hspace{4pt}
\bigoplus_{\mathop{\mathrm{codim}}W=p}
\mathrm{St}_{W^{\perp}}
\hspace{4pt}\cong\hspace{4pt}
\bigoplus_{\mathop{\mathrm{codim}}W=p}\mathrm{St}_{V/W}
$$
(l'isomorphisme ci-dessus provient des isomorphismes $\delta:\mathrm{St}_{V/W}\overset{\cong}{\to}\mathrm{St}_{(V/W)^{*}}$ et $(V/W)^{*}\cong W^{\perp}$). On pose $\mathrm{Lu}_{p,2}:=\bigoplus_{\mathop{\mathrm{codim}}W=p}\mathrm{St}_{V/W}$, on dispose donc d'un {\em second complexe de Lusztig},
$$
\hspace{24pt}
\mathrm{Lu}_{\bullet,2}
\hspace{4pt}:=\hspace{4pt}
(\hspace{2pt}\mathrm{Lu}_{n,2}\to\ldots\to\mathrm{Lu}_{p,2}\overset{\mathrm{d}_{p}}{\to}\mathrm{Lu}_{p-1,2}\to\ldots\to\mathrm{Lu}_{0,2}\hspace{2pt})
\hspace{24pt}.
$$
Les propositions \ref{bordsLu} et \ref{dualrs} conduisent à la suivante~:

\medskip
\begin{pro}\label{bordsLu2} Le $p$-ième opérateur de bord du second complexe de Lusztig
$$
\mathrm{d}_{p}:
\mathrm{Lu}_{p,2}:=\bigoplus_{\mathop{\mathrm{codim}}W=p}\mathrm{St}_{V/W}
\longrightarrow
\bigoplus_{\mathop{\mathrm{codim}}W'=p-1}\mathrm{St}_{V/W'}=:\mathrm{Lu}_{p-1,2}
$$
est donné par
$$
\pi_{W'}\circ\mathrm{d}_{p}\circ\iota_{W}
\hspace{4pt}=\hspace{4pt}
\begin{cases}
\mathrm{s}_{W,W'} & \text{pour\hspace{6pt}$W'\subset W$}, \\
0 & \text{pour\hspace{6pt}$W'\not\subset W$,}
\end{cases}
$$
$\iota_{W}$ et $\pi_{W'}$ désignant respectivement l'inclusion de $\mathrm{St}_{V/W}$ dans $\mathrm{Lu}_{p,2}$ et la projection de $\mathrm{Lu}_{p-1,2}$ sur $\mathrm{St}_{V/W'}$.
\end{pro}

\bigskip
Les complexes de cochaînes $\mathrm{Hom}_{\mathbb{F}_{2}}(\mathrm{Lu}_{\bullet},\mathbb{F}_{2})$ et $\mathrm{Hom}_{\mathbb{F}_{2}}(\mathrm{Lu}_{\bullet,2},\mathbb{F}_{2})$ seront respectivement notés $\mathrm{Lu}^{\bullet}$  et $\mathrm{Lu}^{\bullet}_{2}$. Quand cela nous paraîtra utile voire nécessaire ces notations seront précisées par un $V$ en indice.

\medskip
\subsect{Algèbre homologique dans la catégorie abélienne $\mathcal{E}^{\mathcal{W}}$}\label{resinj}

On rappelle que l'on note respectivement $\mathcal{E}$ et $\mathcal{W}$, la catégorie des $\mathbb{F}_{2}$-espaces vectoriels et l'ensemble ordonné par inclusion des sous-groupes de $V$ (que l'on peut voir comme une catégorie).

\medskip
\footnotesize
\textit{Commentaires.} Dans cette sous-section on change de paradigme~: on remplace $\mathcal{E}^{\mathcal{W}_{0}}$ par $\mathcal{E}^{\mathcal{W}}$. Expliquons la raison de l'apparition de  $\mathcal{W}_{0}$ dans notre mémoire. Celui-ci est influencé par les travaux de Henn \cite{Hennalg}\cite{Henntop} eux-mêmes influencés par ceux de Jackowski et McClure \cite{JMcC}. Soient $G$ un groupe de Lie compact, $p$ un nombre premier fixé, $\widehat{\mathcal{A}}_{G}$ la catégorie (introduite par Quillen \cite{Qui}) dont les objets sont les $p$-groupes abéliens élémentaires $E\subset G$ et dont les morphismes sont les applications linéaires $E\to E'$ induite par une conjugation dans $G$ et enfin $\mathcal{A}_{G}$ la sous-catégorie pleine de $\widehat{\mathcal{A}}_{G}$ obtenue en excluant l'objet $0$~; Jackowski et McClure montrent que l'espace classifiant $\mathrm{B}G$, ``en $p$'', peut être obtenu comme une colimite homotopique, indexée par $\mathcal{A}_{G}$, des $\mathrm{B\hspace{0.5pt}C}(E)$, $\mathrm{C}(E)$ désignant le centralisateur de $E$ dans $G$. Comme l'on a $\mathrm{C}(0)=G$, le résultat en question manquerait de sel si $0$ n'était pas exclu~! L'article \cite{Ol} est un sous-produit de la collaboration d'Oliver avec Jackowski et McClure.
\normalsize

\medskip
Soit $F$ un foncteur de $\mathcal{W}$  dans $\mathcal{E}$. On définit une filtration décroissante de $F$
$$
F={\mathop{\mathrm{filt}}}_{1}^{0}F\supset{\mathop{\mathrm{filt}}}_{1}^{1}F\supset\ldots\supset{\mathop{\mathrm{filt}}}_{1}^{n}F\supset{\mathop{\mathrm{filt}}}_{1}^{n+1}F=0
$$
en prenant
$$
(\mathrm{filt}_{1}^{p}F)(U)
\hspace{4pt}:=\hspace{4pt}
\begin{cases}
F(U) & \text{pour}\hspace{4pt}\dim U\geq p, \\
0 & \text{pour}\hspace{4pt}\dim U<p.
\end{cases}
$$
On pose ${\mathop{\mathrm{Gr}}}_{1}^{p}F:={\mathop{\mathrm{filt}}}_{1}^{p}F/{\mathop{\mathrm{filt}}}_{1}^{p+1}F$.

\medskip
Soit $W$ un sous-groupe de $V$, on note $\mathrm{S}_{W}$ le foncteur de $\mathcal{W}$  dans $\mathcal{E}$ défini par
$$
\mathrm{S}_{W}(U)
\hspace{4pt}:=\hspace{4pt}
\begin{cases}
\mathbb{F}_{2} & \text{pour}\hspace{4pt}U=W, \\
0 & \text{pour}\hspace{4pt}U\not=W.
\end{cases}
$$
Il est clair que $\mathrm{S}_{W}$ est un objet simple de $\mathcal{E}^{\mathcal{W}}$. On constate que l'on a
$$
{\mathop{\mathrm{filt}}}_{1}^{p}F/{\mathop{\mathrm{filt}}}_{1}^{p+1}F
\hspace{4pt}\cong\hspace{4pt}
\bigoplus_{\dim W=p}
F(W)\otimes\mathrm{S}_{W}
$$
(si $E$ est un objet de $\mathcal{E}$ et $S$ un objet de $\mathcal{E}^{\mathcal{W}}$ alors la notation $E\otimes S$ désigne le foncteur $U\mapsto E\otimes S(U)$). Ceci implique en particulier que les $\mathrm{S}_{W}$ sont les seuls objets simples de $\mathcal{E}^{\mathcal{W}}$.

\medskip
On introduit maintenant les projectifs et injectifs ``tautologiques'' de $\mathcal{E}^{\mathcal{W}}$~:

\medskip
-- On note  $\mathbb{P}_{W}$ le foncteur $\mathcal{W}\to\mathcal{E},U\mapsto\mathbb{F}_{2}[\mathrm{Hom}_{\mathcal{W}}(W,U)]$~; le lemme de Yoneda donne $\mathrm{Hom}_{\mathcal{E}^{\mathcal{W}}}(\mathbb{P}_{W},F)=F(W)$, ce qui montre que $\mathbb{P}_{W}$ est un objet projectif de $\mathcal{E}^{\mathcal{W}}$. On observe que l'on dispose pour tout objet $F$ de $\mathcal{E}^{\mathcal{W}}$ d'un épimorphisme canonique $\bigoplus_{W\in\mathcal{W}}F(W)\otimes\mathbb{P}_{W}
\twoheadrightarrow F$~; ceci montre à la fois que~$\mathcal{E}^{\mathcal{W}}$ a assez de projectifs et qu'un projectif de $\mathcal{E}^{\mathcal{W}}$ est facteur direct d'une somme directe de $\mathbb{P}_{W}$.

\medskip
-- On note  $\mathbb{I}^{W}$ le foncteur $\mathcal{W}\to\mathcal{E},U\mapsto\mathbb{F}_{2}^{\mathrm{Hom}_{\mathcal{W}}(U,W)}$~; le lemme de Yoneda donne cette fois $\mathrm{Hom}_{\mathcal{E}^{\mathcal{W}}}(F,\mathbb{I}^{W})=(F(W))^{*}$, ce qui montre que $\mathbb{I}^{W}$ est un objet injectif de $\mathcal{E}^{\mathcal{W}}$. Soit $E$ un $\mathbb{F}_{2}$-espace vectoriel, on constate que le foncteur $E\otimes\mathbb{I}^{W}$ coïncide avec le foncteur $\mathcal{W}\to\mathcal{E},U\mapsto E^{\mathrm{Hom}_{\mathcal{W}}(U,W)}$~; on en déduit que l'on a $\mathrm{Hom}_{\mathcal{E}^{\mathcal{W}}}(F,E\otimes\mathbb{I}^{W})=\mathrm{Hom}_{\mathcal{E}}(F(W),E)$ et donc que $E\otimes\mathbb{I}^{W}$ est un objet injectif de $\mathcal{E}^{\mathcal{W}}$. On dispose pour tout objet $F$ de $\mathcal{E}^{\mathcal{W}}$ d'un monomorphisme canonique $F\hookrightarrow\prod_{W\in\mathcal{W}}F(W)\otimes\mathbb{I}^{W}=\bigoplus_{W\in\mathcal{W}}F(W)\otimes\mathbb{I}^{W}$~; il en résulte que~$\mathcal{E}^{\mathcal{W}}$ a assez d'injectifs et qu'un injectif de $\mathcal{E}^{\mathcal{W}}$ est facteur direct d'une somme directe de~$\mathbb{I}^{W}$.

\vspace{0.75cm}
\textit{Résolution projective de $\mathrm{S}_{0}$ dans la catégorie $\mathcal{E}^{\mathcal{W}}$}

\bigskip
On considère le $\mathcal{E}^{\mathcal{W}}$-épimorphisme canonique $\varepsilon:\mathbb{P}_{0}\to\mathrm{S}_{0}$ (qui est en fait une couverture projective)~; on se propose de ``prolonger'' $\varepsilon$ en une résolution projective.

\medskip
Soit $k$ un entier naturel~; on pose
$$
\mathrm{P}_{k}
\hspace{4pt}:=\hspace{4pt}
\bigoplus_{\dim Z=k}
\mathrm{St}_{Z}\otimes\mathbb{P}_{Z}
$$
et on note
$$
\mathrm{d}_{k+1}
\hspace{4pt}:\hspace{4pt}
\mathrm{P}_{k+1}
\hspace{4pt}=\hspace{4pt}
\bigoplus_{\dim Z'=k+1}
\mathrm{St}_{Z'}\otimes\mathbb{P}_{Z'}
\longrightarrow
\bigoplus_{\dim Z=k}
\mathrm{St}_{Z}\otimes\mathbb{P}_{Z}
\hspace{4pt}=\hspace{4pt}
\mathrm{P}_{k}
$$
 la transformation naturelle dont la composante d'indices $Z,Z'$ (convention ``matricielle''), disons ${\mathrm{d}_{k+1}}_{Z,Z'}$, est donnée par
 $$
{\mathrm{d}_{k+1}}_{Z,Z'}
\hspace{4pt}=\hspace{4pt}
\begin{cases}
\mathrm{r}_{Z',Z}\otimes\nu_{Z',Z} & \text{pour $Z\subset Z'$} \\
0 & \text{pour $Z\not\subset Z'$},
\end{cases}
 $$
 $\nu_{Z',Z}$ désignant l'unique élément non trivial de  $\mathrm{Hom}_{\mathcal{E}^{\mathcal{W}}}(\mathbb{P}_{Z'},\mathbb{P}_{Z})$.
 
 \medskip
 \begin{pro}\label{resprojS0} La suite de $\mathcal{E}^{\mathcal{W}}$-morphismes
 $$
 \hspace{24pt}
 0\longleftarrow\mathrm{S}_{0}
 \overset{\varepsilon}{\longleftarrow}
 \mathrm{P}_{0} \overset{\mathrm{d}_{1}}{\longleftarrow}
 \mathrm{P}_{1} \overset{\mathrm{d}_{2}}{\longleftarrow}
 \ldots\overset{\mathrm{d}_{n}}{\longleftarrow}
  \mathrm{P}_{n} 
  \hspace{24pt},
 $$
est une résolution projective de $\mathrm{S}_{0}$ dans la catégorie $\mathcal{E}^{\mathcal{W}}$.
\end{pro}
 
\medskip
 \textit{Démonstration.} Soit $U$ un sous-groupe non nul de $V$. On constate que
 $$
  \mathrm{P}_{0} (U)\overset{{\mathrm{d}_{1}}_{U}}{\longleftarrow}
 \mathrm{P}_{1}(U)\overset{{\mathrm{d}_{2}}_{U}}{\longleftarrow}
 \ldots\overset{{\mathrm{d}_{n}}_{U}}{\longleftarrow}
  \mathrm{P}_{n}(U)
  $$
est le ``complexe de Lusztig'' $\mathrm{Lu}_{\bullet,U}$ et que $\mathrm{S}_{0}(U)$ est nul~; or $\mathrm{Lu}_{\bullet,U}$ est acyclique. Dans le cas $U=0$, on a $\mathrm{P}_{k}(U)=0$ pour $k>0$ et $\varepsilon_{U}$ est un isomorphisme.
\hfill$\square$

\medskip
\begin{cor}\label{ExtS0SW} Soit $W$ un sous-groupe de $V$~; on a
$$
\mathrm{Ext}_{\mathcal{E}^{\mathcal{W}}}^{k}(\mathrm{S}_{0},\mathrm{S}_{W})
\hspace{4pt}=\hspace{4pt}
\begin{cases}
\mathrm{St}_{W}^{*} & \text{pour $k=\dim W$}, \\
0 & \text{pour $k\not=\dim W$}.
\end{cases}
$$
\end{cor}

\medskip
\textit{Démonstration.} On constate que l'on a $\mathrm{Hom}_{\mathcal{E}^{\mathcal{W}}}(\mathrm{P}_{k},\mathrm{S}_{W})=0$ pour $k\not=\dim W$ et $\mathrm{Hom}_{\mathcal{E}^{\mathcal{W}}}(\mathrm{P}_{k},\mathrm{S}_{W})=\mathrm{St}_{W}^{*}$ pour $k=\dim W$.
\hfill$\square$

\bigskip
\begin{rem}\label{convSt0} On a $\mathrm{Ext}_{\mathcal{E}^{\mathcal{W}}}^{0}(\mathrm{S}_{0},\mathrm{S}_{0})=\mathrm{Hom}_{\mathcal{E}^{\mathcal{W}}}(\mathrm{S}_{0},\mathrm{S}_{0})=\mathbb{F}_{2}$~; cette observation peut être vue comme une justification de notre convention $\mathrm{St}_{0}=\mathbb{F}_{2}$.
\end{rem}

\bigskip
On remplace maintenant le foncteur $\mathrm{S}_{W}$ de l'énoncé \ref{ExtS0SW} par un foncteur $F:\mathcal{W}\to\mathcal{E}$ arbitraire. On pose $\mathrm{C}^{\bullet}_{1}F:=\mathrm{Hom}_{\mathcal{E}^{\mathcal{W}}}(\mathrm{P}_{\bullet},F)$ ($\mathrm{P}_{\bullet}\to\mathrm{S}_{0}$ étant ici la résolution projective de \ref{resprojS0}). Explicitons ce complexe de cochaînes. On a
$$
\mathrm{C}^{k}_{1}F
\hspace{4pt}:=\hspace{4pt}
\bigoplus_{\dim Z=k}\mathrm{St}_{Z}^{*}\otimes F(Z)
$$
et le cobord
$$
\mathrm{d}^{k}
\hspace{2pt}:\hspace{2pt}
\mathrm{C}^{k}_{1}F
\hspace{2pt}=\hspace{2pt}
\bigoplus_{\dim Z=k}
\mathrm{St}_{Z}^{*}\otimes F(Z)
\longrightarrow
\bigoplus_{\dim Z'=k+1}
\mathrm{St}_{Z'}^{*}\otimes F(W')
\hspace{2pt}=\hspace{2pt}
\mathrm{C}^{k+1}_{1}F
$$
est l'homomorphisme dont la composante d'indice $Z,Z'$ (convention matricielle), disons $\mathrm{d}_{Z',Z}^{k}$, est donnée par
$$
\mathrm{d}_{Z',Z}^{k}=
\begin{cases}
(\mathrm{r}_{Z',Z})^{*}\otimes F(Z\subset Z') & \text{pour $Z\subset Z'$},  \\
0 & \text{pour $Z\not\subset Z'$}.
\end{cases}
$$
Par définition même on a~:

\medskip
\begin{pro}\label{Bob} Soient $F$ un foncteur de $\mathcal{W}$ dans $\mathcal{E}$ et $k$ un entier naturel~; on a
$$
\hspace{24pt}
\mathrm{H}^{k}\mathrm{C}^{\bullet}_{1}F
\hspace{4pt}=\hspace{4pt}
\mathrm{Ext}_{\mathcal{E}^{\mathcal{W}}}^{k}(\mathrm{S}_{0},F)
\hspace{24pt}.
$$
\end{pro}

\medskip
Le complexe d'Oliver $\mathrm{C}^{\bullet}_{\mathrm{St}}F$ dont nous avons traité en \ref{speOl}.1 est, à décalage près, le tronqué (brutal) $\sigma^{\geq 1}\mathrm{C}^{\bullet}_{1}F$~; ce phénomène n'est pas surprenant~:

\vspace{0,75cm}
\textit{Relation entre $\mathrm{Ext}_{\mathcal{E}^{\mathcal{W}}}^{p+1}(\mathrm{S}_{0},F)$ et $\lim_{\mathcal{W}_{0}}^{p}(F\circ\mathrm{i})$}

\bigskip
La notation $\mathrm{i}$ qui apparaît dans l'intertitre ci-dessus désigne le ``foncteur inclusion canonique'' de $\mathcal{W}_{0}$ dans $\mathcal{W}$. On pose $\mathrm{K}:=\ker\hspace{1pt}(\varepsilon:\mathbb{P}_{0}\to\mathrm{S}_{0})$.

\begin{pro}\label{Klim1} Soit $F$ un foncteur de $\mathcal{W}$ dans $\mathcal{E}$~; on a un $\mathcal{E}$-isomorphisme canonique
$$
\hspace{24pt}
\mathrm{Hom}_{\mathcal{E}^{\mathcal{W}}}(\mathrm{K},F)
\hspace{4pt}\cong\hspace{4pt}{\lim}_{\mathcal{W}_{0}}\hspace{1pt}(F\circ\mathrm{i})
\hspace{24pt}.
$$
\end{pro}

\textit{Démonstration.} On constate que l'on a
$$
\mathrm{K}(U)
\hspace{4pt}=\hspace{4pt}
\begin{cases}
\mathbb{F}_{2}& \text{pour}\hspace{4pt}U\not=0, \\
0 & \text{pour}\hspace{4pt}U=0.
\end{cases}
$$
Il en résulte $\mathrm{Hom}_{\mathcal{E}^{\mathcal{W}}}(\mathrm{K},F)=
\mathrm{Hom}_{\mathcal{E}^{\mathcal{W}_{0}}}(\Delta_{\mathbb{F}_{2}},F\circ\mathrm{i})$.
\hfill$\square$

\begin{cor}\label{Klim2} Soient $F$ un foncteur de $\mathcal{W}$ dans $\mathcal{E}$ et $k$ un entier naturel~; on a un $\mathcal{E}$-isomorphisme canonique
$$
\hspace{24pt}
\mathrm{Ext}_{\mathcal{E}^{\mathcal{W}}}^{k}(\mathrm{K},F)
\hspace{4pt}\cong\hspace{4pt}{\lim}_{\mathcal{W}_{0}}^{k}\hspace{1pt}(F\circ\mathrm{i})
\hspace{24pt}.
$$\end{cor}

\textit{Démonstration.} On observe que le foncteur  $\mathcal{E}^{\mathcal{W}}\to\mathcal{E}^{\mathcal{W}_{0}},F\mapsto F\circ\mathrm{i}$ est exact et transforme injectifs en injectifs.
\hfill$\square$

\begin{pro}\label{Klim3} Soit $F$ un foncteur de $\mathcal{W}$ dans $\mathcal{E}$~; on a une suite exacte dans $\mathcal{E}$
$$
0\to\mathrm{Hom}_{\mathcal{E}^{\mathcal{W}}}(\mathrm{S}_{0},F)
\to F(0)\to{\lim}_{\mathcal{W}_{0}}\hspace{1pt}(F\circ\mathrm{i})
\to\mathrm{Ext}_{\mathcal{E}^{\mathcal{W}}}^{1}(\mathrm{S}_{0},F)\to 0
$$
et pour $k>0$ un isomorphisme
$$
\hspace{24pt}
{\lim}_{\mathcal{W}_{0}}^{k}\hspace{1pt}(F\circ\mathrm{i})
\hspace{4pt}\cong\hspace{4pt}
\mathrm{Ext}_{\mathcal{E}^{\mathcal{W}}}^{k+1}(\mathrm{S}_{0},F)
\hspace{24pt}.
$$
\end{pro}

\textit{Démonstration.} On considère la longue suite exacte des $\mathrm{Ext}$ associée à la suite exacte dans $\mathcal{E}^{\mathcal{W}}$, $0\to\mathrm{K}\to\mathbb{P}_{0}\to\mathrm{S}_{0}\to 0$.
\hfill$\square$

\vspace{0.75cm}
\textit{Résolution projective de $\mathrm{S}_{W}$ ($\hspace{2pt}W\subset V$) dans la catégorie $\mathcal{E}^{\mathcal{W}}$}

\bigskip
Soient $W$ un sous-groupe de $V$ et $k$ un entier naturel~; on pose
$$
\mathrm{P}_{k,W}
\hspace{4pt}:=\hspace{4pt}
\bigoplus_{\text{$Z\supset W$ et $\dim Z/W=k$}}
\mathrm{St}_{Z/W}\otimes\mathbb{P}_{Z}
$$
et on note
\begin{multline*}
\mathrm{d}_{k+1}
\hspace{4pt}:\hspace{4pt}
\mathrm{P}_{k+1,W}
\hspace{4pt}=\hspace{4pt}
\bigoplus_{\text{$Z'\supset W$ et $\dim Z'/W=k+1$}}
\mathrm{St}_{Z'/W}\otimes\mathbb{P}_{Z'} \\
\longrightarrow
\bigoplus_{\text{$Z\supset W$ et $\dim Z/W=k$}}
\mathrm{St}_{Z/W}\otimes\mathbb{P}_{Z}
\hspace{4pt}=\hspace{4pt}
\mathrm{P}_{k,W}
\end{multline*}
la transformation naturelle dont la composante d'indices $Z,Z'$ (convention ``matricielle''), disons ${\mathrm{d}_{k+1}}_{Z,Z'}$, est donnée par
 $$
{\mathrm{d}_{k+1}}_{Z,Z'}
\hspace{4pt}=\hspace{4pt}
\begin{cases}
\mathrm{r}_{Z'/W,Z/W}\otimes\nu_{Z',Z} & \text{pour $Z\subset Z'$} \\
0 & \text{pour $Z\not\subset Z'$},
\end{cases}
 $$
 $\nu_{Z',Z}$ désignant l'unique élément non trivial de  $\mathrm{Hom}_{\mathcal{E}^{\mathcal{W}}}(\mathbb{P}_{Z'},\mathbb{P}_{Z})$.
 
 \smallskip
 On note encore $\varepsilon:\mathrm{P}_{0,W}=\mathbb{P}_{W}\to\mathrm{S}_{W}$ le $\mathcal{E}^{\mathcal{W}}$-épimorphisme canonique (qui est  à nouveau une couverture projective).
 
 \medskip
 \begin{pro}\label{resprojSW} Soit $W$ un sous-groupe de $V$ avec $\mathop{\mathrm{codim}}W=p$~; la~suite de $\mathcal{E}^{\mathcal{W}}$-morphismes
 $$
 \hspace{24pt}
 0\longleftarrow\mathrm{S}_{W}
 \overset{\varepsilon}{\longleftarrow}
 \mathrm{P}_{0,W} \overset{\mathrm{d}_{1}}{\longleftarrow}
 \mathrm{P}_{1,W} \overset{\mathrm{d}_{2}}{\longleftarrow}
 \ldots\overset{\mathrm{d}_{p}}{\longleftarrow}
  \mathrm{P}_{p,W} 
  \hspace{24pt},
 $$
est une résolution projective de $\mathrm{S}_{W}$ dans la catégorie $\mathcal{E}^{\mathcal{W}}$.
\end{pro}
 
\medskip
 \textit{Démonstration.} Soit $U$ un sous-groupe de $V$. Si $W$ est un sous-groupe strict de~$U$ alors $\mathrm{P}_{\bullet,W}(U)$ s'identifie au complexe $\mathrm{Lu}_{\bullet,U/W}$ qui est acyclique et $\mathrm{S}_{W}(U)$ est nul. Dans le cas $U=W$, on a $\mathrm{P}_{k,W}(U)=0$ pour $k>0$ et~$\varepsilon_{U}$ est un isomorphisme. Si $W$ n'est pas contenu dans $U$ alors $\mathrm{P}_{\bullet,W}(U)$ et $\mathrm{S}_{W}(U)$ sont nuls.
\hfill$\square$

\medskip
\begin{cor}\label{ExtSWSZ} Soient $W$ et $Z$ deux sous-groupes de $V$. On a un isomorphisme canonique
$$
\mathrm{Ext}_{\mathcal{E}^{\mathcal{W}}}^{k}(\mathrm{S}_{W},\mathrm{S}_{Z})
\hspace{4pt}\cong\hspace{4pt}
\begin{cases}
\mathrm{St}_{Z/W}^{*} & \text{pour $W\subset Z$ et $k=\dim Z/W$,} \\
0 & \text{sinon}.
\end{cases}
$$
\end{cor}

\textit{Démonstration.} On constate que l'on a
$$
\hspace{40pt}
\mathrm{Hom}_{\mathcal{E}^{\mathcal{W}}}^{k}(\mathrm{P}_{k,W},\mathrm{S}_{Z})
\hspace{4pt}\cong\hspace{4pt}
\begin{cases}
\mathrm{St}_{Z/W}^{*} & \text{pour $W\subset Z$ et $k=\dim Z/W$,} \\
0 & \text{sinon.$\hspace{143pt}\square$}
\end{cases}
$$

\vspace{0.75cm}
\textit{Résolution injective de $\mathrm{S}_{V}$ dans la catégorie $\mathcal{E}^{\mathcal{W}}$}

\bigskip
On considère le $\mathcal{E}^{\mathcal{W}}$-monomorphisme canonique $\eta:\mathrm{S}_{V}\to\mathbb{I}^{V}$ (qui est en fait une enveloppe injective)~; on se propose de ``prolonger'' $\eta$ en une résolution injective.

\medskip
Soit $k$ un entier naturel~; on pose
$$
\mathrm{I}^{k}
\hspace{4pt}:=\hspace{4pt}
\bigoplus_{\mathop{\mathrm{codim}}W=k}
\mathrm{St}_{V/W}^{*}\otimes\mathbb{I}^{W}
$$
et on note
$$
\mathrm{d}^{k}
\hspace{4pt}:\hspace{4pt}
\mathrm{I}^{k}
\hspace{4pt}=\hspace{4pt}
\bigoplus_{\mathop{\mathrm{codim}}W=k}
\mathrm{St}_{V/W}^{*}\otimes\mathbb{I}^{W}
\longrightarrow
\bigoplus_{\mathop{\mathrm{codim}}W'=k+1}
\mathrm{St}_{V/W'}^{*}\otimes\mathbb{I}^{W'}
\hspace{4pt}=\hspace{4pt}
\mathrm{I}^{k+1}
$$
 la transformation naturelle dont la composante d'indices $W',W$ (convention ``matricielle''), disons $\mathrm{d}^{k}_{W',W}$, est donnée par
 $$
 \mathrm{d}^{k}_{W',W}
 \hspace{4pt}=\hspace{4pt}
\begin{cases}
\mathrm{s}_{V/W',V/W}^{*}\otimes\nu_{W,W'} & \text{pour $W\supset W'$} \\
0 & \text{pour $W\not\supset W'$},
\end{cases}
 $$
 $\nu_{W,W'}$ désignant ici l'unique élément non trivial de  $\mathrm{Hom}_{\mathcal{E}^{\mathcal{W}}}(\mathbb{I}^{W},\mathbb{I}^{W'})$.
 
 \medskip
 \begin{pro}\label{resinjSV} La suite de $\mathcal{E}^{\mathcal{W}}$-morphismes
 $$
 \hspace{24pt}
 0\longrightarrow\mathrm{S}_{V}
 \overset{\eta}{\longrightarrow}
 \mathrm{I}^{0} \overset{\mathrm{d}^{0}}{\longrightarrow}
 \mathrm{I}^{1} \overset{\mathrm{d}^{1}}{\longrightarrow}
 \ldots\overset{\mathrm{d}^{n-1}}{\longrightarrow}
  \mathrm{I}_{n} 
  \hspace{24pt},
 $$
est une résolution injective de $\mathrm{S}_{V}$ dans la catégorie $\mathcal{E}^{\mathcal{W}}$.
\end{pro}
 
\medskip
 \textit{Démonstration.} Soit $U$ un sous-groupe strict de $V$. On constate que
 $$
  \mathrm{I}^{0} (U)\overset{\mathrm{d}^{0}_{U}}{\longrightarrow}
 \mathrm{I}^{1}(U)\overset{\mathrm{d}^{2}_{U}}{\longrightarrow}
 \ldots\overset{\mathrm{d}^{n-1}_{U}}{\longrightarrow}
  \mathrm{I}^{n}(U)
  $$
  est le ``complexe de Lusztig'' $\mathrm{Lu}^{\bullet}_{2,U}$ et que $\mathrm{S}_{V}(U)$ est nul~; or $\mathrm{Lu}^{\bullet}_{2,U}$ est acyclique. Dans le cas $U=V$, on a $\mathrm{I}^{k}(U)=0$ pour $k>0$ et $\eta_{U}$ est un isomorphisme.
\hfill$\square$

\vspace{0.75cm}
\textit{Résolution injective de $\mathrm{S}_{Z}$ ($\hspace{2pt}Z\subset V$) dans la catégorie $\mathcal{E}^{\mathcal{W}}$}

\medskip
Soient $Z$ un sous-groupe de $V$ et $k$ un entier naturel~; on pose
$$
\mathrm{I}^{k,Z}
\hspace{4pt}:=\hspace{4pt}
\bigoplus_{\text{$W\subset Z$ et $\dim Z/W=k$}}
\mathrm{St}_{Z/W}^{*}\otimes\mathbb{I}^{W}
$$
et on note
\begin{multline*}
\mathrm{d}^{k}
\hspace{4pt}:\hspace{4pt}
\mathrm{I}^{k,Z}
\hspace{4pt}=\hspace{4pt}
\bigoplus_{\text{$W\subset Z$ et $\dim Z/W=k$}}
\mathrm{St}_{Z/W}^{*}\otimes\mathbb{I}^{W} \\
\longrightarrow
\bigoplus_{\text{$W'\subset Z$ et $\dim Z/W'=k+1$}}
\mathrm{St}_{Z/W'}^{*}\otimes\mathbb{I}^{W'}
\hspace{4pt}=\hspace{4pt}
\mathrm{I}^{k+1,Z}
\end{multline*}
 la transformation naturelle dont la composante d'indices $W',W$ (convention ``matricielle''), disons $\mathrm{d}^{k}_{W',W}$, est donnée par
 $$
 \mathrm{d}^{k}_{W',W}
 \hspace{4pt}=\hspace{4pt}
\begin{cases}
\mathrm{s}_{Z/W',\hspace{1pt}Z/W}^{*}\otimes\nu_{W,W'} & \text{pour $W\supset W'$}
\\
0 & \text{pour $W\not\supset W'$},
\end{cases}
 $$
 $\nu_{W,W'}$ désignant l'unique élément non trivial de  $\mathrm{Hom}_{\mathcal{E}^{\mathcal{W}}}(\mathbb{I}^{W},\mathbb{I}^{W'})$.
 
\smallskip
On note encore $\eta:\mathrm{S}_{Z}\to\mathbb{I}^{Z}=\mathrm{I}^{0,Z}$ le $\mathcal{E}^{\mathcal{W}}$-monomorphisme canonique (qui est à nouveau une enveloppe injective).
 
\begin{pro}\label{resinjSZ} La suite de $\mathcal{E}^{\mathcal{W}}$-morphismes
 $$
 \hspace{24pt}
 0\longrightarrow\mathrm{S}_{Z}
 \overset{\eta}{\longrightarrow}
 \mathrm{I}^{0,Z} \overset{\mathrm{d}^{0}}{\longrightarrow}
 \mathrm{I}^{1,Z} \overset{\mathrm{d}^{1}}{\longrightarrow}
 \ldots
 \hspace{24pt},
$$
est une résolution injective de $\mathrm{S}_{Z}$ dans la catégorie $\mathcal{E}^{\mathcal{W}}$.
\end{pro}

\medskip
 \textit{Démonstration.} Soit $U$ un sous-groupe de $V$. Si $U$ est un sous-groupe strict de $Z$ alors $\mathrm{I}^{\bullet,Z}(U)$ s'identifie au complexe $\mathrm{Lu}_{2,Z/U}^{\bullet}$ qui est acyclique et $\mathrm{S}_{Z}(U)$ est nul. Dans la cas $U=Z$, on a $\mathrm{I}^{k,Z}(U)=0$ pour $k>0$ et $\eta_{U}$ est un isomorphisme. Si $U$ n'est pas contenu dans $Z$ alors  $\mathrm{I}^{\bullet,Z}(U)$ et $\mathrm{S}_{Z}(U)$ sont nuls. 
\hfill$\square$

\pagebreak

\vspace{0.75cm}
\begin{rems}\label{dualsuite-Kan} Le lecteur aura noté que les quatre démonstrations des propositions \ref{resprojS0}, \ref{resprojSW}, \ref{resinjSV} et \ref{resinjSZ} sont \textit{mutatis mutandis} les mêmes. Nous présentons ci-dessous, en petits caractères, des arguments conceptuels qui auraient pu éviter ces répétitions (tout en rallongeant l'exposition~!).

\footnotesize

\bigskip
1) \textit{Dualité (suite)}

\medskip
On note $\mathcal{E}_{\mathrm{f}}$ la sous-catégorie pleine de $\mathcal{E}$ dont les objets sont les $\mathbb{F}_{2}$-espaces vectoriels de dimension finie. Par ``general nonsense'' on a un isomorphisme fonctoriel
$$
{(\mathcal{E}_{\mathrm{f}}^{\mathcal{W}_{V}})}^{\mathrm{op}}
\hspace{4pt}\cong\hspace{4pt}
{(\mathcal{E}_{\mathrm{f}}^{\mathrm{op}})}^{\mathcal{W}_{V}^{\mathrm{op}}}
$$
(rappelons que lorsque nous sommes amenés à faire varier le $2$-groupe abélien élémentaire $V$, nous précisons, si nécessaire, la notation $\mathcal{W}$ en $\mathcal{W}_{V}$).

\smallskip
Or on dispose~:

\smallskip
-- d'un isomorphisme d'ensembles ordonnés $\mathcal{W}_{V}^{\mathrm{op}}\overset{\cong}{\to}\mathcal{W}_{V^{*}}$, à savoir l'application $U\mapsto U^{\perp}$~;

\smallskip
-- d'un isomorphisme fonctoriel $\mathcal{E}_{\mathrm{f}}^{\mathrm{op}}\overset{\cong}{\to}\mathcal{E}_{\mathrm{f}}$, à savoir le foncteur $E\mapsto E^{*}$.

\smallskip
Il existe donc un isomorphisme fonctoriel canonique, disons
$$
\hspace{24pt}
\mathrm{D}
\hspace{2pt}:\hspace{2pt}
\mathcal{E}_{\mathrm{f}}^{\mathcal{W}_{V^{*}}}
\overset{\cong}{\longrightarrow}{(\mathcal{E}_{\mathrm{f}}^{\mathcal{W}_{V}})}^{\mathrm{op}}
\hspace{24pt}.
$$
Explicitons. Soit $F$ un foncteur de $\mathcal{W}_{V^{*}}$ dans  $\mathcal{E}_{\mathrm{f}}$~; $\mathrm{D}F$ est le foncteur de $\mathcal{W}_{V}$ dans  $\mathcal{E}_{\mathrm{f}}$, $U\mapsto~{F(U^{\perp})}^{*}$. Il est clair que le foncteur $\mathrm{D}$ est exact et qu'il transforme les objets projectifs (resp. injectifs) de la catégorie abélienne $\mathcal{E}_{\mathrm{f}}^{\mathcal{W}_{V^{*}}}$ en objets injectifs (resp. projectifs) de la catégorie abélienne $\mathcal{E}_{\mathrm{f}}^{\mathcal{W}_{V}}$. 
Donnons deux exemples. Considérons le foncteur $\mathbb{P}_{W^{\perp}}:\mathcal{W}_{V^{*}}\to~\mathcal{E}_{\mathrm{f}}$, $W$ désignant un sous-groupe de $V$~; on constate que l'on a $\mathrm{D}\hspace{0.5pt}\mathbb{P}_{W^{\perp}}=\mathbb{I}^{W}$.  Considérons le foncteur $\mathrm{S}_{0_{W^{*}}}:\mathcal{W}_{V^{*}}\to\mathcal{E}_{\mathrm{f}}$~; on constate que l'on a $\mathrm{D}\hspace{0.5pt}\mathrm{S}_{0_{V^{*}}}=\mathrm{S}_{V}$.

\medskip
La discussion précédente montre que l'on peut obtenir la résolution injective  de \ref{resinjSV} en appliquant le foncteur $\mathrm{D}$ à la résolution projective de \ref{resprojS0} ($V$~étant remplacé par $V^{*}$).

\bigskip
2) \textit{Extensions de Kan}

\medskip
Soit $Z$ un sous-groupe de $V$ on note $\mathcal{W}_{\leq Z}$ (resp. $\mathcal{W}_{\geq Z}$) le sous-ensemble (ordonné) de $\mathcal{W}$ constitué des $W$ avec $W\subset Z$ (resp. $W\supset Z$)~; on note $\mathrm{i}_{\leq Z}:\mathcal{W}_{\leq Z}\to\mathcal{W}$ (resp. $\mathrm{i}_{\geq Z}:\mathcal{W}_{\geq Z}\to\mathcal{W}$) l'inclusion canonique. On identifiera respectivement les catégories $\mathcal{W}_{\leq Z}$ et $\mathcal{W}_{\geq Z}$  aux catégories $\mathcal{W}_{Z}$ et $\mathcal{W}_{V/Z}$ auxquelles elles sont canoniquement isomorphes.

\medskip
On note
$$
\mathrm{Ran}_{Z}
\hspace{4pt}:\hspace{4pt}
\mathcal{E}^{\mathcal{W}_{Z}}
\longrightarrow
\mathcal{E}^{\mathcal{W}}
$$
le foncteur extension de Kan à droite le long de $\mathrm{i}_{\leq Z}$, c'est-à-dire l'adjoint à droite du fonteur $\mathcal{E}^{\mathcal{W}}
\to\mathcal{E}^{\mathcal{W}_{Z}},F\mapsto F\circ\mathrm{i}_{\leq Z}$. Soit $G$ un foncteur de $\mathcal{W}_{Z}$ dans~$\mathcal{E}$~; on constate que l'on~a~:
$$
(\mathrm{Ran}_{Z}\hspace{1pt}G)(U)
\hspace{4pt}=\hspace{4pt}
\begin{cases}
G(U) & \text{pour $U\subset Z$}, \\
0 & \text{pour $U\not\subset Z$},
\end{cases}
$$
et que la co-unité d'adjonction $\epsilon_{G}:\mathrm{Ran}_{Z}\hspace{1pt}G\circ\mathrm{i}_{\leq Z}\to G$ est un isomorphisme.

\medskip
On note
$$
\mathrm{Lan}_{Z}
\hspace{4pt}:\hspace{4pt}
\mathcal{E}^{\mathcal{W}_{V/Z}}
\longrightarrow
\mathcal{E}^{\mathcal{W}}
$$
le foncteur extension de Kan à gauche le long de $\mathrm{i}_{\geq Z}$ c'est-à-dire l'adjoint à gauche du fonteur $\mathcal{E}^{\mathcal{W}}
\to\mathcal{E}^{\mathcal{W}_{V/Z}},G\mapsto G\circ\mathrm{i}_{\geq Z}$. Soit $F$ un foncteur de $\mathcal{W}_{V/Z}$ dans~$\mathcal{E}$~; on constate que l'on a~:
$$
(\mathrm{Lan}_{Z}\hspace{1pt}F)(U)
\hspace{4pt}=\hspace{4pt}
\begin{cases}
F(U) & \text{pour $U\supset Z$}, \\
0 & \text{pour $U\not\supset Z$},
\end{cases}
$$
et que l'unité d'adjonction $\eta_{F}:F\to\mathrm{Lan}_{Z}\hspace{1pt}F\circ\mathrm{i}_{\geq Z}$ est un isomorphisme.

\medskip
Les deux propositions suivantes sont  immédiates~:

\bigskip
\textsc{Proposition R.} {\em Soit $Z$ un sous-groupe de $V$.

\medskip
{\em (a)} Les deux foncteurs $F\mapsto F\circ\mathrm{i}_{\leq Z}$ et $\mathrm{Ran}_{Z}$ sont exacts.

\medskip
{\em (b)} Le foncteur $F\mapsto F\circ\mathrm{i}_{\leq Z}$  envoie projectifs sur projectifs.

\medskip
{\em (c)} Le foncteur $\mathrm{Ran}_{Z}$ envoie injectifs sur injectifs.

\medskip
{\em (d)} Pour tout entier $k\geq 0$, tout foncteur $F:\mathcal{W}\to\mathcal{E}$ et tout foncteur $G:\mathcal{W}_{Z}\to\mathcal{E}$, on a un isomorphisme (naturel)~:
$$
\hspace{24pt}
\mathrm{Ext}_{\mathcal{E}^{\mathcal{W}_{Z}}}^{k}(F\circ\mathrm{i}_{\leq Z},G)
\hspace{4pt}\cong\hspace{4pt}
\mathrm{Ext}_{\mathcal{E}^{\mathcal{W}}}^{k}(F,\mathrm{Ran}_{Z}\hspace{1pt}G)
\hspace{24pt}.
$$}

\textsc{Propositon L.} {\em Soit $Z$ un sous-groupe de $V$.

\medskip
{\em (a)} Les deux foncteurs $G\mapsto G\circ\mathrm{i}_{\geq Z}$ et $\mathrm{Lan}_{Z}$ sont exacts.

\medskip
{\em (b)} Le foncteur $G\mapsto G\circ\mathrm{i}_{\geq Z}$  envoie injectifs sur injectifs.

\medskip
{\em (c)} Le foncteur $\mathrm{Lan}_{Z}$ envoie projectifs sur projectifs.

\medskip
{\em (d)} Pour tout entier $k\geq 0$, tout foncteur $F:\mathcal{W}_{V/Z}\to\mathcal{E}$ et tout foncteur $G:\mathcal{W}\to\mathcal{E}$, on a un isomorphisme (naturel)~:
$$
\hspace{24pt}
\mathrm{Ext}_{\mathcal{E}^{\mathcal{W}_{V/Z}}}^{k}(F,G\circ\mathrm{i}_{\geq Z})
\hspace{4pt}\cong\hspace{4pt}
\mathrm{Ext}_{\mathcal{E}^{\mathcal{W}}}^{k}(\mathrm{Lan}_{Z}\hspace{1pt}F,G)
\hspace{24pt}.
$$}

Le fait que $\mathrm{Lan}_{W}$ est exact et envoie projectifs sur projectifs permet de déduire \ref{resprojSW} de \ref{resprojS0} (avec $V$ remplacé par $V/W$). Pareillement, le fait que $\mathrm{Ran}_{Z}$ est exact et envoie injectifs sur injectifs permet de déduire \ref{resinjSZ} de \ref{resinjSV} (avec $V$ remplacé par $Z$). 

\end{rems}

\vspace{0.75cm}
\textit{Résolution injective d'un foncteur arbitraire de $\mathcal{W}$ dans $\mathcal{E}$}

\bigskip
Soit $F$ un foncteur de $\mathcal{W}$ dans $\mathcal{E}$.

\bigskip
--\hspace{8pt}Soient $W$ et $Z$ deux sous-groupes de $V$ avec $W\subset Z$ (alternativement, de façon plus pédante, soit $W\subset Z$ un morphisme de $\mathcal{W}$). On pose
$$
\hspace{24pt}
\mathrm{I}^{\hspace{1pt}W\subset Z}F
\hspace{4pt}:=\hspace{4pt}
\mathrm{St}_{Z/W}^{*}\otimes F(Z)\otimes\mathbb{I}^{W}
\hspace{24pt};
$$
$\mathrm{I}^{\hspace{1pt}W\subset Z}F$ est un objet injectif de $\mathcal{E}^{\mathcal{W}}$.

\bigskip
--\hspace{8pt}Soient $W',W,Z$ trois sous-groupes de $V$ avec $W'\subset W\subset Z$ et\linebreak $\dim W/W'=1$. On note
$$
\mathrm{d}_{\mathrm{h}}^{W'\subset W\subset Z}
\hspace{2pt}:\hspace{2pt}
\mathrm{I}^{\hspace{1pt}W\subset Z}F\to\mathrm{I}^{\hspace{1pt}W'\subset Z}F
$$
le $\mathcal{E}^{\mathcal{W}}$-morphisme induit par

\smallskip
- l'homomorphisme $\mathrm{s}_{Z/W',Z/W}^{*}:\mathrm{St}_{Z/W}^{*}\to\mathrm{St}_{Z/W'}^{*}$,

\smallskip
- l'identité de $F(Z)$,

\smallskip
- l'unique $\mathcal{E}^{\mathcal{W}}$-morphisme non trivial de $\mathbb{I}^{W}$ dans $\mathbb{I}^{W'}$.

\bigskip
--\hspace{8pt}Soient $W,Z,Z'$ trois sous-groupes de $V$ avec $W\subset Z\subset Z'$ et $\dim Z'/Z=1$. On note
$$
\mathrm{d}_{\mathrm{v}}^{W\subset Z\subset Z'}
\hspace{2pt}:\hspace{2pt}
\mathrm{I}^{\hspace{1pt}W\subset Z}F\to\mathrm{I}^{\hspace{1pt}W\subset Z'}F
$$
le $\mathcal{E}^{\mathcal{W}}$-morphisme induit par

\smallskip
-- l'homomorphisme $\mathrm{r}_{Z'/W,Z/W}^{*}:\mathrm{St}_{Z/W}^{*}\to\mathrm{St}_{Z'/W}^{*}$,

\smallskip
-- l'homomorphisme $F(Z\subset Z')$,

\smallskip
-- l'identité de $\mathbb{I}^{W}$.

\bigskip
--\hspace{8pt}On organise la famille $\{\mathrm{I}^{\hspace{1pt}W\subset Z}F\}$ en un objet $(\mathbb{Z}\times\mathbb{Z})$-gradué de $\mathcal{E}^{\mathcal{W}}$. On pose~:
$$
\hspace{24pt}
\mathrm{I}^{p,q}F
\hspace{4pt}:=\hspace{4pt}
\bigoplus_{W\subset Z,\hspace{4pt}\mathop{\mathrm{codim}}W=p,\hspace{4pt}\mathop{\mathrm{codim}}Z=-q}
\mathrm{I}^{\hspace{1pt}W\subset Z}F
\hspace{24pt};
$$
par définition, $\mathrm{I}^{p,q}F$ est nul pour  $(p,q)\not\in\mathrm{Tg}_{n}$ (on rappelle que $\mathrm{Tg}_{n}$ est le sous-ensemble de $\mathbb{Z}\times\mathbb{Z}$ constitué des couples $(p,q)$ vérifiant $p\geq 0$, $p\leq n$, $p+q\geq 0$ et $q\leq 0$).

\medskip
On note $\mathrm{d}_{\mathrm{h}}^{p,q}:\mathrm{I}^{p,q}F\to\mathrm{I}^{p+1,q}F$, le $\mathcal{E}^{\mathcal{W}}$-morphisme induit par les $\mathrm{d}_{\mathrm{h}}^{W'\subset W\subset Z}$ (on a $\mathrm{d}_{\mathrm{h}}^{p,q}=0$ si $(p,q)$ ou $(p+1,q)$ n'appartient pas à $\mathrm{Tg}_{n}$).

\medskip
On note $\mathrm{d}_{\mathrm{v}}^{p,q}:\mathrm{I}^{p,q}F\to\mathrm{I}^{p,q+1}F$, le $\mathcal{E}^{\mathcal{W}}$-morphisme induit par les $\mathrm{d}_{\mathrm{v}}^{W\subset Z\subset Z'}$ (on a $\mathrm{d}_{\mathrm{v}}^{p,q}=0$ si $(p,q)$ ou $(p,q+1)$ n'appartient pas à $\mathrm{Tg}_{n}$).

\bigskip
\begin{pro-def}\label{IpqF} Soit $F$ un foncteur de $\mathcal{W}$ dans $\mathcal{E}$~; soit $(p,q)$ un élément de $\mathbb{Z}\times\mathbb{Z}$.

\medskip
{\em (a)} Le $\mathcal{E}^{\mathcal{W}}$-diagramme
$$
\begin{CD}
\mathrm{I}^{p,q+1}F @>\mathrm{d}_{\mathrm{h}}^{p,q+1}>> \mathrm{I}^{p+1,q+1}F \\
@A\mathrm{d}_{\mathrm{v}}^{p,q}AA @A\mathrm{d}_{\mathrm{v}}^{p+1,q}AA \\
\mathrm{I}^{p,q}F @>\mathrm{d}_{\mathrm{h}}^{p,q}>> \mathrm{I}^{p+1,q}F
\end{CD}
$$
est commutatif.

\pagebreak

\medskip
{\em (b)} Le $\mathcal{E}^{\mathcal{W}}$-morphisme composé
$$
\begin{CD}
\mathrm{I}^{p,q}F @>\mathrm{d}_{\mathrm{h}}^{p,q}>> \mathrm{I}^{p+1,q}F@>\mathrm{d}_{\mathrm{h}}^{p+1,q}>> \mathrm{I}^{p+2,q}F
\end{CD}
$$
est nul.

\medskip
{\em (c)} Le $\mathcal{E}^{\mathcal{W}}$-morphisme composé
$$
\begin{CD}
\mathrm{I}^{p,q}F @>\mathrm{d}_{\mathrm{v}}^{p,q}>> \mathrm{I}^{p,q+1}F@>\mathrm{d}_{\mathrm{v}}^{p,q+1}>> \mathrm{I}^{p,q+2}F
\end{CD}
$$
est nul.

\bigskip
Le $\mathcal{E}^{\mathcal{W}}$-objet bigradué $\{\mathrm{I}^{p,q}F\}_{(p,q)\in\mathbb{Z}\times\mathbb{Z}}$, muni des cobords $\mathrm{d}_{\mathrm{h}}^{p,q}$ et  $\mathrm{d}_{\mathrm{v}}^{p,q}$, est un bicomplexe de cochaînes dans la catégorie $\mathcal{E}^{\mathcal{W}}$ que l'on note $\mathrm{I}^{\bullet,\bullet}F$.
\end{pro-def}

\bigskip
\textit{Démonstration du (a).} Elle résulte de \ref{Stcom}.
\hfill$\square$

\medskip
\textit{Démonstration du (b).}  On considère la résolution injective $\mathrm{S}_{Z}\to\mathrm{I}^{\bullet,Z}$ de la proposition \ref{resinjSZ}. Soit $q$ un entier avec $-n\leq q\leq 0$~; on a tout fait pour que la suite de $\mathcal{E}^{\mathcal{W}}$-morphismes
$$
\mathrm{I}^{0,q}F\longrightarrow\ldots\longrightarrow
\mathrm{I}^{p,q}F\overset{\mathrm{d}_{\mathrm{h}}^{p,q}}{\longrightarrow}
\mathrm{I}^{p+1,q}F\longrightarrow
\ldots\longrightarrow\mathrm{I}^{n,q}F
$$
coïncide avec $\bigoplus_{\mathop{\mathrm{codim}}Z=-q}\hspace{4pt}
F(Z)\otimes\mathrm{I}^{\bullet,Z}[-q]$.
\hfill$\square$

\medskip
\textit{Démonstration du (c).} On considère la résolution projective $\mathrm{P}_{\bullet,W}\to\mathrm{S}_{W}$ de la proposition \ref{resprojSW} et le $\mathcal{E}$-complexe de cochaînes $\mathrm{Hom}_{\mathcal{E}^{\mathcal{W}}}(\mathrm{P}_{\bullet,W},F)$. Soit $p$ un entier avec $0\leq p\leq n$~; on constate cette fois que la suite de $\mathcal{E}^{\mathcal{W}}$-morphismes
$$
\mathrm{I}^{p,-n}F\longrightarrow\ldots\longrightarrow
\mathrm{I}^{p,q}F\overset{\mathrm{d}_{\mathrm{v}}^{p,q}}{\longrightarrow}
\mathrm{I}^{p,q+1}F\longrightarrow
\ldots\longrightarrow\mathrm{I}^{p,0}F
$$
s'identifie à $\bigoplus_{\mathop{\mathrm{codim}}W=p}\hspace{4pt}
\mathrm{Hom}_{\mathcal{E}^{\mathcal{W}}}(\mathrm{P}_{\bullet,W},F)[-p]\otimes\mathbb{I}^{W}$.
\hfill$\square$

\bigskip
On note $\mathrm{Tot}\hspace{1pt}\mathrm{I}^{\bullet,\bullet}F$ le totalisé du bicomplexe défini ci-dessus~;  $\mathrm{Tot}\hspace{1pt}\mathrm{I}^{\bullet,\bullet}F$ est donc un complexe de cochaînes dans la catégorie $\mathcal{E}^{\mathcal{W}}$ dont tous les termes sont injectifs. Soit $k\geq 0$ un entier~; on constate que l'on a
$$
\mathrm{Tot}^{k}\hspace{2pt}\mathrm{I}^{\bullet,\bullet}F
\hspace{4pt}=\hspace{4pt}
\bigoplus_{W\subset Z,\hspace{4pt}\dim Z/W=k}
\mathrm{St}_{Z/W}^{*}\otimes F(Z)\otimes\mathbb{I}^{W}
$$
et donc en particulier
$$
\hspace{4pt}
\mathrm{Tot}^{0}\hspace{2pt}\mathrm{I}^{\bullet,\bullet}F
\hspace{2pt}=\hspace{2pt}
\bigoplus_{W}
F(W)\otimes\mathbb{I}^{W}
\hspace{12pt},\hspace{12pt}
\mathrm{Tot}^{1}\hspace{2pt}\mathrm{I}^{\bullet,\bullet}F
\hspace{2pt}=\hspace{2pt}
\bigoplus_{W\subset Z,\hspace{4pt}\dim Z/W=1}
F(Z)\otimes\mathbb{I}^{W}
\hspace{4pt}.
$$
On note $\eta:F\to\bigoplus_{W}F(W)\otimes\mathbb{I}^{W}=\mathrm{Tot}^{0}\hspace{1pt}\mathrm{I}^{\bullet,\bullet}F$ le $\mathcal{E}^{\mathcal{W}}$-monomorphisme canonique.

\bigskip
\begin{pro}\label{resinjF} Soit $F$ un foncteur de $\mathcal{W}$ dans $\mathcal{E}$. Le $\mathcal{E}^{\mathcal{W}}$-morphisme composé
$$
\begin{CD}
F@>\eta>>
\mathrm{Tot}^{0}\hspace{2pt}\mathrm{I}^{\bullet,\bullet}F
@>\mathrm{d}^{0}_{\mathrm{Tot}}>>
\mathrm{Tot}^{1}\hspace{2pt}\mathrm{I}^{\bullet,\bullet}F
\end{CD}
$$
est nul (la notation $\mathrm{d}^{0}_{\mathrm{Tot}}$ est transparente) et $\hspace{1pt}F\overset{\eta}{\to}\mathrm{Tot}\hspace{1pt}\mathrm{I}^{\bullet,\bullet}F\hspace{1pt}$ est une résolution injective de $F$ dans la catégorie $\mathcal{E}^{\mathcal{W}}$.
\end{pro}

\bigskip
\textit{Démonstration.} On vérifie tout d'abord l'égalité $\mathrm{d}^{0}_{\mathrm{Tot}}\circ\eta=0$ . Soient $W$ et $Z$ deux sous-groupes de $V$ avec $W\subset Z$ et $\dim Z/W=1$ ; soit $\pi$ la projection
 de $\mathrm{Tot}^{1}\hspace{2pt}\mathrm{I}^{\bullet,\bullet}F$ sur $F(Z)\otimes\mathbb{I}^{W}$. Par construction $\pi\circ(\mathrm{d}^{0}_{\mathrm{Tot}}\circ\eta)$ est la somme dans $\mathrm{Hom}_{\mathcal{E}^{\mathcal{W}}}(F,F(Z)\otimes\mathbb{I}^{W})$ des deux homomorphismes suivants~:
 
\smallskip
 -- le composé de l'homomorphisme $F\to F(W)\otimes\mathbb{I}^{W}$  ``donné par Yoneda''  et de l'homomorphisme  $F(W)\otimes\mathbb{I}^{W}\to F(Z)\otimes\mathbb{I}^{W}$ induit par $F(W\subset Z)$,

 \smallskip
 -- le composé de l'homomorphisme $F\to F(Z)\otimes\mathbb{I}^{Z}$ ``donné par Yoneda''  et de l'homomorphisme $F(Z)\otimes\mathbb{I}^{Z}\to F(Z)\otimes\mathbb{I}^{W}$ induit par l'unique $\mathcal{E}^{\mathcal{W}}$-morphisme non trivial de $\mathbb{I}^{Z}$ dans $\mathbb{I}^{W}$.
 
 \smallskip
L'égalité $\mathrm{d}^{0}_{\mathrm{Tot}}\circ\eta=0$ est donc équivalente à la commutativité des diagrammes suivants
$$
\begin{CD}
F(Z)\otimes\mathbb{I}^{Z} @>>> F(Z)\otimes\mathbb{I}^{W} \\
@AAA @AAA \\
F @>>> F(W)\otimes\mathbb{I}^{W}
\end{CD}
$$
dont la vérification est immédiate.

\medskip
On montre maintenant que $\hspace{1pt}F\overset{\eta}{\to}\mathrm{Tot}\hspace{1pt}\mathrm{I}^{\bullet,\bullet}F\hspace{1pt}$ est bien une résolution injective de $F$ dans la catégorie $\mathcal{E}^{\mathcal{W}}$. On pose $\Gamma^{\bullet}F:=\mathrm{Tot}\hspace{1pt}\mathrm{I}^{\bullet,\bullet}F$ et on note $\widetilde{\Gamma}^{\bullet}F$ le complexe coaugmenté défini par $\eta$~; on considère ci-après $F\mapsto\widetilde{\Gamma}^{\bullet}F$ comme un foncteur défini sur la catégorie $\mathcal{E}^{\mathcal{W}}$ et à valeurs dans la catégorie des $\mathcal{E}^{\mathcal{W}}$-complexes de cochaînes coaugmentés. Le foncteur $\widetilde{\Gamma}^{\bullet}$ possède les deux propriétés suivantes~:

\smallskip
1) Il est exact. En effet, il est manifeste que le foncteur $F\mapsto\mathrm{I}^{W\subset Z}F$ est exact.

\smallskip
2) Pour tout sous-groupe $Z$ de $V$ le complexe $\widetilde{\Gamma}^{\bullet}\hspace{1pt}\mathrm{S}_{Z}$ est acyclique. En effet, on constate que $\widetilde{\Gamma}^{\bullet}\hspace{1pt}\mathrm{S}_{Z}$ coïncide avec la résolution injective de \ref{resinjSZ}.

\smallskip
On en déduit que le complexe $\widetilde{\Gamma}^{\bullet}F$ est acyclique pour tout $F$ en utilisant la filtration décroissante de $F$ par les $\mathrm{filt}_{1}^{k}F,\hspace{4pt}0\leq k\leq n+1$. En effet, on a $\mathrm{filt}_{1}^{n+1}F=0$, $\mathrm{filt}_{1}^{0}F=F$ et le quotient $\mathrm{Gr}_{1}^{k}F:=\mathrm{filt}_{1}^{k}F/\mathrm{filt}_{1}^{k+1}F$ est  isomorphe à une somme directe de $\mathrm{S}_{Z}$ (avec $\dim Z=k$) pour tout $k$ avec $0\leq k\leq n$.
\hfill$\square$

\pagebreak

\bigskip
\begin{rem}\label{Keller} Par construction le complexe $\mathrm{I}^{\bullet,q}F\hspace{1pt}[q]$ est une résolution injective du foncteur  $\mathrm{Gr}_{1}^{q+n}F$ dans la catégorie $\mathcal{E}^{\mathcal{W}}$ (voir la démonstration du point (b) de \ref{IpqF}).  Nous approfondissons cette remarque (en petits caractères) ci-dessous.

\footnotesize
\medskip
Soit $k$ un entier avec $0\leq k\leq n-1$. Soit $\epsilon_{k}$ l'élément de $\mathrm{Ext}_{\mathcal{E}^{\mathcal{W}}}^{1}(\mathrm{Gr}_{1}^{k}F,\mathrm{Gr}_{1}^{k+1}F)$ associé à la suite exacte
$$
\hspace{24pt}
0\to\mathrm{Gr}_{1}^{k+1}F\to\mathrm{filt}_{1}^{k}F/\mathrm{filt}_{1}^{k+2}F\to
\mathrm{Gr}_{1}^{k}F\to 0
\hspace{24pt};
$$
un argument formel montre que le produit de Yoneda $\epsilon_{k+1}\smile\epsilon_{k}$ est nul pour $k\leq n-2$. Soit $\mathrm{e}_{k}:\mathrm{I}^{\bullet,k-n}F\hspace{1pt}[k-n]\hspace{1pt}[1]\to\mathrm{I}^{\bullet,k+1-n}F\hspace{1pt}[k+1-n]$ l'homomorphisme de complexes induit par le cobord vertical du bicomplexe $\mathrm{I}^{\bullet,\bullet}F$. On observe que l'on a par définition
$$
\mathrm{I}^{\bullet,q}\hspace{1pt}(\mathrm{filt}_{1}^{k}F/\mathrm{filt}_{1}^{k+2}F)
\hspace{4pt}=\hspace{4pt}
\begin{cases}
\mathrm{I}^{\bullet,q}F & \text{pour $q=k-n,k+1-n$}  \\
0 & \text{pour $q\not=k-n,k+1-n$}
\end{cases}
$$
et que le cobord vertical du bicomplexe $\mathrm{I}^{\bullet,\bullet}\hspace{1pt}(\mathrm{filt}_{1}^{k}F/\mathrm{filt}_{1}^{k+2}F)$ est $\mathrm{e}_{k}\hspace{1pt}[n-k-1]$. Comme  la proposition \ref{resinjF} dit que $\mathrm{Tot}\hspace{1pt}\mathrm{I}^{\bullet,\bullet}\hspace{1pt}(\mathrm{filt}_{1}^{k}F/\mathrm{filt}_{1}^{k+2}F)$ est une résolution injective du foncteur $\mathrm{filt}_{1}^{k}F/\mathrm{filt}_{1}^{k+2}F$, un argument presque aussi formel que celui évoqué précédemment montre que $\mathrm{e}_{k}$ représente $\epsilon_{k}$. Ceci implique que le composé $\mathrm{e}_{k+1}\circ\mathrm{e}_{k}$ est  homotopiquement nul. En fait ce composé est nul~; cette annulation est cruciale dans notre construction du bicomplexe $\mathrm{I}^{\bullet,\bullet}F$ (point (c) de \ref{IpqF}). Bernhard Keller nous a signalé que ce phénomène n'était par contre pas formel. Précisons. Soient $\mathcal{A}$ une catégorie abélienne avec assez d'injectifs et $A$ un objet de $\mathcal{A}$, muni d'une filtration décroissante $A=A^{0}\supset\ldots\supset A^{k}\supset\ldots\supset~A^{n+1}=0$. Soient $\epsilon_{k}$ les éléments de $\mathrm{Ext}_{\mathcal{A}}^{1}(A^{k}/A^{k+1},A^{k+1}/A^{k+2})$ définis comme ci-dessus (qui vérifient encore $\epsilon_{k+1}\smile\epsilon_{k}=0$ pour $k\leq n-2$)~; en général, il n'est pas possible d'exhiber des résolutions injectives $I^{\bullet}_{k}$ de $A^{k}/A^{k+1}$ et des homomorphismes de complexes $e_{k}:I^{\bullet}_{k}[1]\to\nolinebreak I^{\bullet}_{k+1}$ qui représentent $\epsilon_{k}$ et vérifient $e_{k+1}\circ e_{k}=0$. Bernard Keller nous a également indiqué que l'on peut trouver des informations sur cette question dans \cite[3.1]{BBDG}.

\end{rem}
\normalsize

\subsect{Sur les foncteurs de $\mathcal{W}$ dans $V\hspace{-3pt}\text{-}\hspace{1.5pt}\mathcal{U}$ vus comme certaines suites de foncteurs de $\mathcal{W}$ dans $ \mathcal{E}$}\label{FV-U}

\medskip
Le début de cette sous-section est vraiment sans surprises d'où les petits caractères.

\medskip
\footnotesize
Soit $M$ un objet de $V\hspace{-1.5pt}\text{-}\hspace{1pt}\mathcal{U}$. La  donnée de $M$ est équivalente aux données suivantes~:

\medskip
-- une suite $(M^{k})_{k\in\mathbb{N}}$ de $\mathcal{E}$-objets ($M^{k}$ est constitué des éléments de degré $k$ de $M$)~;

\medskip
- des $\mathcal{E}$-morphismes, disons $\mu(k,h):M^{k}\to M^{k+\vert h\vert}$, $h$ désignant un élément (homogène) de $\mathrm{H}^{*}V$ dont le degré est noté $\vert h\vert$ ($\mu(k,h)$ est la multiplication par $h$),

\medskip
- des $\mathcal{E}$-morphismes, disons $\alpha(k,\theta):M^{k}\to M^{k+\vert\theta\vert}$, $\theta$ désignant un élément (homogène) de $\mathrm{A}$ dont le degré est noté $\vert\theta\vert$ ($\alpha(k,\theta)$ est fourni par l'action de $\theta$ sur $M^{k}$),

\medskip
les $\mathcal{E}$-morphismes ci-dessus satisfaisant la liste, disons $(\mathcal{R})$, des relations ci-dessous~:

\smallskip
--\hspace{8pt} $\mu(k,h_{1}+h_{2})=\mu(k,h_{1})+\mu(k,h_{2})$ ($h_{1}$ et $h_{2}$ de même degré)

\smallskip
--\hspace{8pt} $\mu(k,h_{2}h_{1})=\mu(k+\vert h_{1}\vert,h_{2})\circ\mu(k,h_{1})$

\smallskip
--\hspace{8pt} $\alpha(k,\theta_{1}+\theta_{2})=\alpha(k,\theta_{1})+\alpha(k,\theta_{2})$ ($\theta_{1}$ et $\theta_{2}$ de même degré)

\smallskip
--\hspace{8pt} $\alpha(k,\theta_{2}\theta_{1})=\alpha(k+\vert \theta_{1}\vert,\theta_{2})\circ\alpha(k,\theta_{1})$

\smallskip
--\hspace{8pt} $\alpha(k+\vert h\vert,\theta)\circ\mu(k,h)=\sum\mu(k+\vert \theta''\vert,\theta'h)\circ\alpha(k,\theta'')$ avec $\Delta\theta=\sum\theta'\otimes\theta''$, $\Delta:\mathrm{A}\to\mathrm{A}\otimes\mathrm{A}$ désignant la diagonale de l'algèbre de Steenrod.

\bigskip
La discussion (un peu lourde~!) ci-dessus montre que la donnée d'un foncteur $F$ de $\mathcal{W}$ dans $V\hspace{-2pt}\text{-}\hspace{1pt}\mathcal{U}$ est équivalente à celle d'une suite de foncteurs $F^{k}$ de $\mathcal{W}$ dans~$\mathcal{E}$ ($F^{k}(U)$, $U\subset V$, est l'ensemble des éléments de degré $k$ de $F(U)$)  et de transformations naturelles $\mathrm{m}(k,h):F^{k}\to F^{k+\vert h\vert}$,  $h$ désignant un élément (homogène) de $\mathrm{H}^{*}V$,  et $\mathrm{a}(k,\theta):F^{k}\to F^{k+\vert\theta\vert}$, $\theta$ désignant un élément (homogène) de $\mathrm{A}$, transformations naturelles qui satisfont une liste de relations, disons $(\mathcal{S})$, calquée sur la liste $(\mathcal{R})$ ($\mu(-;-)$ est remplacé par $\mathrm{m}(-,-)$ et $\alpha(-,-)$ par $\mathrm{a}(-,-)$).
\normalsize

\medskip
\begin{defi}\label{HomSF} Soient $S$ un foncteur de $\mathcal{W}$ dans $\mathcal{E}$ et $F$ un foncteur de~$\mathcal{W}$ dans $V\hspace{-3pt}\text{-}\hspace{2pt}\mathcal{U}$. On note $\mathrm{Hom}_{\mathcal{E}^{\mathcal{W}}}(S,F)$ le $\mathrm{H}^{*}V$-$\mathrm{A}$-module instable dont le $\mathbb{F}_{2}$-espace vectoriel des éléments de degré $k$ et $\mathrm{Hom}_{\mathcal{E}^{\mathcal{W}}}(S,F^{k})$ et dont la structure est donnée par les $\mathcal{E}$-morphismes $\mathrm{Hom}_{\mathcal{E}^{\mathcal{W}}}(S, \mathrm{m}(k,h))$ et $\mathrm{Hom}_{\mathcal{E}^{\mathcal{W}}}(S, \mathrm{a}(k,\theta))$.
\end{defi}

\medskip
\begin{exple}\label{HomPWF} Soit $W$ un sous-groupe de $V$ et  $F$ un foncteur de $\mathcal{W}$ dans $V\hspace{-3pt}\text{-}\hspace{2pt}\mathcal{U}$~; on a $\mathrm{Hom}_{\mathcal{E}^{\mathcal{W}}}(\mathbb{P}_{W},F)=F(W)$.
\end{exple}

\medskip
\begin{pro}\label{HomStf} Soient $S$ un foncteur de $\mathcal{W}$ dans $\mathcal{E}$ et $F$ un foncteur de $\mathcal{W}$ dans $V_{\mathrm{tf}}\text{-}\mathcal{U}$. Si $S(W)$ est de dimension finie pour tout sous-groupe $W$ de $V$ alors $\mathrm{Hom}_{\mathcal{E}^{\mathcal{W}}}(S,F)$ est de type fini comme $\mathrm{H}^{*}V$-module.
\end{pro}

\medskip
\textit{Démonstration.} On considère le $\mathcal{E}^{\mathcal{W}}$-épimorphisme canonique
$$
\bigoplus_{W\in\mathcal{W}}S(W)\otimes\mathbb{P}_{W}
\twoheadrightarrow S
$$
(on observera incidemment que l'existence de cet épimorphisme entraîne que la catégorie $\mathcal{E}^{\mathcal{W}}$ a assez de projectifs). Cet épimorphisme induit un $(V\hspace{-3pt}\text{-}\hspace{2pt}\mathcal{U})$-monomorphisme
$$
\mathrm{Hom}_{\mathcal{E}^{\mathcal{W}}}(S,F)
\hookrightarrow
\bigoplus_{W\in\mathcal{W}}(S(W))^{*}\otimes F(W)
$$
(observer que $F(W)$ est de dimension finie en chaque degré). On conclut en invoquant le fait que $\mathrm{H}^{*}V$ est noethérien.
\hfill$\square$

\medskip
La proposition-définition ci-dessous étend la définition \ref{HomSF}. La première partie de la proposition est immédiate~; la seconde résulte du fait que si $I$ est un objet injectif de $(V\hspace{-3pt}\text{-}\hspace{2pt}\mathcal{U})^{\mathcal{W}}$ alors, pour tout entier $k\geq 0$, $I^{k}$ ($I$ en degré $k$) est un objet injectif de $\mathcal{E}^{\mathcal{W}}$.

\medskip
\begin{pro-def}\label{ExtSF} Soient  $p$ un entier naturel, $S$ un foncteur de $\mathcal{W}$ dans $\mathcal{E}$ et $F$ un foncteur de~$\mathcal{W}$ dans $V\hspace{-3pt}\text{-}\hspace{2pt}\mathcal{U}$.

\smallskip
On note $\mathrm{Ext}_{\mathcal{E}^{\mathcal{W}}}^{p}(S,F)$ le $\mathrm{H}^{*}V$-$\mathrm{A}$-module instable dont le $\mathbb{F}_{2}$-espace vectoriel des éléments de degré $k$ est $\mathrm{Ext}_{\mathcal{E}^{\mathcal{W}}}^{p}(S,F^{k})$ et dont la structure est donnée par les $\mathcal{E}$-morphismes $\mathrm{Ext}_{\mathcal{E}^{\mathcal{W}}}^{p}(S, \mathrm{m}(k,h))$ et $\mathrm{Ext}_{\mathcal{E}^{\mathcal{W}}}^{p}(S, \mathrm{a}(k,\theta))$.

\smallskip
Soit $P_{\bullet}\to S$ une résolution projective dans la catégorie $\mathcal{E}^{\mathcal{W}}$ alors on a un isomorphisme canonique
$$
\hspace{24pt}
\mathrm{Ext}_{\mathcal{E}^{\mathcal{W}}}^{p}(S,F)
\hspace{4pt}\cong\hspace{4pt}
\mathrm{H}^{p}\hspace{2pt}
\mathrm{Hom}_{\mathcal{E}^{\mathcal{W}}}(P_{\bullet},F)
\hspace{24pt}.
$$

\smallskip
Soit $F\to I^{\bullet}$ une résolution injective dans la catégorie $(V\hspace{-3pt}\text{-}\hspace{2pt}\mathcal{U})^{\mathcal{W}}$ alors on a un isomorphisme canonique
$$
\hspace{24pt}
\mathrm{Ext}_{\mathcal{E}^{\mathcal{W}}}^{p}(S,F)
\hspace{4pt}\cong\hspace{4pt}
\mathrm{H}^{p}\hspace{2pt}
\mathrm{Hom}_{\mathcal{E}^{\mathcal{W}}}(S,I^{\bullet})
\hspace{24pt}.
$$
\end{pro-def}

\medskip
La  proposition \ref{Pf-bis}  et le corollaire \ref{derivePf-ter} (avatar de \ref{derivePf-bis}) ci-après, illustrent les définitions \ref{HomSF} et \ref{ExtSF}. Soit $M$ un $\mathrm{H}^{*}V$-$\mathrm{A}$-module instable~; on rapelle que la notation $\widehat{\Psi}_{M}$ désigne le foncteur de $\mathcal{W}$ dans $V\hspace{-2.5pt}\text{-}\hspace{0.5pt}\mathcal{U}$, $W\mapsto\mathrm{EFix}_{(V,W)}M$.

\medskip
\begin{pro}\label{Pf-bis} Soit $M$ un $\mathrm{H}^{*}V_{\mathrm{tf}}$-$\mathrm{A}$-module instable. On a un isomorphisme canonique de $\mathrm{H}^{*}V_{\mathrm{tf}}$-$\mathrm{A}$-modules instables
$$
\hspace{24pt}
\mathrm{Pf}\hspace{1pt}M
\hspace{4pt}\cong\hspace{4pt}
\mathrm{Hom}_{\mathcal{E}^{\mathcal{W}}}(\mathrm{S}_{0},\widehat{\Psi}_{M})
\hspace{24pt}.
$$
\end{pro}

\textit{Démonstration.}  Soit $F$ un foncteur de $\mathcal{W}$ dans $\mathcal{E}$~; on a
$$
\hspace{24pt}
\mathrm{Hom}_{\mathcal{E}^{\mathcal{W}}}(\mathrm{S}_{0},F)
\hspace{4pt}\cong\hspace{4pt}
\mathrm{H}^{0}\hspace{2pt}
\mathrm{Hom}_{\mathcal{E}^{\mathcal{W}}}(\mathrm{P}_{\bullet},F)
\hspace{24pt},
$$
$\mathrm{P}_{\bullet}\to\mathrm{S}_{0}$ étant la résolution projective de \ref{resprojS0}, c'est-à dire
$$
\hspace{24pt}
\mathrm{Hom}_{\mathcal{E}^{\mathcal{W}}}(\mathrm{S}_{0},F)
\hspace{4pt}\cong\hspace{4pt}
\bigcap_{\dim W=1}
\ker\hspace{1pt}(F(0)\to F(W))
\hspace{24pt}.
$$
Soit maintenant $M$ un $\mathrm{H}^{*}V$-$\mathrm{A}$-module instable. L'isomorphisme ci-dessus donne
$$
\hspace{24pt}
\mathrm{Hom}_{\mathcal{E}^{\mathcal{W}}}(\mathrm{S}_{0},\widehat{\Psi}_{M})
\hspace{4pt}\cong\hspace{4pt}
\bigcap_{\dim W=1}
\ker\hspace{1pt}(M\to \mathrm{EFix}_{(V,W)}M)
\hspace{24pt}.
$$
Par définition le membre de droite est le terme $\mathrm{F}^{n}M$ de la filtration décrois\-sante de $M$ introduite au début de la section 5. Or on a $\mathrm{F}^{n}M=\mathrm{Pf}M$ si $M$ est de type fini comme $\mathrm{H}^{*}V$-module (Proposition \ref{Fn=Pf}).
\hfill
$\square$

\medskip
\begin{cor}\label{derivePf-ter} Soit $M$ un $\mathrm{H}^{*}V_{\mathrm{tf}}$-$\mathrm{A}$-module instable. On a pour tout entier $k\geq 0$,  un isomorphisme canonique de $\mathrm{H}^{*}V_{\mathrm{tf}}$-$\mathrm{A}$-modules instables
$$
\hspace{24pt}
\mathrm{R}^{k}\mathrm{Pf}\hspace{1pt}M
\hspace{4pt}\cong\hspace{4pt}
\mathrm{Ext}_{\mathcal{E}^{\mathcal{W}}}^{k}(\mathrm{S}_{0},\widehat{\Psi}_{M})
\hspace{24pt}.
$$
\end{cor}

\textit{Démonstration.} Soit $M\to I^{\bullet}$ une résolution injective dans la catégorie $V_{\mathrm{tf}}\text{-}\mathcal{U}$ alors $,\widehat{\Psi}_{M}\to\widehat{\Psi}_{I^{\bullet}}$ est une résolution injective dans la catégorie $(V_{\mathrm{tf}}\text{-}\mathcal{U})^{\mathcal{W}}$ (remplacer $\Psi$ par $\widehat{\Psi}$ dans la démonstration des points (a) et (b) de \ref{derivePf}). On conclut à l'aide de \ref{ExtSF} (dernière partie de la proposition) et \ref{Pf-bis}.
\hfill$\square$

\bigskip
\begin{rem}\label{derivPf-quatro} Les énoncés \ref{Klim3} et \ref{derivePf-ter} redonnent \ref{derivePf}.
\end{rem}

\subsect{Construction du bicomplexe $\mathrm{B}^{p,q}M$}\label{Bpq}

\medskip
Soient $M$ un $\mathrm{H}^{*}V$-$\mathrm{A}$-module instable et $\mathrm{F}^{p}M$ le $p$-ième terme de la filtration décroissante de $M$ introduite en section 5. Par définition on a~:
$$
\hspace{24pt}
\mathrm{F}^{p}M
\hspace{4pt}=\hspace{4pt}
\bigcap_{\mathop{\mathrm{codim}}W<p}
\ker\hspace{1pt}(M=\widehat{\Psi}_{M}(0)\to\widehat{\Psi}_{M}(W))
\hspace{24pt};
$$
ceci motive la définition ci-après.

\medskip
\begin{defi}\label{filt2} Soit $F$ un foncteur de $\mathcal{W}$ dans $\mathcal{E}$ ou $V\hspace{-2.5pt}\text{-}\hspace{0.5pt}\mathcal{U}$. On définit une filtration décroissante de $F(0)$~:
$$
F(0)=\mathrm{filt}_{2}^{0}\hspace{1pt}F(0)\supset\ldots\supset\mathrm{filt}_{2}^{p}\hspace{1pt}F(0)\supset\ldots\supset\mathrm{filt}_{2}^{n}\hspace{1pt}F(0)\supset\mathrm{filt}_{2}^{n+1}\hspace{1pt}F(0)=0
$$
en prenant
$$
\hspace{24pt}
\mathrm{filt}_{2}^{p}\hspace{1pt}F(0)
\hspace{4pt}:=\hspace{4pt}
\bigcap_{\mathop{\mathrm{codim}}W<p}
\ker\hspace{1pt}(F(0)\to F(W))
\hspace{24pt}.
$$
On pose $\mathrm{Gr}_{2}^{p}\hspace{1pt}F(0):=\mathrm{filt}_{2}^{p}\hspace{1pt}F(0)/\mathrm{filt}_{2}^{p+1}\hspace{1pt}F(0)$.
\end{defi}

\bigskip
\begin{rem}\label{filt2-bis}
La filtration de $F(0)$ que nous venons d'introduire se prolonge en fait en une filtration décroissante de~$F$~:
$$
F={\mathop{\mathrm{filt}}}_{2}^{0}F\supset\ldots\supset{\mathop{\mathrm{filt}}}_{2}^{p}F\supset\ldots\supset{\mathop{\mathrm{filt}}}_{2}^{n}F\supset{\mathop{\mathrm{filt}}}_{2}^{n+1}F=0
$$
en prenant
$$
\hspace{24pt}
({\mathop{\mathrm{filt}}}_{2}^{p}F)(U)
\hspace{4pt}:=\hspace{4pt}
\bigcap_{\mathop{\mathrm{codim}}W<p}
\ker\hspace{1pt}(F(U)\to F(U+W))
\hspace{24pt}.
$$
Il est clair que l'on a tout fait pour avoir $\mathrm{filt}_{2}^{p}\hspace{1pt}F(0)=({\mathop{\mathrm{filt}}}_{2}^{p}F)(0)$. Si  $M$ est un $\mathrm{H}^{*}V_{\mathrm{tf}}$-$\mathrm{A}$-module instable alors ${\mathop{\mathrm{filt}}}_{2}^{p}\hspace{1pt}\widehat{\Psi}_{M}$ s'identifie à $\widehat{\Psi}_{\mathrm{F}^{p}M}$. Pour s'en convaincre observer que le foncteur $M\mapsto\widehat{\Psi}_{M}$ est exact et que \ref{compoEFix} implique l'identification $\widehat{\Psi}_{\mathrm{EFix}_{(V,W)}M}(U)=\widehat{\Psi}_{M}(U+W)$.
\end{rem}

\bigskip
Le formalisme ci-dessus permet de définir pour tout  foncteur $F:\mathcal{W}\to\mathcal{A}$, avec $\mathcal{A}=\mathcal{E}$ ou $\mathcal{A}=V\hspace{-2.5pt}\text{-}\hspace{0.5pt}\mathcal{U}$,  un $\mathcal{A}$-complexe de cochaînes, que nous notons  $\mathrm{C}_{2}^{\bullet}F$. La définition de $\mathrm{C}_{2}^{\bullet}F$ imite celle du complexe $\mathrm{C}^{\bullet}M$, associé en section 5 à tout $\mathrm{H}^{*}V$-$\mathrm{A}$-module instable $M$, de telle sorte que l'on ait $\mathrm{C}_{2}^{\bullet}\hspace{1pt}\widehat{\Psi}_{M}=\mathrm{C}^{\bullet}M$.

\medskip
Explicitons. Soit $F\to I^{\bullet}$ une résolution injective de $F$ dans la catégorie $\mathcal{A}^\mathcal{W}$~; $\mathrm{C}_{2}^{\bullet}F$ est  le complexe de cochaînes dans la catégorie $\mathcal{A}$ dont le  $p$-ième terme est $\mathrm{H}^{p}\mathrm{Gr}_{2}^{p}\hspace{1pt}I^{\bullet}(0)$ et dont le cobord est donné par le connectant associé à la suite exacte de $\mathcal{A}$-complexes
$$
\hspace{24pt}
0\to\mathrm{Gr}_{2}^{p+1}\hspace{1pt}I^{\bullet}(0)\to\mathrm{filt}_{2}^{p}\hspace{1pt}I^{\bullet}(0)/\mathrm{filt}_{2}^{p+2}\hspace{1pt}I^{\bullet}(0)\to\mathrm{Gr}_{2}^{p}\hspace{1pt}I^{\bullet}(0)\to 0
\hspace{24pt}.
$$
Là encore, les mantras de la théorie des résolutions injectives montrent que le complexe $\mathrm{C}_{2}^{\bullet}F$ est indépendant du choix de $I^{\bullet}$.

\medskip
\begin{pro}\label{pre-C2Psihat} Si $F$ un foncteur de $\mathcal{W}$ dans $V\hspace{-3pt}\text{-}\hspace{1pt}\mathcal{U}$ alors on a
$$
\hspace{24pt}
\mathrm{C}_{2}^{\bullet}F
\hspace{4pt}=\hspace{4pt}
\mathrm{C}^{\bullet}F(0)
\hspace{24pt}.
$$
\end{pro}

\textit{Démonstration.} Observer que si $F\to I^{\bullet}$ est une résolution injective de $F$ dans la catégorie $(V\hspace{-3pt}\text{-}\hspace{1pt}\mathcal{U})^\mathcal{W}$ alors  $F(0)\to I^{\bullet}(0)$ est une résolution injective de~$F$ dans la catégorie $V\hspace{-3pt}\text{-}\hspace{1pt}\mathcal{U}$ (voir le point (d) de \ref{co-induit}).
\hfill$\square$

\medskip
\begin{cor}\label{C2Psihat} Soit $M$ un $\mathrm{H}^{*}V$-$\mathrm{A}$-module instable~; on a
$$
\hspace{24pt}
\mathrm{C}_{2}^{\bullet}\hspace{1pt}\widehat{\Psi}_{M}
\hspace{4pt}=\hspace{4pt}
\mathrm{C}^{\bullet}M
\hspace{24pt}.
$$
\end{cor}

Soient $M$ un $\mathrm{H}^{*}V$-$\mathrm{A}$-module instable et $M\to I^{\bullet}$ une résolution injective de~$M$ dans la catégorie $V\hspace{-2pt}\text{-}\mathcal{U}$. Nous avons mis en oeuvre dans la démonstration de l'implication (i)$\Rightarrow$(ii) du théorème \ref{pendantalg} la suite spectrale définie par la filtration du complexe $I^{\bullet}$ par les $\mathrm{F}^{p}I^{\bullet}$~;  nous notons ci-après $\mathrm{E}_{1}^{\bullet,\bullet}M$ son terme~$\mathrm{E}_{1}$. La ligne $\mathrm{E}_{1}^{\bullet,q}M$, munie de la différentielle~$\mathrm{d}_{1}$, est un complexe de cochaines dans la catégorie $V\hspace{-2pt}\text{-}\mathcal{U}$ et~$\mathrm{C}^{\bullet}M$ n'est rien d'autre que $\mathrm{E}_{1}^{\bullet,0}M$.

\smallskip
Soit maintenant $F$ un foncteur de $\mathcal{W}$ dans $\mathcal{A}$, avec $\mathcal{A}=\mathcal{E}$ ou $\mathcal{A}=V\hspace{-2pt}\text{-}\mathcal{U}$, et $F\to I^{\bullet}$ une résolution injective de $F$ dans la catégorie $\mathcal{A}^\mathcal{W}$~; nous notons~$\mathrm{E}_{1}^{\bullet,\bullet}F$ le terme $\mathrm{E}_{1}$ de la suite spectrale associée à la filtration du complexe~$I^{\bullet}(0)$ par les ${\mathop{\mathrm{filt}}}_{2}^{p}\hspace{2pt}I^{\bullet}(0)$. A nouveau $\mathrm{E}_{1}^{\bullet,\bullet}F$ ne dépend pas du choix de~$I^{\bullet}$ et  $\mathrm{C}_{2}^{\bullet}F$ n'est rien d'autre que $\mathrm{E}_{1}^{\bullet,0}F$. La proposition \ref{pre-C2Psihat} et le corollaire \ref{C2Psihat} se généralisent immédiatement~:

\medskip
\begin{pro}\label{pre-E1Psihat} Si $F$ est un foncteur de $\mathcal{W}$ dans $V\hspace{-2pt}\text{-}\mathcal{U}$ alors on a
$$
\hspace{24pt}
\mathrm{E}_{1}^{\bullet,\bullet}F
\hspace{4pt}=\hspace{4pt}
\mathrm{E}_{1}^{\bullet,\bullet}F(0)
\hspace{24pt}.
$$
\end{pro}

\begin{cor}\label{E1Psihat} Soit $M$ un $\mathrm{H}^{*}V$-$\mathrm{A}$-module instable~; on a
$$
\hspace{24pt}
\mathrm{E}_{1}^{\bullet,\bullet}\hspace{1pt}\widehat{\Psi}_{M}
\hspace{4pt}=\hspace{4pt}
\mathrm{E}_{1}^{\bullet,\bullet}M
\hspace{24pt}.
$$
\end{cor}

\begin{exple}\label{C2SV} Le $\mathcal{E}$-complexe $\mathrm{C}_{2}^{\bullet}\hspace{1pt}\mathrm{S}_{V}$ est canoniquement isomorphe au complexe de Lusztig $\mathrm{Lu}^{\bullet}_{2,V}$ et $\mathrm{E}_{1}^{\bullet,q}\hspace{1pt}$ est nul pour $q\not=0$ . On vérifie ces affirmations ci-après. On considère la résolution injective $\mathrm{S}_{V}\to\mathrm{I}^{\bullet}$ de \ref{resinjSV}. Comme l'on a $\mathbb{I}^{W}(0)=\mathbb{F}_{2}$ pour tout sous-groupe $W$ de $V$, on constate que l'on a $\mathrm{I}^{\bullet}(0)=\mathrm{Lu}^{\bullet}_{2,V}$. On a d'autre part
$$
{\mathop{\mathrm{filt}}}_{2}^{p}\hspace{2pt}\mathbb{I}^{W}(0)
\hspace{4pt}=\hspace{4pt}
\begin{cases}
\mathbb{I}^{W}(0)& \text{pour $\mathop{\mathrm{codim}}W\geq p$,}  \\
0 & \text{pour $\mathop{\mathrm{codim}}W< p$,}
\end{cases}
$$
si bien que la filtration de $\mathrm{I}^{\bullet}(0)$ par les ${\mathop{\mathrm{filt}}}_{2}^{p}\hspace{2pt}\mathrm{I}^{\bullet}(0)$ coïncide avec celle par les troncations brutales $\sigma^{\geq p}\hspace{2pt}\mathrm{I}^{\bullet}(0)$ (pour la définition de la troncation brutale on pourra se reporter à la démonstration du point (a) de \ref{efinispecial1}). On conclut à l'aide du lemme \ref{E1troncbrutal}.
\end{exple}

\bigskip
Nous en arrivons à présent à la définition du bicomplexe $\mathrm{B}^{\bullet,\bullet}M$, $M$ désignant un $\mathrm{H}^{*}V$-$\mathrm{A}$-module instable, évoqué dans la proposition \ref{bicomplexe}~; nous suivons toujours la même stratégie~:

\smallskip
-- nous définissons d'abord $\mathrm{B}^{\bullet,\bullet}F$ pour $F$ un foncteur de $\mathcal{W}$ dans $\mathcal{E}$, 

\smallskip
-- nous étendons ensuite la définition ci-dessus aux foncteurs $F$ de $\mathcal{W}$ dans~$V\hspace{-2pt}\text{-}\mathcal{U}$ grâce au yoga de la sous-section 8.4,

\smallskip
-- nous faisons enfin $F=\widehat{\Psi}_{M}$.

\pagebreak
\bigskip
\begin{defi}\label{Bpq-1} Soit $F$ un foncteur de $\mathcal{W}$ dans $\mathcal{E}$. On pose
$$
\hspace{24pt}
\mathrm{B}^{\bullet,\bullet}F
\hspace{4pt}:=\hspace{4pt}
(\mathrm{I}^{\bullet,\bullet}F)(0)
\hspace{24pt},
$$
$\mathrm{I}^{\bullet,\bullet}F$ désignant  le $\mathcal{E}^{\mathcal{W}}$-bicomplexe introduit en \ref{IpqF}~;  $\mathrm{B}^{\bullet,\bullet}F$ est donc un $\mathcal{E}$-bicomplexe de cochaînes avec
$$
\hspace{24pt}
\mathrm{B}^{p,q}F
\hspace{4pt}=\hspace{4pt}
\bigoplus_{W\subset Z,\hspace{4pt}\mathop{\mathrm{codim}}W=p,\hspace{4pt}\mathop{\mathrm{codim}}Z=-q}
\mathrm{St}_{Z/W}^{*}\otimes F(Z)
\hspace{24pt},
$$
dont les cobords horizontaux et verticaux sont obtenus ``en évaluant en $0$'' ceux de \ref{IpqF}.
\end{defi}

\medskip
\begin{pro-def}\label{Bpq-2} Soit $F$ un foncteur de $\mathcal{W}$ dans $\mathcal{E}$.

\smallskip
Soit ${}_{\mathrm{B}}\mathrm{E}_{1}^{\bullet,\bullet}F$ le terme $\mathrm{E}_{1}^{\bullet,\bullet}$ de la suite spectrale associée à $\mathrm{B}^{\bullet,\bullet}F$ obtenue en dérivant ``verticalement''. Alors on a un isomorphisme canonique
$$
\hspace{24pt}
{}_{\mathrm{B}}\mathrm{E}_{1}^{\bullet,\bullet}F
\hspace{4pt}\cong\hspace{4pt}
\mathrm{E}_{1}^{\bullet,\bullet}F
\hspace{24pt},
$$
préservant les différentielles $\mathrm{d}_{1}$ (en clair, on considère ici ${}_{\mathrm{B}}\mathrm{E}_{1}^{\bullet,\bullet}F$ et $\mathrm{E}_{1}^{\bullet,\bullet}F$ comme des bicomplexes avec cobord vertical nul, suivant en cela {\em \cite{CE}}\hspace{1pt}).
\end{pro-def}

\bigskip
\textit{Démonstration}

\medskip
1) D'une part la proposition \ref{resinjF} dit que $F\overset{\eta}{\to}\mathrm{Tot}\hspace{1pt}\mathrm{I}^{\bullet,\bullet}F$ est une résolution injective de $F$ dans la catégorie $\mathcal{E}^{\mathcal{W}}$ si bien que $\mathrm{E}_{1}^{\bullet,\bullet}F$ peut être vu comme le terme $\mathrm{E}_{1}$ de la suite spectrale définie par la filtration du $\mathcal{E}$-complexe $(\mathrm{Tot}\hspace{1pt}\mathrm{I}^{\bullet,\bullet}F)(0)$ par les ${\mathop{\mathrm{filt}}}_{2}^{p}\hspace{2pt}(\mathrm{Tot}\hspace{1pt}\mathrm{I}^{\bullet,\bullet}F)(0)$.

\medskip
2) D'autre part ${}_{\mathrm{B}}\mathrm{E}_{1}^{\bullet,\bullet}F$ est obtenu de la façon suivante~:

\smallskip
Soit $p\geq 0$ un entier~; on considère le sous-bicomplexe $\mathrm{F}^{p}_{\mathrm{II}}(\mathrm{I}^{\bullet,\bullet}F)(0)$ de  $(\mathrm{I}^{\bullet,\bullet}F)(0)$ (notation de \cite{CE}) défini par
$$
{(\mathrm{F}^{p}_{\mathrm{II}}(\mathrm{I}^{\bullet,\bullet}F)(0))}^{s,t}
=
\begin{cases}
(\mathrm{I}^{s,t}F)(0) & \text{pour $s\geq p$}, \\
0 & \text{pour $s<p$}.
\end{cases}
$$
La filtration (décroissante) du bicomplexe $(\mathrm{I}^{\bullet,\bullet}F)(0)$ par les $\mathrm{F}^{p}_{\mathrm{II}}(\mathrm{I}^{\bullet,\bullet}F)(0)$ induit une filtration du complexe $\mathrm{Tot}\hspace{1pt}(\mathrm{I}^{\bullet,\bullet}F)(0)=(\mathrm{Tot}\hspace{1pt}\mathrm{I}^{\bullet,\bullet}F)(0)$ par des sous-complexes $\mathrm{F}^{p}_{\mathrm{II}}(\mathrm{Tot}\hspace{1pt}\mathrm{I}^{\bullet,\bullet}F)(0):=\mathrm{Tot}\hspace{1pt}(\mathrm{F}^{p}_{\mathrm{II}}\mathrm{I}^{\bullet,\bullet}F)(0)$~; ${}_{\mathrm{B}}\mathrm{E}_{1}^{\bullet,\bullet}F$ est le terme $\mathrm{E}_{1}$ de la suite spectrale définie par la filtration du $\mathcal{E}$-complexe $(\mathrm{Tot}\hspace{1pt}\mathrm{I}^{\bullet,\bullet}F)(0)$ par les $\mathrm{F}^{p}_{\mathrm{II}}(\mathrm{Tot}\hspace{1pt}\mathrm{I}^{\bullet,\bullet}F)(0)$.

\medskip
Or on constate que l'on a
$$
\hspace{24pt}
{\mathop{\mathrm{filt}}}_{2}^{p}\hspace{2pt}(\mathrm{Tot}\hspace{1pt}\mathrm{I}^{\bullet,\bullet}F)(0)
\hspace{4pt}=\hspace{4pt}
\mathrm{F}^{p}_{\mathrm{II}}(\mathrm{Tot}\hspace{1pt}\mathrm{I}^{\bullet,\bullet}F)(0)
\hspace{24pt}.
$$
Ceci résulte de l'argument déjà utilisé dans l'exemple \ref{C2SV}.~:
$$
{\mathop{\mathrm{filt}}}_{2}^{p}\hspace{2pt}\mathbb{I}^{W}(0)
\hspace{4pt}=\hspace{4pt}
\begin{cases}
\mathbb{I}^{W}(0)& \text{pour $\mathop{\mathrm{codim}}W\geq p,$}  \\
0 & \text{pour $\mathop{\mathrm{codim}}W< p$}.
\end{cases}
$$
\hfill$\square$

\bigskip
Soit $F$ un foncteur de $\mathcal{W}$ dans $V\hspace{-2pt}\text{-}\mathcal{U}$~;  on définit $\mathrm{B}^{\bullet,\bullet}F$, qui est maintenant un bicomplexe de cochaînes dans la catégorie $V\hspace{-2pt}\text{-}\mathcal{U}$, en posant
$$
\mathrm{B}^{p,q}F
\hspace{4pt}:=\hspace{4pt}
\bigoplus_{W\subset Z,\hspace{4pt}\mathop{\mathrm{codim}}W=p,\hspace{4pt}\mathop{\mathrm{codim}}Z=-q}
\mathrm{St}_{Z/W}^{*}\otimes F(Z)
$$
et en explicitant les cobords, horizontaux et verticaux, ``par les mêmes formules'' que précédemment. Alternativement, on peut définir $\mathrm{B}^{\bullet,\bullet}F$ en adoptant le point de vue de 8.4. Donnons quelques détails. Soit $k$ un entier naturel~; on note $\mathrm{deg}^{k}:V\hspace{-2pt}\text{-}\mathcal{U}\to\mathcal{E}$ le foncteur  $M\mapsto M^{k}$ et $\mathrm{Deg}^{k}:(V\hspace{-2pt}\text{-}\mathcal{U})^{\mathcal{W}}\to\mathcal{E}^{\mathcal{W}}$ le foncteur $F\mapsto\mathrm{deg}^{k}\circ F$. On prend $\mathrm{deg}^{k\hspace{1pt}}\mathrm{B}^{\bullet,\bullet}F:=\mathrm{B}^{\bullet,\bullet}\hspace{1pt}\mathrm{Deg}^{k}F$, la $(V\hspace{-2pt}\text{-}\mathcal{U})$-structure de $\mathrm{B}^{\bullet,\bullet}F$ est obtenue en appliquant le foncteur $\mathrm{B}^{\bullet,\bullet}$, défini sur la catégorie $\mathcal{E}^{\mathcal{W}}$ et à valeurs dans la catégorie des $\mathcal{E}$-bicomplexes de cochaînes, aux $\mathcal{E}^{\mathcal{W}}$-morphismes $\mathrm{Deg}^{k}F\to\mathrm{Deg}^{k'}F$ donnés par la $(V\hspace{-2pt}\text{-}\mathcal{U})$-structure des $F(U)$, $U\subset V$ (voir la sous-section 8.4). Ce point de vue permet d'étendre l'énoncé \ref{Bpq-2} aux foncteurs de $\mathcal{W}$ dans $V\hspace{-2pt}\text{-}\mathcal{U}$~:

\medskip
\begin{pro-def}\label{Bpq-3} Soit $F$ un foncteur de $\mathcal{W}$ dans $V\hspace{-2pt}\text{-}\mathcal{U}$.

\smallskip
Soit ${}_{\mathrm{B}}\mathrm{E}_{1}^{\bullet,\bullet}F$ le terme $\mathrm{E}_{1}^{\bullet,\bullet}$ de la suite spectrale associée à $\mathrm{B}^{\bullet,\bullet}F$ obtenue en dérivant ``verticalement''. Alors on a un isomorphisme canonique
$$
\hspace{24pt}
{}_{\mathrm{B}}\mathrm{E}_{1}^{\bullet,\bullet}F
\hspace{4pt}\cong\hspace{4pt}
\mathrm{E}_{1}^{\bullet,\bullet}F
\hspace{24pt},
$$
préservant les différentielles $\mathrm{d}_{1}$.
\end{pro-def}

\medskip
\textit{Démonstration.} On a par définition $\mathrm{deg}^{k}\hspace{1pt}\mathrm{B}^{\bullet,\bullet}F=\mathrm{B}^{\bullet,\bullet}\hspace{1pt}\mathrm{Deg}^{k}F$. D'autre part, comme le foncteur $\mathrm{Deg}^{k}$ est exact et envoit injectifs sur injectifs, on a un isomorphisme canonique $\mathrm{deg}^{k}\hspace{1pt}\mathrm{E}_{1}^{\bullet,\bullet}F\cong\mathrm{E}_{1}^{\bullet,\bullet}\hspace{1pt}\mathrm{Deg}^{k}F$, naturel en $F$. 
\hfill$\square$

\bigskip
\begin{defi}\label{Bpq-4} Soit $M$ un $\mathrm{H}^{*}V$-$\mathrm{A}$-module instable. On pose
$$
\hspace{24pt}
\mathrm{B}^{\bullet,\bullet}M
\hspace{4pt}:=\hspace{4pt}
\mathrm{B}^{\bullet,\bullet}\hspace{2pt}\widehat{\Psi}_{M}
\hspace{24pt};
$$
$\mathrm{B}^{\bullet,\bullet}M$ est donc un bicomplexe de cochaînes dans la catégorie $V\hspace{-2pt}\text{-}\mathcal{U}$ dont les  termes sont donnés par l'égalité
$$
\hspace{24pt}
\mathrm{B}^{p,q}M
\hspace{4pt}:=\hspace{4pt}
\bigoplus_{W\subset Z,\hspace{4pt}\mathop{\mathrm{codim}}W=p,\hspace{4pt}\mathop{\mathrm{codim}}Z=-q}
\mathrm{St}_{Z/W}^{*}\otimes\mathrm{EFix}_{(V,W)}M
\hspace{24pt}.
$$
\end{defi}

\bigskip
\textit{Vérification des affirmations de la proposition \ref{bicomplexe}}

\medskip
-- Si $M$ est un $\mathrm{H}^{*}V_{\mathrm{tf}}$-$\mathrm{A}$-module instable alors il en est de même pour les $\mathrm{B}^{p,q}M$ d'après \ref{typefiniFix}.

\medskip
-- La propriété $(\mathcal{B}_{0})$ est satisfaite par définition.

\medskip
-- Compte tenu de l'égalité ${}_{\mathrm{B}}\mathrm{E}_{1}^{\bullet,\bullet}M={}_{\mathrm{B}}\mathrm{E}_{1}^{\bullet,\bullet}\hspace{2pt}\widehat{\Psi}_{M}$, la propriété $(\mathcal{B}_{1})$ est un cas particulier de \ref{Bpq-3}.

\medskip
 -- L'isomorphisme $\mathrm{EFix}_{(V,W)}M\cong\mathrm{H}^{*}V\otimes_{\mathrm{H}^{*}V/W}\mathrm{Fix}_{(V,W)}M$ de  \ref{EFix} implique immédiatement la propriété $(\mathcal{B}_{2})$.

\bigskip
\begin{rem}\label{Bpq-5} Soit $M$ un $\mathrm{H}^{*}V$-$\mathrm{A}$-module instable, l'isomorphisme
$$
{}_{\mathrm{B}}\mathrm{E}_{1}^{\bullet,\bullet}\hspace{2pt}\widehat{\Psi}_{M}
\hspace{4pt}\cong\hspace{4pt}
\mathrm{E}_{1}^{\bullet,\bullet}M
$$
apporte un éclairage concret sur le complexe $\mathrm{C}^{\bullet}M$ introduit en section 5. Expliquons pourquoi. Comme l'on a, par définition, $\mathrm{C}^{\bullet}M=\mathrm{E}_{1}^{\bullet,0}M$, l'isomorphisme ci-dessus dit en particulier que l'on dispose d'un isomorphisme 
$$
\hspace{24pt}
\mathrm{C}^{\bullet}M
\hspace{4pt}\cong\hspace{4pt}
\mathop{\mathrm{coker}}\hspace{2pt}(\mathrm{d}_{\mathrm{h}}^{\bullet,-1}:
\mathrm{B}^{\bullet,-1}M\to\mathrm{B}^{\bullet,0}M)
\hspace{24pt}.
$$
Toujours par définition, on a
$$
\mathrm{B}^{\bullet,0}M
\hspace{4pt}=\hspace{4pt}
\mathrm{Lu}_{2,V}^{\bullet}
\otimes\mathrm{EFix}_{(V,V)}M
$$
et
$$
\hspace{24pt}
\mathrm{B}^{\bullet,-1}M
\hspace{4pt}=\hspace{4pt}
\bigoplus_{\mathop{\mathrm{codim}}Z=1}
\mathrm{Lu}_{2,Z}^{\bullet}[-1]
\otimes\mathrm{EFix}_{(V,Z)}M
\hspace{24pt}.
$$
L'homomorphisme de $(V\hspace{-2pt}\text{-}\mathcal{U})$-complexes $\mathrm{d}_{\mathrm{h}}^{\bullet,-1}$ est induit par

\smallskip
--  les homomorphismes de $\mathcal{E}$-complexes $\mathrm{Lu}_{2,Z}^{\bullet}[-1]\to\mathrm{Lu}_{2,V}^{\bullet}$ définis en termes des homomorphismes $\mathrm{r}_{V/W,Z/W}:\mathrm{St}_{V/W}\to\mathrm{St}_{Z/W}$ ($W\subset Z\subset V$),

\smallskip
-- les $(V\hspace{-2pt}\text{-}\mathcal{U})$-morphismes $\mathrm{EFix}_{(V,Z)}M\to\mathrm{EFix}_{(V,V)}M$.

\smallskip
On observera que dans le cas $M=\mathrm{H}^{*}_{V}X$, $X$ désignant un $V$-CW-complexe fini, c'est-à-dire dans le cas où l'on considère le complexe $\mathrm{C}_{\mathrm{alg}}^{\bullet}X$, l'homomorphisme $\mathrm{EFix}_{(V,Z)}M\to\mathrm{EFix}_{(V,V)}M$ s'identifie à l'homomomorphisme restriction $\mathrm{H}^{*}_{V}X^{Z}\to\mathrm{H}^{*}_{V}X^{V}$.

\medskip
Illustrons cette remarque par un cas particulier. On prend $M=\mathrm{c}_{V}\hspace{1pt}\mathrm{H}^{*}V$ avec~$V\not=0$. On a vu lors de la démonstration du point (b) de la proposition \ref{MV} que l'on a $(\mathrm{Deg}^{0}\hspace{1pt}\widehat{\Psi}_{\mathrm{c}_{V}\hspace{1pt}\mathrm{H}^{*}V})(W)=0$ pour $W\not=V$ (le cas $W=0$, exclu en section $7$, est trivial)~;  il en résulte l'égalité $\mathrm{deg}^{0}\hspace{1pt}\mathrm{B}^{\bullet,-1}(\mathrm{c}_{V}\hspace{1pt}\mathrm{H}^{*}V)=0$ et l'isomorphisme $\mathrm{deg}^{0}\hspace{1pt}\mathrm{C}^{\bullet}(\mathrm{c}_{V}\hspace{1pt}\mathrm{H}^{*}V)\cong\mathrm{Lu}_{2,V}^{\bullet}$. Une méthode alternative pour démontrer le point (b) de \ref{MV} est de se convaincre que cet isomorphisme est équivariant pour les actions à droite de $\mathrm{GL}(V)$ et de procéder ensuite par récurrence sur la dimension de $V$.

\end{rem}

\pagebreak

\sect{Filtration par la codimension du support et foncteurs $\mathrm{Fix}$}\label{Appendice}

\subsect{Filtration par la codimension du support}\label{fpcds}

 \medskip
Soit $R$ un anneau (commutatif unitaire)~;  on note $\mathop{\mathrm{Spec}}R$ son spectre~:  l'ensem\-ble des idéaux premiers $\mathfrak{p}$ de $R$, muni de la topologie de Zariski. 

\medskip
Soit $M$ un $R$-module~; soit $k\geq 0$ un entier. On pose
$$
\hspace{24pt}
\mathrm{F}^{k}M
\hspace{4pt}:=\hspace {4pt}
\bigcap_{\mathrm{ht}(\mathfrak{p})<k}\ker\hspace{2pt}(M\to M_{\mathfrak{p}})
\hspace{24pt}.
$$
Les notations $\mathrm{ht}(\mathfrak{p})$ et $M_{\mathfrak{p}}$ ci-dessus désignent respectivement  la hauteur de l'idéal premier $\mathfrak{p}$ \cite[Chap. VIII, \S 1, \numero\hspace{1pt}3, Définition 6]{Bo2} et le localisé de~$M$ en~$\mathfrak{p}$. Pour le confort du lecteur  rappelons la définition ``concrète'' \cite[Chap. VIII, \S 1, \numero\hspace{1pt}3, Proposition 7-a)]{Bo2} de $\mathrm{ht}(\mathfrak{p})$~:  c'est la borne supérieure des entiers $n$ tels que l'on a $\mathfrak{p}_{0}\subsetneq\mathfrak{p}_{1}\subsetneq\ldots\subsetneq\mathfrak{p}_{n}=\mathfrak{p}$ avec $\mathfrak{p}_{i}$ premier pour $0\leq i\leq n$. Les sous-modules $\mathrm{F}^{k}M$ fournissent une filtration décroissante naturelle de $M$~:
$$
M=\mathrm{F}^{0}M\supset\mathrm{F}^{1}M\supset\ldots\supset\mathrm{F}^{k}M\supset\ldots
$$
que l'on appelle la {\em filtration par la codimension du support}. On fait ci-après quelques observations naïves qui sont censées justifier cette terminologie.

\medskip
-- Rappelons que le {\em support} d'un $R$-module $M$ est le sous-ensemble de $\mathop{\mathrm{Spec}}R$ constitué des idéaux premier $\mathfrak{p}$ tels que l'on a $M_{\mathfrak{p}}\not=0$~; on le note $\mathrm{Supp}\hspace{1pt}(M)$.

\medskip
-- Notons $\ell_{\mathfrak{p}}:M\to M_{\mathfrak{p}}$ l'homomorphisme canonique. Soit $x$ un élément de~$M$, les deux conditions suivantes sont équivalentes (on utilise l'exactitude du foncteur localisation en $\mathfrak{p}$)~:

\smallskip
(i) $\ell_{\mathfrak{p}}(x)=0$~;

\smallskip
(ii) $\mathfrak{p}\not\in\mathrm{Supp}\hspace{1pt}(R\hspace{1pt}x)$.

\smallskip
Il en résulte que $x\in\mathrm{F}^{k}M$ et $\mathfrak{p}\in\mathrm{Supp}\hspace{1pt}(R\hspace{1pt}x)$ implique $\mathrm{ht}(\mathfrak{p})\geq k$.

\medskip
-- Soit $\mathfrak{a}$ un idéal de $R$, nous notons $\mathrm{Var}(\mathfrak{a})$ le sous-ensemble $\{\mathfrak{p};\mathfrak{a}\subset\mathfrak{p}\}$ de $\mathop{\mathrm{Spec}}R$ (ce sous-ensemble est fermé par définition, il est classiquement noté $\mathrm{V}(\mathfrak{a})$ mais le lecteur devinera aisément pourquoi nous avons écarté cette notation dans notre mémoire). Soient $M'$ un sous-module de $M$ et $\mathrm{Ann}(M')$ son idéal annulateur~; si $M'$ est de type fini alors on a $\mathrm{Supp}\hspace{1pt}(M')=\mathrm{Var}\hspace{1pt}(\mathrm{Ann}(M'))$ \cite[Chap. II, \S4, \numero\hspace{1pt}4, Proposition 17]{Bo1} si bien que le support de $M'$ est fermé dans $\mathop{\mathrm{Spec}}R$ et que l'on peut parler de sa codimension (dans $\mathop{\mathrm{Spec}}R$), que l'on note $\mathop{\mathrm{codim}}\mathrm{Supp}\hspace{1pt}(M')$ \cite[Chap. 8, \S1, \numero\hspace{1pt}2]{Bo2}. Ce qui précède conduit à l'énoncé suivant (utiliser notamment \cite[Chap. II, \S 4, \numero\hspace{1pt}4, Corollaire de la proposition 16]{Bo1} et \cite[Chap. 8, \S1, \numero\hspace{1pt}3, Définition 6 et Proposition 7-a)]{Bo2})~:

\begin{pro}\label{cds1} Soient $M$ un $R$-module et $M'$ un sous-module de type fini de $M$. Les trois conditions suivantes sont équivalentes~:
\begin{itemize}
\item[(i)] $M'\subset\mathrm{F}^{k}M$~;
\item[(ii)] $M'=\mathrm{F}^{k}M'$~;
\item[(iii)] $\mathop{\mathrm{codim}}\mathrm{Supp}\hspace{1pt}(M')\geq k$.
\end{itemize}

\end{pro}

\begin{cor}\label{cds2} Soient $R$ un anneau noethérien et $M$ un $R$-module de type fini de $M$. Alors pour tout sous-module $M'$ de $M$ et tout entier  $k\geq 0$ les trois conditions suivantes sont équivalentes~:
\begin{itemize}
\item[(i)] $M'\subset\mathrm{F}^{k}M$~;
\item[(ii)] $M'=\mathrm{F}^{k}M'$~;
\item[(iii)] $\mathop{\mathrm{codim}}\mathrm{Supp}\hspace{1pt}(M')\geq k$.
\end{itemize}
(En particulier on a $\mathop{\mathrm{codim}}\mathrm{Supp}\hspace{1pt}(\mathrm{F}^{k}M)\geq k$.)

\end{cor}

\bigskip
\textit{Exemple.} Soit $K$ un corps~; on prend pour $R$ l'anneau $K[X_{1},X_{2},\ldots,X_{n}]$ des polynômes en $n$ indéterminées à coefficients dans $K$. Dans ce cas~:

\medskip
(1) La dimension de Krull de $R$ \cite[Chap. VIII, \S 1, \numero\hspace{1pt}3, Définition 5]{Bo2} est égale à $n$ \cite[Chap. VIII, \S 2, \numero\hspace{1pt}4, Corollaire 1]{Bo2}.

\medskip
(2) Soit $\mathfrak{p}$ un idéal premier de $R$, alors on a les équivalences suivantes~:

\smallskip
(2.1) $\mathrm{ht}(\mathfrak{p})=0\iff\mathfrak{p}=\{0\}$~;

\smallskip
(2.2) $\mathrm{ht}(\mathfrak{p})=1\iff\mathfrak{p}=R\hspace{1pt}a\hspace{4pt}\text{avec}\hspace{4pt}$a$\hspace{4pt}\text{irréductible}$~;

\smallskip
(2.3) $\mathrm{ht}(\mathfrak{p})=n\iff\mathfrak{p}\hspace{4pt}\text{maximal}$.

\medskip
Soit $M$ un $R$-module.

\smallskip
-- Le point (1) implique que l'on a $\mathrm{F}^{k}M=0$ pour $k>n$.

\smallskip
-- Le point (2.1) montre que  $\mathrm{F}^{1}M$ est la torsion de $M$.

\smallskip
-- Si $M$ est de type fini alors  $\mathrm{F}^{n}M$ est de dimension finie en tant qu'espace vectoriel sur $K$ (utiliser  \cite[Chap. IV, \S 1, \numero\hspace{1pt}4, Théorèmes 1 et 2]{Bo1} et le nullstellensatz).

\subsect{Comparaison entre la filtration d'un $\mathrm{H}^{*}V$-$\mathrm{A}$-module instable définie en termes des foncteurs $\mathrm{Fix}$ et la filtration par la codimension du support du $\mathrm{H}^{*}V$-module sous-jacent}

\medskip
Soit $M$ un $\mathrm{H}^{*}V$-$\mathrm{A}$-module instable. Nous avons défini dans l'introduction (et mis en {\oe}uvre dans la section 5) une première filtration décroissante naturelle de $M$, par des sous-$\mathrm{H}^{*}V$-$\mathrm{A}$-modules instables,
$$
M=\mathrm{F}^{0}M\supset\mathrm{F}^{1}M\supset\ldots\supset\mathrm{F}^{k}M\supset\ldots\supset\mathrm{F}^{n}M\supset\mathrm{F}^{n+1}M=0
$$
en posant
$$
\hspace{24pt}
\mathrm{F}^{k}M
\hspace{4pt}=\hspace{4pt}
\bigcap_{\mathop{\mathrm{codim}}W<k}
\ker\hspace{1pt}(\eta_{(V,W),M}:M\to\mathrm{H}^{*}V\otimes_{\mathrm{H}^{*}V/W}\mathrm{Fix}_{(V,W)}M)
\hspace{24pt}.
$$
D'autre part le $\mathrm{H}^{*}V$-module sous-jacent à $M$ admet une seconde filtration décroissante naturelle, par des sous-$\mathrm{H}^{*}V$-modules, à savoir la filtration par la codimension du support qui est du même type, compte tenu du fait que  $\mathrm{H}^{*}V$ est isomorphe à un anneau de polynômes, à coefficients dans $\mathbb{F}_{2}$, en $n$ indéterminées (voir l'exemple à la fin de \ref{fpcds}.1).

\begin{pro}\label{comparaison1} Soit $M$ un $\mathrm{H}^{*}V$-$\mathrm{A}$-module instable qui est de type fini comme  $\mathrm{H}^{*}V$-module. La filtration de $M$ définie en termes des foncteurs $\mathrm{Fix}$ et la filtration par la codimension du support du $\mathrm{H}^{*}V$-module sous-jacent coïncident (en tant que filtrations par des sous-$\mathrm{H}^{*}V$-modules).
\end{pro}

\bigskip
\textit{Démonstration.} Notons (localement) $\mathrm{F}_{1}^{k}M$ les sous-$\mathrm{H}^{*}V$-$\mathrm{A}$-modules instables de la première filtration et $\mathrm{F}_{2}^{k}M$ les sous-$\mathrm{H}^{*}V$-modules de la seconde. Soit $M\subset I$ une inclusion de $\mathrm{H}^{*}V$-$\mathrm{A}$-modules instables alors on a
$$
\mathrm{F}_{j}^{k}M
\hspace{4pt}=\hspace{4pt}
M\cap\mathrm{F}_{j}^{k}I\leqno{(1)}
$$
pour $j=1,2$. Cette égalité résulte pour $j=1$ de l'exactitude des foncteurs $\mathrm{EFix}$ (voir \ref{incfilt}) et pour $j=2$ de celle des foncteurs de localisation en $\mathfrak{p}$. Comme la catégorie $V_{\mathrm{tf}}\text{-}\mathcal{U}$ a assez d'injectifs, l'égalité $(1)$ montre qu'il suffit de vérifier \ref{comparaison1} pour $M$ injectif. Compte tenu de la classification des $(V_{\mathrm{tf}}\text{-}\mathcal{U})$-injectifs (voir \ref{injectifs3}) on est ramené à vérifier \ref{comparaison1} pour $M=\mathrm{H}^{*}V\otimes_{\mathrm{H}^{*}V/W}N$ avec $W\subset V$ et $N$ un $\mathrm{H}^{*}V/W$-$\mathrm{A}$-module instable fini. On conclut en constatant que si $M$ est de cette forme alors on a
$$
\mathrm{F}_{j}^{k}M
\hspace{4pt}=\hspace{4pt}
\begin{cases}
M & \text{pour $k\leq\mathop{\mathrm{codim}}W$}, \\
0 & \text{pour $k>\mathop{\mathrm{codim}}W$}.
\end{cases}
\leqno{(2)}
$$
Le cas $j=1$ a déjà été traité dans la section 5 (Proposition \ref{filtparticulier}). Reste à traiter le cas $j=2$.

\medskip
On pose $\mathcal{P}_{W}:=\ker(\mathrm{H}^{*}V\overset{i^{*}}{\to}\mathrm{H}^{*}W)$, $i$ désignant l'inclusion de $W$ dans $V$~; $\mathcal{P}_{W}$ est un idéal homogène de $\mathrm{H}^{*}V$ stable sous l'action de l'algèbre de Steenrod. Cet idéal est premier, puisque $\mathrm{H}^{*}W$ est intègre, sa hauteur est la codimension de $W$ (invoquer par exemple \cite[Chap. VIII, \S 2, \numero\hspace{1pt}4, Corollaire 2]{Bo2}). On observera que $i^{*}$ induit un isomorphisme canonique de $\mathrm{A}$-algèbres instables
$$
\mathrm{H}^{*}V\otimes_{\mathrm{H}^{*}V/W}\mathbb{F}_{2}
\hspace{4pt}\cong\hspace{4pt}
\mathrm{H}^{*}W
$$
($\mathbb{F}_{2}$ étant considéré ci-dessus comme un $\mathrm{H}^{*}V/W$-$\mathrm{A}$-module instable).

\medskip
Le $\mathrm{H}^{*}V/W$-$\mathrm{A}$-module instable fini $N$ admet une filtration finie
$$
0=N_{-1}\subset N_{0}\subset N_{1}\subset\ldots\subset N_{r}=N
$$
avec $N_{s}/N_{s-1}\cong\Sigma^{k}\mathbb{F}_{2}$, $k\in\mathbb{N}$, pour~$0\leq s\leq r$. Du coup  le $\mathrm{H}^{*}V$-$\mathrm{A}$-module instable $M$ admet une filtration finie
$$
0=M_{-1}\subset M_{0}\subset M_{1}\subset\ldots\subset M_{r}=M
\leqno{(3)}
$$
avec $M_{s}/M_{s-1}\cong\Sigma^{k}\hspace{1pt}\mathrm{H}^{*}W\cong\Sigma^{k}\hspace{1pt}(\mathrm{H}^{*}V/\mathcal{P}_{W})$, $k\in\mathbb{N}$, pour~$0\leq s\leq r$. Le cas $j=2$ de l'égalité (2) en résulte. Détaillons lourdement~:

\medskip
Soient $R$ un anneau et $\mathfrak{p}$ un idéal premier de $R$~; on a $\mathrm{Supp}\hspace{1pt}(R/\mathfrak{p})=\mathrm{Var}(\mathfrak{p})$ et  $(R/\mathfrak{p})_{\mathfrak{p}}$ est le corps des fractions de l'anneau intègre $R/\mathfrak{p}$, en particulier $\ell_{\mathfrak{p}}:R/\mathfrak{p}\to(R/\mathfrak{p})_{\mathfrak{p}}$ est injective. Spécialisons~: on a $\mathrm{Supp}\hspace{1pt}(\mathrm{H}^{*}V/\mathcal{P}_{W})=\mathrm{Var}(\mathcal{P}_{W})$ et $\ell_{\mathcal{P}_{W}}:\mathrm{H}^{*}V/\mathcal{P}_{W}\to(\mathrm{H}^{*}V/\mathcal{P}_{W})_{\mathcal{P}_{\mathrm{W}}}$ est injective. On en déduit par récurrence sur la longueur de la filtration (3) que l'on a $\mathrm{Supp}\hspace{1pt}(M)=\mathrm{Var}(\mathcal{P}_{W})$(\cite[Chapitre II, \S4, \numero 4, Proposition 16 (i)]{Bo1}) et que $\ell_{\mathcal{P}_{W}}:M\to M_{\mathcal{P}_{\mathrm{W}}}$ est injective. Les égalités $\mathop{\mathrm{codim}}\mathrm{Supp}\hspace{1pt}(M)=\mathrm{ht}(\mathcal{P}_{W})$ (\cite[Chap. 8, \S1, \numero\hspace{1pt}3, Définition 6 et Proposition 7-a)]{Bo2}), $\mathrm{ht}(\mathcal{P}_{W})=\mathop{\mathrm{codim}}W$, et l'implication $(iii)\Rightarrow(ii)$ de \ref{cds1} montrent que l'on a $\mathrm{F}_{2}^{k}M=M$ pour $k\leq\mathop{\mathrm{codim}}W$. L'injectivité de  $\ell_{\mathcal{P}_{W}}:M\to M_{\mathcal{P}_{\mathrm{W}}}$ (et à nouveau l'égalité $\mathrm{ht}(\mathcal{P}_{W})=\mathop{\mathrm{codim}}W$) implique $\mathrm{F}_{2}^{k}M=0$ pour $k>\mathop{\mathrm{codim}}W$.
\hfill$\square$

\medskip
\begin{scho}\label{PWtorsions} Soit $M$ un $\mathrm{H}^{*}V$-$\mathrm{A}$-module instable qui est de type fini comme  $\mathrm{H}^{*}V$-module~; soit $\mathrm{F}^{k}M$ le $k$-ième sous-$\mathrm{H}^{*}V$-$\mathrm{A}$-module instable de la filtration de $M$ définie en termes des foncteurs $\mathrm{Fix}$. Alors on a
$$
\mathrm{F}^{k}M
\hspace{4pt}=\hspace{4pt}
\underset{\mathop{\mathrm{codim}}W=k}{+}\hspace{4pt}\Gamma_{\mathcal{P}_{W}}
$$
($\hspace{2pt}\Gamma_{\mathcal{P}_{W}}$ désigne ci-dessus la $\mathcal{P}_{W}$-torsion de $M$, voir le début de la section 3).
\end{scho}

\medskip
\begin{scho}\label{comparaison2} Soit $M$ un $\mathrm{H}^{*}V$-$\mathrm{A}$-module instable qui est de type fini comme  $\mathrm{H}^{*}V$-module. Alors les sous-$\mathrm{H}^{*}V$-modules de la filtration par la codimension du support du $\mathrm{H}^{*}V$-module sous-jacent sont stables sous l'action de l'algèbre de Steenrod.
\end{scho}

\pagebreak

\medskip
Compte tenu de \ref{Filttop2}, on a également~:

\begin{scho}\label{comparaison3} Soit $X$ un $V$-CW-complexe fini. Alors la filtration décrois\-sante de $\mathrm{H}_{V}^{*}X$ définie par
$$
\hspace{24pt}
\mathrm{F}^{k}\hspace{1pt}\mathrm{H}_{V}^{*}X
\hspace{4pt}=\hspace{4pt}
\bigcap_{\mathop{\mathrm{codim}}W<k}
\ker\hspace{1pt}(\mathrm{H}_{V}^{*}X\to\mathrm{H}_{V}^{*}X^{W})
\hspace{24pt}
$$
coïncide (en tant que filtration par des sous-$\mathrm{H}^{*}V$-modules) avec la filtration par la codimension du support du $\mathrm{H}^{*}V$-module sous-jacent.
\end{scho}

\vspace{0.75cm}
\textit{Sur l'hypothèse de finitude qui figure dans les énoncés \ref{comparaison1}, \ref{PWtorsions}} et \ref{comparaison2}

\bigskip
Nous avons supposé dans les énoncés en question que le $\mathrm{H}^{*}V$-$\mathrm{A}$-module instable $M$ est de type fini comme  $\mathrm{H}^{*}V$-module. Nous avions deux raisons pour cela. La première était que dans notre contexte (sections 3, 4, 5 et 6) cette condition était satisfaite, la seconde que cela simplifiait notre exposition. En fait cette restriction peut être levée. Pour cela on doit apporter de légères modifications à la démonstration que nous avons donnée de \ref{comparaison1}, modifications que nous décrivons brièvement ci-après.

\bigskip
1) On remplace le théorème de classification des $V_{\mathrm{tf}}\text{-}\mathcal{U}$-injectifs par celui de classification des $V\text{-}\mathcal{U}$-injectifs (respectivement \ref{injectifs3} et \ref{injectifs1}).

\bigskip
2) On observe que les foncteurs localisation en $\mathfrak{p}$ et les foncteurs $\mathrm{EFix}$, qui sont des adjoints à gauche, commutent aux sommes directes arbitraires si bien que l'on est ramené à nouveau au cas où $M$ est un $V\text{-}\mathcal{U}$-injectif indécomposable et plus généralement au cas $M=S\otimes(\mathrm{H}^{*}V\otimes_{\mathrm{H}^{*}V/W}N)$ avec $S$ un $\mathrm{A}$-module instable et $N$ un $\mathrm{H}^{*}V/W$-$\mathrm{A}$-module instable fini.

\bigskip
3) On revisite le point (a) de \ref{suspensionFix}. Soit $P$ un $\mathrm{H}^{*}V$-$\mathrm{A}$-module instable~; on vérifie que l'unité d'adjonction
$$
\eta_{(V,W),S\otimes P}:S\otimes P
\longrightarrow
\mathrm{H}^{*}V\otimes_{\mathrm{H}^{*}V/W}\mathrm{Fix}_{(V,W)}(S\otimes P)
$$
s'identifie au produit tensoriel de l'homomorphisme $\iota_{S}:S=\mathrm{T}_{0}S\to\mathrm{T}_{W}S$ et de l'unité d'adjonction $\eta_{(V,W),P}:P\to
\mathrm{H}^{*}V\otimes_{\mathrm{H}^{*}V/W}\mathrm{Fix}_{(V,W)}P$. Comme $\iota_{S}$ admet une rétraction naturelle, à savoir l'homomorphisme $\mathrm{T}_{W}S\to\mathrm{T}_{0}S=S$, on constate que l'on dispose de la variante suivante de la proposition \ref{proclefbis}~:

\pagebreak

\bigskip
\begin{pro}\label{proclefter} Soient $V$ un $\ell$-groupe abélien élémentaire et $W$, $U$ des~sous-groupes~; soient $S$ un $\mathrm{A}$-module instable et $N$ un $\mathrm{H}^{*}V/U$-$\mathrm{A}$-module instable fini. Alors l'unité d'adjonction
$$
\begin{CD}
S\otimes(\mathrm{H}^{*}V\otimes_{\mathrm{H}^{*}V/U}N)
@>\eta_{(V,W)}>>
\mathrm{H}^{*}V\otimes_{\mathrm{H}^{*}V/W}\mathrm{Fix}_{(V,W)}(S\otimes(\mathrm{H}^{*}V\otimes_{\mathrm{H}^{*}V/U}N))
\end{CD}
$$
est un monomorphisme (scindé) si l'on a $W\subset U$  et est nulle si l'on a $W\not\subset U$.
\end{pro}

Cette proposition implique la généralisation ci-dessous de la proposition \ref{filtparticulier}~:

\begin{pro}\label{filtparticulierbis} Soient $U$ un sous-groupe de $V$, $S$ un $\mathrm{A}$-module instable et $N$ un $\mathrm{H}^{*}V/U$-$\mathrm{A}$-module instable fini. Alors on a
$$
\mathrm{F}^{k}\hspace{2pt}(S\otimes(\mathrm{H}^{*}V\otimes_{\mathrm{H}^{*}V/U}N))
\hspace{4pt}=\hspace{4pt}
\begin{cases}
S\otimes(\mathrm{H}^{*}V\otimes_{\mathrm{H}^{*}V/U}N) & \text{pour $k\leq\mathop{\mathrm{codim}}U$}, \\
0 & \text{pour $k>\mathop{\mathrm{codim}}U$}.
\end{cases}
$$
\end{pro}

(Ci-dessus $\mathrm{F}^{k}$ est défini en termes des foncteurs $\mathrm{Fix}$.)

\medskip
\subsect{Compléments}\label{Serre1}

\medskip
La spécificité de la filtration par la codimension du support d'un $\mathrm{H}^{*}V$-module gradué muni d'une action ``compatible et instable'' de l'algèbre de Steenrod tient essentiellement au résultat ci-dessous dû à Jean-Pierre Serre \cite[\S2, Proposition (1)]{Ser}.  On rappelle que la {\em racine} (ou le {\em radical}\hspace{1pt}) d'un idéal $\mathfrak{a}$ d'un anneau $R$, est le sous-ensemble de $R$ constitué des éléments dont une puissance appartient à $\mathfrak{a}$ ; ce sous-ensemble est aussi un idéal de $R$, il est noté~$\sqrt{\mathfrak{a}}$. On dit que $\mathfrak{a}$ est {\em radiciel} si l'on a $\mathfrak{a}=\sqrt{\mathfrak{a}}$~; l'idéal $\sqrt{\mathfrak{a}}$ est radiciel.

\begin{theo}\label{Serre2} Soient $\mathfrak{a}$ un idéal homogène de $\mathrm{H}^{*}V$ stable sous l'action de l'algèbre de Steenrod. Soit $\mathcal{W}_{\mathfrak{a}}$ le sous-ensemble de $\mathcal{W}$ constitué des sous- groupes $W\subset V$ tels que l'on a $\mathfrak{a}\subset\mathcal{P}_{W}$ et qui sont maximaux pour cette propriété. Alors on a~:
 $$
 \hspace{24pt}
 \sqrt{\mathfrak{a}}
 \hspace{4pt}=\hspace{4pt}
 \bigcap_{W\in\mathcal{W}_{\mathfrak{a}}}
 \mathcal{P}_{W}
 \hspace{24pt}.
 $$

\end{theo}

Compte tenu d'un lemme bien connu sur les idéaux premiers d'un anneau (voir par exemple \cite[Proposition 1.11 (ii)]{AM}) ce théorème admet le corollaire suivant~:

\begin{cor}\label{Serre3} Soient $\mathfrak{p}$ un idéal premier homogène de $\mathrm{H}^{*}V$ stable sous l'action de l'algèbre de Steenrod.  Alors il existe un sous-groupe $W$ de $V$, uniquement déterminé, tel que l'on a $\mathfrak{p}=\mathcal{P}_{W}$.

\end{cor}

On se propose de montrer que la théorie des $V_{\mathrm{tf}}\text{-}\mathcal{U}$-injectifs conduit à une démonstration de \ref{Serre2}.

\begin{lem}\label{racine} Soit $\mathfrak{a}$ un idéal de $\mathrm{H}^{*}V$. Si $\mathfrak{a}$ est homogène et stable sous l'action de l'algèbre de Steenrod alors il en est de même pour $\sqrt{\mathfrak{a}}$.

\end{lem}

\textit{Démonstration.} Le fait que $\sqrt{\mathfrak{a}}$ est homogène est un cas particulier d'un résultat général sur les idéaux homogènes d'un anneau gradué (voir par exemple \cite[Chap. VII, \S2, Theorem 8 (b)]{ZS}). Dans le cas de l'anneau gradué $\mathrm{H}^{*}V$ on peut se convaincre de ce que $\sqrt{\mathfrak{a}}$ est homogène grâce à l'observation triviale ci-après. Soit $x$ un élément, pas nécessairement  homogène, de $\mathrm{H}^{*}V$ alors les deux conditions suivantes sont équivalentes~:

 \begin{itemize}
 \item[(i)] il existe $m$ dans $\mathbb{N}-\{0\}$ avec $x^{m}\in\mathfrak{a}$;
 \item[(ii)] il existe $m$ dans $\mathbb{N}$ avec $x^{2^{m}}\in\mathfrak{a}$.
 \end{itemize}
 
 \smallskip
 Cette observation conduit aussi à une preuve du fait que $\sqrt{\mathfrak{a}}$ est stable sous l'action de l'algèbre de Steenrod. Soient $x$ un élément, cette fois homogène, de $\sqrt{\mathfrak{a}}$ et $i$ un entier. Il existe un entier $m$ tel que l'on a $x^{2^{m}}\in\mathfrak{a}$. Comme $\mathfrak{a}$ est stable sous l'action de l'algèbre de Steenrod on a aussi $\mathrm{Sq}^{2^{m}i}\hspace{2pt}x^{2^{m}}\in\mathfrak{a}$~; l'égalité $\mathrm{Sq}^{2^{m}i}\hspace{2pt}x^{2^{m}}=(\mathrm{Sq}^{i}x)^{2^{m}}$ montre $\mathrm{Sq}^{i}x\in\mathfrak{a}$.
 \hfill$\square$
 
 \bigskip
 On va démontrer l'énoncé d'algèbre commutative \ref{Serre2} en utilisant l'énoncé d'algèbre linéaire \ref{enveloppe-reduit}.
  
\medskip
Soit $M$ un $\mathrm{A}$-module instable. Rappelons que l'on dit que $M$ est {\em réduit} si pour tout entier $i$ l'opération $\mathrm{Sq}^{i}: M^{i}\to M^{2i}$, souvent notée $\mathrm{Sq}_{0}$, est injective. Rappelons également l'origine de cette terminologie. Si $M$ est une $A$-algèbre instable alors les deux conditions suivantes sont équivalentes~:
 \begin{itemize}
 \item[(i)] la $\mathbb{F}_{2}$-algèbre graduée sous-jacente à $M$ est réduite (au sens de l'algèbre commutative)~;
 \item[(ii)] le $\mathrm{A}$-module instable sous-jacent à $M$ est réduit (au sens ci-dessus).
  \end{itemize}

\medskip
Cette équivalence provient du fait que l'on a $\mathrm{Sq}_{0}\hspace{1pt}x=x^{2}$ pour tout $x$ dans $M$ (la définition de la notion de $A$-algèbre instable s'inspire des propriétés de la cohomologie modulo $2$ d'un espace).

\bigskip
Ces rappels étant faits, nous pouvons énoncer :

\begin{pro}\label{enveloppe-reduit} Soient $M$ un $\mathrm{H}^{*}V$-$\mathrm{A}$-module instable qui est de type fini comme $\mathrm{H}^{*}V$-module et $i:M\to E$ une enveloppe injective (dans la catégorie $V_{\mathrm{tf}}\text{-}\mathcal{U}$ ou $V\text{-}\mathcal{U}$, voir \ref{smith0.13bis}). Les deux conditions suivantes sont équivalentes~:
 \begin{itemize}
 
 \medskip
 \item[(i)] le $A$-module sous-jacent à $M$ est réduit ;
 
\medskip
 \item[(ii)] $E$ est isomorphe à une somme directe de la forme
$$
\bigoplus_{W\in\mathcal{W}}
\hspace{4pt}
(\mathrm{H}^{*}W)^{\hspace{1pt}a_{W}}
$$
avec $a_{W}\in\mathbb{N}$ (\hspace{2pt}$\mathrm{H}^{*}W$ étant un $\mathrm{H}^{*}V$-$\mathrm{A}$-module instable \textit{via} l'homomorphisme canonique de $\mathrm{A}$-algèbres instables  $\mathrm{H}^{*}V\to\mathrm{H}^{*}W$).
 \end{itemize}

\end{pro}

\medskip
\textit{Démonstration de $(ii)\Rightarrow(i)$.} Si $(ii)$ est vérifiée alors le $\mathrm{A}$-module instable sous-jacent à $E$ est réduit (par exemple parce que les $\mathbb{F}_{2}$-algèbres $\mathrm{H}^{*}W$ sont  réduites) et il en est donc de même pour $M$.

\medskip
\textit{Démonstration de $(i)\Rightarrow(ii)$.} On sait \textit{a priori}, d'après \ref{injectifs3}, que $E$ est isomorphe à une somme directe de la forme
$$
\bigoplus_{(W,k)\in\mathcal{W}\times\mathbb{N}}
\hspace{4pt}
(\mathrm{H}^{*}V\otimes_{\mathrm{H}^{*}V/W}\mathrm{J}_{V/W}(k))^{\hspace{1pt}a_{W,k}}
$$
avec $a_{W,k}\in\mathbb{N}$ et $a_{W,k}=0$ pour $k$ assez grand. Comme l'on a
$$
\hspace{24pt}
\mathrm{H}^{*}V\otimes_{\mathrm{H}^{*}V/W}\mathrm{J}_{V/W}(0)
\hspace{4pt}\cong\hspace{4pt}
\mathrm{H}^{*}V\otimes_{\mathrm{H}^{*}V/W}\mathbb{F}_{2}
\hspace{4pt}\cong\hspace{4pt}
\mathrm{H}^{*}W
\hspace{24pt},
$$
il suffit de montrer que $(i)$ implique $a_{W,k}=0$ pour $k>0$.

\medskip
On dispose d'une inclusion canonique de $\mathrm{H}^{*}V$-$\mathrm{A}$-modules instables, disons $\iota_{V,k}:\Sigma^{k}\mathbb{F}_{2}\hookrightarrow\mathrm{J}_{V}(k)$ (en fait $\mathrm{J}_{V}(k)$ est une enveloppe injective de $\Sigma^{k}\mathbb{F}_{2}$). En effet, on a par définition $\mathrm{Hom}_{V\text{-}\mathcal{U}}(\Sigma^{k}\mathbb{F}_{2},\mathrm{J}_{V}(k))=\mathrm{Hom}_{\mathbb{F}_{2}}(\mathbb{F}_{2},\mathbb{F}_{2})$ et l'unique homomorphisme non nul de $\Sigma^{k}\mathbb{F}_{2}$ dans $\mathrm{J}_{V}(k)$ est nécessairement injectif. En appliquant le foncteur $\mathrm{H}^{*}V\otimes_{\mathrm{H}^{*}V/W}-$ à l'inclusion canonique $\iota_{V/W,k}$, on obtient une inclusion de $\mathrm{H}^{*}V$-$\mathrm{A}$-modules instables, tout aussi canonique,  $\Sigma^{k}\mathrm{H}^{*}W\hookrightarrow\mathrm{H}^{*}V\otimes_{\mathrm{H}^{*}V/W}\mathrm{J}_{V/W}(k)$ (en fait $\mathrm{H}^{*}V\otimes_{\mathrm{H}^{*}V/W}\mathrm{J}_{V/W}(k)$ est une enveloppe injective de $\Sigma^{k}\mathrm{H}^{*}W$). Ce qui précède permet de conclure. Si l'on a~$a_{W,k}\not=0$, alors $E$ contient un sous-module $P$ isomorphe à $\Sigma^{k}\mathrm{H}^{*}W$~; si $M$ est réduit on a $i^{-1}(P)=0$  pour $k>0$ (puisque dans ce cas l'opération $\mathrm{Sq}_{0}$ est nulle sur $i^{-1}(P)$) et donc $P=0$ (propriété d'une enveloppe injective) ce qui force $a_{W,k}=0$ pour $k>0$.
\hfill$\square$

\pagebreak

\bigskip
Fin de l'interlude linéaire. Revenons à la démonstration de \ref{Serre2}. Soit $\mathfrak{a}$ un idéal homogène de $\mathrm{H}^{*}V$ stable sous l'action de l'algèbre de Steenrod. Les trois conditions suivantes sont équivalentes~:
\begin{itemize}
\item[(i)] l'idéal $\mathfrak{a}$ est radiciel ($\mathfrak{a}=\sqrt{\mathfrak{a}}$)~;
\item[(ii)] la $\mathbb{F}_{2}$-algèbre $\mathrm{H}^{*}V/\mathfrak{a}$ est réduite~;
\item[(iii)] le $\mathrm{A}$-module instable $\mathrm{H}^{*}V/\mathfrak{a}$ est réduit.
\end{itemize}
(Pour $(ii)\iff(iii)$ voir la discussion qui précède \ref{enveloppe-reduit}.)

\bigskip
Compte tenu du lemme \ref{racine}, le cas particulier où $\mathfrak{a}$ est radiciel du théorème \ref{Serre2} implique le cas général~; on suppose donc $\mathfrak{a}$ radiciel. Dans ce cas, d'après la proposition \ref{enveloppe-reduit}, $\mathrm{H}^{*}V/\mathfrak{a}$ admet une enveloppe injective, dans la catégorie $V_{\mathrm{tf}}$-$\mathcal{U}$, de la forme 
$$
i:\mathrm{H}^{*}V/\mathfrak{a}\to\mathrm{H}^{*}W_{1}\oplus\mathrm{H}^{*}W_{2}\oplus\ldots\oplus\mathrm{H}^{*}W_{r}
$$
avec $W_{1},W_{2},\ldots,W_{r}$ une suite finie de sous-groupes de $V$.

\medskip
Soient $q:\mathrm{H}^{*}V\to\mathrm{H}^{*}V/\mathfrak{a}$ l'homorphisme canonique de $\mathrm{A}$-algèbres instables  et $\pi_{k}:\mathrm{H}^{*}V/\mathfrak{a}\to\mathrm{H}^{*}W_{k}$ l'homomorphisme de  $\mathrm{H}^{*}V$-$\mathrm{A}$-modules instables composé de $i$ et de la projection sur le facteur $\mathrm{H}^{*}W_{k}$. Le composé $\pi_{k}\circ q:\mathrm{H}^{*}V\to\mathrm{H}^{*}W_{k}$ coïncide avec l'homomorphisme de $\mathrm{A}$-algèbres instables induit par l'inclusion, disons $j_{k}:W_{k}\hookrightarrow V$. Détaillons lourdement. Posons  $\varphi_{k}:=\pi_{k}\circ q$, puisque $\varphi_{k}$ est un homomorphisme de $\mathrm{H}^{*}V$-$\mathrm{A}$-modules instables, on a $\varphi_{k}(a)=j_{k}^{*}(a)\varphi_{k}(1)$ pour tout $a$ dans $\mathrm{H}^{*}V$ et donc $\varphi_{k}=\lambda\hspace{1pt}j_{k}^{*}$ en posant $\lambda=\varphi_{k}(1)$~; $\lambda=0$ est interdit car $\varphi_{k}$ est non nul (propriété d'enveloppe injective).

\medskip
Arrivé là, on a  déjà démontré le théorème de Serre. En effet, on a $\mathfrak{a}=\ker( i\circ q)$ et donc, d'après ce qui précède,
$$
\hspace{24pt}
\mathfrak{a}
\hspace{4pt}=\hspace{4pt}
\bigcap_{k=1}^{r}\mathcal{P}_{W_{k}}
\hspace{24pt}.
$$
On peut ``nettoyer"  l'égalité ci-dessus de la façon suivante. Soit $M\hspace{-1.5pt}ax$ le sous-ensemble de $\{W_{1},W_{2},\ldots,W_{r}\}$ constitué des éléments maximaux pour l'inclusion alors on aussi
$$
\hspace{24pt}
\mathfrak{a}
\hspace{4pt}=\hspace{4pt}
\bigcap_{W\in M\hspace{-1pt}ax}\mathcal{P}_{W}
\hspace{24pt}.
$$
En fait ce nettoyage est inutile dans notre contexte. Considérons à nouveau l'enveloppe injective
$$
\hspace{24pt}
i:\mathrm{H}^{*}V/\mathfrak{a}\to\mathrm{H}^{*}W_{1}\oplus\mathrm{H}^{*}W_{2}\oplus\ldots\oplus\mathrm{H}^{*}W_{r}
\hspace{4pt}:=\hspace{4pt}E
\hspace{24pt}.
$$
Les propriétés d'enveloppe injective font que l'on a $W_{k}\not\subset W_{l}$ et $W_{l}\not\subset W_{k}$ pour $k\not=l$ (en termes plus pédants que la relation d'ordre sur $\{W_{1},W_{2},\ldots,W_{r}\}$  définie par l'inclusion est l'égalité). Détaillons. Supposons $W_{k}\subset W_{l}$, alors on a $\pi_{k}=\rho\circ\pi_{l}$, $\rho$ désignant  l'homomorphisme $\mathrm{H}^{*}W_{l}\to\mathrm{H}^{*}W_{k}$  induit par cette inclusion. Identifions $\mathrm{H}^{*}W_{k}$ à un sous-module de $E$~; l'égalité $\pi_{k}=\rho\circ\pi_{l}$ entraîne $i^{-1}(\mathrm{H}^{*}W_{k})=0$, contradiction.
\hfill$\square$

\bigskip
\textit{Commentaires}

\medskip
1) La démonstration que nous avons donnée du théorème \ref{Serre2} est bien moins ``géodésique'' que la démonstration originale  de Serre~! Signalons cependant que notre méthode de ``linéarisation" fournit des informations supplémen\-taires sur les idéaux satisfaisant les conditions du théorème.

\medskip
2) On peut donner une démonstration du théorème \ref{Serre2} à l'aide de la partie II de \cite{HLS}. En effet celle-ci implique que $\mathrm{H}^{*}V/\sqrt{\mathfrak{a}}$ se plonge comme $\mathrm{A}$-algèbre instable dans un produit fini $\prod_{k}\mathrm{H}^{*}V_{k}$ et on conclut comme ci-dessus. Il est à signaler que la partie II de \cite{HLS} se déduit, là encore par un ``principe de linéarisation'',  de la partie I de cette référence qui traite des $\mathrm{A}$-modules instables ``aux nilpotents près''.

\medskip
3) Nous avons supposé $\ell=2$ afin que $\mathrm{H}^{*}V$ soit commutative et rester dans le confort douillet de l'algèbre commutative. Si l'on veut rester dans ce cadre pour $\ell>2$, il faut remplacer $\mathrm{H}^{*}V$ par sa sous-algèbre, disons $\mathrm{P}V$, engendrée par le Bockstein des classes de degré $1$ ($\mathrm{P}V$ est isomorphe à un anneau de polynômes, à coefficients dans $\mathbb{F}_{\ell}$, en $n$ indéterminées de degré $2$) et l'algèbre de Steenrod $\mathrm{A}$ par sa sous-algèbre, disons $\mathrm{A}'$, engendrée par les opérations de Steenrod $\mathrm{P}^{i}$.  C'est ce que fait Serre dans \cite{Ser}.

\medskip
4) Serre utilise le théorème \ref{Serre2} pour obtenir l'énoncé suivant~:

\begin{cor}\label{Serre4} Soit $\mathfrak{a}$ un idéal de $\mathrm{H}^{*}V$ vérifiant les hypothèses du théorème \ref{Serre2}. Si $\mathfrak{a}$ est non nul alors il contient un produit d'éléments de~$\mathrm{H}^{1}V-\{0\}$.
\end{cor}

\bigskip
\textit{Démonstration.} On reprend les notations de \ref{Serre2}. Puisque $\mathfrak{a}$ est non nul, on~a $W\not=V$ pour tout $W$ dans $\mathcal{W}_{\mathfrak{a}}$~; il existe donc, pour tout $W$ dans $\mathcal{W}_{\mathfrak{a}}$, une forme linéaire $u_{W}$ non nulle sur $V$ et nulle sur $W$. On considère $u_{W}$ comme un élément de $\mathrm{H}^{1}V-\{0\}$~; par construction $\prod_{W\in\mathcal{W}_{\mathfrak{a}}}u_{W}$ appartient à $\sqrt{\mathfrak{a}}$ . Une puissance de ce produit appartient à $\mathfrak{a}$.
\hfill$\square$

\bigskip
On pose
$\mathrm{c}_{V}:=\prod_{u\in\mathrm{H}^{1}V-\{0\}}u$. Il est bien connu que l'énoncé ci-dessus conduit au suivant~:

\pagebreak

\begin{pro}\label{Serre5} Soit $M$ un $\mathrm{H}^{*}V$-$\mathrm{A}$-module instable. Soit $x$ un élément de $M$, les deux conditions suivantes sont équivalentes~:
\begin{itemize}

\medskip
\item[(i)] $x$ est un élément de torsion du $\mathrm{H}^{*}V$-module sous-jacent à $M$~;

\medskip
\item[(ii)] il existe un entier naturel $m$ tel que l'on a $\mathrm{c}_{V}^{m}x=0$.
\end{itemize}
\end{pro}

\textit{Démonstration de $(i)\Rightarrow(ii)$.} Cette implication résulte de \ref{Serre4} et du lemme suivant dû à Landweber et Stong (\cite[Proposition 3]{LSt})~:

\begin{lem}\label{annulateur} Soit $\mathrm{Ann}(x)\subset\mathrm{H}^{*}V$ l'idéal annulateur de $x$. Alors l'idéal $\sqrt{\mathrm{Ann}(x)}$ est stable sous l'action de $\mathrm{A}$.
\end{lem}

\bigskip
\textit{Démonstration.} Soient $a$ un élément de $\mathrm{H}^{*}V$, $x$ un élément de $M$ et $m$ un entier avec $2^{m}>\vert x\vert$ ($\vert x\vert$ désigne le degré de $x$). On a $\mathrm{Sq}^{2^{m}i}(a^{2^{m}}x)=(\mathrm{Sq}^{i}a)^{2^{m}}x$ pour tout entier $i\geq 0$. Cette identité montre que $a\in\sqrt{\mathrm{Ann}(x)}$ implique $\mathrm{Sq}^{i}a\in\sqrt{\mathrm{Ann}(x)}$.
\hfill$\square\square$

\pagebreak

\renewcommand{\refname}{Références}

\end{document}